\documentclass{svmult}

\usepackage{makeidx}         
\usepackage{graphicx}        
\usepackage{multicol}        
\usepackage[bottom]{footmisc}
\usepackage{amsfonts,amsmath,amscd,color}

\makeindex             
\def\noz{m}
\def\noi{n}
\def\SVc{c(\cC)}
\def\toread#1{\subparagraph{#1: more serious reading.}}
\def\Cocycle{B}
\def\Re{\operatorname{Re}}
\def\Im{\operatorname{Im}}
\def\Vol{\operatorname{Vol}}
\def\Area{\operatorname{Area}}

\def\ind{\operatorname{ind}}
\def\id{\operatorname{id}}
\def\Id{\operatorname{Id}}
\def\N{\mathbb N}

\def\Z#1{{\mathbb Z}^{#1}}

\def\Q{\mathbb Q}

\def\R#1{{\mathbb R}^{#1}}

\def\C#1{{\mathbb C}^{#1}}

\def\torus{{\mathbb T}^2}
\def\T#1{{\mathbb T}^{#1}}

\def\Hyp{{\mathbb H}^2}

\def\cB{{\mathcal B}}
\def\cC{{\mathcal C}}

\def\cG{{\mathcal G}}
\def\cH{{\mathcal H}}
\def\cK{{\mathcal K}}
\def\cL{{\mathcal L}}
\def\cM{{\mathcal M}}
\def\cN{{\mathcal N}}
\def\cO{{\mathcal O}}
\def\cQ{{\mathcal Q}}
\def\cS{{\mathcal S}}
\def\cT{{\mathcal T}}
\def\cU{{\mathcal U}}
\def\cV{{\mathcal V}}

\newlength{\halfbls}

\spnewtheorem*{GeneralProblem}{General~Problem}{\bfseries}{\itshape}
\spnewtheorem*{CountingProblem}{Counting~Problem}{\bfseries}{\itshape}
\spnewtheorem*{Problem}{Problem}{\bfseries}{\itshape}

\spnewtheorem*{NNLemma}{Lemma}{\bfseries}{\itshape}
\spnewtheorem*{NNProposition}{Proposition}{\bfseries}{\itshape}
\spnewtheorem*{NNCorollary}{Corollary}{\bfseries}{\itshape}
\spnewtheorem*{NNTheorem}{Theorem}{\bfseries}{\itshape}
\spnewtheorem*{NNConjecture}{Conjecture}{\bfseries}{\itshape}
\spnewtheorem*{KeyTheorem}{Key Theorem}{\bfseries}{\itshape}
\spnewtheorem*{ErgodicTheorem}{Ergodic
Theorem}{\bfseries}{\itshape}

\spnewtheorem{Convention}{Convention}{\itshape}{\rmfamily}
\spnewtheorem*{Exercise}{Exercise}{\itshape}{\rmfamily}
\spnewtheorem*{Example}{Example}{\itshape}{\rmfamily}
\spnewtheorem*{NNRemark}{Remark}{\itshape}{\rmfamily}
\spnewtheorem{FProblem}{Problem}{\itshape}{\rmfamily}

\spnewtheorem*{KacLemma}{Kac Lemma}{\itshape}{\rmfamily}

\begin{document}
\title*{Flat Surfaces}

\titlerunning{Flat Surfaces}
\author{Anton Zorich}
\authorrunning{Anton Zorich}

\institute{IRMAR, Universit\'e de Rennes 1, Campus de Beaulieu,
35042 Rennes, France
\texttt{Anton.Zorich@univ-rennes1.fr}
   %
}

\maketitle

\begin{abstract}
Various problems of  geometry,  topology and dynamical systems on
surfaces as  well  as  some  questions concerning one-dimensional
dynamical systems lead to the  study  of  closed surfaces endowed
with a flat metric with  several  cone-type  singularities.  Such
flat surfaces are naturally organized into  families which appear
to be isomorphic to moduli spaces of holomorphic one-forms.

One can obtain much information about the geometry and dynamics of
an individual flat surface by studying both its orbit  under  the
Teichm\"uller geodesic flow and
under the linear group action. In
particular, the Teichm\"uller geodesic flow plays  the  role  of a time
acceleration machine (renormalization procedure) which allows  to
study the asymptotic  behavior  of  interval exchange transformations
and of surface foliations.

This survey is an attempt to present some selected ideas, concepts
and facts in Teichm\"uller dynamics in a playful way.
\end{abstract}

\begin{picture}(-5,-280)(-5,-280)
\put(-20,60){To appear in}
\put(-20,50){{\sc Frontiers in Number Theory, Physics, and Geometry Vol.I},}
\put(-20,40){P.~Cartier; B.~Julia; P.~Moussa; P.~Vanhove (Editors),}
\put(-20,30){Springer Verlag, 2006.}
\end{picture}

\noindent
57M50,  
32G15   
(37D40, 
37D50,  
30F30)  

\keywords{Flat  surface,   billiard  in  polygon,   Teichm\"uller
geodesic flow,  moduli space of Abelian differentials, asymptotic
cycle,  Lyapunov  exponent,   interval  exchange  transformation,
renormalization, Teichm\"uller disc, Veech surface}

\setcounter{minitocdepth}{2}
\dominitoc

\section{Introduction}
\label{zorich:s:Introduction}

These  notes  correspond  to    lectures   given  first  at
Les~Houches and later, in an extended version, at ICTP (Trieste).
As a result they keep all blemishes of oral presentations. I rush
to announce  important theorems as  facts, and then I deduce from
them numerous corollaries (which in reality  are  used  to  prove
these very  keystone theorems). I  omit proofs or replace them by
conceptual ideas  hiding  under  the  carpet  all  technicalities
(which sometimes constitute the main value of the proof). Even in
the choice  of the subjects I poach the most  fascinating issues,
ignoring  those   which are difficult to present  no  matter  how
important  the  latter ones are. These notes  also  contain  some
philosophical  discussions  and hopes which  some
emotional speakers like me include in their talks  and which one,
normally, never dares to put into a written text.

I am telling  all  this   to warn the reader  that  this  playful
survey of  some selected ideas, concepts  and facts in  this area
cannot replace any serious introduction in the subject and should
be taken with reservation.

As a much more serious accessible introduction I  can recommend a
collection        of        introductory        surveys        of
A.~Eskin~\cite{zorich:Eskin:Handbook},   G.~Forni~\cite{zorich:Forni:Handbook},
P.~Hubert      and      T.~Schmidt~\cite{zorich:Hubert:Schmidt:Handbook}
and H.~Masur~\cite{zorich:Masur:Handbook:1B},
organized  as  a  chapter  of  the
Handbook of Dynamical  Systems.  I also recommend recent surveys
of  H.~Masur  and S.~Tabachnikov~\cite{zorich:Masur:Tabachnikov}  and of
J.~Smillie~\cite{zorich:Smillie:billiards}.    The    part    concerning
renormalization   and  interval   exchange   transformations   is
presented          in          the           article           of
J.-C.~Yoccoz~\cite{zorich:Yoccoz:Les:Houches} of the current volume in a
much more responsible way than  my  introductory  exposition  in
Sec.~\ref{zorich:s:Renormalization:Rauzy:Veech:Induction}.

\subsection{Flat Surfaces}
\label{zorich:ss:Flat:Surfaces}

There is a common prejudice which makes us implicitly associate a
metric of  constant positive curvature to  a sphere, a  metric of
constant zero curvature  to  a torus,  and  a metric of  constant
negative curvature to a surface  of  higher  genus. Actually, any
surface can be endowed with  a  flat metric, no matter  what  the
genus of  this surface is...  with the only reservation that this
flat metric will have  several  singular points. Imagine that our
surface is made  from  plastic. Then we can  flatten  it from the
sides pushing all curvature  to  some small domains; making these
domains  smaller  and smaller  we  can  finally  concentrate  all
curvature at several points.

\begin{figure}[htb]
\centering
\includegraphics{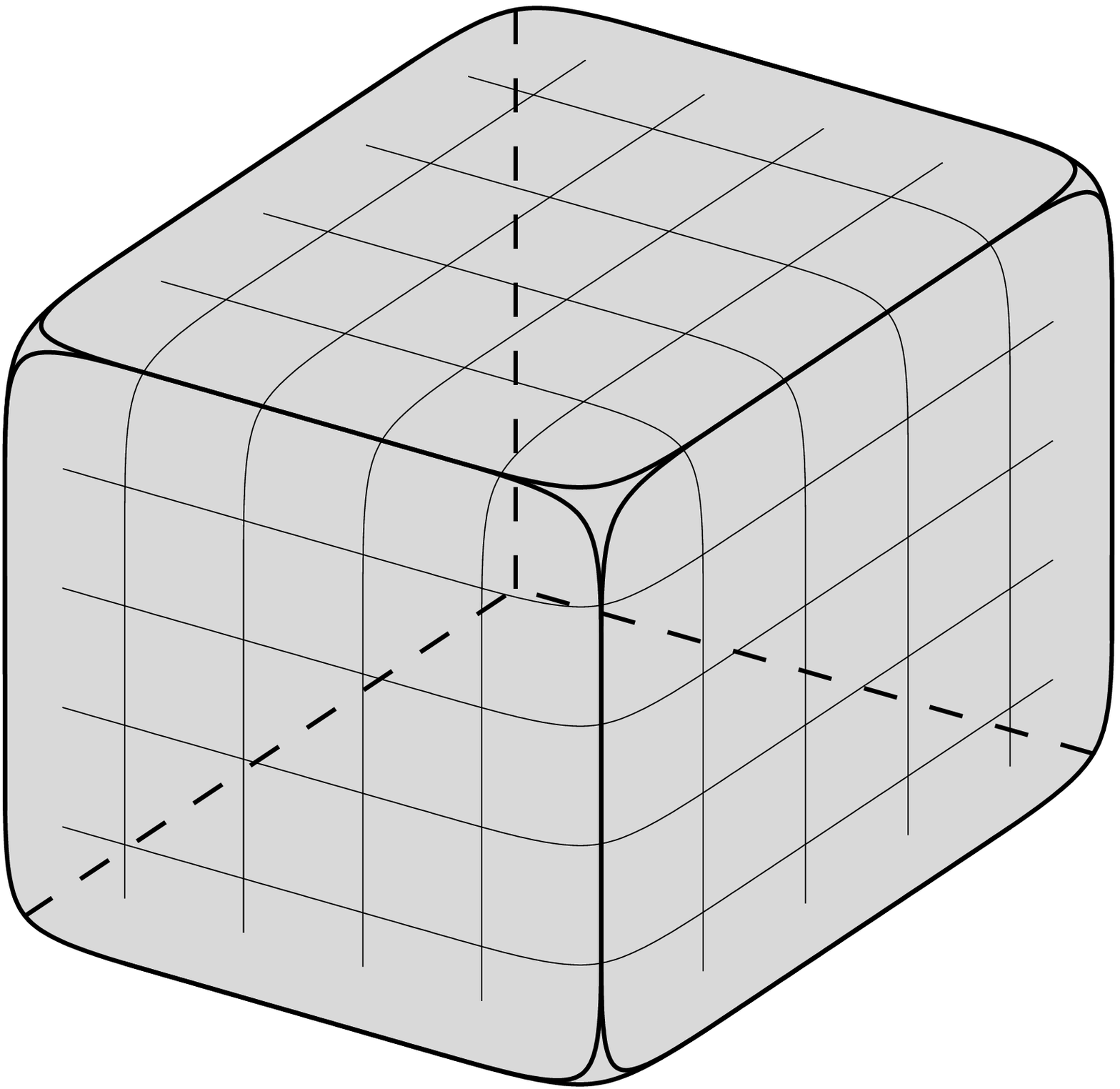}
\includegraphics{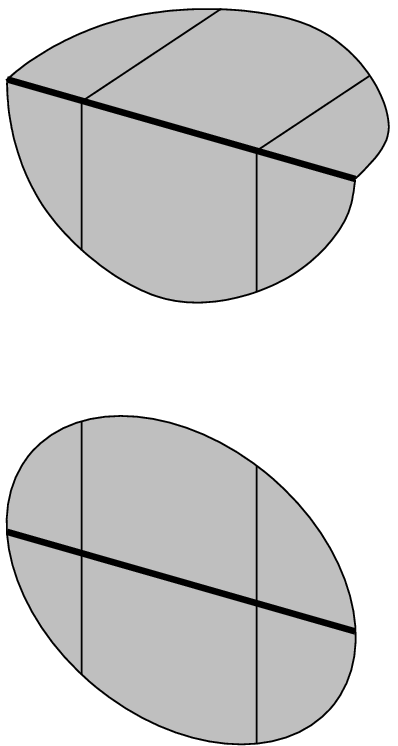}
\includegraphics{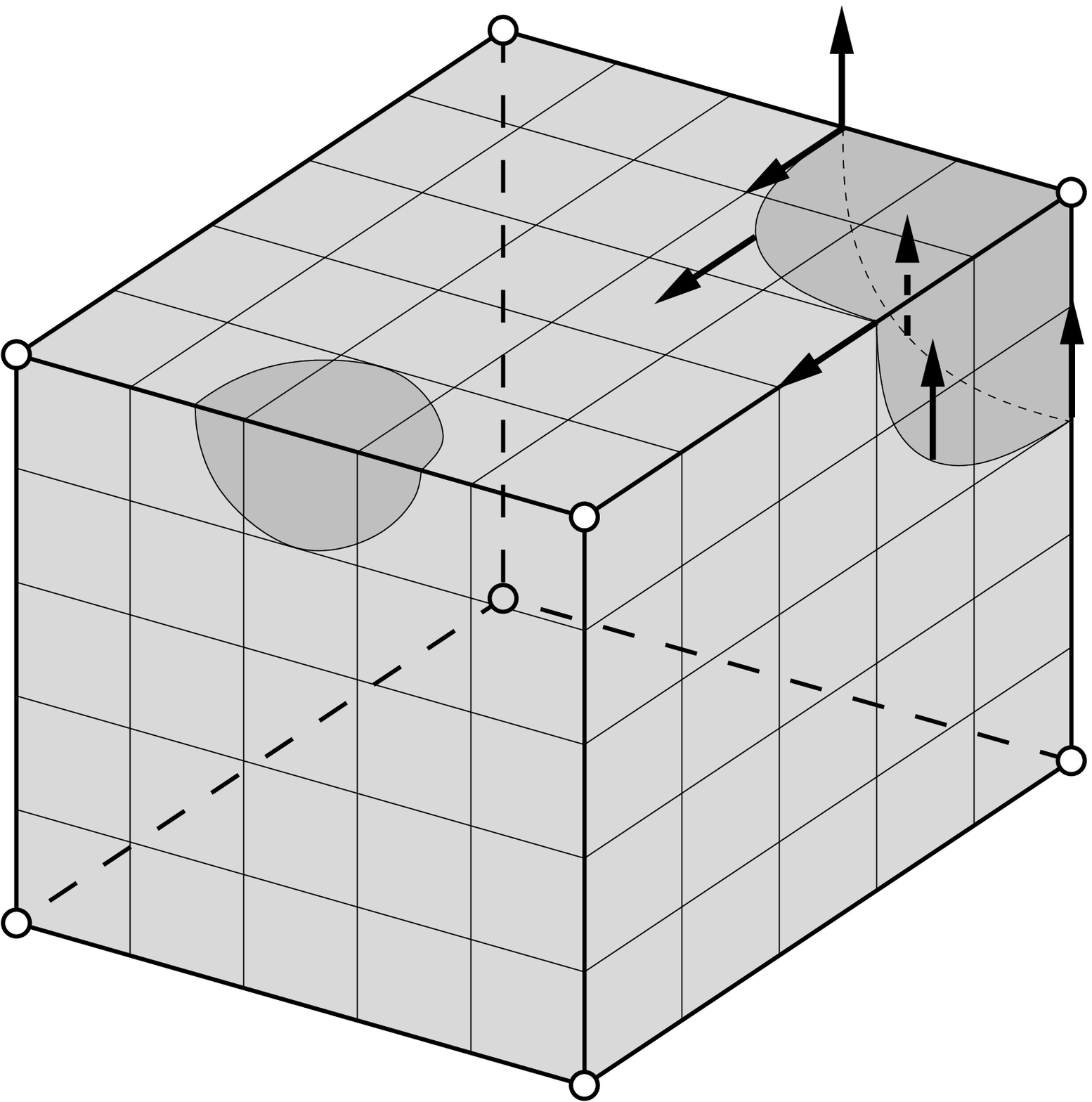}
\includegraphics{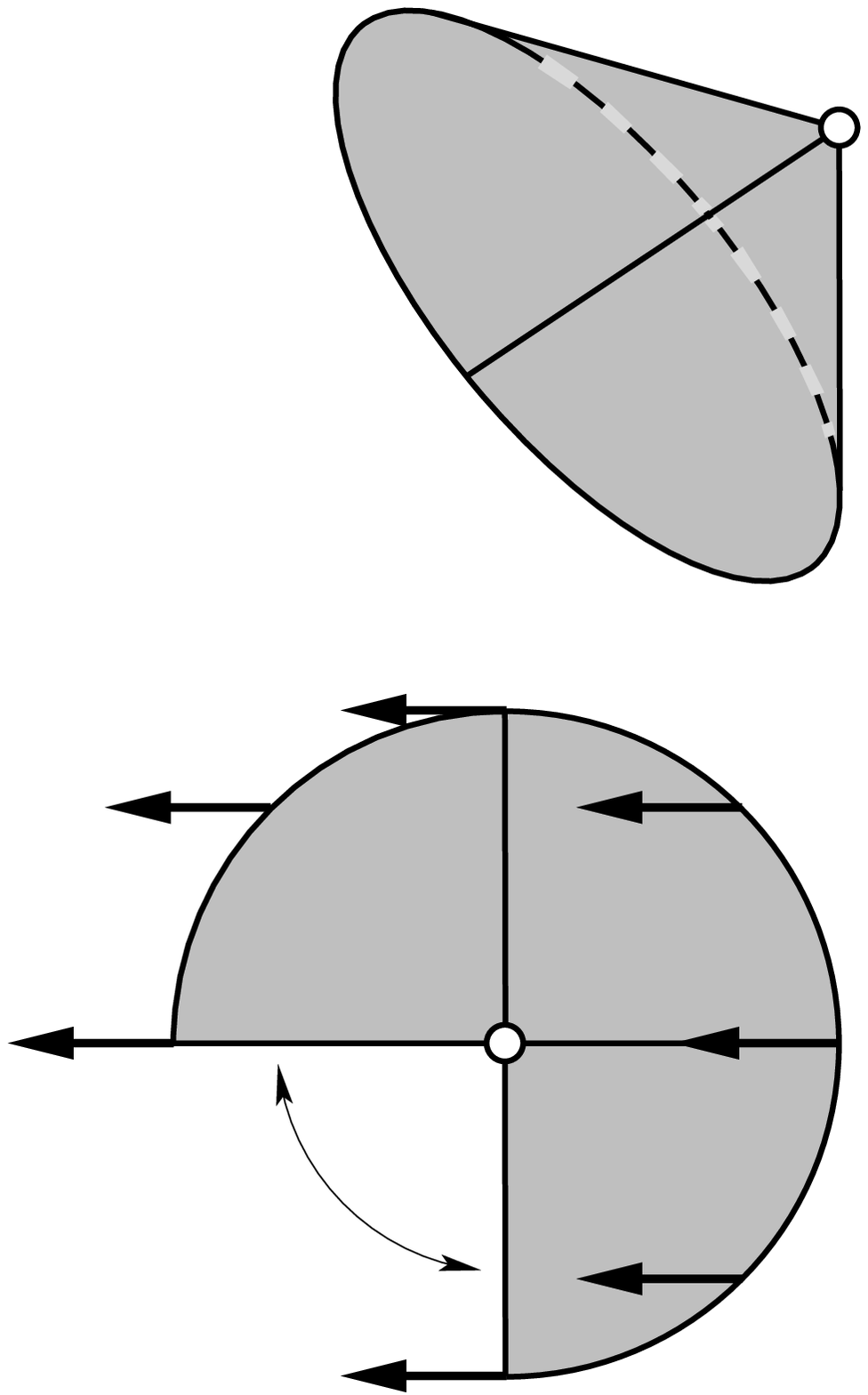}
\includegraphics{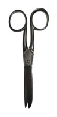}
\vspace{120bp}
\caption{
The surface of  the cube represents  a flat sphere with eight conical
singularities.
The metric \emph{does not} have singularities on the edges.
After parallel  transport  around a conical  singularity  a vector
comes back pointing to a direction  different  from  the  initial
one, so this flat metric has nontrivial \emph{holonomy}\index{Holonomy}.
\label{zorich:fig:cube}
}
\end{figure}

Consider the surface of a  cube.  It gives an
example  of  a  perfectly  flat   sphere   with   eight   conical
singularities corresponding to  eight vertices of the cube. Note
that  our  metric   is  nonsingular  on  edges:  taking  a  small
neighborhood of an interior point of an edge and unfolding  it we
get  a  domain in a Euclidean plane, see~Fig.~\ref{zorich:fig:cube}.  The
illusion of degeneration  of the metric  on the edges  comes  from the
singularity of the embedding of  our  flat  sphere into the Euclidean
space $\R{3}$.

However, the vertices of the cube correspond to
actual
\emph{conical singularities}\index{Conical!singularity}
of the metric.  Taking  a small neighborhood of  a  vertex we see
that it is isometric to a neighborhood of the vertex of a cone. A
flat cone is characterized by the
\emph{cone angle}\index{Cone angle}:
we  can  cut  the cone along a straight ray with an origin at the
vertex  of the  cone,  place the resulting  flat  pattern in  the
Euclidean plane and measure the  angle  between  the  boundaries,
see~Fig.~\ref{zorich:fig:cube}. Say, any vertex of  the  cube  has  cone
angle $3\pi/2$ which is easy to see since there are three squares
adjacent to any vertex, so  a  neighborhood of a vertex is  glued
from three right angles.

Having a  manifold (which is in our case  just a surface) endowed
with  a  metric it is quite natural  to  study
\emph{geodesics}\index{Geodesic},
which in a flat metric are locally isometric to straight lines.

\begin{GeneralProblem}
Describe the behavior of a generic
geodesic\index{Geodesic!generic}
on a flat surface. Prove
(or disprove) that the geodesic flow is
\index{Ergodic}
ergodic\footnote{In this context \emph{``ergodic''}
means  that  a typical  geodesic  will visit  any  region in  the
\emph{phase space} and, moreover, that in average it will spend a
time   proportional   to   the   volume  of  this   region;   see
Appendix~\ref{zorich:s:Ergodic:Theorem} for details.}
on a typical (in any reasonable sense) flat surface.

Does  any  (almost any)  flat  surface has  at  least one
closed geodesic\index{Geodesic!closed}
which does not pass through singular points?

If yes,  are there many  closed geodesics like that? Namely, find
the asymptotics for the  number  of\index{Geodesic!counting of periodic geodesics} closed geodesics shorter than
$L$ as the bound $L$ goes to infinity.
\end{GeneralProblem}

Believe it or  not  there has been no  (even  partial) advance in
solving  this  problem.  The  problem  remains  open even in  the
simplest case, when a surface is a sphere with only three conical
singularities; in particular, it is  not  known,  whether any (or
even almost any) such flat sphere has \emph{at  least one} closed
geodesic. Note that in this particular case, when  a flat surface
is a flat sphere with three conical singularities  the problem is
a  reformulation  of the corresponding billiard problem which  we
shall discuss in Sect.~\ref{zorich:ss:Billiards:in:Polygons}.

\subsection{Very Flat Surfaces}
\label{zorich:ss:Very:Flat:Surfaces}

A  general  flat  surface  with conical singularities  much  more
resembles a general Riemannian manifold  than  a  flat torus. The
reason is that it has nontrivial
\emph{holonomy}\index{Holonomy}.

Locally a  flat surface is isometric  to a Euclidean  plane which
defines  a  \emph{parallel transport} along paths on the  surface
with punctured conical  points. A parallel transport along a path
homotopic to  a trivial path on  this punctured surface  brings a
vector tangent to the surface to itself. However, if the  path is
not homotopic to  a trivial one,  the resulting vector  turns  by
some angle.  Say, a parallel  transport along a small closed path
around a  conical singularity makes  a vector turn exactly by the
cone  angle,  see  Fig.~\ref{zorich:fig:cube}.   (Exercise:   perform  a
parallel transport of a vector around a vertex of a cube.)

Nontrivial linear holonomy forces a generic geodesic to come back
and to intersect  itself again and again in different directions;
geodesics on a flat torus (which  has  trivial  linear  holonomy)
exhibit radically different behavior. Having chosen a direction to
the North, we can transport it to any other point of the
torus; the  result would not depend on the path. A geodesic on the
torus emitted  in some direction  will forever keep going in this
direction.  It  will either close up producing  a  regular  closed
geodesic, or  will never intersect  itself. In the latter case it
will produce a dense irrational winding line on the torus.

Fortunately,  the  class  of  flat surfaces with  trivial  linear
holonomy\index{Holonomy} is not reduced to flat tori. Since we cannot  advance in
the General Problem  from  the previous  section, from now on  we
confine ourselves to the study of these
\emph{very flat} surfaces\index{Surface!very flat}
(often called \emph{translation surfaces}):
that is, to closed orientable surfaces endowed with  a flat metric
having a finite number of conical singularities and having trivial
linear holonomy.

\begin{figure}[htb]
\centering
\includegraphics{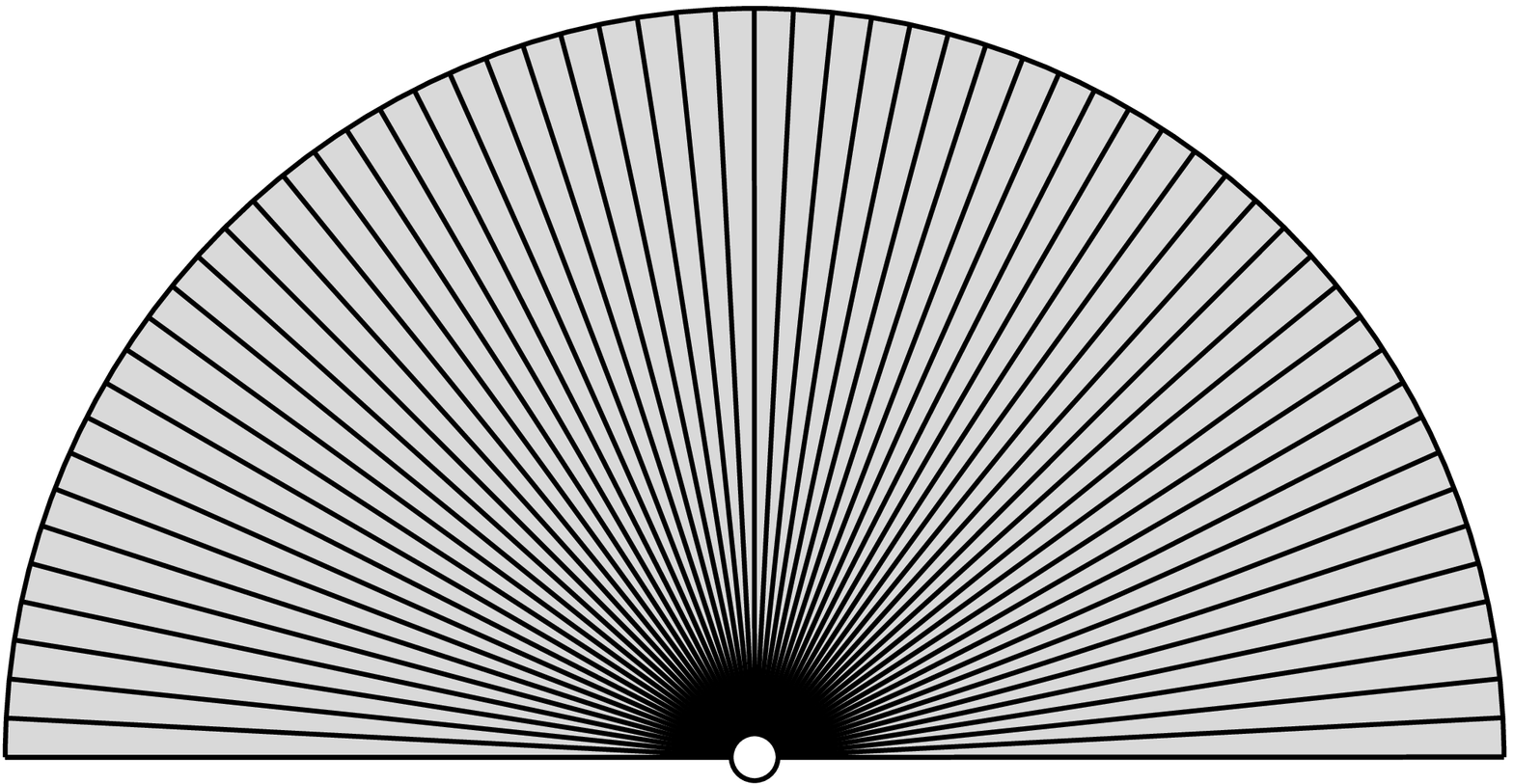}
\includegraphics{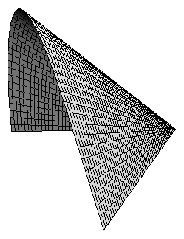}
\includegraphics{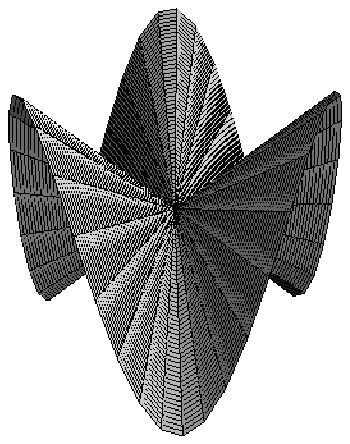}
   %
%
   %
   %
\vspace{140bp} 
\caption{
\label{zorich:fig:monkey:saddle}
A neighborhood of a conical point with a cone angle $6\pi$ can be
glued from six metric half discs
}
\end{figure}

Triviality  of  linear holonomy implies, in particular, that  all
cone angles\index{Cone angle} at
conical singularities\index{Conical!singularity}
are  integer  multiples  of
$2\pi$. Locally a neighborhood of such a conical point looks like
a ``\emph{monkey saddle}''\index{Conical!singularity},
see Fig.~\ref{zorich:fig:monkey:saddle}.

As  a  first  example  of a nontrivial very  flat  surface
consider a regular octagon with identified  opposite sides. Since
identifications  of  the  sides  are   isometries,   we   get   a
well-defined flat metric on  the  resulting surface. Since in our
identifications  we  used  only  parallel  translations  (and  no
rotations), we, actually, get  a very flat (translation) surface.
It is easy to see (check  it!) that our gluing rules identify all
vertices of the  octagon  producing a single conical singularity.
The cone  angle at this singularity is  equal to  the sum of  the
interior angles of the octagon, that is to $6\pi$.

\begin{figure}[htb]
\centering
\includegraphics{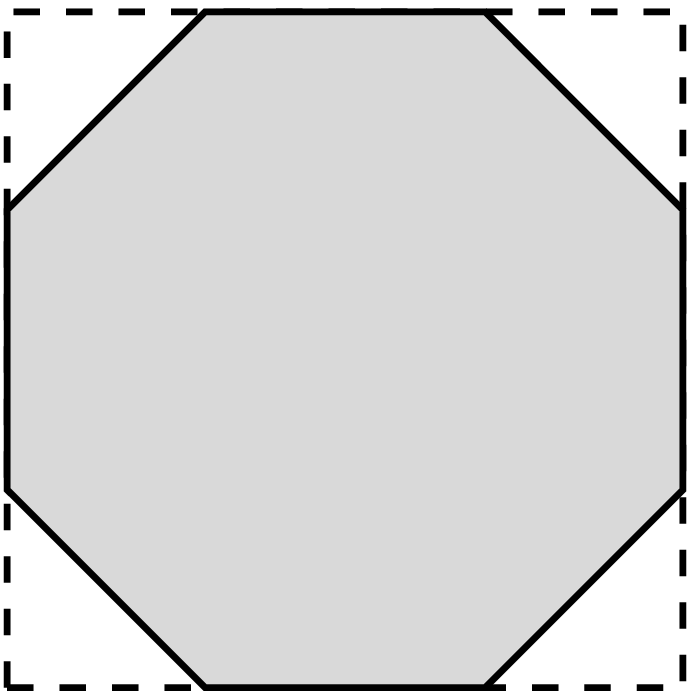}
\includegraphics{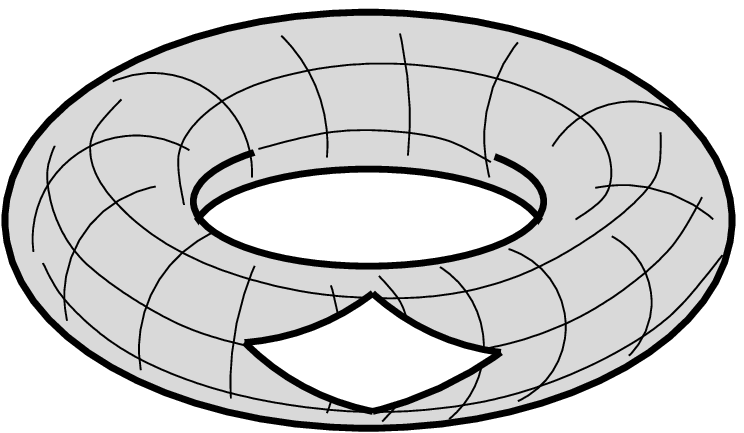}
\includegraphics{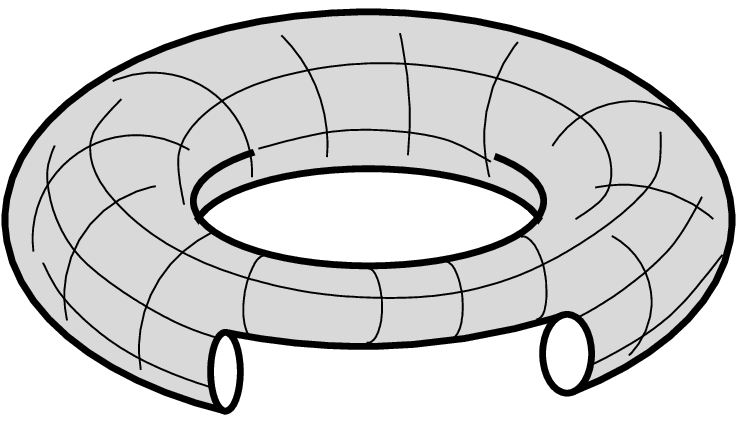}
\includegraphics{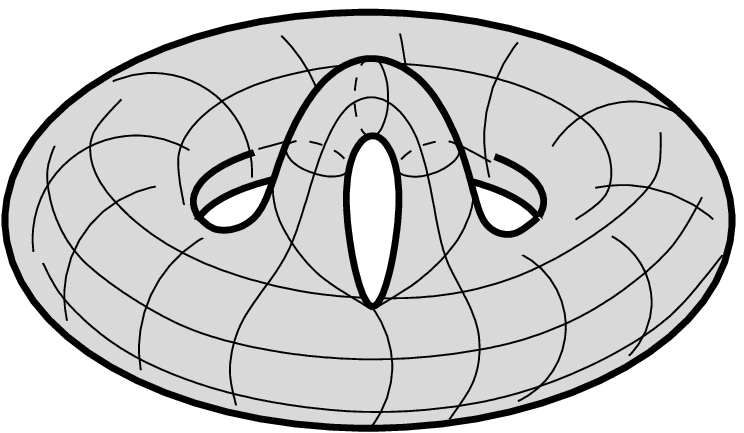}
\vspace{150bp}
\caption{
\label{zorich:fig:octagon:to:pretzel}
Gluing a pretzel from a regular octagon
}
\end{figure}

Figure~\ref{zorich:fig:octagon:to:pretzel} is an attempt to convince the
reader  that  the  resulting  surface  has  genus  two.  We first
identify  the  vertical  sides  and the horizontal sides  of  the
octagon obtaining a torus with a hole of the form of a square. To
simplify  the  drawing  we  slightly cheat: namely,  we  consider
another  torus  with  a hole of the form of a square, but the new
square hole is turned by $\pi/4$ with respect to the initial one.
Identifying a pair of horizontal sides of the hole by an isometry
we get a  torus with two  holes (corresponding to  the  remaining
pair  of  sides,  which  are  still   not  identified).  Finally,
isometrically identifying  the pair of  holes we get a surface of
genus two.

\begin{Convention}
\label{zorich:conv:flat:surface}
From now on  by  a \emph{flat  surface}\index{Flat surface}\index{Surface!flat}
we mean a closed  oriented
surface  with  a  flat  metric having a finite number  of  conical
singularities, such that the metric has  trivial linear holonomy\index{Holonomy}.
Moreover, we always assume that the flat surface  is endowed with
a distinguished  direction; we refer  to this direction as the
``direction to the North'' or as the
``\emph{vertical direction}''\index{Direction!vertical}.
\end{Convention}

The convention above implies, in  particular,  that  if we rotate
the octagon from Fig.~\ref{zorich:fig:octagon:to:pretzel} (which changes
the ``direction to the North'') and glue a flat surface from this
rotated octagon, this will give us a different flat surface.

We make three  exceptions to Convention~\ref{zorich:conv:flat:surface}
in this paper:
billiards in  general  polygons\index{Billiard!in polygon}
considered at the
beginning Sec.~\ref{zorich:ss:Billiards:in:Polygons}  give rise to  flat
metrics     with      nontrivial     linear     holonomy.      In
Sec.~\ref{zorich:ss:Toy:Example:Family:of:Flat:Tori}  we  consider  flat
tori forgetting the direction to the North.

Finally, in Sec.~\ref{zorich:ss:Extremal:Quasiconformal:Map} we consider
\emph{half-translation  surfaces}\index{Half-translation surface}\index{Surface!half-translation}
corresponding to  flat metrics
with holonomy group $\Z{}/2\Z{}$. Such  flat  metric  is a slight
generalization of a \emph{very flat} metric: a parallel transport
along a loop  may change the direction of a vector,
that is a  vector $\vec  v$ might  return  as $-\vec  v$ after  a
parallel transport.

\subsection{Synopsis and Reader's Guide}
\label{zorich:ss:Readers:Guide}

These lectures are an attempt to give some idea of what  is known
(and what is not known) about  flat surfaces, and to show what an
amazing and marvellous object a  flat  surface  is: problems from
dynamical  systems,   from  solid  state  physics,  from  complex
analysis,  from  algebraic  geometry,  from  combinatorics,  from
number theory, ... (the  list  can be considerably extended) lead
to the study of flat surfaces.

\paragraph{Section~\ref{zorich:s:Some:Motivations}. Motivations}
To give an  idea  of  how  flat surfaces  appear  in  different
guises we give some motivations in Sec.~\ref{zorich:s:Some:Motivations}.
Namely,  we  consider
billiards in polygons\index{Billiard!in polygon},
and, in  particular,
billiards in
rational polygons\index{Billiard!in rational polygon}\index{Polygon!rational}
(Sec.~\ref{zorich:ss:Billiards:in:Polygons}) and show that the consideration
of
billiard  trajectories\index{Billiard!trajectory}
is equivalent to the consideration  of
geodesics\index{Geodesic!in flat metric}
on   the   corresponding   flat  surface.  As  another
motivation                we                show               in
Sec.~\ref{zorich:ss:Electron:Transport:on:Fermi:Surfaces}    how     the
electron transport on
Fermi-surfaces\index{Fermi-surface}\index{Surface!Fermi-surface}
leads to  study of
foliation\index{Foliation!measured foliation}
defined    by    a    closed    1-form   on   a    surface.    In
Sec.~\ref{zorich:ss:Flows:on:Surfaces:and:Surface:Foliations}  we  show
that under some conditions on  the  closed  1-form such a foliation
can  be  ``straightened  out''  into  an  appropriate  flat  metric.
Similarly,  a  Hamiltonian  flow  defined  by  the  corresponding
multivalued Hamiltonian  on  a  surface  follows  geodesics in an
appropriate flat metric.

\paragraph{Section~\ref{zorich:s:Families:Of:Flat:Surfaces:and:Moduli:Spaces:of:Abelian:Differentials}.
Basic Facts}
A reader who is not interested in motivations can proceed
directly to
Sec.~\ref{zorich:s:Families:Of:Flat:Surfaces:and:Moduli:Spaces:of:Abelian:Differentials}
which describes the basic facts  concerning  flat  surfaces.  For
most of applications  it  is important  to  consider not only  an
individual flat surface, but an  entire  family  of flat surfaces
sharing the same topology: genus, number  and  types  of  conical
singularities.   In   Sec.~\ref{zorich:ss:Families:of:Flat:Surfaces}  we
discuss deformations of flat metric  inside  such  families. As a
model          example           we          consider          in
Sec.~\ref{zorich:ss:Toy:Example:Family:of:Flat:Tori} the  family of flat
tori. In Sec.~\ref{zorich:ss:Dictionary:of:Complex:Analytic:Language} we
show  that  a  flat  structure  naturally  determines  a  complex
structure  on   the   surface   and   a
holomorphic   one-form\index{Holomorphic!1-form}.
Reciprocally, a  holomorphic one-form naturally determines a flat
structure.  The dictionary  establishing  correspondence  between
geometric  language  (in  terms  of   the   flat   metrics)   and
complex-analytic language  (in terms of holomorphic one-forms) is
very  important   for   the   entire   presentation;   it   makes
Sec.~\ref{zorich:ss:Dictionary:of:Complex:Analytic:Language}        more
charged        than       an        average        one.        In
Sec.~\ref{zorich:ss:Moduli:Space:of:Holomorphic:One:Forms}  we  continue
establishing  correspondence  between   families   of  flat
surfaces and strata of
moduli spaces of holomorphic one-forms\index{Moduli space!of holomorphic 1-forms}\index{Stratum in the moduli space}.
In
Sec.~\ref{zorich:ss:Action:of:SL2R:on:the:Moduli:Space} we describe  the
action  of  the  linear  group $SL(2,\R{})$\index{Action on the moduli space!ofSL@of $SL(2,\R{})$}\index{0SL@$SL(2,\R{})$-action on the moduli space}
on flat  surfaces  --
another key issue of this theory.

We complete
Sec.~\ref{zorich:s:Families:Of:Flat:Surfaces:and:Moduli:Spaces:of:Abelian:Differentials}
with an attempt to present the following general principle in the
study of flat surfaces. In order to get some information about an
individual flat surface it is often very convenient  to find (the
closure  of)  the  orbit  of   corresponding   element   in   the
\emph{family}  of  flat surfaces under the action  of  the  group
$SL(2,\R{})$ (or, sometimes, under the  action  of  its  diagonal
subgroup). In many  cases  the structure  of  this orbit gives
comprehensive  information  about  the   initial   flat  surface;
moreover, this information might be  not  accessible  by a direct
approach.      These       ideas      are      expressed       in
Sec.~\ref{zorich:ss:General:Philosophy}.  This  general   principle   is
illustrated in Sec.~\ref{zorich:ss:Implementation:of:General:Philosophy}
presenting  Masur's   criterion   of
unique  ergodicity\index{Ergodic!uniquely ergodic}
of  the
directional flow  on a flat  surface. (A reader not familiar with
the ergodic theorem can  either  skip this  last  section or read  an
elementary       presentation        of       ergodicity       in
Appendix~\ref{zorich:s:Ergodic:Theorem}.)

\paragraph{Section~\ref{zorich:s:How:Do:Generic:Geodesics:Wind:Around:Flat:Surfaces}.
Topological Dynamics of Generic Geodesics}
This section  is independent from the  others; a reader  can pass
directly to  any of the further ones. However,  it gives a strong
motivation      for      renormalization       discussed       in
Sec.~\ref{zorich:s:Renormalization:Rauzy:Veech:Induction}  and  in   the
lectures   by   J.-C.~Yoccoz~\cite{zorich:Yoccoz:Les:Houches}  in   this
volume. It can  also be  used as  a  formalism for  the study  of
electron            transport            mentioned             in
Sec.~\ref{zorich:ss:Electron:Transport:on:Fermi:Surfaces}.

In  Sec.~\ref{zorich:ss:Asymptotic:Cycle}  we  discuss  the  notion   of
\emph{asymptotic cycle}\index{Asymptotic!cycle}
generalizing the \emph{rotation  number}
of an  irrational winding line on a torus.  It describes how
an ``irrational winding line'' on a surface of higher genus winds
around      a      flat     surface      in      average.      In
Sec.~\ref{zorich:ss:Deviation:from:Asymptotic:Cycle}  we   heuristically
describe the further terms of approximation, and we complete with
a    formulation     of     the     corresponding    result    in
Sec.~\ref{zorich:ss:Asymptotic:Flag:and:Dynamical:Hodge:Decomposition}.

In fact, this description is equivalent  to  the  description  of
the deviation  of a directional flow on a flat  surface  from the ergodic
mean.
Sec.~\ref{zorich:ss:Asymptotic:Flag:and:Dynamical:Hodge:Decomposition}
involves  some  background in ergodic theory; it  can  be  either
omitted  in the  first  reading, or can  be  read accompanied  by
Appendix~\ref{zorich:s:Multiplicative:ergodic:theorem}        presenting
the multiplicative ergodic theorem\index{Multiplicative ergodic theorem}\index{Ergodic!multiplicative ergodic theorem}\index{Theorem!multiplicative ergodic}.

\paragraph{Section~\ref{zorich:s:Renormalization:Rauzy:Veech:Induction}.
Renormalization}

This section  describes  the  relation  between the
Teichm\"uller geodesic flow\index{Teichm\"uller!geodesic flow}
and \emph{renormalization}\index{Renormalization}
for
\emph{interval exchange   transformations}\index{Interval exchange transformation}
discussed  in   the   lectures   of
J.-C.~Yoccoz~\cite{zorich:Yoccoz:Les:Houches} in  this volume. It
is slightly more technical than other sections and can be omitted
by a  reader who is not interested  in the  proof of the  Theorem
from
Sec.~\ref{zorich:ss:Asymptotic:Flag:and:Dynamical:Hodge:Decomposition}.

In Sec.~\ref{zorich:ss:First:Return:Map} we show that interval  exchange
transformations naturally arise as the \emph{first return map} of
a directional flow on a flat surface to a transversal segment. In
Sec.~\ref{zorich:ss:Evaluation:of:the:asymptotic:cycle:using:an:iet}  we
perform an  explicit  computation  of the
\emph{asymptotic cycle}\index{Asymptotic!cycle}
(defined   in   Sec.~\ref{zorich:ss:Asymptotic:Cycle})  using   interval
exchange                   transformations.                    In
Sec.~\ref{zorich:ss:Time:acceleration:machine:Renormalization}        we
present a  conceptual  idea  of
\emph{renormalization}\index{Renormalization},
a powerful
technique of acceleration of  motion along trajectories
of  the   directional   flow.   This   idea   is  illustrated  in
Sec.~\ref{zorich:ss:Euclidean:algorithm} in the simplest  case  where we
interpret the  Euclidean algorithm as a renormalization procedure
for rotation of a circle.

We  develop these  ideas  in  Sec.~\ref{zorich:ss:Rauzy:Veech:induction}
describing a concrete geometric renormalization procedure (called
\emph{Rauzy--Veech   induction}\index{Rauzy--Veech induction})
applicable  to   general  flat
surfaces  (and  general interval  exchange  transformations).  We
continue in Sec.~\ref{zorich:ss:space:of:iet}  with the elementary formalism
of
\emph{multiplicative      cocycles}\index{Cocycle!multiplicative}
(see       also
Appendix~\ref{zorich:s:Multiplicative:ergodic:theorem}).       Following
W.~Veech we describe  in  Sec.~\ref{zorich:ss:Teichmuller:Geodesic:Flow}
\emph{zippered  rectangles}\index{Zippered rectangle}
coordinates  in  a  family  of  flat
surfaces  and  describe the action of the
Teichm\"uller  geodesic flow\index{Teichm\"uller!geodesic flow}
in  a fundamental domain  in these coordinates. We show that
the first  return map of the  Teichm\"uller geodesic flow  to the
boundary  of   the   fundamental   domain   corresponds   to  the
Rauzy--Veech induction\index{Rauzy--Veech induction}.
In Sec.~\ref{zorich:ss:Spectrum:of:Lyapunov:exponents} we present  a  short
overview of recent results of  G.~Forni, M.~Kontsevich, A.~Avila  and  M.~Viana
concerning the  spectrum  of
\emph{Lyapunov  exponents}\index{Lyapunov exponent}\index{Exponent!Lyapunov exponent}
of  the
corresponding cocycle (completing  the  proof of the Theorem from
Sec.~\ref{zorich:ss:Asymptotic:Flag:and:Dynamical:Hodge:Decomposition}).
As   an    application    of    the    technique   developed   in
Sec.~\ref{zorich:s:Renormalization:Rauzy:Veech:Induction} we show  in
Sec.~\ref{zorich:ss:Encoding:a:Continued:Fraction} that  in the simplest
case of tori  it  gives the  well-known  encoding of a  continued
fraction\index{Continued fraction}  by  a cutting  sequence  of  a  geodesic  on  the upper
half-plane.

\paragraph{Section~\ref{zorich:s:Closed:Geodesics:and:Saddle:Connections:on:Flat:Surfaces}.
Closed geodesics}

This  section  is  basically  independent of other  sections;  it
describes the  relation  between  closed  geodesics on individual
flat surfaces and  ``cusps''\index{Cusp of the moduli space@Cusp of the moduli space, \emph{see also} Moduli space; principal boundary}
on the corresponding moduli spaces.
It might be  useful  for those who are  interested  in the global
structure of the moduli spaces.

Following    A.~Eskin    and    H.~Masur    we    formalize    in
Sec.~\ref{zorich:ss:Counting:Closed:Geodesics:and:Saddle:Connections}
the  counting  problems  for  closed  geodesics  and  for  saddle
connections of bounded length  on  an individual flat surface. In
Sec.~\ref{zorich:ss:Siegel:Veech:Formula} we present  the
Siegel--Veech Formula\index{Siegel--Veech!formula}\index{Formula!of Siegel--Veech}
and explain a  relation  between the counting problem and
evaluation of the volume of a tubular neighborhood  of a ``cusp''
in       the       corresponding      moduli       space.      In
Sec.~\ref{zorich:ss:Simplest:Cusps:of:the:Moduli:Space} we describe  the
structure  of  a  simplest  cusp.  We  describe the structure  of
general ``cusps''\index{Cusp of the moduli space@Cusp of the moduli space, \emph{see also} Moduli space; principal boundary}
(the structure of
\emph{principal boundary}\index{Moduli space!principal boundary of the moduli space@principal boundary of the moduli space, \emph{see also} Cusp of the moduli space}
of the moduli space) in
Sec.~\ref{zorich:ss:Multiple:Isometric:Geodesics:and:Principal:Boundary}.
As  an  illustration of  possible  applications  we  consider  in
Sec.~\ref{zorich:ss:Application:Billiards:in:Rectangular:Polygons}
billiards in rectangular polygons.

\paragraph{Section~\ref{zorich:s:Volume:of:the:Moduli:Space}
Volume of the Moduli Space}

In Sec.~\ref{zorich:ss:Square:tiled:surfaces}  we consider very  special
flat surfaces,  so called
\emph{square-tiled surfaces}\index{Square-tiled surface}\index{Surface!square-tiled}
which play
a  role  of
\emph{integer  points}\index{Lattice!in the moduli space}
in  the  moduli  space.  In
Sec.~\ref{zorich:ss:Approach:of:Eskin:and:Okounkov}   we   present   the
technique  of  A.~Eskin   and   A.~Okounkov  who  have  found  an
asymptotic formula for the number of  square-tiled surfaces glued
from a bounded  number of squares  and applied these  results  to
evaluation of volumes of moduli spaces\index{Moduli space!volume of the moduli space}.

As usual,  Sec.~\ref{zorich:s:Volume:of:the:Moduli:Space} is independent
of others; however,  the notion of a square-tiled surface appears
later   in   the   discussion   of
\emph{Veech surfaces}\index{Veech!surface}\index{Surface!Veech}
in Sec.~\ref{zorich:ss:Veech:Surfaces}--\ref{zorich:ss:Teichmuller:Discs}.

\paragraph{Section~\ref{zorich:s:Crash:Course:in:Teichmuller:Theory}.
Crash Course in Teichm\"uller Theory}
We  proceed  in
Sec.~\ref{zorich:s:Crash:Course:in:Teichmuller:Theory} with  a  very
brief  overview of some elementary  background  in Teichm\"uller
theory.      Namely,     we     discuss      in
Sec.~\ref{zorich:ss:Extremal:Quasiconformal:Map}   the
\emph{extremal quasiconformal   map}\index{Quasiconformal!extremal quasiconformal map}
and  formulate   the
\emph{Teichm\"uller theorem}\index{Teichm\"uller!theorem}\index{Theorem!Teichm\"uller},
which           we use           in
Sec.~\ref{zorich:ss:Teichmuller:Metric:and:Teichmuller:Geodesic:Flow}
we to    define    the    distance   between   complex
structures (\emph{Teichm\"uller metric}\index{Teichm\"uller!metric}\index{Metric!Teichm\"uller metric}).
We finally explain
why the action of the diagonal  subgroup in $SL(2;\R{})$  on the
space  of  flat surfaces  should  be  interpreted   as   the
\emph{Teichm\"uller geodesic flow}\index{Teichm\"uller!geodesic flow}.

\paragraph{Section~\ref{zorich:s:Hope:for:a:Magic:Stick:and:Recent:Results}.
Main Conjecture and Recent Results}
In this last section we  discuss  one of the central problems  in
the  area  --  a  conjectural structure of \emph{all}  orbits  of
$GL^+(2,\R{})$. The main hope  is  that the closure of \emph{any}
such orbit  is a nice complex subvariety, and  that in this sense
the moduli spaces of holomorphic 1-forms\index{Moduli space!of holomorphic 1-forms}
and the moduli spaces of
quadratic  differentials\index{Moduli space!of quadratic differentials}
resemble homogeneous  spaces  under  an
action of a unipotent group.  In  this section we also present  a
brief  survey  of  some  very  recent  results  related  to  this
conjecture obtained by K.~Calta, C.~McMullen and others.

We start  in Sec.~\ref{zorich:ss:SL2R:action:in:geometric:terms} with  a
geometric   description    of   the
$GL^+(2,\R{})$-action\index{Action on the moduli space!ofGL@of $GL^+(2,\R{})$}\index{0GL@$GL^+(2,\R{})$-action on the moduli space}
in Sec.~\ref{zorich:ss:SL2R:action:in:geometric:terms}  and  show  why  the
projections of the orbits
(so-called \emph{Teichm\"uller  discs}\index{Teichm\"uller!disc})
to  the  moduli  space
$\cM_g$\index{0M10@$\cM_g$ -- moduli space of complex structures}\index{Moduli space!of complex structures} of complex structures  should  be
considered       as
\emph{complex      geodesics}\index{Geodesic!complex geodesic}.       In
Sec.~\ref{zorich:ss:Geometric:Counterparts:of:Ratners:Theorem}        we
present some results telling that analogous ``complex geodesics''
in a homogeneous  space  have a very nice  behavior.  It is known
that  the   moduli  spaces  are  \emph{not}  homogeneous  spaces.
Nevertheless, in  Sec.~\ref{zorich:ss:Main:Hope}  we announce one of the
main hopes in  this  field  telling that in the  context  of  the
closures of ``complex  geodesics'' the moduli spaces behave as if
they were.

We  continue  with  a  discussion  of  two  extremal  examples of
$GL^+(2,\R{})$-invariant submanifolds. In
Sec.~\ref{zorich:ss:Classification:of:Connected:Components:of:the:Strata}
we describe the ``largest'' ones: the connected components of the
strata. In Sec.~\ref{zorich:ss:Veech:Surfaces} we consider flat surfaces
$S$ (called \emph{Veech surfaces}\index{Veech!surface}\index{Surface!Veech})
with the ``smallest'' possible
orbits: the ones  which are closed.  Since recently the  list  of
known  Veech  surfaces was  very  short.  However,  K.~Calta  and
C.~McMullen have discovered an infinite family  of Veech surfaces
in genus two and  have  classified them. Developing these results
C.~McMullen has proved the main  conjecture  in  genus two. These
results  of   K.~Calta   and   C.~McMullen   are   discussed   in
Sec.~\ref{zorich:ss:Revolution:in:Genus:Two}.  Finally,  we consider  in
Sec.~\ref{zorich:ss:Teichmuller:Discs}     the     classification     of
Teichm\"uller discs  in  $\cH(2)$ due to P.~Hubert, S.~Leli\`evre
and to C.~McMullen.

\paragraph{Section~\ref{zorich:s:Open:Problems}.
Open Problems}

In this section we collect  open  problems  dispersed through the
text.

\paragraph{Appendix~\ref{zorich:s:Ergodic:Theorem}.
Ergodic Theorem}
In  Appendix~\ref{zorich:s:Ergodic:Theorem}   we  suggest  a   two-pages
exposition of some key facts and constructions in ergodic theory.

\paragraph{Appendix~\ref{zorich:s:Multiplicative:ergodic:theorem}.
Multiplicative Ergodic Theorem}

Finally,  in  Appendix~\ref{zorich:s:Multiplicative:ergodic:theorem}  we
discuss the
Multiplicative Ergodic Theorem\index{Multiplicative ergodic theorem}\index{Ergodic!multiplicative ergodic theorem}\index{Theorem!multiplicative ergodic}
which  is mentioned in
Sec.~\ref{zorich:s:How:Do:Generic:Geodesics:Wind:Around:Flat:Surfaces}
and used in Sec.~\ref{zorich:s:Renormalization:Rauzy:Veech:Induction}.

We start  with  some  elementary  linear-algebraic motivations
in Sec.~\ref{zorich:ss:A:Crash:Course:of:Linear:Algebra} which  we
apply in Sec.~\ref{zorich:ss:Multiplicative:ergodic:theorem:torus}
to     the simplest case of a ``linear''  map  of  a
multidimensional torus. This  examples  give  us  intuition
necessary  to  formulate  in
Sec.~\ref{zorich:ss:Multiplicative:ergodic:theorem}
the \emph{multiplicative ergodic theorem}. Morally,  we
associate to an ergodic dynamical system a matrix  of \emph{mean
differential} (or of \emph{mean monodromy} in some  cases).  We
complete  this section  with  a  discussion  of   some   basic
properties   of \emph{Lyapunov  exponents}\index{Lyapunov exponent}\index{Exponent!Lyapunov exponent}
playing  a  role
of  logarithms  of eigenvalues of  the ``mean differential''
(``mean monodromy'') of the dynamical system.

\subsection{Acknowledgments}

I would like to thank organizers and participants of the workshop
\emph{``Frontiers in  Number Theory, Physics and Geometry''} held
at Les Houches, organizers and participants of the activity on
\emph{Algebraic and Topological Dynamics} held at  MPI, Bonn, and
organizers and  participants  of  the workshop on \emph{Dynamical
Systems}   held   at   ICTP,   Trieste,   for   their   interest,
encouragement, helpful  remarks  and for fruitful discussions. In
particular,  I   would   like   to   thank  A.~Avila,  C.~Boissy,
J.-P.Conze,   A.~Eskin,   G.~Forni,   P.~Hubert,   M.~Kontsevich,
F.~Ledrappier, S.~Leli\`evre, H.~Masur,   C.~McMullen,  Ya.~Pesin,  T.~Schmidt,
M.~Viana, Ya.~Vorobets and J.-C.~Yoccoz.

These notes would  be never written without tactful, friendly and
persistent pressure and help of B.~Julia and P.~Vanhove.

I would  like to  thank MPI f\"ur Mathematik at  Bonn and IHES at
Bures-sur-Yvette for their hospitality while preparation of these
notes. I highly  appreciate the help of V.~Solomatina and of M.-C.~Vergne who prepared
several most complicated pictures. I am grateful to M.~Duchin
and to G.~Le Floc'h for their kind permission to use the photographs.
I would like to thank the The State Hermitage Museum
and Succession H.~Matisse/VG Bild--Kunst for their kind permission
to use ``La Dance'' of H.~Matisse as an illustration in
Sec.~\ref{zorich:s:Closed:Geodesics:and:Saddle:Connections:on:Flat:Surfaces}.

\section{Eclectic Motivations}
\label{zorich:s:Some:Motivations}

In this section  we show how  flat surfaces appear  in  different
guises:  we  consider billiards in polygons, and, in  particular,
billiards         in         rational        polygons.         In
Sec.~\ref{zorich:ss:Billiards:in:Polygons} we show that consideration of
billiard trajectories  is  equivalent  to  the  consideration  of
geodesics  on   the   corresponding   flat  surface.  As  another
motivation                we                show               in
Sec.~\ref{zorich:ss:Electron:Transport:on:Fermi:Surfaces}    how     the
electron  transport  on Fermi-surfaces\index{Fermi-surface}\index{Surface!Fermi-surface}
leads  to  the  study  of
foliation  defined  by  a  closed   1-form   on   a  surface.  In
Sec.~\ref{zorich:ss:Flows:on:Surfaces:and:Surface:Foliations}  we  show,
that under some conditions on  the  closed  1-form such foliation
can  be  ``straightened  up''  in  an  appropriate  flat  metric.
Similarly,  a  Hamiltonian  flow  defined  by  the  corresponding
multivalued Hamiltonian  on  a  surface  follows  geodesics in an
appropriate flat metric.

\subsection{Billiards in Polygons}
\label{zorich:ss:Billiards:in:Polygons}

\paragraph{Billiards in General Polygons}

Consider a
\emph{polygonal billiard table}\index{Billiard!table}
and an  ideal billiard ball which reflects from  the walls of the
table by the ``optical'' rule: the angle of  incidence equals the
angle after the reflection. We assume that the mass of  our ideal
ball is concentrated at one point; there is no friction, no spin.

We mostly consider
\emph{regular trajectories}\index{Billiard!trajectory}\index{Trajectory!billiard trajectory},
which do  not pass through the  corners of the  polygon. However,
one  can  also  study  trajectories emitted from one  corner  and
trapped after several reflections in some  other  (or  the  same)
corner.  Such  trajectories  are  called  the
\emph{generalized diagonals}\index{Generalized diagonal}\index{Diagonal!generalized}.

To  simplify  the  problem  let us start our  consideration  with
billiards in  \emph{triangles}.  A  triangular  billiard table\index{Billiard!triangular}
is
defined by  angles $\alpha,\beta,\gamma$ (proportional  rescaling
of the  triangle does not change  the dynamics of  the billiard).
Since $\alpha+\beta+\gamma=\pi$ the family of triangular billiard
tables is described by two real parameters.

It is difficult to believe that the following Problem is open for
many decades.

\begin{Problem}[Billiard in a Polygon]\index{Billiard!in polygon}

\begin{enumerate}
\item  Describe  the  behavior  of  a  generic  regular  billiard
trajectory\index{Billiard!trajectory!generic}
in  a generic  triangle,  in  particular,  prove  (or
disprove) that the billiard flow is ergodic\footnote{On behalf of
the Center for Dynamics and Geometry  of Penn  State  University,
A.~Katok promised a prize of 10.000 euros for a solution  of this
problem.}\index{Ergodic};

\item
Does any (almost  any) billiard table  has at least  one  regular
periodic trajectory?\index{Billiard!trajectory!periodic}
If  the  answer  is  affirmative, does this
trajectory survive under deformations of the billiard table?

\item
If  a  periodic  trajectory  exists,  are   there  many  periodic
trajectories like  that?  Namely,  find  the  asymptotics for the
number of periodic trajectories\index{Billiard!counting of periodic trajectories}
of length shorter than $L$ as the
bound $L$ goes to infinity. \end{enumerate}
\end{Problem}

\begin{figure}[htb]
\centering
\includegraphics{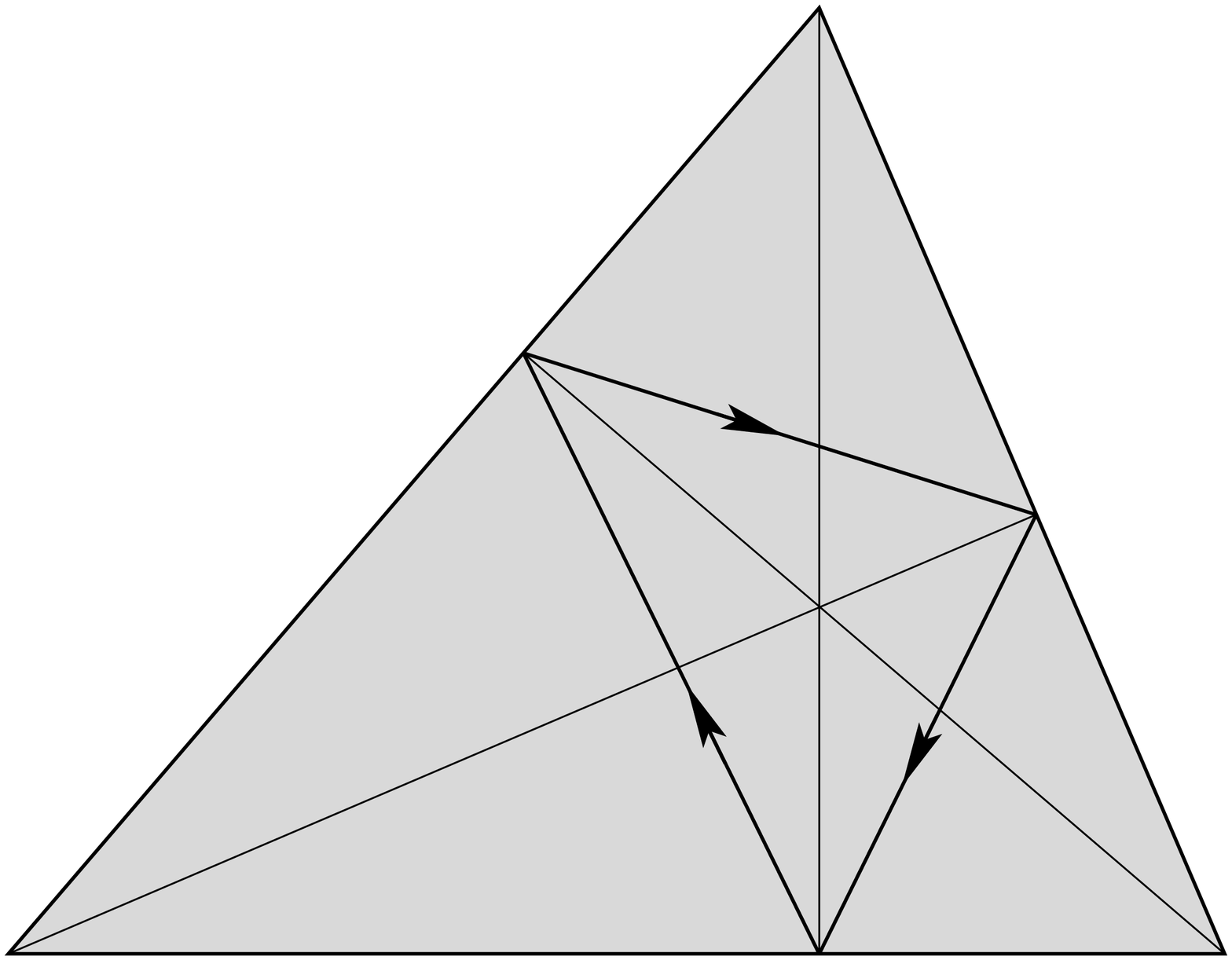}
\includegraphics{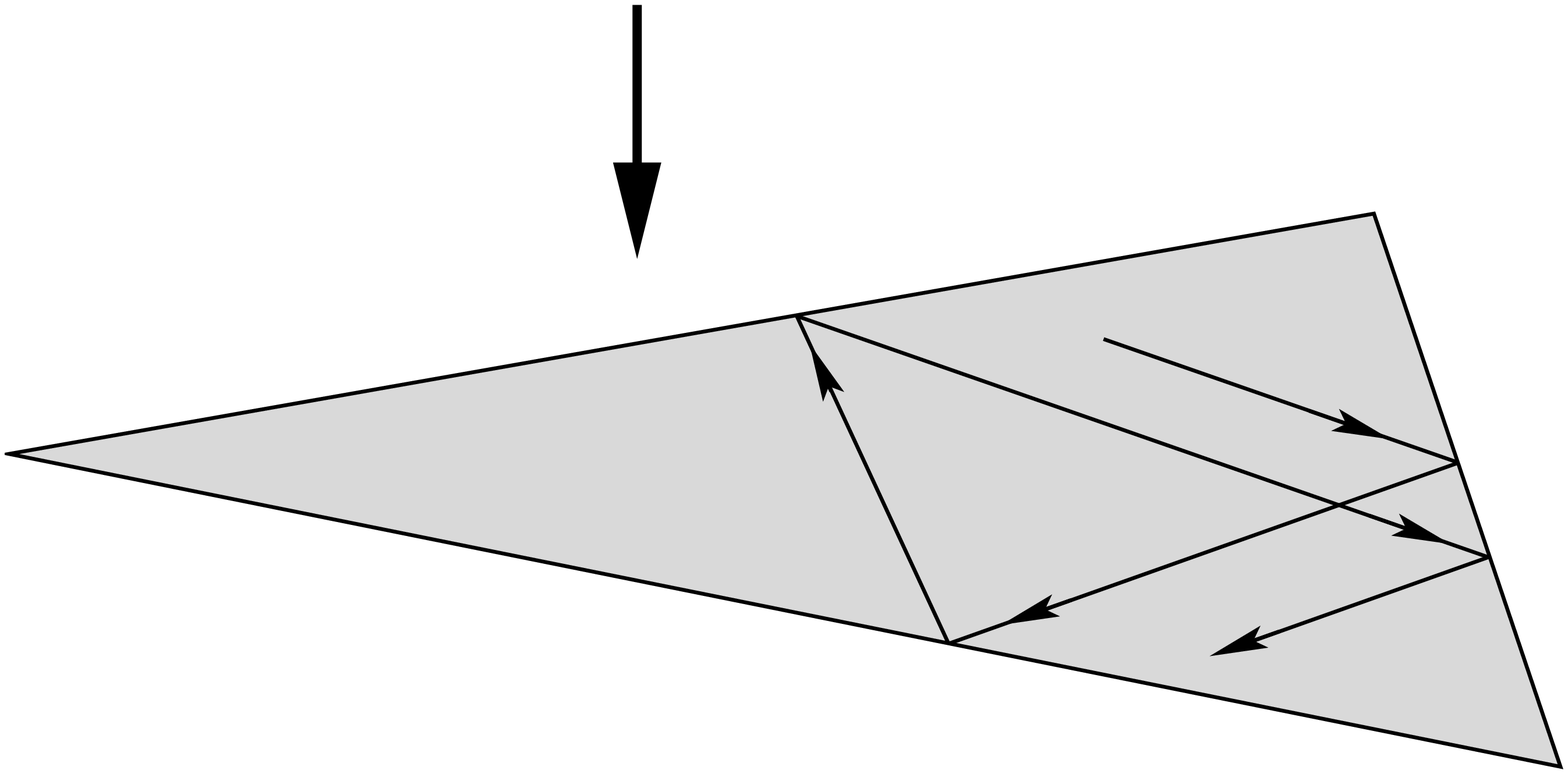}
\includegraphics{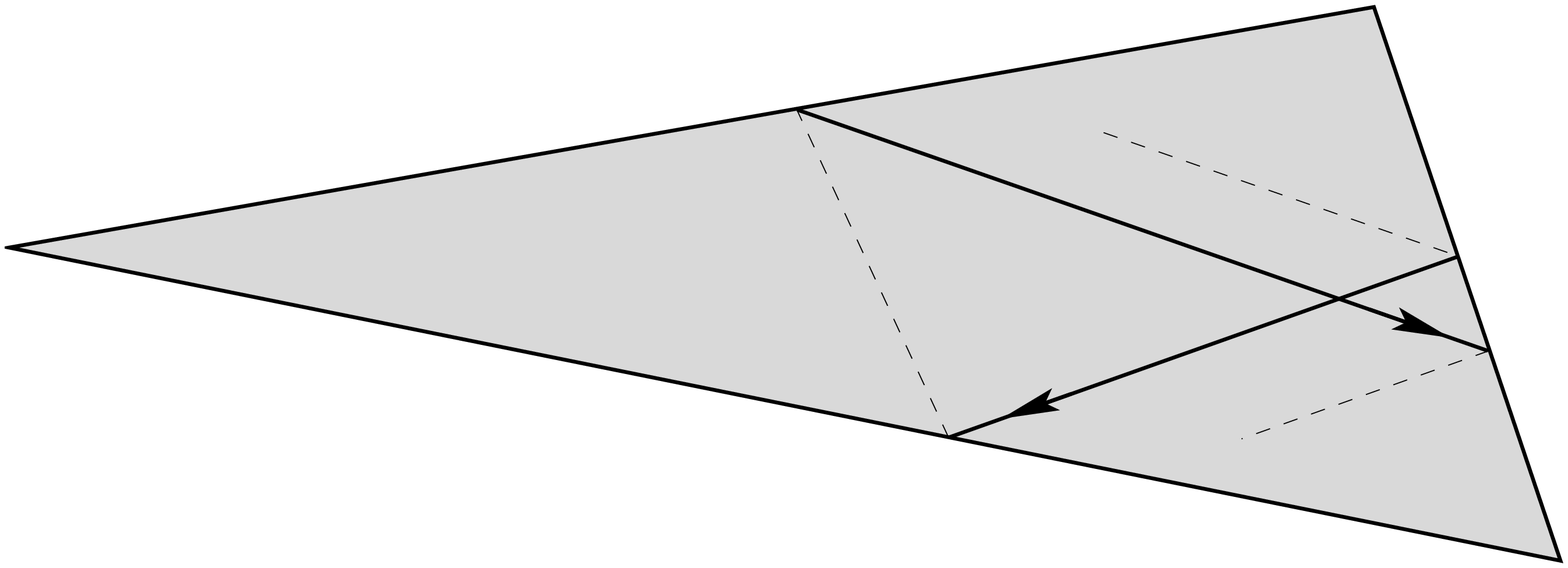}
\includegraphics{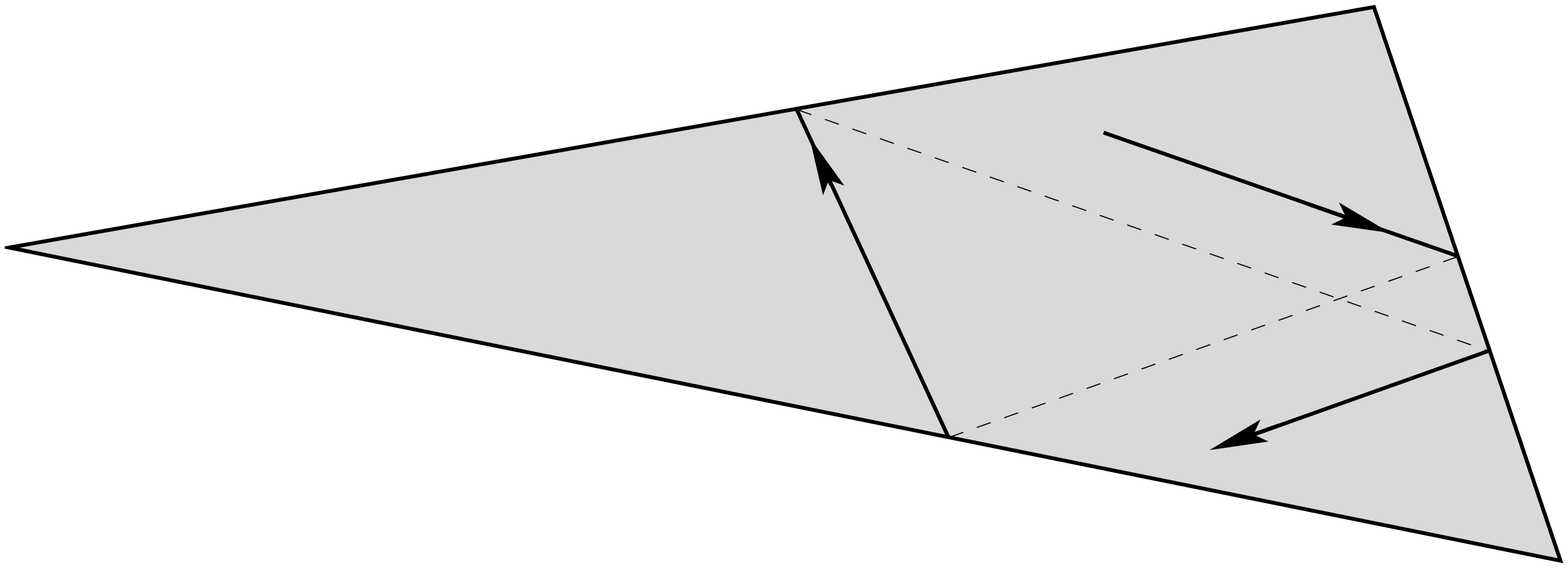}
\vspace{120bp}
\caption{Fagnano  trajectory and Fox--Kershner construction
\label{zorich:fig:faniano}
}
\end{figure}

It  is  easy to find a  special  closed regular trajectory in  an
acute triangle:  see  the  left picture at Fig.~\ref{zorich:fig:faniano}
presenting   the
\emph{Fagnano   trajectory}\index{Fagnano trajectory}\index{Trajectory!Fagnano trajectory}.
This   periodic
trajectory is known for at least two centuries.
However, it is not known whether any (or at least almost any)
obtuse triangle has a periodic billiard trajectory.

Obtuse triangles  with  the angles $\alpha\le\beta<\gamma$ can be
parameterized by a point of a  ``simplex''  $\Delta$  defined  as
$\alpha+\beta<   \pi/2$,  $\alpha,\beta>0$.   For   some   obtuse
triangles the existence  of a regular periodic trajectory is  known.
Moreover,   some   of   these   periodic   trajectories,   called
\emph{stable   periodic   trajectories},  survive   under   small
deformations of the  triangle, which proves existence of periodic
trajectories for some regions in  the  parameter  space  $\Delta$
(see the  works  of  G.~Galperin,  A.~M.~Stepin, Ya.~Vorobets and
A.~Zemliakov~\cite{zorich:Galperin:Stepin:Vorobets:1},
\cite{zorich:Galperin:Stepin:Vorobets:2},     \cite{zorich:Galperin:Zemliakov},
\cite{zorich:Vorobets:billiards:in:triangles}). It remains to prove that
such regions cover the entire parameter space $\Delta$. Currently
R.~Schwartz is in progress of extensive computer search of stable
periodic trajectories hoping to cover $\Delta$ with corresponding
computer-generated regions.

Now, following  Fox  and Kershner~\cite{zorich:Fox:Kershner}, let us see
how billiards in polygons lead naturally  to  geodesics  on  flat
surfaces.

Place  two copies  of  a polygonal billiard  table  one atop  the
other. Launch a billiard trajectory on one of the copies  and let
it jump from one copy to the other after each reflection (see the
right   picture   at  Fig.~\ref{zorich:fig:faniano}).   Identifying  the
boundaries of the two copies  of  the polygon we get a  connected
path $\rho$  on  the corresponding topological sphere. Projecting
this path to any of the two ``polygonal hemispheres'' we  get the
initial billiard trajectory.

It remains to note that our topological sphere is endowed  with a
flat metric (coming from the  polygon).  Analogously  to the flat
metric on the surface of a cube which is nonsingular on the edges
of    the     cube    (see    Sec.~\ref{zorich:ss:Flat:Surfaces}     and
Fig.~\ref{zorich:fig:cube}),  the   flat   metric   on   our  sphere  is
nonsingular on  the  ``equator''  obtained  from  the  identified
boundaries of  the two equal ``polygonal hemispheres''. Moreover,
let $x$ be a point where the path $\rho$ crosses the ``equator''.
Unfolding a neighborhood of a point $x$ on the ``equator'' we see
that the  corresponding fragment of  the path $\rho$ unfolds to a
straight segment in  the flat metric.  In other words,  the  path
$\rho$ is a geodesic in the corresponding flat metric.

The  resulting  flat  metric  is  \emph{not  very  flat}:  it has
nontrivial linear holonomy\index{Holonomy}.
The conical singularities\index{Conical!singularity}
of the flat
metric correspond to the vertices of the polygon;
the cone angle\index{Cone angle}
of a  singularity is twice  the angle at the corresponding vertex
of the polygon.

We have proved that every geodesic\index{Geodesic!in flat metric}
on our flat sphere projects to
a billiard trajectory\index{Billiard!trajectory}
and every  billiard  trajectory  lifts to a
geodesic.     This     is    why     General     Problem     from
Sec.~\ref{zorich:ss:Flat:Surfaces}    is    so   closely    related   to
the Billiard Problem.

\paragraph{Two Beads on a Rod and Billiard in a Triangle}

\begin{figure}[htb]
\centering
\includegraphics{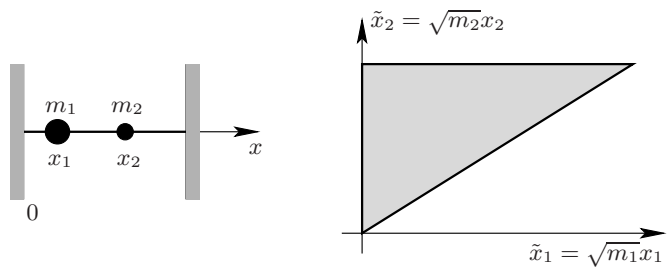}
\begin{picture}(0,0)(0,0)
\put(-108,-46){$m_1$}
\put(-83,-46){$m_2$}
\put(-107,-66){$x_1$}
\put(-81,-66){$x_2$}
\put(-31,-62){$x$}
\put(-115,-87){$0$}
\put(75,-102){$\tilde x_1=\sqrt{m_1} x_1$}
\put(15,-14){$\tilde x_2=\sqrt{m_2} x_2$}
\end{picture}
\vspace{100bp}
\caption{
Gas of two molecules in a one-dimensional chamber
\label{zorich:fig:two:balls:on:a:rod}
}
\end{figure}

It would be  unfair  not to  mention  that billiards in  polygons
attracted a lot of attention  as  (what initially seemed to be)  a
simple  model of a  Boltzman  gas.  To  give a  flavor  of  this
correspondence we consider a system of two elastic beads confined
to    a     rod     placed     between     two     walls,     see
Fig.~\ref{zorich:fig:two:balls:on:a:rod}.  (Up   to   the   best  of  my
knowledge   this   construction   originates   in   lectures   of
Ya.~G.~Sinai~\cite{zorich:Sinai:73}.)

The  beads  have different masses $m_1$ and  $m_2$  they  collide
between themselves, and also with  the  walls.  Assuming that the
size of the beads is negligible we can describe the configuration
space of our system using coordinates $0<x_1\le x_2\le  a$ of the
beads, where $a$ is the distance between the walls. Rescaling the
coordinates as
$$
\begin{cases}
\tilde{x}_1 &= \sqrt{m_1} x_1\\
\tilde{x}_2 &= \sqrt{m_2} x_2
\end{cases}
$$
we see that the configuration space in the new coordinates is given
by a right triangle $\Delta$, see
Fig.~\ref{zorich:fig:two:balls:on:a:rod}. Consider now a trajectory of
our dynamical system. We leave to the reader the pleasure to prove
the following elementary Lemma:

\begin{NNLemma}
In coordinates $(\tilde{x}_1,\tilde{x}_2)$ trajectories of the
system of two beads on a rod correspond to billiard trajectories\index{Billiard!trajectory}
in the triangle $\Delta$.
\end{NNLemma}

\paragraph{Billiards in Rational Polygons}

We have  seen that  taking two copies of a  polygon we can reduce
the study  of a billiard  in a  general polygon to  the study  of
geodesics  on   the  corresponding  flat  surface.  However,  the
resulting flat  surface has \emph{nontrivial} linear holonomy\index{Holonomy}, it
is \emph{not} ``very flat''.

Nevertheless, a  more  restricted  class  of  billiards,  namely,
billiards in
\emph{rational  polygons}\index{Polygon!rational}\index{Billiard!in rational polygon},
lead  to
``very flat''\index{Surface!very flat}
(translation) surfaces.

A polygon is called
\emph{rational}\index{Polygon!rational}
if all its  angles are rational  multiples of $\pi$.  A  billiard
trajectory  emitted  in some direction will change direction
after the first reflection,  then  will change direction once
more after the second reflection, etc.  However,  for  any  given
billiard trajectory in a rational  billiard  the  set of possible
directions is  \emph{finite},  which  make  billiards in rational
polygons so different from general ones.

As a basic example consider  a  billiard in a rectangle. In  this
case a generic trajectory  at any moment goes in one of four
possible directions. Developing the idea with the general polygon
we can  take \emph{four} (instead  of two) copies of our billiard
table (one  copy for each direction).  As soon as  our trajectory
hits the  wall and  changes the direction we make  it jump to the
corresponding copy of the billiard  --  the  one representing the
corresponding direction.

\begin{figure}[htb]
\centering
\includegraphics{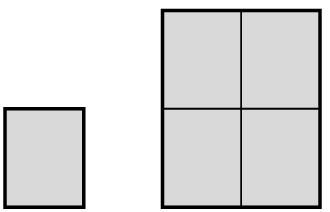}
\includegraphics{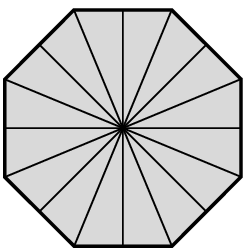}
\includegraphics{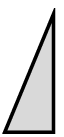}
\includegraphics{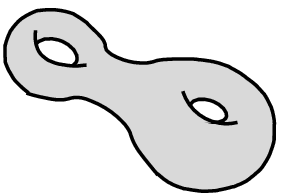}
\includegraphics{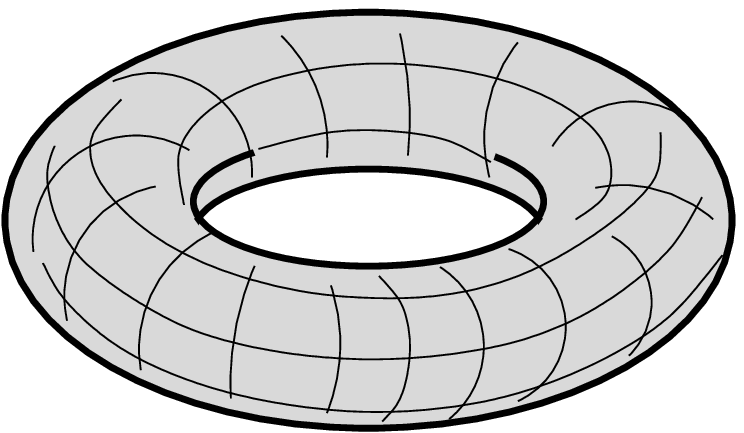}
\vspace{165bp} 
\caption{Billiard in a rectangle corresponds to directional
flow on a flat torus. Billiard in a right triangle $(\pi/8,
3\pi/8, \pi/2)$ leads to directional flow on a flat surface
obtained from the regular octagon.
\label{zorich:fig:billiard:in:a:polygon}
}
\end{figure}

By  construction  each side  of  every copy  of  the billiard  is
identified with exactly one side of another copy of the billiard.
Upon  these  identifications the  four  copies  of  the  billiard
produce a closed surface and  the  unfolded  billiard  trajectory
produces a connected  line  on this  surface.  We suggest to  the
reader to  check that \emph{the  resulting surface is a torus and
the unfolded  trajectory is a geodesic  on this flat  torus}, see
Fig.~\ref{zorich:fig:billiard:in:a:polygon}.

A similar unfolding construction (often called
\emph{Katok--Zemliakov construction})\index{Katok--Zemliakov construction}
works for a billiard in any rational polygon\index{Billiard!in rational polygon}.
Say, for a
billiard in a right triangle with angles $(\pi/8, 3\pi/8,
\pi/2)$ one has to take $16$ copies (corresponding to $16$
possible directions of a  given  billiard trajectory).
Appropriate  identifications  of these $16$  copies  produce  a
regular  octagon with identified opposite  sides  (see
Fig.~\ref{zorich:fig:billiard:in:a:polygon}).  We know       from
Sec.~\ref{zorich:ss:Very:Flat:Surfaces}        and from
Fig.~\ref{zorich:fig:octagon:to:pretzel}  that  the  corresponding flat
surface is a ``very flat'' surface\index{Surface!very flat}
of  genus two having a single
conical singularity with the cone angle $6\pi$.

\begin{figure}[htb]
\centering
\includegraphics{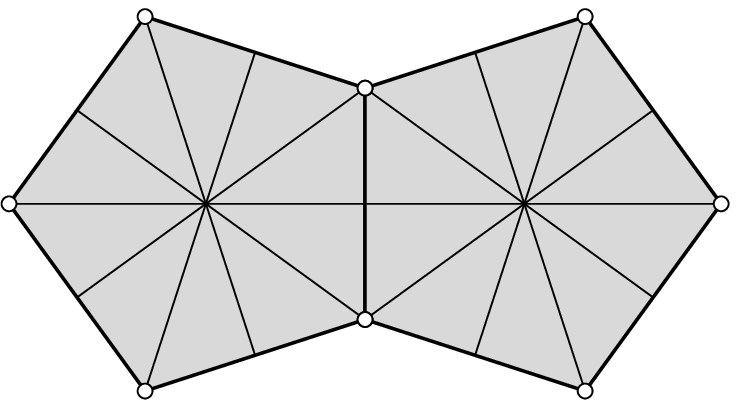}
\includegraphics{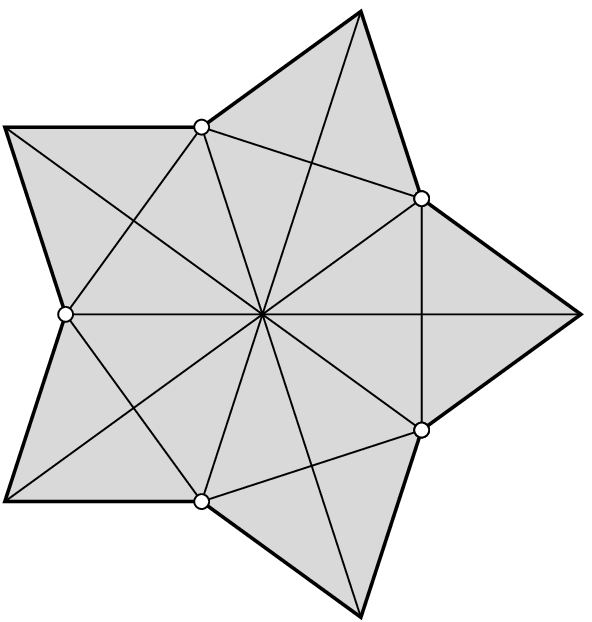}
\vspace{145bp} 
\caption{   We   can  unfold the billiard  in  the  right triangle
$(\pi/2,\pi/5,3\pi/10)$  into different  polygons. However,  the
resulting     very     flat     surfaces     are the     same
(see also~\cite{zorich:Hubert:Schmidt:Handbook}).
\label{zorich:fig:different:unfoldings}
}
\end{figure}

\begin{Exercise}
What is the genus of the   surface  obtained  by  Katok--Zemliakov
construction  from  an  isosceles   triangle   $(3\pi/8,  3\pi/8,
\pi/4)$? How many conical points\index{Conical!singularity}
does it have? What are  the cone
angles\index{Cone angle}
at  these points? Hint: this surface  \emph{can  not}  be
glued from a regular octagon.

It is quite common to unfold a rational billiard in two steps. We
first unfold  the billiard table  to a polygon, and then identify
the appropriate pairs of sides of  the  resulting  polygon.  Note
that  the  polygon  obtained  in  this  intermediate step is  not
canonical.

Show that a generic billiard trajectory in the right triangle
with angles $(\pi/2,\pi/5,3\pi/10)$ has $20$ directions. Show that
both  polygons  at   Fig.~\ref{zorich:fig:different:unfoldings}  can  be
obtained by  Katok--Zemliakov  construction  from  $20$ copies of
this triangle. Verify that after identification of parallel sides
of these polygons we obtain isometric  very  flat  surfaces  (see
also~\cite{zorich:Hubert:Schmidt:Handbook}).   What   genus,  and   what
conical points do they have?  What  are the cone angles at  these
points? \end{Exercise}

Note that in comparison  with  the initial construction, where
we had  only  two  copies  of  the  billiard  table  we get  a
more complicated surface.  However, what we gain is that in this
new construction our flat surface is actually ``very flat''\index{Surface!very flat}:
it
has trivial linear holonomy.  It has a lot   of   consequences;
say,   due    to    a    Theorem    of
H.~Masur~\cite{zorich:Masur:existence:of:closed:regular:geodesic} it is
possible to find a regular periodic geodesic on \emph{any}
``very flat''  surface.  If  the  flat  surface  was constructed
from  a billiard, the corresponding closed geodesic\index{Geodesic!closed}
projects to
a regular periodic trajectory\index{Billiard!trajectory!periodic}
of  the  corresponding billiard\index{Billiard!in rational polygon}
which solves part  of the Billiard Problem
for billiards  in rational polygons.

We did  not intend to present  in this section  any comprehensive
information  about  billiards,  our  goal  was  just  to  give  a
motivation for the study of flat surfaces. A reader interested in
billiards can  get a good idea on  the subject  from a very  nice
book   of    S.~Tabachnikov~\cite{zorich:Tabachnikov}.   Details   about
billiards in polygons (especially rational polygons) can be found
in  the  surveys of  E.~Gutkin~\cite{zorich:Gutkin:survey:on:billiards},
P.~Hubert and T.~Schmidt~\cite{zorich:Hubert:Schmidt:Handbook},
H.~Masur    and    S.~Tabachnikov~\cite{zorich:Masur:Tabachnikov}    and
J.~Smillie~\cite{zorich:Smillie:billiards}.

\subsection{Electron Transport on Fermi-Surfaces}
\label{zorich:ss:Electron:Transport:on:Fermi:Surfaces}

\begin{figure}[htb]
\centering
\includegraphics{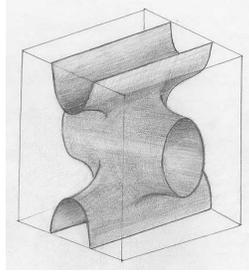}
\vspace{100bp} 
\caption{
Riemann surface of genus $3$ embedded into a torus $\T{3}$
\label{zorich:fig:periodic:surface}
}
\end{figure}

Consider  a  periodic  surface  $\tilde M^2$ in $\R{3}$  (i.e.  a
surface invariant under translations by  any  integer  vector  in
$\Z{3}$). Such a surface can be constructed in a fundamental
domain of a cubic  lattice,  see Fig.~\ref{zorich:fig:periodic:surface},
and  then  reproduced repeatedly in the lattice.  Choose  now  an
affine plane in $\R{3}$ and consider an intersection  line of the
surface  by  the plane. This intersection line  might  have  some
closed components and it may also have some unbounded components.
The question is  \emph{how  does an unbounded component propagate
in $\R{3}$}?

The study of this subject  was  suggested  by S.~P.~Novikov about
1980 (see~\cite{zorich:Novikov:82}) as a mathematical formulation of the
corresponding problem  concerning electron transport in metals. A
periodic surface represents a
\emph{Fermi-surface}\index{Fermi-surface}\index{Surface!Fermi-surface},
affine  plane
is a  plane orthogonal to  a magnetic field, and the intersection
line  is  a  trajectory  of  an   electron   in   the   so-called
\emph{inverse lattice}.

\begin{figure}[hbt]
\centering
\includegraphics{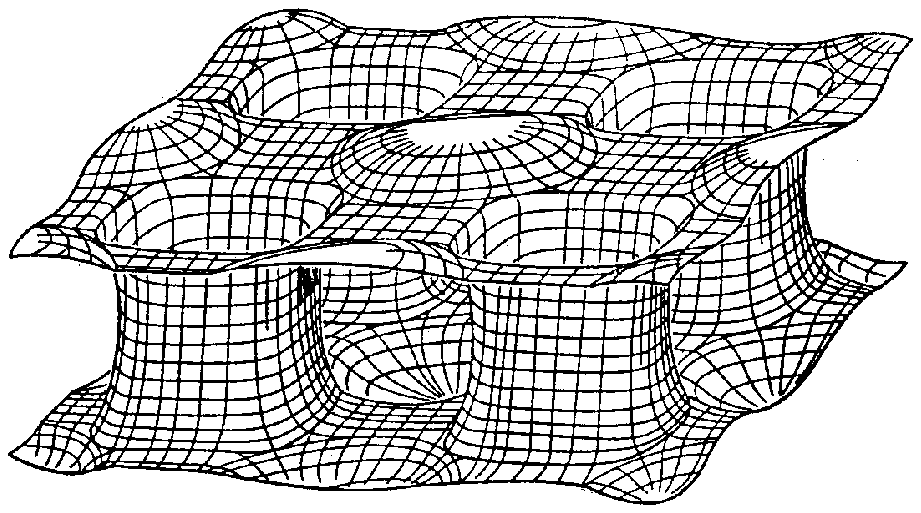}
\includegraphics{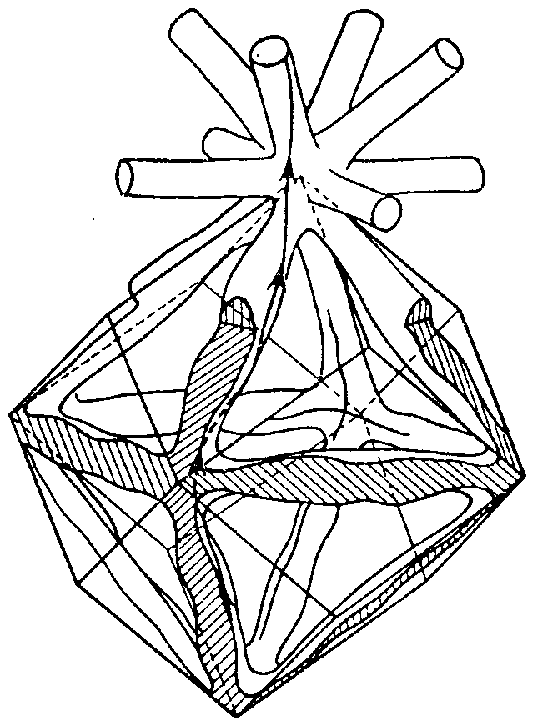}
\includegraphics{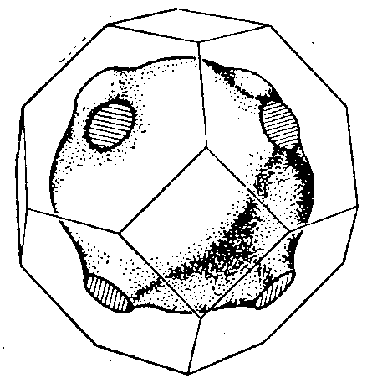}
\vspace{105bp} 
\caption{
\label{zorich:fig:fermi:surfaces}
Fermi  surfaces  of  tin,  iron  and  gold correspond to  Riemann
surfaces  of  high  genera.  (Reproduced  from~\cite{zorich:LAK}   which
cites~\cite{zorich:41} and \cite{zorich:75} as the source)
}
\end{figure}

\begin{figure}[htb]
\centering
\includegraphics{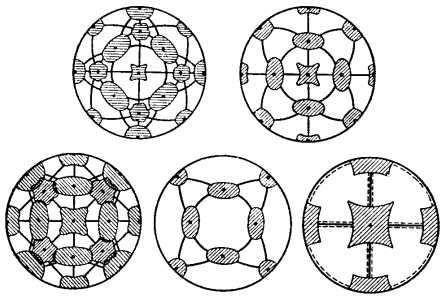}
\vspace{105bp}
\caption{
\label{zorich:fig:opentraj:iron}
Stereographic projection of the magnetic field directions (shaded
regions   and   continuous  curves)  which  give  rise  to   open
trajectories  for   some  Fermi-surfaces  (experimental   results
in~\cite{zorich:LAK}).
}
\end{figure}

It was known since extensive experimental research in the 50s and 60s
that  Fermi-surfaces\index{Fermi-surface}\index{Surface!Fermi-surface}
may  have  fairly  complicated  shape,  see
Fig.~\ref{zorich:fig:fermi:surfaces}; it  was also known that open (i.e.
unbounded) trajectories  exist, see Fig.~\ref{zorich:fig:opentraj:iron},
however, up to the beginning of the 80s there were no general  results in
this area.

In particular, it  was not known whether open trajectories follow
(in a large  scale) the same  direction, whether there  might  be
some scattering (trajectory comes from infinity  in one direction
and then after  some scattering goes  to infinity in  some  other
direction, whether the trajectories may even exhibit some chaotic
behavior?

Let us see now how this problem is related to flat surfaces.

First note that  passing to a quotient $\R{3}/\Z{3}=\T{3}$ we get
a closed  orientable  surface  $M^2\subset\T{3}$ from the initial
periodic surface  $\tilde  M^2$.  Say,  identifying  the opposite
sides of a  unit cube at Fig.~\ref{zorich:fig:periodic:surface} we get a
closed surface $M^2$ of genus $g=3$.

We  are  interested  in  plane sections of the  initial  periodic
surface $\tilde M^2$. This plane sections can be  viewed as level
curves of  a  linear  function  $f(x,y,z)=ax+by+cz$ restricted to
$\tilde M^2$.

Consider      now     a      closed      differential      1-form
$\tilde\omega=a\,dx+b\,dy+c\,dz$ in  $\R{3}$ and its  restriction
to  $\tilde  M^2$. A  closed  1-form  defines  a  codimension-one
foliation  on  a  manifold:  locally one can represent  a  closed
one-form as a differential of a  function, $\tilde\omega=df$; the
foliation is defined by the levels of the function $f$. We prefer
to use  the 1-form $\tilde\omega=a\,dx+b\,dy+c\,dz$ to the linear
function $f(x,y,z)=ax+by+cz$  because we cannot push the function
$f(x,y,z)$   into    a    torus    $\T{3}$   while   the   1-form
$\omega=a\,dx+b\,dy+c\,dz$ is well-defined  in $\T{3}$. Moreover,
after    passing    to    a    quotient    over    the    lattice
$\R{3}\to\R{3}/\Z{3}$ the  plane sections of $\tilde M^2$ project
to the leaves  of  the foliation  defined  by restriction of  the
closed 1-form $\omega$ in $\T{3}$ to the surface $M^2$.

Thus, our initial problem can be reformulated as follows.

\begin{Problem}[Novikov's Problem on Electron Transport]
   %
Consider   a   foliation\index{Foliation!measured foliation}
defined  by  a  linear  closed   1-form
$\omega=a\,dx+b\,dy+c\,dz$ on a  closed surface $M^2\subset\T{3}$
embedded into a three-dimensional torus. How  do  the  leaves  of
this foliation get unfolded, when we unfold the  torus $\T{3}$ to
its universal cover $\R{3}$?
\end{Problem}

The  foliation  defined  by a closed  1-form  on  a  surface is a
subject of discussion  of  the next  section.  We shall see  that
under   some   natural   conditions   such a foliation   can   be
``straightened up'' to a geodesic  foliation  in  an  appropriate
flat metric.

The way  in which  a geodesic on a flat  surface gets unfolded in
the  universal  Abelian  cover  is   discussed   in   detail   in
Sec.~\ref{zorich:s:How:Do:Generic:Geodesics:Wind:Around:Flat:Surfaces}.

To   be   honest,   we   should   admit   that   1-forms   as  in the
Problem above  usually  do not  satisfy these
conditions.       However,      a       surface       as       in this
Problem can be  decomposed into several
components, which  (after some surgery) already satisfy the necessary
requirements.

\paragraph{References for Details}

There  is a lot  of  progress  in  this area,  basically  due  to
S.~Novikov's  school,  and especially to I.~Dynnikov, and in  many
cases Novikov's Problem on electron transport is solved.

In~\cite{zorich:Zorich:84}  the   author   proved   that   for  a  given
Fermi-surface\index{Fermi-surface}\index{Surface!Fermi-surface}
and  for an open  dense set of directions of planes
any open trajectory is bounded by a pair of parallel lines inside
the corresponding plane.

In a series  of papers I.~Dynnikov applied a different approach: he
fixed the direction of the plane  and  deformed  a
Fermi-surface\index{Fermi-surface}\index{Surface!Fermi-surface}
inside  a family  of  level surfaces of  a  periodic function  in
$\R{3}$. He proved  that  for  all but at most one level any  open
trajectory is also bounded between two lines.

However,  I.~Dynnikov   has   constructed   a  series  of  highly
elaborated examples showing that in some cases an open trajectory
can ``fill'' the  plane. In particular, the following question is
still open. Consider the  set  of directions of those hyperplanes
which give ``nontypical'' open trajectories. Is it true that this
set  has  measure  zero  in  the  space  $\R{}\text{P}^2$  of all
possible directions? What can be said about Hausdorff dimension of
this set?

For    more     details     we     address    the    reader    to
papers~\cite{zorich:Dynnikov:98},                     \cite{zorich:Dynnikov:99}
and~\cite{zorich:Novikov:Maltsev}.

\subsection{Flows on Surfaces and Surface Foliations}
\label{zorich:ss:Flows:on:Surfaces:and:Surface:Foliations}

Consider a closed  1-form on a closed orientable surface. Locally
a closed 1-form $\omega$ can be represented as the differential of
a  function  $\omega=df$. The level curves of  the  function  $f$
locally define the leaves of the  closed  1-form  $\omega$.  (The
fact that the function $f$ is  defined only up to a constant does
not affect the structure of the level curves.) We get a
foliation\index{Foliation!measured foliation}
on a surface.

In this section we  present  a necessary and sufficient condition
which tells when  one can find  an appropriate flat  metric  such
that  the  foliation  defined  by  the  closed  1-form  becomes a
foliation of parallel  geodesics going in some fixed direction on
the surface. This criterion was given  in  different  context  by
different     authors:     \cite{zorich:Calabi},      \cite{zorich:Katok:1973},
\cite{zorich:Hubbard:Masur}. We present here one more formulation of the
criterion. Morally, it says that \emph{the foliation defined by a
closed  1-form  $\omega$   can   be  ``straightened  up''  in  an
appropriate flat metric if and only if the form $\omega$ does not
have closed  leaves homologous to zero}.  In the remaining  part of
this section we present a rigorous formulation of this statement.

Note  that  a  closed   1-form   $\omega$  on  a  closed  surface
necessarily  has  some critical  points:  the  points  where  the
function $f$  serving  as  a ``local antiderivative'' $\omega=df$
has critical points.

\begin{figure}[htb]
\centering
\includegraphics{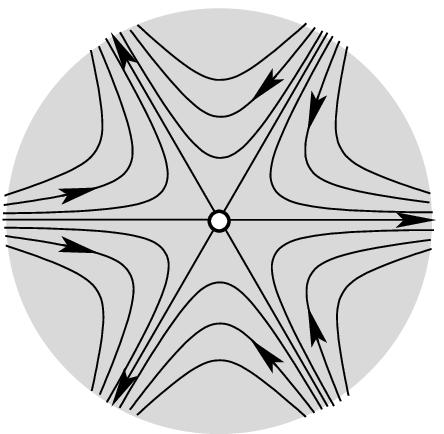}
\vspace{95bp}
\caption{
\label{zorich:fig:saddle:point}
Horizontal foliation\index{Foliation!horizontal}
in a neighborhood of a saddle point\index{Conical!singularity}.
Topological (nonmetric) picture}
\end{figure}

The  first  obstruction for ``straightening'' is the presence  of
minima and  maxima:  such  critical  points  should be forbidden.
Suppose now that the closed 1-form  has  only  isolated  critical
points and  all of them are  ``saddles'' (i.e. $\omega$  does not
have minima and maxima). Say, a form defined in local coordinates
as  $df$, where  $f=x^3+y^3$  has a saddle  point\index{Conical!singularity}
in the  origin
$(0,0)$ of our coordinate chart, see Fig.~\ref{zorich:fig:saddle:point}.

There are several singular leaves of  the  foliation  landing  at
each saddle; say, a saddle point from Fig.~\ref{zorich:fig:saddle:point}
has six  prongs (separatrices\index{Separatrix})
representing  the  critical  leaves. Sometimes a
critical leaf emitted  from one saddle  can land to  another  (or
even to  the same) saddle  point. In  this case we  say that  the
foliation has a
\emph{saddle connection}\index{Saddle!connection}.

Note that the foliation defined by a closed 1-form on an oriented
surface     gets     natural      orientation     (defined     by
$\operatorname{grad}(f)$ and  by the orientation of the surface).
Now  we  are ready  to  present  a  rigorous  formulation  of the
criterion.

\begin{NNTheorem}
\label{zorich:th:Calabi}
Consider a foliation\index{Foliation!measured foliation}
defined on a closed orientable  surface by a
closed  1-form  $\omega$. Assume  that  $\omega$  does  not  have
neither minima nor maxima but  only  isolated  saddle points. The
foliation defined by $\omega$ can  be  represented  as a geodesic
foliation in an appropriate flat metric if and only if  any cycle
obtained as a  union of closed  paths following in  the  positive
direction a sequence of  saddle  connections\index{Saddle!connection}
is not homologous to
zero.
\end{NNTheorem}

In  particular,  if  there  are  no  saddle  connections  at  all
(provided  there are no  minima   and   maxima)  it  can  always  be
straightened  up.  (In  slightly  different terms it  was  proved
in~\cite{zorich:Katok:1973} and in~\cite{zorich:Hubbard:Masur}.)

Note that saddle points of the closed 1-form $\omega$ correspond
to conical points\index{Conical!singularity}
of the resulting flat metric\index{Surface!very flat}.

One  can  consider  a  closed  1-form  $\omega$ as a  multivalued
Hamiltonian and consider corresponding Hamiltonian flow along the
leaves of the foliation defined by $\omega$. On the torus $\T{2}$
is was  studied by V.~I.~Arnold~\cite{zorich:Arnold:91} and by K.~Khanin
and Ya.~G.~Sinai~\cite{zorich:Khanin:Sinai}.

\section{Families of Flat Surfaces and Moduli Spaces of Abelian Differentials}
\label{zorich:s:Families:Of:Flat:Surfaces:and:Moduli:Spaces:of:Abelian:Differentials}

In this section we present the generalities on
flat surfaces\index{Surface!flat}.
We start   in    Sec.~\ref{zorich:ss:Families:of:Flat:Surfaces}   with   an
elementary  construction  of a  flat  surface  from  a  polygonal
pattern. This construction explicitly shows that any flat surface
can  be  deformed  inside  an appropriate \emph{family}  of  flat
surfaces.    As     a    model    example    we    consider    in
Sec.~\ref{zorich:ss:Toy:Example:Family:of:Flat:Tori} the  family of flat
tori. In Sec.~\ref{zorich:ss:Dictionary:of:Complex:Analytic:Language} we
show  that  a  flat  structure  naturally  determines  a  complex
structure  on   the   surface   and   a   holomorphic   one-form.
Reciprocally, a  holomorphic one-form naturally determines a flat
structure.  The dictionary  establishing  correspondence  between
geometric  language  (in  terms  of   the   flat   metrics)   and
complex-analytic language  (in terms of holomorphic one-forms) is
very  important   for   the   entire   presentation;   it   makes
Sec.~\ref{zorich:ss:Dictionary:of:Complex:Analytic:Language}        more
charged        than       an        average        one.        In
Sec.~\ref{zorich:ss:Moduli:Space:of:Holomorphic:One:Forms}  we  continue
with  establishing  correspondence  between   families   of  flat
surfaces and strata of moduli spaces of holomorphic one-forms. In
Sec.~\ref{zorich:ss:Action:of:SL2R:on:the:Moduli:Space} we describe  the
action  of  the  linear  group $SL(2,\R{})$ on flat  surfaces  --
another key issue of this theory.

We complete
Sec.~\ref{zorich:s:Families:Of:Flat:Surfaces:and:Moduli:Spaces:of:Abelian:Differentials}
with an attempt to present the following general principle in the
study of flat surfaces. In order to get some information about an
individual flat surface it is often very convenient  to find (the
closure  of)  the  orbit  of  the corresponding element  in   the
\emph{family}  of  flat surfaces under the action  of  the  group
$SL(2,\R{})$ (or, sometimes, under the  action  of  its  diagonal
subgroup). In many  cases  the structure  of  this orbit gives  a
comprehensive  information  about  the   initial   flat  surface;
moreover, this information might be  not  accessible  by a direct
approach.      These       ideas      are      expressed       in
Sec.~\ref{zorich:ss:General:Philosophy}.  This  general   principle   is
illustrated in Sec.~\ref{zorich:ss:Implementation:of:General:Philosophy}
presenting  Masur's   criterion   of  unique  ergodicity  of  the
directional flow  on a flat  surface. (A reader not familiar with
ergodic theorem can  either  skip this  last  section or read  an
elementary       presentation        of       ergodicity       in
Appendix~\ref{zorich:s:Ergodic:Theorem}.)

\subsection{Families of Flat Surfaces}
\label{zorich:ss:Families:of:Flat:Surfaces}

In  this section  we  present a construction  which
allows  to  obtain  a  large  variety  of flat surfaces,  and,
moreover,  allows  to  continuously  deform  the  resulting  flat
structure. Later  on we shall  see that this construction is even
more general than it may seem  at the beginning: it allows to get
\emph{almost all}  flat  surfaces  in  any  \emph{family} of flat
surfaces sharing the  same geometry (i.e. genus, number and types
of conical singularities). The construction is strongly motivated
by    an    analogous    construction    in    the    paper    of
H.~Masur~\cite{zorich:Masur:Annals:82}.

\begin{figure}[hbt]
%
%
\includegraphics{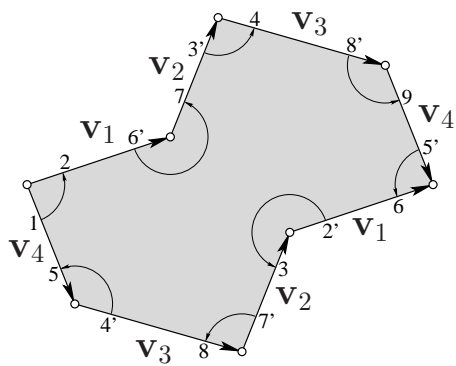}
\begin{picture}(35,2)(35,2)
\put(142,-48){\Large $\vec v_1$}
\put(168,-24){\Large $\vec v_2$}
\put(222,-7){\Large $\vec v_3$}
\put(270,-42){\Large $\vec v_4$}
\put(115,-93){\Large $\vec v_4$}
\put(163,-131){\Large $\vec v_3$}
\put(216,-111){\Large $\vec v_2$}
\put(245,-85){\Large $\vec v_1$}
\end{picture}
\vspace{140bp}
\caption{
\label{zorich:fig:suspension}
Identifying corresponding pairs of sides of this polygon by
parallel translations we obtain a flat surface\index{Surface!flat}.
}
\end{figure}

Consider a collection of vectors $ \vec v_1,  \dots, \vec v_\noi$
in $\R{2}$ and construct from  these  vectors a broken line in  a
natural way:  a $j$-th edge of the broken  line is represented by
the vector $\vec{v}_j$. Construct another broken line starting at
the same point as the initial  one by taking the same vectors but
this time in the order $v_{\pi(1)},  \dots, v_{\pi(\noi)}$, where
$\pi$ is some permutation of $\noi$ elements.

By construction the two broken  lines  share  the same endpoints;
suppose    that     they     bound     a     polygon     as    in
Fig.~\ref{zorich:fig:suspension}.   Identifying   the  pairs   of  sides
corresponding to the same  vectors  $v_j$, $j=1, \dots, \noi$, by
parallel translations we obtain a flat surface.

The  polygon  in  our  construction depends continuously  on  the
vectors  $\vec{v}_i$.  This means that the topology of the  resulting
flat  surface (its  genus  $g$, the number $\noz$  and the types of  the
resulting  conical  singularities)  do  not  change  under  small
deformations  of the vectors  $\vec{v}_i$.  Say, we suggest  to  the
reader to  check that the  flat surface obtained from the polygon
presented in Fig.~\ref{zorich:fig:suspension} has genus two and a single
conical singularity with cone angle $6\pi$.

\subsection{Toy Example: Family of Flat Tori}
\label{zorich:ss:Toy:Example:Family:of:Flat:Tori}

In the previous section we have seen that a flat structure can be
deformed. This allows to consider a flat surface as an element of
a \emph{family}  of  flat  surfaces  sharing  a common geometry
(genus, number of  conical points). In  this section we  study the
simplest example of  such  family: we  study  the family of  flat
tori.  This  time   we  consider  the  family  of  flat  surfaces
\emph{globally}. We shall see that it has a  natural structure of
a noncompact  complex-analytic  manifold  (to  be  more honest --
\emph{orbifold}).  This   ``baby   family''   of  flat  surfaces,
actually, exhibits all  principal features of any other family of
flat surfaces, except that the family of flat  tori constitutes a
homogeneous  space  endowed with a nice hyperbolic metric,  while
general  families  of  flat  surfaces  \emph{do   not}  have  the
structure of a homogeneous space.

\begin{figure}[htb]
%
\centering
\includegraphics{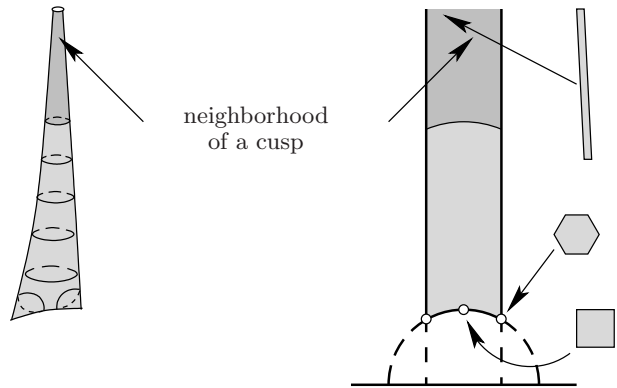}
\begin{picture}(0,0)(0,0)
\put(-50,-40){neighborhood}
\put(-50,-50){\ \ \ of a cusp}
\end{picture}
\vspace{140bp}
\caption{
\label{zorich:fig:space:of:flat:tori}
Space of flat tori
}
\end{figure}

To simplify consideration  of flat tori  as much as  possible  we
make two exceptions from the usual way in which we  consider flat
surfaces. Temporarily (only  in this section) we forget about the
choice of  the direction to  North: in this section two isometric
flat tori define the same element of the family of all flat tori.
Another exception concerns normalization. Almost everywhere below
we consider the area of any  flat surface to be normalized to one
(which can be achieved by a simple homothety). In this section it
would be more  convenient  for us to apply  homothety  in the way
that the  shortest closed geodesic on  our flat torus  would have
length  1.  Find the  closed  geodesic which  is  next after  the
shortest one in the length  spectrum.  Measure  the angle $\phi$,
where $0\le \phi \le \pi$  between  these  two geodesics; measure
the length $r$ of the second geodesic and mark a point with polar
coordinates  $(r,\phi)$  on  the  upper  half-plane.  This  point
encodes our torus.

Reciprocally, any point of the  upper  half-plane  defines a flat
torus in the  following way. (Following the tradition we consider
the upper half-plane as a complex one.) A point $x+iy$  defines a
parallelogram generated  by  vectors  $\vec  v_1=(1,0)$ and $\vec
v_2=(x+iy)$,  see  Fig.~\ref{zorich:fig:space:of:flat:tori}. Identifying
the opposite sides of the parallelogram we get a flat torus.

To make this correspondence bijective we have to be sure that the
vector $\vec v_1$  represents  the shortest closed geodesic. This
means that the point  $(x+iy)$  representing $\vec v_2$ cannot be
inside the unit disc. The  condition  that  $\vec v_2$ represents
the geodesic  which is next after the shortest  one in the length
spectrum   implies    that    $-1/2\le    x    \le    1/2$,   see
Fig.~\ref{zorich:fig:space:of:flat:tori}.  Having  mentioned  these  two
hints we suggest to the reader to prove the following Lemma.

\begin{lemma}
\label{zorich:lm:space:of:flat:tori}
The  family  of   flat  tori  is  parametrized  by  the  shadowed
fundamental domain  from Fig.~\ref{zorich:fig:space:of:flat:tori}, where
the parts of the boundary of  the  fundamental  domain  symmetric
with respect to the vertical axis $(0,iy)$ are identified.
\end{lemma}

Note  that  topologically  we  obtain a sphere punctured  at  one
point: the resulting surface has a
\emph{cusp}\index{Cusp of the moduli space@Cusp of the moduli space, \emph{see also} Moduli space; principal boundary}.
Tori represented
by points close to the cusp  are  ``disproportional'':  they  are
very narrow and very long. In other words they have an abnormally
short geodesic.

Note also that there are two special points on our modular curve:
they correspond to points with coordinates $(0+i)$ and $\pm 1/2 +
i\sqrt{3}/2$. The corresponding tori has extra symmetry, they can
be  represented  by  a  square  and  by a  regular  hexagon  with
identified opposite sides correspondingly. The surface glued from
the fundamental domains has ``corners'' at these two points.

There is an alternative  more  algebraic approach to our problem.
Actually, the fundamental  domain constructed above is known as a
\emph{modular  curve},  it  parameterizes   the   \emph{space  of
lattices of  area one} (which is isomorphic to  the space of flat
tori). It can be seen as a double quotient
$$
\begin{array}{rcl}
&\backslash\, SL(2,\R{}) /&\\
[-\halfbls]
SO(2,\R{})&& SL(2,\Z{})
\end{array}
\quad
=
\quad
\begin{array}{rl}
{\mathbb H} /&\\
[-\halfbls]  & SL(2,\Z{})
\end{array}
$$

\subsection{Dictionary of Complex-Analytic Language}
\label{zorich:ss:Dictionary:of:Complex:Analytic:Language}

We have  seen  in  Sec.~\ref{zorich:ss:Families:of:Flat:Surfaces} how to
construct    a     flat    surfaces    from    a    polygon    as
on~Fig.~\ref{zorich:fig:suspension}.

Note that  the polygon is embedded  into a complex  plane $\C{}$,
where the  embedding is defined  up to a parallel translation. (A
rotation of the polygon changes the
vertical direction\index{Direction!vertical}
and hence,
according  to  Convention~\ref{zorich:conv:flat:surface} it  changes the
corresponding flat surface.)

Consider the natural coordinate $z$ in the complex plane. In this
coordinate the parallel translations which we use to identify the
sides of the polygon are represented as
\begin{equation}   
\label{zorich:eq:z:equals:z:plus:const}
z'=z+const
\end{equation}
Since this correspondence is holomorphic, it means that our flat
surface $S$ with punctured conical points inherits the complex
structure. It is an exercise in complex analysis to check that
the complex structure extends to the punctured points.

Consider now a holomorphic 1-form $dz$  in  the  initial  complex
plane.  When  we pass to the  surface  $S$ the coordinate $z$ is  not
globally defined anymore. However, since  the  changes  of  local
coordinates          are          defined         by          the
rule~\eqref{zorich:eq:z:equals:z:plus:const} we see that $dz=dz'$. Thus,
the  holomorphic  1-form  $dz$  on $\C{}$ defines  a  holomorphic
1-form $\omega$  on $S$ which in  local coordinates has  the form
$\omega=dz$. Another exercise in complex analysis  shows that the
form $\omega$ has
\emph{zeroes}\index{Zero}
exactly  at those  points  of $S$ where  the  flat structure  has
conical singularities.

In an appropriate local coordinate $w$ in a  neighborhood of zero
(different from the  initial  local coordinate $z$) a holomorphic
1-form can  be represented as  $w^d\,dw$, where $d$ is called the
\emph{degree} of zero\index{Degree of zero}\index{Zero!degree of zero}.
The form $\omega$ has a zero of  degree $d$ at a
\emph{conical point}\index{Conical!singularity}\index{Cone angle}
with cone angle $2\pi(d+1)$.

Recall the formula for  the sum of degrees of zeroes of a
holomorphic 1-form on a Riemann surface of genus $g$:
$$
\sum_{j=1}^\noz d_j = 2g-2
$$
This relation can be interpreted as the formula of
Gauss--Bonnet\index{Gauss--Bonnet formula}\index{Formula!of Gauss--Bonnet}
for the flat metric.

Vectors $\vec  v_j$ representing the  sides of the polygon can be
considered as complex numbers. Let $\vec v_j$ be joining vertices
$P_j$  and  $P_{j+1}$  of  the  polygon.  Denote by $\rho_j$  the
resulting path on  $S$ joining the points $P_j,P_{j+1}\in S$. Our
interpretation of  $\vec v_j$ as  of a complex number implies the
following obvious relation:
\begin{equation}
\label{zorich:eq:vj:is:a:period}
\vec v_j=\int_{P_j}^{P_{j+1}}dz=\int_{\rho_j}\omega
\end{equation}
Note that the path $\rho_j$ represents a
\emph{relative cycle}\index{Relative!cycle}\index{Cycle!relative}\index{Relative!homology group}:
an element of the relative homology group
$H_1(S,\{P_1,\dots,P_\noz\};\Z{})$ of the surface $S$ relative to
the finite collection of  conical  points $\{P_1,\dots,P_\noz\}$.
Relation~\eqref{zorich:eq:vj:is:a:period}   means    that   $\vec   v_j$
represents a
\emph{period}\index{Period}
of $\omega$: an integral of $\omega$ over a relative cycle $\rho_j$.

Note also that the  flat area of the surface $S$ equals  the area
of the original polygon, which can be measured as an  integral of
$dx\wedge dy$ over the  polygon.  Since in the complex coordinate
$z$ we have $dx\wedge dy=\frac{i}{2}dz\wedge d\bar{z}$ we get the
following formula for the flat area of $S$:
\begin{equation}
\label{zorich:eq:Riemann:bilinear:relation}
area(S)=\frac{i}{2}\int_S \omega\wedge\bar{\omega}
=\frac{i}{2}\sum_{j=1}^g (A_j\bar{B}_j-\bar{A}_j B_j)
\end{equation}
Here we also used  the  Riemann bilinear relation which expresses
the  integral  $\int_S  \omega\wedge\bar{\omega}$  in  terms   of
\emph{absolute  periods}\index{Period!absolute}
$A_j,  B_j$  of  $\omega$,  where  the
absolute periods $A_j, B_j$ are  the  integrals  of $\omega$ with
respect to some symplectic basis of cycles.

An  individual  flat surface defines a pair: (complex  structure,
holomorphic 1-form). A family of  flat  surfaces  (where the flat
surfaces  are as  usual  endowed with a  choice  of the  vertical
direction)       corresponds       to      a
\emph{stratum}\index{Stratum!in the moduli space}
$\cH(d_1,\dots,d_\noz)$\index{0H20@$\cH(d_1,\dots,d_\noz)$ -- stratum in the moduli space}\index{Stratum!in the moduli space}
in the
\emph{moduli space of  holomorphic 1-forms\/}\index{Moduli space!of complex structures}\index{Moduli space!of holomorphic 1-forms}.
Points of the  stratum  are represented  by  pairs (point in  the
moduli space of complex structures,  holomorphic  1-form  in  the
corresponding  complex   structure   having   zeroes  of  degrees
$d_1,\dots,d_\noz$).

The notion ``stratum'' has the following origin. The moduli space
of pairs (holomorphic 1-form, complex structure)  forms a natural
vector   bundle   over  the  moduli  space  $\cM_g$\index{0M10@$\cM_g$ -- moduli space of complex structures}\index{Moduli space!of complex structures}  of   complex
structures.  A fiber  of  this vector bundle  is  a vector  space
$\C{g}$ of holomorphic 1-forms  in  a given complex structure. We
already  mentioned  that  the  sum  of  degrees of  zeroes  of  a
holomorphic  1-form  on  a  Riemann surface of genus  $g$  equals
$2g-2$. Thus, the  total space
$\cH_g$\index{0H10@$\cH_g$ -- moduli space of holomorphic 1-forms}
of our vector  bundle  is
stratified  by  subspaces  of  those forms which have  zeroes  of
degrees         exactly         $d_1,\dots,d_\noz$,         where
$d_1+\dots+d_\noz=2g-2$.  Say,   for   $g=2$  we  have  only  two
partitions  of number  $2$,  so we get  two  strata $\cH(2)$  and
$\cH(1,1)$.   For    $g=3$   we   have   five   partitions,   and
correspondingly   five  strata   $\cH(4),   \cH(3,1),   \cH(2,2),
\cH(2,1,1), \cH(1,1,1,1)$.
\begin{equation}
\label{zorich:eq:projection:H:to:M}
\begin{array}{ccc}
\cH(d_1,\dots,d_\noz)&\subset&\cH_g\\
&& \downarrow\\
&&\, \cM_g
\end{array}
\end{equation}
Every stratum $\cH(d_1,\dots,d_\noz)$\index{0H20@$\cH(d_1,\dots,d_\noz)$ -- stratum in the moduli space}\index{Stratum!in the moduli space}
is a complex-analytic
orbifold of dimension\index{Dimension of a stratum}
\begin{equation}
\label{zorich:eq:dim:of:stratum}
\dim_{\C{}} \cH(d_1,\dots,d_\noz)=2g+\noz-1
\end{equation}
Note, that an individual stratum
$\cH(d_1,\dots,d_\noz)$\index{0H20@$\cH(d_1,\dots,d_\noz)$ -- stratum in the moduli space}\index{Stratum!in the moduli space}
does not
form a  fiber bundle over  $\cM_g$\index{0M10@$\cM_g$ -- moduli space of complex structures}\index{Moduli space!of complex structures}. For example, according to our
formula,     $\dim_{\C{}}\cH(2g-2)=2g$,     while    $\dim_{\C{}}
\cM_g=3g-3$.

We showed how the geometric structures related to  a flat surface
define  their   complex-analytic  counterparts.  Actually,   this
correspondence goes in two  directions.  We suggest to the reader
to make the  inverse  translation: to start with complex-analytic
structure  and to  see  how it defines  the  geometric one.  This
correspondence can be summarized in the following dictionary.

\begin{table}[htb]
\centering
\caption{Correspondence of geometric and complex-analytic notions}
\label{zorich:tab:dictionary:geometry:to:complex:analysis}
\index{Conical!singularity}
\index{Zero}
\index{Period!relative}
\index{Holomorphic!1-form}
%
%
\begin{tabular}{|c|c|}
\hline&\\[-\halfbls]
       Geometric language              &      Complex-analytic language\\
[-\halfbls] &\\ \hline & \\ [-\halfbls]
flat structure (including a choice     &      complex structure + a choice\\
of the vertical direction)             &   of a holomorphic 1-form $\omega$\\
[-\halfbls] &\\ \hline & \\ [-\halfbls]
conical point                          & zero of degree $d$\\
with a cone angle $2\pi(d+1)$          & of the holomorphic 1-form $\omega$\\
                                       &(in local coordinates $\omega=w^d\,dw$)\\
[-\halfbls] &\\ \hline & \\ [-\halfbls]
side $\vec v_j$ of a polygon           & relative period $\int_{P_j}^{P_{j+1}}\omega=\int_{\vec v_j}\omega$ \\
                                       & of the 1-form $\omega$ \\
[-\halfbls] &\\ \hline & \\ [-\halfbls]
area of the flat surface $S$           & $\frac{i}{2}\int_S \omega\wedge\bar{\omega}=$ \\
                                       & $=\frac{i}{2}\sum_{j=1}^g (A_j\bar{B}_j-\bar{A}_j B_j)$ \\
[-\halfbls] &\\ \hline & \\ [-\halfbls]
$\ $ family of flat surfaces sharing the same $\,$& stratum $\cH(d_1,\dots,d_\noz)$ in the \\
types $2\pi(d_1+1),\dots,2\pi(d_\noz+1)$& $\ $ moduli space of Abelian differentials $\ $\\
of cone angles                         & \\
[-\halfbls] &\\ \hline & \\ [-\halfbls]
coordinates in the family:             & coordinates in $\cH(d_1,\dots,d_\noz)$ : \\
vectors $\vec v_i$                     & collection of relative periods of $\omega$,  \\
defining the polygon                   & i.e. cohomology class \\
                                       & $[\omega]\in H^1(S, \{P_1,\dots,P_\noz\};\C{})$ \\
[-\halfbls] &\\ \hline
\end{tabular}
\end{table}
%

\subsection{Volume Element in the Moduli Space of Holomorphic One-Forms}
\label{zorich:ss:Moduli:Space:of:Holomorphic:One:Forms}

In  the  previous  section  we  have   considered  vectors  $\vec
v_1,\dots,\vec v_\noi$ determining the polygon from which we glue
a  flat  surface  $S$  as on Fig.~\ref{zorich:fig:suspension}.  We  have
identified these vectors $\vec v_j\in\R{2}\sim\C{}$ with  complex
numbers   and   claimed   (without   proof)   that   under   this
identification $\vec v_1,\dots,\vec v_\noi$ provide us with local
coordinates in  the  corresponding  family  of  flat surfaces. We
identify every such family with a stratum
$\cH(d_1,\dots,d_\noz)$\index{0H20@$\cH(d_1,\dots,d_\noz)$ -- stratum in the moduli space}\index{Stratum!in the moduli space}
in the moduli  space  of holomorphic 1-forms. In complex-analytic
language we have locally identified a neighborhood of a ``point''
(complex   structure,   holomorphic   1-form  $\omega$)  in   the
corresponding stratum with a neighborhood of the cohomology class
$[\omega]\in H^1(S,\{P_1, \dots, P_\noz\};\C{})$.

Note   that   the   cohomology    space    $H^1(S,\{P_1,   \dots,
P_\noz\};\C{})$ contains a natural integer lattice  $H^1(S,\{P_1,
\dots, P_\noz\};\Z{}\oplus  \sqrt{-1}\,\Z{})$. Consider a  linear
volume  element\index{Volume element!in the moduli space}
$d\nu$\index{0d10@$d\nu$ -- volume element in the moduli space}
in the vector space $H^1(S,\{P_1,  \dots,
P_\noz\};\C{})$ normalized in such a way that  the  volume  of  the
fundamental domain in the ``cubic'' lattice
$$
H^1(S,\{P_1, \dots, P_\noz\};\Z{}\oplus \sqrt{-1}\,\Z{})\subset
H^1(S,\{P_1, \dots, P_\noz\};\C{})
$$
is equal to one. In other terms
$$
d\nu=\frac{1}{J}\ \frac{1}{(2\sqrt{-1})^\noi}\ d{\vec
v_1}d\bar{\vec v}_1 \dots d{\vec v_\noi}d\bar{\vec v}_\noi,
$$
where $J$ is the determinant of a change of the basis $\{\vec v_1,
\dots, \vec v_\noi\}$ considered as a basis in the first
relative \emph{homology} to some ``symplectic'' basis in the
first relative homology.

Consider now the real hypersurface\index{0H30@$\cH_1(d_1,\dots,d_\noz)$ -- ``unit hyperboloid''}\index{Stratum!in the moduli space}\index{Unit hyperboloid}
$$
\cH_1(d_1,\dots,d_\noz)\subset\cH(d_1,\dots,d_\noz)
$$
defined  by the equation  $area(S)=1$.  Taking  into   consideration
formula~\eqref{zorich:eq:Riemann:bilinear:relation}  for  the   function
$area(S)$ we see that the hypersurface
$\cH_1(d_1,\dots,d_\noz)$\index{0H30@$\cH_1(d_1,\dots,d_\noz)$ -- ``unit hyperboloid''}\index{Stratum!in the moduli space}\index{Unit hyperboloid}
defined  as  $area(S)=1$  can  be   interpreted   as   a   ``unit
hyperboloid''  defined  in local coordinates as a  level  of  the
indefinite quadratic form~\eqref{zorich:eq:Riemann:bilinear:relation} .

The  volume  element  $d\nu$  can  be  naturally restricted to  a
hyperplane  defined  as a level hypersurface of  a  function.  We
denote      the      corresponding
volume      element\index{Volume element!on the ``unit hyperboloid''}
on
$\cH_1(d_1,\dots,d_\noz)$\index{0H30@$\cH_1(d_1,\dots,d_\noz)$ -- ``unit hyperboloid''}\index{Stratum!in the moduli space}\index{Unit hyperboloid}
by $d\nu_1$\index{0d20@$d\nu_1$ -- volume element on the ``unit hyperboloid''}.

\begin{NNTheorem}[H.~Masur. W.~A.~Veech]
The total volume\index{Moduli space!volume of the moduli space}
$$
\int_{\cH_1(d_1,\dots,d_\noz)} d\nu_1
$$
of every stratum is finite.
\end{NNTheorem}

The  values  of  these  volumes  were  computed only recently  by
A.~Eskin  and   A.~Okounkov~\cite{zorich:Eskin:Okounkov},  twenty  years
after  the  Theorem  above was proved  in~\cite{zorich:Masur:Annals:82},
\cite{zorich:Veech:Annals:82} and~\cite{zorich:Veech:Flat:surfaces}. We discuss
this computation in Sec.~\ref{zorich:s:Volume:of:the:Moduli:Space}.

\subsection{Action of $SL(2,\R{})$ on the Moduli Space}
\label{zorich:ss:Action:of:SL2R:on:the:Moduli:Space}

\index{Action on the moduli space!ofGL@of $GL^+(2,\R{})$|(}
\index{0GL@$GL^+(2,\R{})$-action on the moduli space|(}
\index{Action on the moduli space!ofSL@of $SL(2,\R{})$|(}
\index{0SL@$SL(2,\R{})$-action on the moduli space|(}

In this  section we discuss a property of flat surfaces which
is, probably,  the most important in our study:  we show that the
linear  group  acts  on  every  family  of  flat  surfaces,  and,
moreover,          acts          \emph{ergodically}          (see
Append.~\ref{zorich:s:Ergodic:Theorem} for discussion of  the  notion of
ergodicity).  This  enables us  to  apply  tools  from  dynamical
systems and from ergodic theory.

\begin{figure}[htb]
%
\centering
\includegraphics{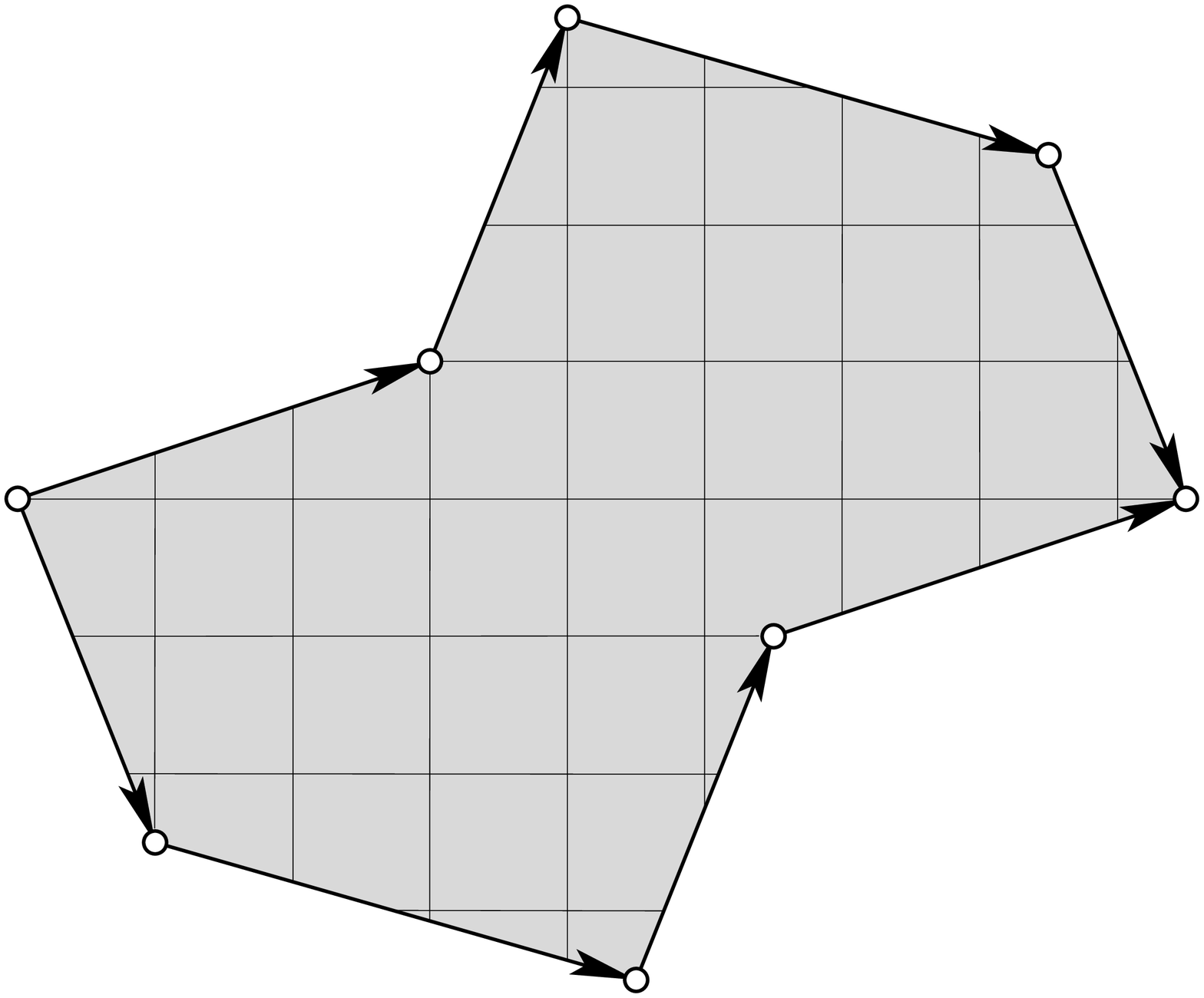}
\includegraphics{zorich_octagon_houches.eps}
\includegraphics{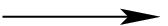}
\vspace{100bp}
\caption{
\label{zorich:fig:GL2R:action}
Action of the linear group on flat surfaces
}
\end{figure}

Consider  a  flat surface $S$ and consider  a  polygonal  pattern
obtained by unwrapping  it along some geodesic cuts. For example,
one can assume that our flat  surface $S$ is glued from a polygon
$\Pi\subset \R{2}$  as  on  Fig.~\ref{zorich:fig:suspension}. Consider a
linear transformation $g\in GL^+(2,\R{})$ of the  plane $\R{2}$. It
changes the  shape of the polygon. However, the  sides of the new
polygon $g \Pi$ are again arranged into pairs, where the sides in
each  pair  are parallel and have equal  length  (different  from
initial one),  see Fig.~\ref{zorich:fig:GL2R:action}. Thus,  identifying
the sides in each pair by  a parallel translation we obtain a new
flat surface $g S$.

It is easy to check that the surface $g S$ does not depend on the
way in  which $S$ was unwrapped to a  polygonal pattern $\Pi$. It
is clear that all topological characteristics  of  the  new  flat
surface  $g  S$  (like   genus,   number  and  types  of  conical
singularities) are the same as those of the  initial flat surface
$S$. Hence, we get a continuous action of  the group $GL^+(2,\R{})$
on each stratum
$\cH(d_1,\dots,d_\noz)$\index{0H20@$\cH(d_1,\dots,d_\noz)$ -- stratum in the moduli space}\index{Stratum!in the moduli space}.

Considering the  subgroup  $SL(2,\R{})$ of area preserving linear
transformations we get the  action  of $SL(2,\R{})$ on the ``unit
hyperboloid'' $\cH_1(d_1,\dots,d_\noz)$\index{0H30@$\cH_1(d_1,\dots,d_\noz)$ -- ``unit hyperboloid''}\index{Stratum!in the moduli space}\index{Unit hyperboloid}. Considering the diagonal
subgroup $\begin{pmatrix}  e^t  &  0\\  0  & e^{-t} \end{pmatrix}
\subset   SL(2,\R{})$   we  get  a  continuous  action  of   this
one-parameter subgroup  on each stratum  $\cH(d_1,\dots,d_\noz)$.
This action induces a  natural flow on the stratum, which
is called the
\emph{Teichm\"uller geodesic flow}\index{Teichm\"uller!geodesic flow}.
\begin{KeyTheorem}[H.~Masur. W.~A.~Veech]
The action of the groups $SL(2,\R{})$ and $\begin{pmatrix} e^t & 0\\
0 & e^{-t}
\end{pmatrix}$ preserves the measure $d\nu_1$\index{0d20@$d\nu_1$ -- volume element on the ``unit hyperboloid''}.
Both actions are
ergodic\index{Ergodic} with respect to this measure on each connected
component of every stratum
$\cH_1(d_1,\dots,d_\noz)$\index{0H30@$\cH_1(d_1,\dots,d_\noz)$ -- ``unit hyperboloid''}\index{Stratum!in the moduli space}\index{Unit hyperboloid}.
\end{KeyTheorem}

This theorem might seem quite  surprising.  Consider  almost  any
flat surface  $S$  as  in Fig.~\ref{zorich:fig:suspension}. ``Almost any
flat  surface''  is  understood  as ``corresponding to a  set  of
parameters $\vec v_1, \dots, \vec v_4$ of full  measure; here the
vectors   $\vec    v_i$    define    the   polygon   $\Pi$   from
Fig.~\ref{zorich:fig:suspension}.

Now  start  contracting  the  polygon  $\Pi$  it in the  vertical
direction and expanding it  in  the horizontal direction with the
same coefficient $e^t$. The theorem says, in particular, that for
an appropriate $t\in\R{}$ the deformed  polygon  will  produce  a
flat surface  $g_t S$ which would be arbitrary  close to the flat
surface  $S_0$  obtained  from  the   regular   octagon   as   on
Fig.~\ref{zorich:fig:octagon:to:pretzel}  since  a trajectory  of almost
any point under  an ergodic flow  is everywhere dense  (and  even
``well distributed'').  However,  it  is  absolutely  clear  that
acting     on     our     initial     polygon     $\Pi$      from
Fig.~\ref{zorich:fig:suspension} with expansion-contraction we never get
close to a regular octagon... Is there a contradiction?..

\begin{figure}[htb]
%
\centering
\includegraphics{zorich_octagon_houches.eps}
\includegraphics{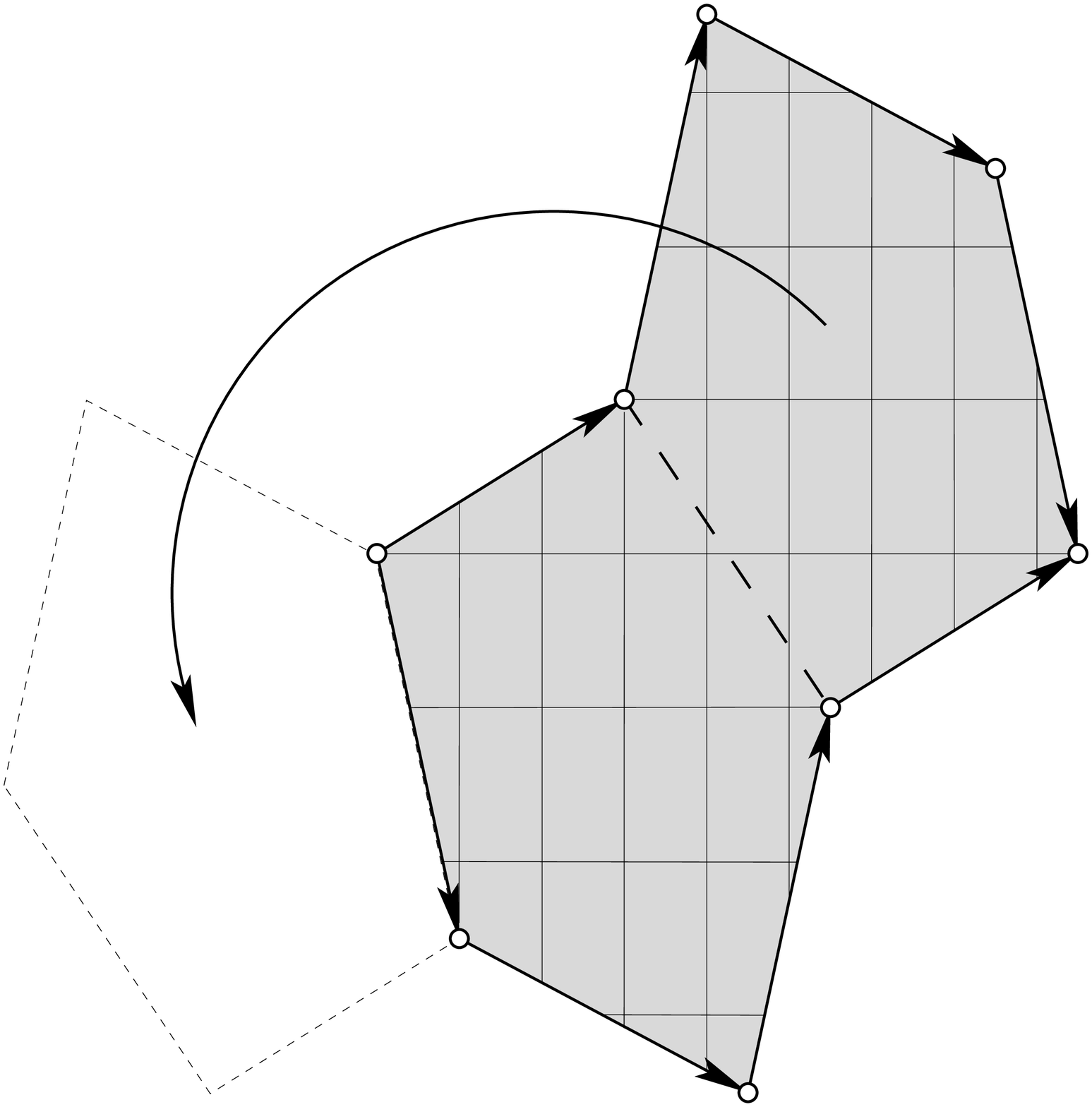}
\includegraphics{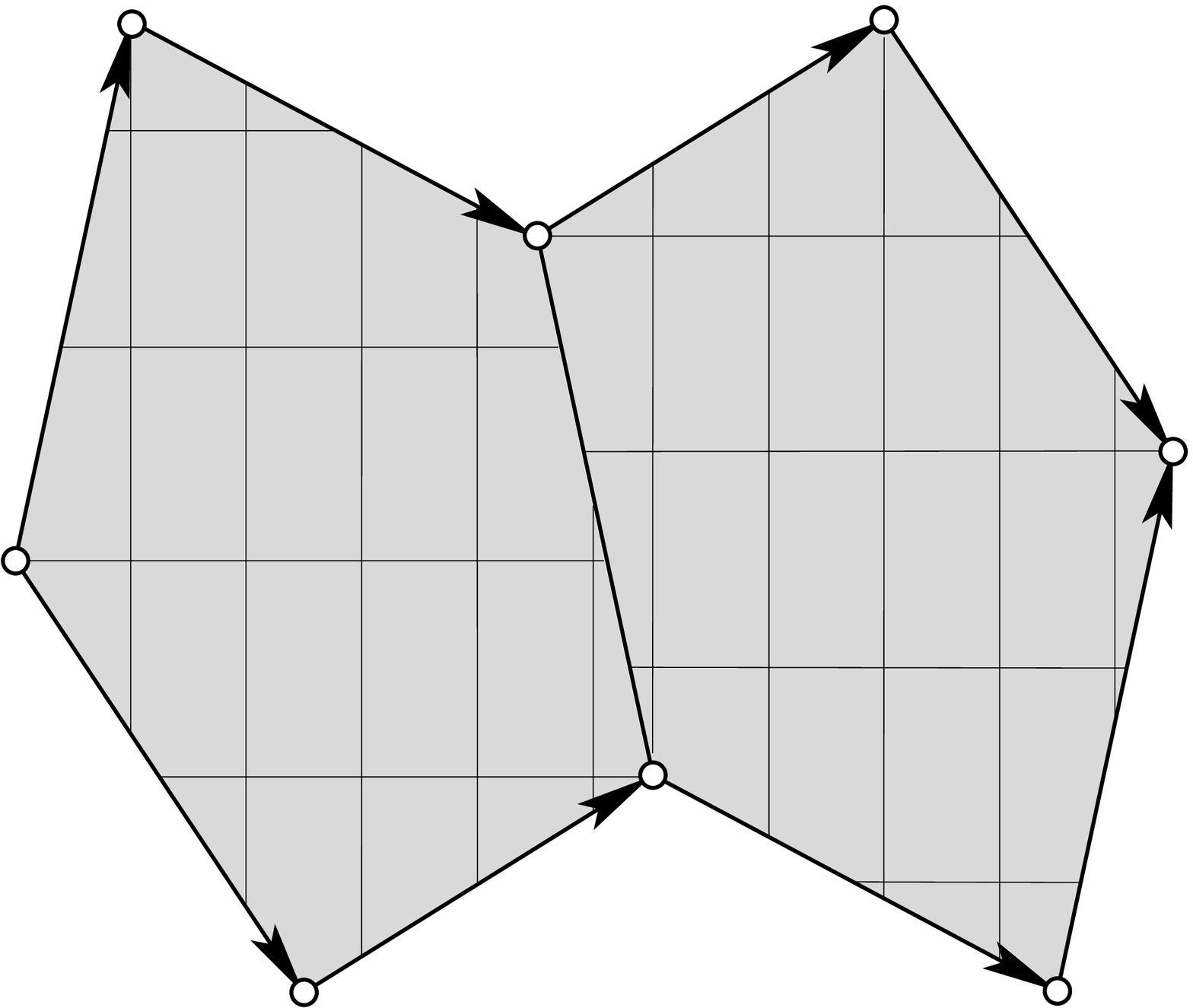}
\includegraphics{zorich_ciseaux6_1cm.eps}
\begin{picture}(0,0)(0,0)
\put(-95,-70){\LARGE$\longrightarrow$}
\put(40,-70){\LARGE$\ =$}
\end{picture}
\vspace{130bp}
\caption{
\label{zorich:fig:different:ways:to:unwrap}
The first modification of the polygon changes the
flat structure while the second one just changes the way in which we unwrap the flat surface
}
\end{figure}

There  is  no contradiction since the statement  of  the  theorem
concerns flat surfaces and  not  polygons. In practice this means
that we can apply expansion-contraction  to  the  polygon  $\Pi$,
which  does  not change too much  the  shape of the polygon,  but
radically changes the flat structure. Then we can  change the way
in   which   we   unwrap   the   flat   surface   $g_t   S$  (see
Fig.~\ref{zorich:fig:different:ways:to:unwrap}). This radically  changes
the shape  of the polygon, but \emph{does not  change at all} the
flat structure!

\index{Action on the moduli space!ofGL@of $GL^+(2,\R{})$|)}
\index{0GL@$GL^+(2,\R{})$-action on the moduli space|)}
\index{Action on the moduli space!ofSL@of $SL(2,\R{})$|)}
\index{0SL@$SL(2,\R{})$-action on the moduli space|)}

\subsection{General Philosophy}
\label{zorich:ss:General:Philosophy}

Now we are  ready  to describe informally the  basic  idea of our
approach to the  study  of  flat surfaces. Of course  it  is  not
universal; however, in many  cases  it appears to be surprisingly
powerful.

Suppose that we  need some information about geometry or dynamics
of an  individual flat surface $S$.  Consider the element  $S$ in
the  corresponding  family  of  flat  surfaces  $\cH(d_1,  \dots,
d_\noz)$\index{0H20@$\cH(d_1,\dots,d_\noz)$ -- stratum in the moduli space}.
Denote   by  $\cC(S)=\overline{GL^+(2,\R{})\,  S}\subset
\cH(d_1, \dots, d_\noz)$ the closure of the $GL^+(2,\R{})$-orbit of
$S$  in   $\cH(d_1,  \dots,  d_\noz)$.  \emph{In  numerous  cases
knowledge about the  structure  of $\cC(S)$ gives a comprehensive
information  about  geometry  and  dynamics of the  initial  flat
surface $S$. Moreover, some delicate numerical characteristics of
$S$ can be  expressed as averages of simpler characteristics over
$\cC(S)$.}

The  remaining  part of this survey  is  an attempt to show  some
implementations  of   this  general  philosophy.  The  first  two
illustrations would be presented in the next section.

We  have  to  confess  that  we do not tell all the truth in  the
formulation above. Actually, there is a hope that this philosophy
extends  much  further. A  closure  of an  orbit  of an  abstract
dynamical  system  might have  extremely  complicated  structure.
According to the most  optimistic  hopes, the closure $\cC(S)$ of
the $GL^+(2,\R{})$-orbit  of  \emph{any}  flat  surface $S$ is  a  nice
complex-analytic variety.  Moreover, according to the most daring
conjecture  it  would  be  possible   to   classify   all   these
$GL^+(2,\R{})$-invariant  subvarieties.  For  genus two the  latter
statements     were    recently     proved     by     C.~McMullen
(see~\cite{zorich:McMullen:Hilbert}   and~\cite{zorich:McMullen:genus:2})   and
partly by K.~Calta~\cite{zorich:Calta}.

We    discuss     this     hope     in     more    details    in
Sec.~\ref{zorich:s:Hope:for:a:Magic:Stick:and:Recent:Results},        in
particular, in Sec.~\ref{zorich:ss:Main:Hope}. We complete this  section
by a Theorem which supports  the  hope for some nice and  simple
description of orbit closures.

\begin{NNTheorem}[M.~Kontsevich]
Suppose that a closure $\cC(S)$ in
$\cH(d_1, \dots, d_\noz)$\index{0H20@$\cH(d_1,\dots,d_\noz)$ -- stratum in the moduli space}
of a
$GL^+(2,\R{})$-orbit of some flat surface $S$ is a complex-analytic
subvariety.      Then      in      cohomological      coordinates
$H^1(S,\{P_1,\dots,P_\noz\};\C{})$ it is represented by an affine
subspace.
\end{NNTheorem}

\subsection{Implementation of General Philosophy}
\label{zorich:ss:Implementation:of:General:Philosophy}

In this section we present two illustrations showing how the
``general philosophy'' works in practice.

Consider a directional flow on a flat surface  $S$. It is called
\emph{minimal}\index{Flow!minimal}
when the closure  of any trajectory  gives the entire  surface.
When  a directional flow  on  a  flat  torus is  minimal,  it  is
necessarily
ergodic\index{Flow!ergodic},
in particular, any trajectory in average  spends  in  any  subset
$U\subset\T{2}$ a time  proportional to the area (measure) of the
subset  $U$.  Surprisingly,  for  surfaces  of  higher  genera  a
directional flow can be minimal but not ergodic!  Sometimes it is
possible  to  find some  special  direction  with  the  following
properties. The flow in  this  direction is minimal. However, the
flat surface $S$  might be decomposed  into a disjoint  union  of
several subsets $V_i$ of positive measure in such a way that some
trajectories of the directional flow  prefer  one  subset to the
others. In other words, the average time spent by a trajectory in
the  subset  $V_i$ is  not  proportional  to  the  area  of $V_i$
anymore.  (The   original   ideas   of   such   examples   appear
in~\cite{zorich:Veech:1969}, \cite{zorich:Katok:1973}, \cite{zorich:Sataev:1975}
\cite{zorich:Keane:nonergodic};  see
also~\cite{zorich:Masur:Tabachnikov} and
especially~\cite{zorich:Masur:Handbook:1B}   for   a   very   accessible
presentation of such examples.)

Suppose  that  we managed  to  find a  direction  on the  initial
surface $S_0$ such that the flow in this direction is minimal but
not ergodic (with respect  to  the natural Lebesgue measure). Let
us apply  a rotation to  $S_0$ which would make the corresponding
direction vertical. Consider the resulting flat  surface $S$ (see
Convention~\ref{zorich:conv:flat:surface}                             in
Sec.~\ref{zorich:ss:Very:Flat:Surfaces}).  Consider  the   corresponding
``point''  $S\in\cH(d_1,\dots,d_\noz)$   and  the  orbit   $\{g_t
S\}_{t\in\R{}}$ of $S$  under the action of the diagonal subgroup
$g_t=\begin{pmatrix} e^t & 0\\ 0 & e^{-t} \end{pmatrix}$.

Recall  that  the stratum (or, more precisely, the  corresponding
``unit hyperboloid'')
$\cH_1(d_1,\dots,d_\noz)$\index{0H30@$\cH_1(d_1,\dots,d_\noz)$ -- ``unit hyperboloid''}\index{Stratum!in the moduli space}\index{Unit hyperboloid}
is never compact,
it always contains
``cusps''\index{Cusp of the moduli space@Cusp of the moduli space, \emph{see also} Moduli space; principal boundary}:
regions where the corresponding flat  surfaces  have  very  short
saddle  connections   or   very   short   closed  geodesics  (see
Sec.~\ref{zorich:ss:Toy:Example:Family:of:Flat:Tori}).

\begin{NNTheorem}[H.~Masur]
Consider a flat surface $S$. If the vertical flow is  minimal but
not ergodic  with respect to  the natural Lebesgue measure on the
flat  surface  then  the  trajectory  $g_t  S$  of
the Teichm\"uller geodesic flow\index{Teichm\"uller!geodesic flow}
is divergent,  i.e.  it eventually leaves any fixed
compact subset $K\subset\cH_1(d_1,\dots,d_\noz)$ in the stratum.
\end{NNTheorem}

Actually, this theorem has an even stronger form.

A  stratum   $\cH_1(d_1,\dots,d_\noz)$\index{0H30@$\cH_1(d_1,\dots,d_\noz)$ -- ``unit hyperboloid''}\index{Stratum!in the moduli space}\index{Unit hyperboloid}
has   ``cusps''  of  two
different origins. A flat surface  may  have  two distinct zeroes
get very close to each other. In this case $S$ has a short
saddle connection\index{Saddle!connection}
(or, what is the same, a short
relative period\index{Period!relative}).
However, the  corresponding  Riemann  surface  is  far from being
degenerate. The cusps of this type  correspond  to  ``simple
noncompactness'':   any   stratum  $\cH_1(d_1,\dots,d_\noz)$\index{0H30@$\cH_1(d_1,\dots,d_\noz)$ -- ``unit hyperboloid''}\index{Stratum!in the moduli space}\index{Unit hyperboloid}   is
adjacent         to         all        ``smaller''         strata
$\cH_1(d_1+d_2,d_3,\dots,d_\noz)$, $\dots$.

Another type of degeneration of a flat surface is the appearance of a
short  closed  geodesic. In  this  case  the  underlying  Riemann
surface  is  close to  a  degenerate  one;  the cusps of this second type
correspond to ``essential noncompactness''.

To formulate a stronger version of the above Theorem consider
the natural projection  of the stratum
$\cH(d_1,\dots,d_\noz)$\index{0H20@$\cH(d_1,\dots,d_\noz)$ -- stratum in the moduli space}\index{Stratum!in the moduli space}
to the moduli      space      $\cM_g$\index{0M10@$\cM_g$ -- moduli space of complex structures}\index{Moduli space!of complex structures}     of      complex
structures (see~\eqref{zorich:eq:projection:H:to:M}
in Sec.~\ref{zorich:ss:Dictionary:of:Complex:Analytic:Language}).
Consider the image of the orbit  $\{g_t  S\}_{t\in\R{}}$  in
$\cM_g$\index{0M10@$\cM_g$ -- moduli space of complex structures}\index{Moduli space!of complex structures} under this  natural  projection.  By  the reasons which
we  explain  in Sec.~\ref{zorich:s:Crash:Course:in:Teichmuller:Theory}
it  is natural to call this image a \emph{Teichm\"uller
geodesic}\index{Teichm\"uller!geodesic flow}.

\begin{NNTheorem}[H.~Masur]
Consider a flat surface $S$. If the vertical flow is  minimal but
not ergodic  with respect to  the natural Lebesgue measure on the
flat  surface  then  the  ``Teichm\"uller geodesic'' $g_t  S$  is
divergent, i.e. it eventually leaves  any  fixed  compact  subset
$K\subset  \cM$  in the moduli space of  complex  structures  and
never visits it again.
\end{NNTheorem}

This statement (in  a  slightly different formulation) is usually
called  \emph{Masur's  criterion  of   unique   ergodicity}  (see
Sec.~\ref{zorich:s:Ergodic:Theorem}   for   discussion  of   the  notion
\emph{unique ergodicity}).

As  a  second illustrations  of  the  ``general  philosophy''  we
present a combination  of \emph{Veech criterion} and of a Theorem
of J.~Smillie.

Recall that closed regular geodesics on a flat  surface appear in
families of parallel closed geodesics. When the flat surface is a
flat  torus, any such  family  covers  all  the torus. However,  for
surfaces of higher genera  such  families usually cover a cylinder
filled  with  parallel closed  geodesic  of  equal  length.  Each
boundary of such a cylinder contains a  conical  point.  Usually  a
geodesic emitted  in the same direction  from a point  outside of
the cylinder is  dense  in the complement to  the  cylinder or at
least in some nontrivial part of the complement. However, in some
rare cases, it may happen that the entire surface decomposes into
several cylinders filled with parallel closed  geodesics going in
some fixed direction. This is  the  case for the vertical or  for
the horizontal direction on the flat surface glued from a regular
octagon,  see  Fig.~\ref{zorich:fig:octagon:to:pretzel} (please  check).
Such direction is called \emph{completely periodic}\index{Direction!completely periodic}\index{Ergodic!uniquely ergodic}.

\begin{NNTheorem}[J.~Smillie; W.~A.~Veech]
Consider a flat  surface $S$. If its $GL^+(2,\R{})$-orbit is closed
in  $\cH(d_1,  \dots,  d_\noz)$  then a directional flow  in  any
direction  on  $S$ is  either  completely  periodic  or  uniquely
ergodic.
\end{NNTheorem}
(see Sec.~\ref{zorich:s:Ergodic:Theorem} for the   notion  of  \emph{unique  ergodicity}).
Note   that   unique   ergodicity
implies, in particular, that any orbit which is not a
\emph{saddle connection}\index{Saddle!connection},
(i.e. which does not hit the singularity both in forward and in backward direction)
is  everywhere
dense.    We     shall     return     to    this    Theorem    in
Sec.~\ref{zorich:ss:Veech:Surfaces}   where   we   discuss   \emph{Veech
surfaces}.

\section{How Do Generic Geodesics Wind Around Flat Surfaces}
\label{zorich:s:How:Do:Generic:Geodesics:Wind:Around:Flat:Surfaces}

In this section we study geodesics on a flat surface $S$ going in
generic\index{Geodesic!generic}
directions on  $S$.  Such  geodesics  are dense on  $S$;
moreover, it is possible  to show that they wind around $S$  in a
relatively regular manner. Namely, it is possible to find a cycle
$c\in  H_1(S;\R{})$  such that  in  some sense  a  long piece  of
geodesic  pretends  to wind around $S$ repeatedly following  this
\emph{asymptotic cycle} $c$.

\index{Asymptotic!cycle|(}

In Sec.~\ref{zorich:ss:Asymptotic:Cycle} we study the model  case of the
torus and give a rigorous definition  of  the  asymptotic  cycle.
Then we study the asymptotic cycles on general  flat surfaces. In
Sec.~\ref{zorich:ss:Deviation:from:Asymptotic:Cycle} we study  ``further
terms of approximation''.  The asymptotic  cycle describes the way in
which  a  geodesic  winds  around  the  surface  in  average.  In
Sec.~\ref{zorich:ss:Deviation:from:Asymptotic:Cycle}   we   present   an
empirical  description   of   the  deviation  from  average.  The
corresponding    rigorous   statements    are    formulated    in
Sec.~\ref{zorich:ss:Asymptotic:Flag:and:Dynamical:Hodge:Decomposition}.
Some  ideas  of the  proof  of this  statement  are presented  in
Sec.~\ref{zorich:s:Renormalization:Rauzy:Veech:Induction}.

\subsection{Asymptotic Cycle}
\label{zorich:ss:Asymptotic:Cycle}

\paragraph{Asymptotic Cycle on a Torus}

As usual  we start from the model  case of  the torus. We  assume
that our  flat torus  is glued from a square  in the natural way.
Consider an  irrational  direction  on  the  torus; any geodesic
going in this direction is dense in the torus.

Fix a point $x_0$ on the  torus and emit a geodesic in the chosen
direction. Wait  till it winds for some time  around the torus and
gets close  to the initial point $x_0$. Join  the endpoints of the
resulting piece of geodesic by a short path. We get a closed loop
on the torus which defines  a  cycle $c_1$ in the first  homology
group  $H_1(\T{2};\Z{})$  of  the  torus.  Now  let  the  initial
geodesic wind around the torus for some longer time; wait till it
get  close  enough  to  the  initial  point $x_0$  and  join  the
endpoints of the longer piece of geodesic by a short path. We get
a new  cycle  $c_2\in  H_1(\T{2};\Z{})$.  Considering  longer and
longer  geodesic  segments we get a sequence  of  cycles  $c_i\in
H_1(\T{2};\Z{})$.

\begin{figure}[htb]
\centering
\includegraphics{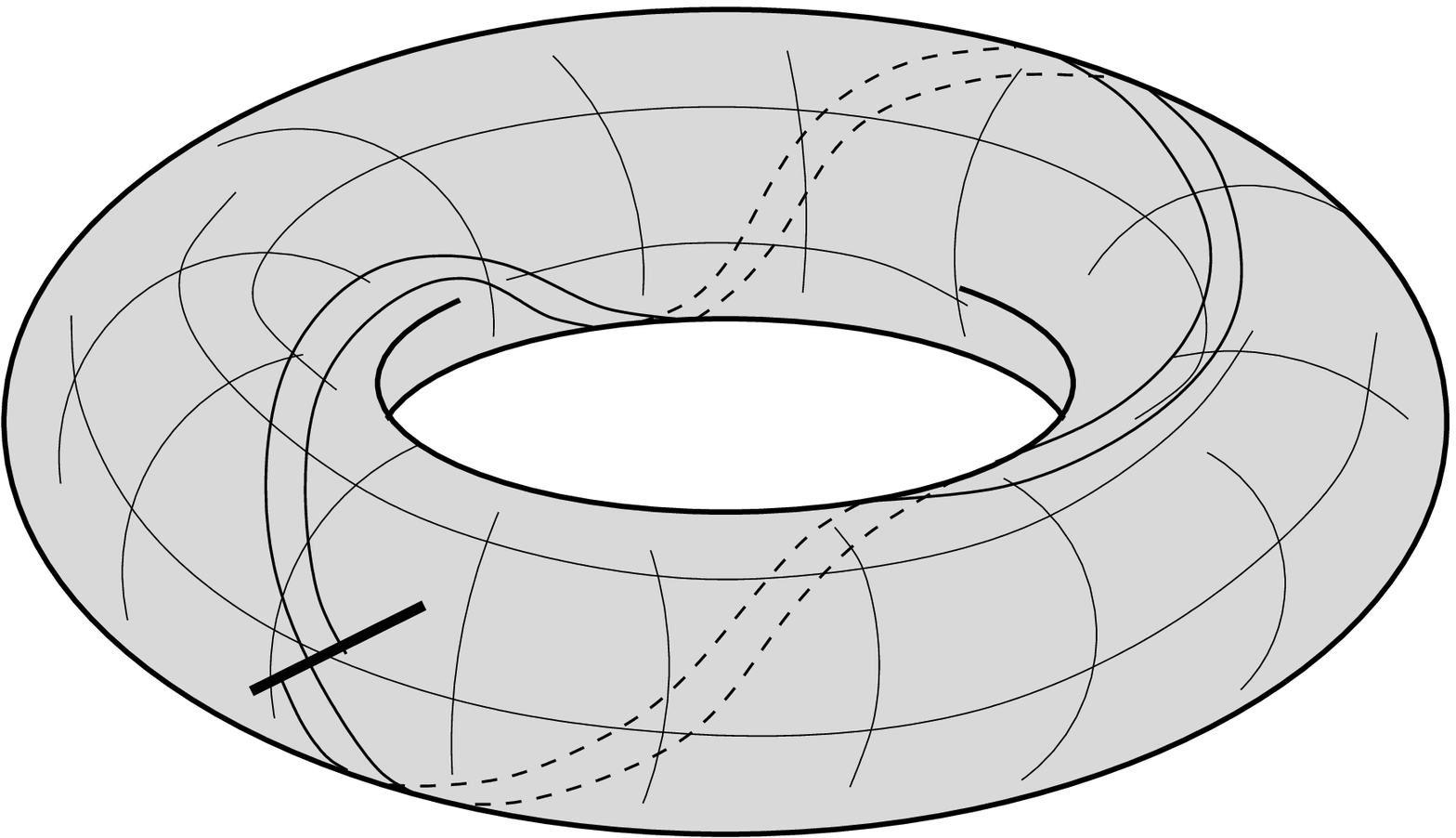}
\includegraphics{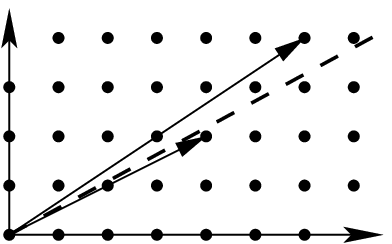}
\vspace{80bp} 
\caption{
\label{zorich:fig:asymptotic:cycle:on:T2}
A sequence of cycles approximating a dense geodesic on
a torus
}
\end{figure}

For example, we  can  choose a  short  segment $X$ going  through
$x_0$ orthogonal (or just transversal)  to  the  direction of the
geodesic. Each time  when the geodesic  crosses $X$ we  join  the
crossing point with the point $x_0$ along $X$  obtaining a closed
loop.  Consecutive  return points  $x_1,  x_2,  \dots$  define  a
sequence     of     cycles     $c_1,     c_2,     \dots$,     see
Fig.~\ref{zorich:fig:asymptotic:cycle:on:T2}.

For the torus case  we  can naturally identify  the  universal
covering space  $\R{2}\to\T{2}$  with  the  first  homology group
$H_1(\T{2};\R{})\simeq\R{2}$. Our irrational  geodesic unfolds to
an irrational straight line  $\cV_1$  in $\R{2}$ and the sequence
of cycles $c_1,  c_2,\dots$ becomes a sequence of integer vectors
$\vec v_1,  \vec  v_2,  \dots \in\Z{2}\subset\R{2}$ approximating
$\cV_1$, see Fig.~\ref{zorich:fig:asymptotic:cycle:on:T2}.

In particular, it is not surprising that there exists the limit
\begin{equation}
\label{zorich:eq:asymptotic:cycle:for:T2:normalized:by:norm}
\lim_{n\to\infty} \frac{c_n}{\|c_n\|} = c
\end{equation}
Under our  identification $H_1(\T{2};\R{})\simeq\R{2}$ the  cycle
$c$ represents a unit vector in direction $\cV_1$.

Let  the area of  the torus  be normalized to  one. Let  the
interval $X$,  which we use  to construct the sequence $c_1, c_2,
\dots$, be orthogonal to the direction of the geodesic. Denote by
$|X|$  its  length.  The  following  limit  also  exists  and  is
proportional to the previous one:
\begin{equation}
\label{zorich:eq:asymptotic:cycle:for:T2:normalized:by:n}
\lim_{n\to\infty} \frac{1}{n}\ c_n = \frac{1}{|X|}\cdot c
\end{equation}
The cycles obtained as
limits~\eqref{zorich:eq:asymptotic:cycle:for:T2:normalized:by:norm}
and~\eqref{zorich:eq:asymptotic:cycle:for:T2:normalized:by:n} are
called
\emph{asymptotic cycles}\index{Asymptotic!cycle}.
They show how the  corresponding irrational geodesic
winds around  the torus \emph{in average}. It is  easy to see that
they do not depend on the starting point $x$.

The    notion    ``asymptotic    cycle''   was   introduced    by
S.~Schwartzman~\cite{zorich:Schwartzman}.

\paragraph{Asymptotic Cycle on a Surface of Higher Genus}

We  can  apply  the same construction  to  a  geodesic  on a flat
surface $S$ of  higher genus. Having a geodesic segment $X\subset
S$ and some point $x\in X$ we emit from $x$ a geodesic orthogonal
to  $X$. From  time  to time the  geodesic  would intersect  $X$.
Denote the corresponding points as $x_1, x_2, \dots $. Closing up
the corresponding pieces of the geodesic by joining the endpoints
$x_0, x_j$ with a path going along $X$ we again get a sequence of
cycles $c_1, c_2, \dots$.

\begin{proposition}
For any flat surface $S$ of area one and for almost any direction
$\alpha$ on it any geodesic going in direction  $\alpha$ is dense
on  $S$  and  has  an  asymptotic  cycle which  depends  only  on
$\alpha$.

In other words, for almost any direction the limit
$$
\lim_{n\to\infty} \frac{1}{n}\ c_n = \frac{1}{|X|}\cdot c
$$
exists and the corresponding asymptotic cycle $c$ does not depend
on the starting point $x_0\in S$.
\end{proposition}

This  proposition  is an elementary corollary from the  following
theorem        of        S.~Kerckhoff,        H.~Masur        and
J.~Smillie~\cite{zorich:Kerckhoff:Masur:Smillie} (which is  another  Key
Theorem in this area).

\begin{NNTheorem}[S.~Kerckhoff, H.~Masur, J.~Smillie]
For  any flat  surface  $S$ the directional  flow  in almost  any
direction is ergodic.
\end{NNTheorem}

In  this  case  the  asymptotic  cycle  has  the  same  dynamical
interpretation as for  the torus: it  shows how  a  geodesic
going  in  the  chosen  direction  winds  around  the  surface $S$
\emph{in average}.

\begin{NNRemark}
Note  that  the asymptotic cycle $c\in H_1(S,\R{})$  also  has  a
topological  interpretation.  Assume  for  simplicity  that   $c$
corresponds  to  the vertical  direction.  Let  $\omega$  be  the
holomorphic 1-form corresponding to the  flat  structure  on  $S$
(see Sec.~\ref{zorich:ss:Dictionary:of:Complex:Analytic:Language}). Then
the  closed  1-form  $\omega_0=\Re(\omega)$ defines the  vertical
foliation  and   $c=D[\omega_0]$   is   Poincar\'e  dual  to  the
cohomology class of $\omega_0$. Choosing other ergodic directions
on  the  flat  surface  $S$  we  get  asymptotic  cycles  in  the
two-dimensional  subspace   $\langle D[\omega_0],D[\omega_1]\rangle_{\R{}}\subset
H_1(S,\R{})$  spanned   by  homology  classes  dual  to  cocycles
$\omega_0=\Re(\omega)$ and $\omega_1=\Im(\omega)$.
\end{NNRemark}

\subsection{Deviation from Asymptotic Cycle}
\label{zorich:ss:Deviation:from:Asymptotic:Cycle}

We have seen in the  previous  section that a sequence of  cycles
$c_1, c_2, \dots$ approximating long pieces  of an ``irrational''
geodesic on  a flat torus $\T{2}$ and  on a  flat surface $S$  of
higher  genus   exhibit   similar   behavior:  their  norm  grows
(approximately) linearly  in $n$ and their direction approaches
the direction of the asymptotic  cycle  $c$.  Note, however, that
for the  torus the cycles $c_n$  live in the two-dimensional space
$H_1(\T{2};\R{})\simeq\R{2}$,  while  for  the surface of  higher
genus    $g\ge    2$   the   cycles   live   in the  larger    space
$H_1(S;\R{})\simeq\R{2g}$. In particular, they have ``more room''
for deviation from the asymptotic direction.

Namely,       observing       the       right       part       of
Fig.~\ref{zorich:fig:asymptotic:cycle:on:T2}  we see  that  all  vectors
$c_n$ follow the line $\cV_1$ spanned by the asymptotic cycle $c$
rather  close:  the norm  of  projection  of  $c_n$  to  the line
orthogonal to $\cV_1$ is  uniformly  bounded (with respect to $n$
and to the choice of the starting point $x_0$).

The situation is  different for surfaces of higher genera. Choose
a hyperplane $\cS\perp c$ in $H_1(S,\R{})$ as a screen orthogonal
(transversal)  to  the  asymptotic  cycle  $c$   and  consider  a
projection to this  screen parallel to $c$. Projections of cycles
$c_n$  would  not  be  uniformly  bounded  anymore.  There  is no
contradiction  since  if the  norms  of  these  projections  grow
sublinearly, then the directions  of  the cycles $c_n$ still tend
to direction of the asymptotic cycle $c$.

\begin{figure}[htb]
\centering
\includegraphics{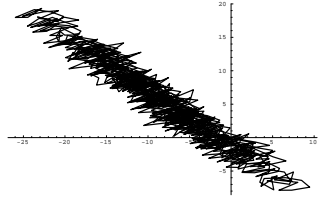}
\vspace{90bp}
\caption{
\label{zorich:fig:asymptotic:cycle:old:experiments}
Projection of a broken line joining the endpoints of $c_1, c_2,
..., c_{100000}$ to a screen orthogonal to the asymptotic cycle.
Genus  $g=3$}
\end{figure}

Let us observe how the projections are distributed  in the screen
$\cS$.  Figure~\ref{zorich:fig:asymptotic:cycle:old:experiments}   shows
results of numerical experiments where we take a  projection of a
broken   line   joining  the  endpoints  of  $c_1,  c_2,   \dots,
c_{100000}$ and we take a two-dimensional screen orthogonal to $c$
to make the picture more explicit.

\begin{figure}[hbt]
\special
{
   %
psfile=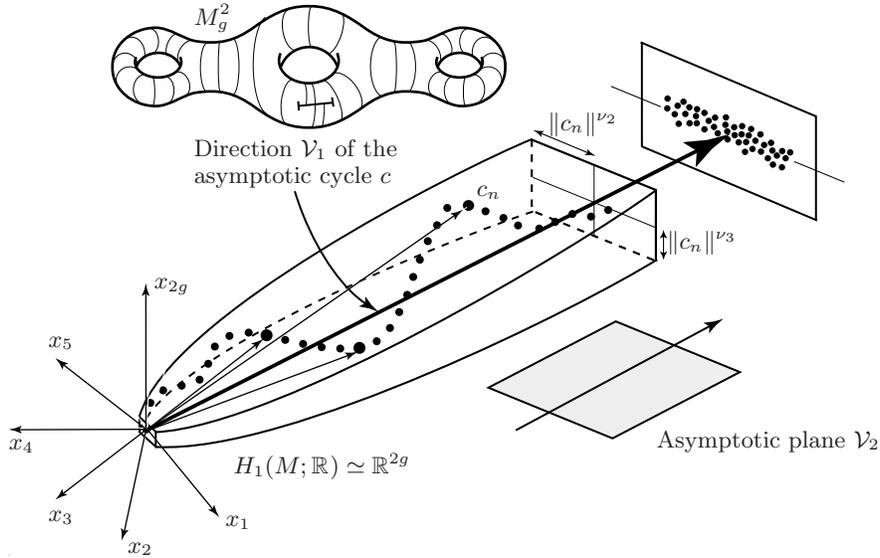
hscale= 100
vscale= 100
hoffset=5
voffset=-210
}
\begin{picture}(15,-10)(15,-10)
\put(0,0) 
 {\begin{picture}(0,0)(0,0)
 \put(197,-85){$c_n$}
 \put(224,-60){$\|c_n\|^{\nu_2}$}
 \put(269,-105){$\|c_n\|^{\nu_3}$}
\end{picture}}
\put(0,0)
 {\begin{picture}(0,0)(0,0)
   \put(105,-190){$H_1(M;\R{})\simeq\R{2g}$}
   \put(102,-210){$x_1$}
   \put(65,-220){$x_2$}
   \put(35,-208){$x_3$}
   \put(20,-181){$x_4$}
   \put(35,-142){$x_5$}
   \put(75,-120){$x_{2g}$}
 \end{picture}}
\put(0,0)
 {\begin{picture}(0,0)(0,0)
   \put(265,-180){\small Asymptotic plane $\cV_2$}
 \end{picture}}
\put(0,0)
 {\begin{picture}(0,0)(0,0)
   \put(90,-70){\small Direction $\cV_1$ of the}
   \put(90,-80){\small asymptotic cycle $c$}
 \end{picture}}
\put(0,0)
 {\begin{picture}(0,0)(0,0)
   \put(90,-20){$M^2_g$}
 \end{picture}}
\end{picture}
\vspace{220bp}
\caption{
\label{zorich:fig:asymptotic:cycle:deviation}
Deviation from the asymptotic direction
}
\end{figure}

We see  that the distribution of  projections of cycles  $c_n$ in
the screen $\cS$ is anisotropic: the projections accumulate along
some line. This means that in the original space $\R{2g}$ the vectors
$c_n$  deviate   from   the   asymptotic  direction  $\cV_1$  not
arbitrarily but along some two-dimensional subspace $\cV_2\supset
\cV_1$, see Fig.~\ref{zorich:fig:asymptotic:cycle:deviation}.

Moreover, measuring the  norms  of the projections $proj(c_n)$ to
the screen $\cS$ orthogonal to $\cL_1=\langle c\rangle_{\R{}}$,
we get\index{0nu10@$\nu_1,\dots,\nu_g$ -- Lyapunov exponents related to the Teichm\"uller geodesic flow}\index{Exponent!Lyapunov exponent}\index{Lyapunov exponent}
$$
\limsup_{n\to\infty} \frac{\log\|proj(c_n)\|}{\log n} = \nu_2<1
$$
In  other  words the vector  $c_n$ is located approximately  in  the
subspace $\cV_2$, and the distance from its endpoint  to the line
$\cV_1\subset \cV_2$  is  bounded  by $const\cdot|c_n\|^{\nu_2}$, see
Fig.~\ref{zorich:fig:asymptotic:cycle:deviation}.

Consider now a new screen  $\cS_2\perp\cV_2$  orthogonal  to  the
plane  $\cV_2$.  Now the screen $\cS_2$ has  codimension  two  in
$H_1(S,\R{})\simeq\R{2g}$.  Considering the projections of  $c_n$ to
$\cS_2$  we  eliminate  the  asymptotic  directions  $\cV_1$  and
$\cV_2$ and we see how do the vectors $c_n$ deviate from $\cV_2$.
On  the  screen   $\cS_2$  we  see   the  same  picture   as   in
Fig.~\ref{zorich:fig:asymptotic:cycle:old:experiments}: the  projections
are located along a one-dimensional subspace.

Coming back  to  the ambient space $H_1(S,\R{})\simeq\R{2g}$, this
means that in the first term of approximation  all vectors $c_n$
are aligned along the one-dimensional subspace $\cV_1$ spanned by
the asymptotic cycle. In  the  second term of approximation, they
can deviate from $\cV_1$, but the deviation occurs  mostly in the
two-dimensional subspace $\cV_2$, and  has  order $\|c\|^{\nu_2}$
where $\nu_2<1$.  In the third  term of approximation we see that
the vector  $c_n$  may  deviate  from  the  plane  $\cV_2$,  but  the
deviation occurs mostly in a three-dimensional  space $\cV_3$ and
has order $\|c\|^{\nu_3}$ where $\nu_3<\nu_2$.

Going on we  get further terms of approximation. However, getting
to a subspace $\cV_g$ which has  half the dimension  of  the  ambient
space we shall see  that, in a sense, there is no  more deviation
from $\cV_g$: the distance from any $c_n$ to $\cV_g$ is uniformly
bounded.

Note   that    the    intersection    form   endows   the   space
$H_1(S,\R{})\simeq\R{2g}$ with a natural symplectic structure. It
can  be  checked  that  the  resulting  $g$-dimensional  subspace
$\cV_g$ is a \emph{Lagrangian} subspace for this symplectic form.

\subsection{Asymptotic Flag and ``Dynamical Hodge Decomposition''}
\label{zorich:ss:Asymptotic:Flag:and:Dynamical:Hodge:Decomposition}

A rigorous formulation of phenomena  described  in  the  previous
section is given by the following Theorem
proved\footnote{Actually, the theorem was initially proved  under
certain hypothesis on the Lyapunov exponents of the Teichm\"uller
geodesic  flow.  These conjectures were proved later by  G.~Forni
and in the most complete form by A.~Avila and M.~Viana;  see the
end       of       this        section       and       especially
Sec.~\ref{zorich:ss:Spectrum:of:Lyapunov:exponents}}
by         the         author          in~\cite{zorich:Zorich:Deviation}
and~\cite{zorich:Zorich:How:do}.

Following Convention~\ref{zorich:conv:flat:surface} we always consider a
flat  surface  together  with  a  choice  of  direction  which by
convention   is   called  the   \emph{vertical   direction},   or
\emph{direction to the North}. Using an  appropriate homothety we
normalize   the   area   of   $S$   to   one,   so   that   $S\in
\cH_1(d_1,\dots,d_\noz)$.

We chose a point $x_0\in S$ and a horizontal segment  $X$ passing
through $x_0$; by $|X|$ we denote the length of $X$.  We consider
a  geodesic  ray  $\gamma$  emitted  from  $x_0$ in the  vertical
direction.  (If  $x_0$  is  a  saddle  point,  there  are several
outgoing vertical geodesic rays;  choose  any of them.) Each time
when  $\gamma$  intersects   $X$  we  join  the  point  $x_n$  of
intersection and the starting point  $x_0$  along  $X$ producing a
closed path. We denote the  homology  class  of the corresponding
loop by $c_n$.

Let $\omega$ be  the holomorphic 1-form representing $S$; let $g$
be   genus   of   $S$.  Choose  some  Euclidean   metric   in
$H_1(S;\R{})\simeq\R{2g}$ which would allow to measure a distance
from  a  vector   to  a  subspace.  Let  by  convention  $\log(0)
=-\infty$.

\begin{NNTheorem}
\label{zorich:th:mf}
For almost any flat surface $S$ in any stratum
$\cH_1(d_1,\dots,d_\noz)$\index{0H30@$\cH_1(d_1,\dots,d_\noz)$ -- ``unit hyperboloid''}\index{Stratum!in the moduli space}\index{Unit hyperboloid} there exists a
flag of subspaces\index{Flag!asymptotic}\index{Asymptotic!flag}
$$
\cV_1\subset \cV_2\subset\dots\subset \cV_g\subset H_1(S;\R{})
$$
in the first homology group of the surface with the following
properties.

Choose any starting point $x_0\in X$ in the horizontal segment
$X$. Consider the corresponding sequence $c_1, c_2, \dots $ of
cycles.

--- The following limit exists
$$
|X| \lim_{n\to\infty} \frac{1}{n}\, c_n = c,
$$
where the nonzero   asymptotic   cycle   $c\in  H_1(M^2_g;\R{})$  is
Poincar\'e     dual     to    the     cohomology     class     of
$\omega_0=\Re[\omega]$,   and    the   one-dimensional   subspace
$\cV_1=\langle c\rangle_{\R{}}$ is spanned by $c$.

--- For any $j=1,\dots, g-1$ one has
$$
\limsup_{n\to\infty} \frac{\log\,dist(c_n,\cV_j)}{\log n}  =
\nu_{j+1}
$$
and
$$
|dist(c_n,\cV_g)|  \le const,
$$
where the constant depends  only on $S$ and on the choice  of the
Euclidean structure in the homology space.

The numbers  $2,1+\nu_2,\dots,1+\nu_g$\index{0nu10@$\nu_1,\dots,\nu_g$ -- Lyapunov exponents related to the Teichm\"uller geodesic flow}
are  the top $g$
Lyapunov exponents\index{Lyapunov exponent}\index{Exponent!Lyapunov exponent}
of the Teichm\"uller geodesic flow\index{Teichm\"uller!geodesic flow}
on the corresponding connected
component\index{Moduli space!connected components of the strata}
of  the  stratum  $\mathcal{H}(d_1,\dots,d_\noz)$;  in
particular, they do not depend on  the  individual  generic  flat
surface $S$ in the connected component.
\end{NNTheorem}

\index{Asymptotic!cycle|)}

A reader who is  not  familiar with \emph{Lyapunov exponents} can
either           read           about           them           in
Appendix~\ref{zorich:s:Multiplicative:ergodic:theorem} or just  consider
the numbers
$\nu_j$\index{0nu10@$\nu_1,\dots,\nu_g$ -- Lyapunov exponents related to the Teichm\"uller geodesic flow}\index{Exponent!Lyapunov exponent}\index{Lyapunov exponent}
as some abstract constants which depend only
on   the   connected   component   $\cH^{comp}(d_1,\dots,d_\noz)$
containing the flat surface $S$.

It  should  be  stressed,  that the theorem above  was  initially
formulated  in~\cite{zorich:Zorich:How:do}  as a  conditional statement:
under the  conjecture  that  $\nu_g>0$  there  exist a Lagrangian
subspace $\cV_g$  such that the  cycles are in a bounded distance
from $\cV_g$; under the further conjecture that all the exponents
$\nu_j$,  for  $j=1,2,\dots,g$,  are  distinct, there is  a  {\it
complete}
Lagrangian flag\index{Flag!Lagrangian}\index{Lagrangian!flag}
(i.e.  the dimensions of the subspaces
$\cV_j$, where $j=1,2,\dots,g$, rise each time by one). These two
conjectures were  later proved by G.~Forni~\cite{zorich:Forni:02} and by
A.~Avila  and   M.~Viana~\cite{zorich:Avila:Viana}.  We  discuss   their
theorems in Sec.~\ref{zorich:ss:Spectrum:of:Lyapunov:exponents}.

Another remark concerns the choice of the horizontal segment $X$.
By convention it  is  chosen in  such  way that the  trajectories
emitted in the vertical direction\index{Direction!vertical} (in  direction  to  the  North)
from the endpoints of $X$ hit the conical points before the first
return to $X$. Usually we just  place the left endpoint of $X$ at
the conical point.

Omitting  this  condition and considering a continuous family  of
horizontal subintervals  $X_t$  of  variable  length (say, moving
continuously one of the  endpoints),  the theorem stays valid for
a subset of $X_t$ of full measure.

\section{Renormalization for Interval Exchange Transformations.
Rauzy--Veech Induction}
\label{zorich:s:Renormalization:Rauzy:Veech:Induction}

In this section we elaborate a powerful time acceleration machine
which  allows   to   study the  asymptotic   cycles   described   in
Sec.~\ref{zorich:s:How:Do:Generic:Geodesics:Wind:Around:Flat:Surfaces}.
Following the  spirit of this  survey we put emphasis on geometric
ideas  and  omit proofs.  This  section can  be  considered as  a
geometric      counterpart      of       the      article      of
J.-C.~Yoccoz~\cite{zorich:Yoccoz:Les:Houches} in the current volume.

I use this opportunity to thank  M.~Kontsevich for numerous ideas
and conjectures  which  were  absolutely crucial for my impact in
this  theory:  without  numerous  discussions  with M.~Kontsevich
papers~\cite{zorich:Zorich:Gauss:map}     and   \cite{zorich:Zorich:Deviation},
probably, would be never written.

\subsection{First Return Maps and Interval Exchange Transformations}
\label{zorich:ss:First:Return:Map}

Our  goal  is  to  study  cycles  obtained from  long  pieces  of
``irrational''  geodesic  on a  flat  surface  by  joining  their
endpoints along a  transversal segment $X$. To perform this study
we elaborate some  simple machine which generates the cycles, and
then we  accelerate this machine to obtain very  long cycles in a
rather short time.

Consider \emph{all}  geodesics  emitted from a transverse segment
$X$ in the same generic direction  and let each of them come back
to $X$ for  the  first  time. We get a  \emph{first  return  map}
$T:X\to X$ which  is  interesting by  itself  and which deserves  a
separate discussion.  Its  properties  play a crucial  role in our
study.   (See    Appendix~\ref{zorich:s:Ergodic:Theorem}   for   general
properties of the first return map.)

\begin{figure}[htb]
\includegraphics{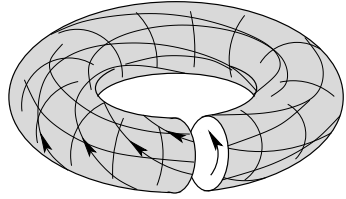}
\vspace{80bp}
\caption{
\label{zorich:fig:rotation:of:the:circle}
The first return map of a meridian to itself induced by a directional
flow is just a twist
}
\end{figure}

As usual  let  us start with a model case of a flat torus. Take a
meridian of  the torus  as a transversal $X$ and  emit from $X$ a
directional flow. Every geodesic  comes  back to $X$ inducing the
\emph{first   return   map}  $T:X\to  X$,  which  in  this   case
isometrically rotates the meridian  $X$  along itself by an angle
which   depends   on    the    direction   of   the   flow,   see
Fig.~\ref{zorich:fig:rotation:of:the:circle}.

Assume that  our flat torus $\T{2}$ is glued  from a unit square.
Let us replace now a meridian  of the torus by a generic geodesic
segment $X$ orthogonal to the direction of the  flow. From every
point of $X$ we emit  a  geodesic in direction orthogonal to  $X$
and wait till it hits $X$ for the first time.  We  again obtain a
first return map $T:X\to X$,  but  this time the map $T$ is slightly  more
complicated.

\begin{figure}[htb]
%
\includegraphics{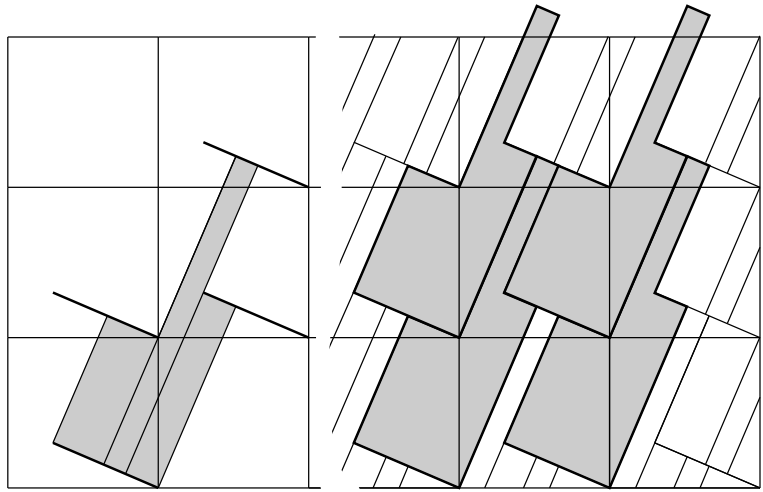}
\vspace{180bp}
\caption{
\label{zorich:fig:iet}
Directional flow on  a  torus. The first return map  of  a segment to
itself  is   an   interval   exchange   transformation  of  three
subintervals
}
\end{figure}

To study  this map it is convenient  to unfold  the torus into  a
plane. The map $T$  is  presented at Fig.~\ref{zorich:fig:iet}. It chops
$X$ into three  pieces and then  shuffles them sending  the  left
subinterval to  the right, the  right subinterval to the left and
keeping the middle one in the middle but shifting it a bit.
The  map $T$ gives an example of  an
{\it interval exchange transformation}\index{Interval exchange transformation|(}.

Note that when the direction  of  the  flow is \emph{irrational},
the  geodesics  emitted  from  $X$ again cover the  entire  torus
$\T{2}$ before  coming back to $X$. The torus  is get ripped into
three rectangles based on the three  subintervals  in  which  $T$
chops $X$. The corresponding building of three rectangles gives a
new fundamental domain representing the torus:  one  can  see  at
Fig.~\ref{zorich:fig:iet}  that  it tiles the plane. Initially we  glued
our flat torus from a  square;  the  building under consideration
gives another way to unwrap $\T{2}$ into a  polygon. We recommend
to  the  reader  to  check  that  identifying the  two  pairs  of
corresponding vertical sides of  the  building and gluing the top
horizontal  sides  of the  rectangles  to the  bottom  of $X$  as
prescribed by the interval exchange transformation $T$ we get the
initial torus.

\begin{figure}[htb]
%
\includegraphics{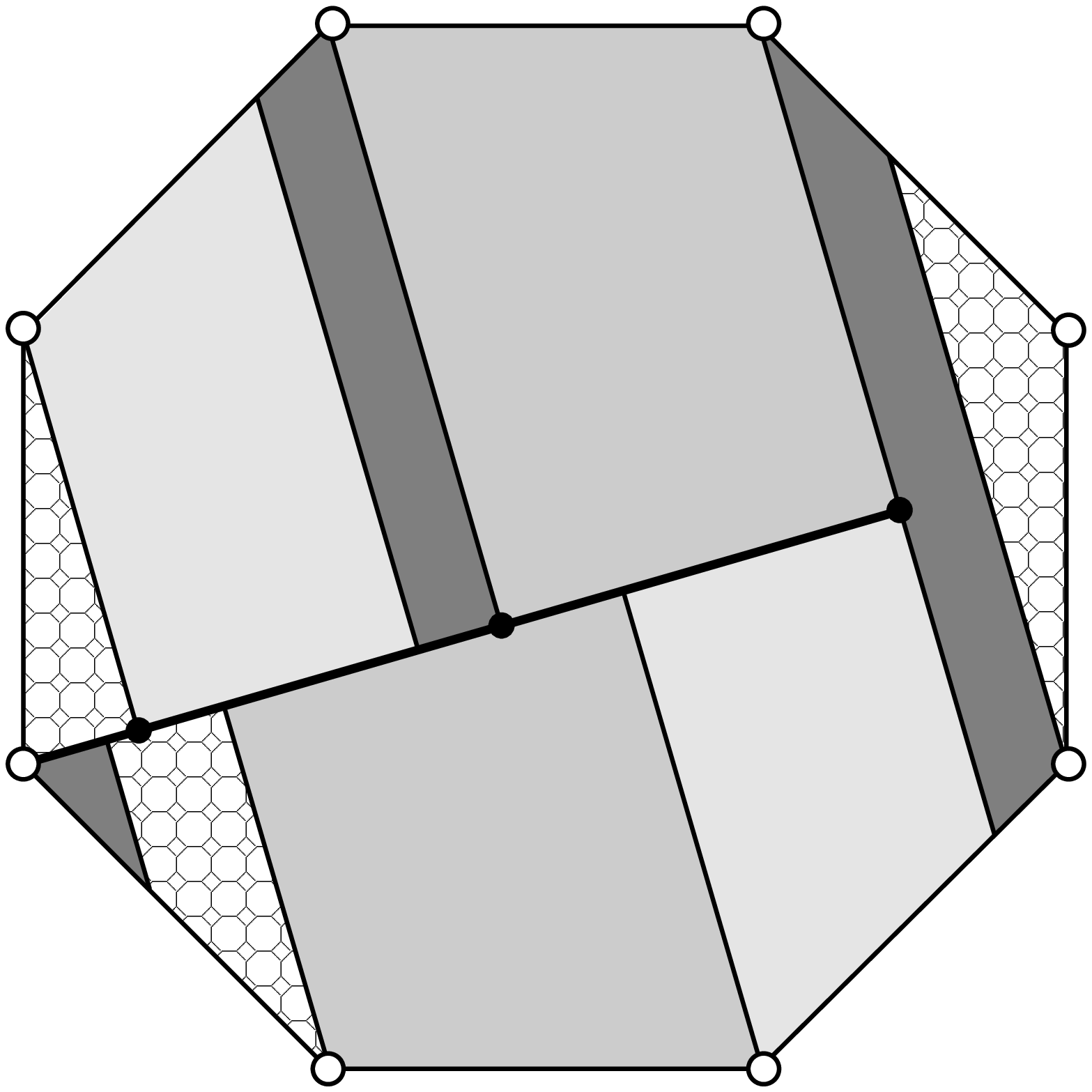}
\includegraphics{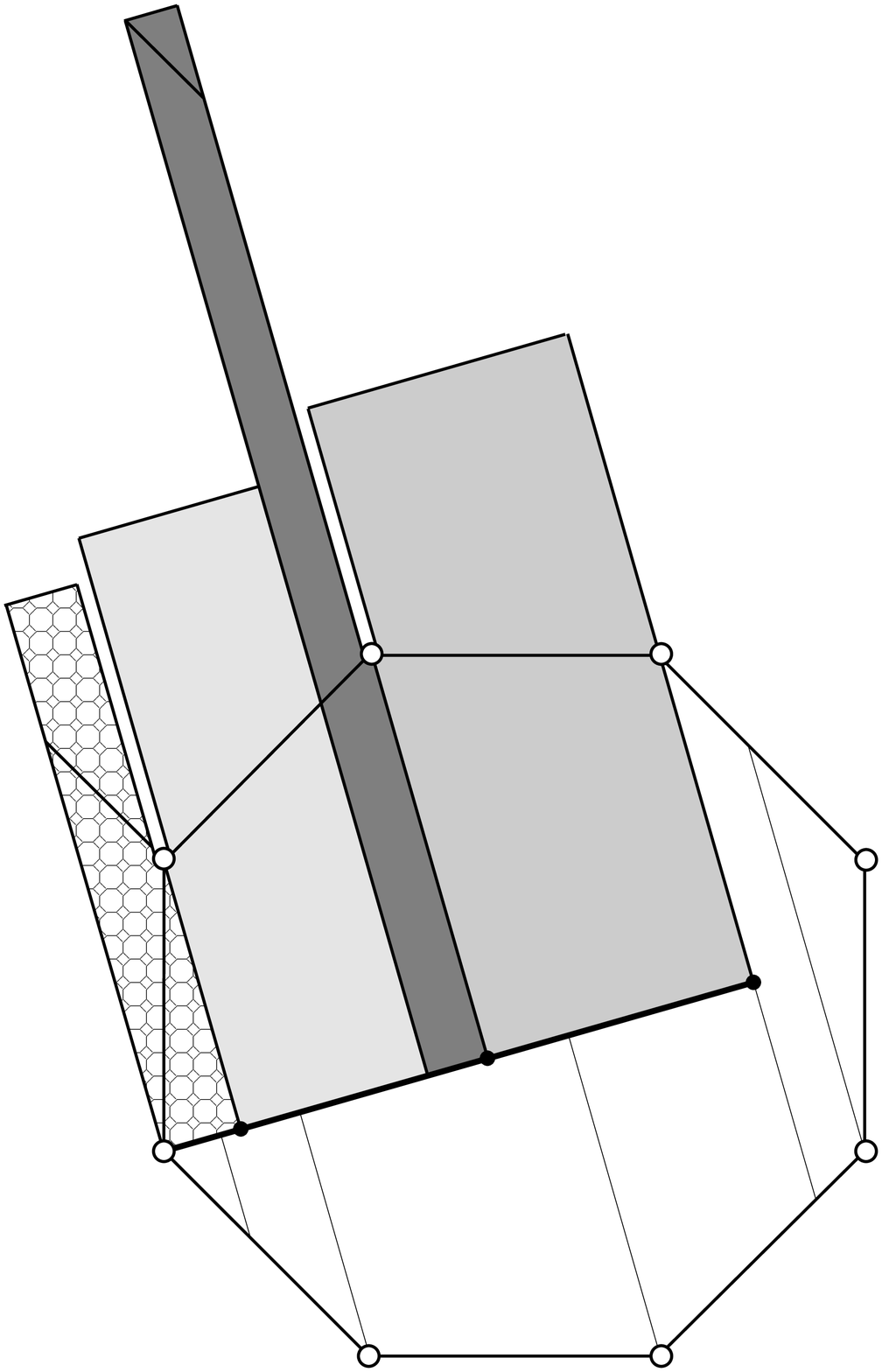}
\begin{picture}(100,50)(100,50)
\put(111,-48){\Large $\vec v_1$}
\put(144,-9){\Large $\vec v_2$}
\put(196,-5){\Large $\vec v_3$}
\put(232,-37){\Large $\vec v_4$}
\put(116,-95){\Large $\vec v_4$}
\put(157,-126){\Large $\vec v_3$}
\put(204,-119){\Large $\vec v_2$}
\put(236,-84){\Large $\vec v_1$}
\end{picture}
\vspace{180bp} 
\caption{
\label{zorich:fig:iet:octagon}
The first return map $T:X\to X$ of a geodesic segment $X$ defined by
a directional flow decomposes the surface into four rectangles
``zippered''\index{Zippered rectangle} along singular trajectories
}
\end{figure}

Consider  now  a flat  surface  of genus  higher  than one.  Say,
consider   a    flat    surface    of    genus    $g=2$   as   on
Fig.~\ref{zorich:fig:iet:octagon}. We  suggest  to  the  reader to check
that this  flat surface has a  single conical singularity  with a
cone angle  $6\pi$ (see Fig.~\ref{zorich:fig:monkey:saddle}). To study a
directional flow choose as before a geodesic segment $X\subset S$
orthogonal to  the direction of the  flow and consider  the first
return    map    $T:X\to   X$   induced   by   the   flow;    see
Fig.~\ref{zorich:fig:iet:octagon}. We  see  that  $X$  is  chopped  into
a larger  number  of  subintervals  (in comparison with  the  torus
case),  namely,  for our choice of  $X$  it is chopped into  four
subintervals.

Now we observe a  new  phenomenon: trajectories emitted from some
points of $X$  hit  the conical  point  and our directional  flow
splits at this point. Since in our particular case the cone angle
at the conical point is $6\pi=3\cdot 2\pi$ there  are {\it three}
trajectories  in   direction   $\vec   v$   which   hit  it.  The
corresponding points at which $X$ is chopped are marked with bold
dots. The remaining discontinuity  point  of $X$ corresponds to a
trajectory which hits the endpoint of $X$.

Our construction  with a segment  $X$ transversal to the flow and
with trajectories of the flow emitted from $X$  and followed till
their first return to $X$  trims  a braid from the flow.  Conical
points play the role  of a comb which splits the flow  into several
locks and  then trims them in  a different order.  Note, however,
that if we follow the flow  till the second return to $X$ it will
pass through  the comb twice,  and thus will be generically split
already  into  seven locks. (If you are  interested  in  details,
think why this second return has {\it seven} and not  {\it eight}
locks and what sort of genericity we need).

Similarly, the interval  exchange  transformation $T$ of the base
interval $X$ can be compared to a shuffling machine. Imagine that
$X$ represents  a stock of cards. We split  the stock into $\noi$
parts of fixed widths and shuffle the parts in a  different order
(given  by  some permutation $\pi$ of $\noi$  elements).  At  the
second iteration we again split the new stock in the parts of the
same widths as before and shuffle the parts according to the same
permutation $\pi$, etc.  Note, that even if the permutation $\pi$
is  such   that  $\pi^2=id$,  say,  $\pi=(4,3,2,1)$,  the  second
iteration $T^2$ is  not  an identical transformation provided the
widths $\lambda_1, \dots, \lambda_4$ are  not  symmetric:  for  a
generic choice  of  $\lambda_1,  \dots,  \lambda_4$  the interval
exchange transformation  $T^2$ has $6$ discontinuities (and hence
$7$ subintervals under exchange).

\begin{Exercise}
Consider  an interval  exchange  transformation  $T(\lambda,\pi)$
corresponding to the permutation $\pi=(4,3,2,1)$. Choose some generic
values  of   the   lengths   $\lambda_1,   \dots,  \lambda_4$  of
subintervals and construct $T^2$ and $T^3$.
\end{Exercise}

\subsection{Evaluation of the Asymptotic Cycle Using an Interval
Exchange Transformation}
\label{zorich:ss:Evaluation:of:the:asymptotic:cycle:using:an:iet}

Now we can return to our  original problem. We want to study long
pieces  of  leaves of the vertical foliation.  Fix  a  horizontal
segment $X$ and emit a vertical trajectory from  some point $x\in
X$. When  the trajectory intersects $X$  for the first  time join
the corresponding point $T(x)$ to  the  original  point $x$ along
$X$ to obtain a closed  loop.  Here $T:X\to X$ denotes the  first
return map  to the transversal  $X$ induced by the vertical flow.
Denote by  $c(x,1)$  the  corresponding  cycle  in $H_1(S;\Z{})$.
Following the vertical trajectory  further  on we shall return to
$X$ once again. Joining $x$ and the point $T(T(x))$ of the second
return to  $X$ along $X$ we obtain the  second cycle $c(x,2)$. We
want to  describe the cycle  $c(x,N)$ obtained after a very large
number $N$ of returns.

\begin{figure}[htb]
\centering
%
   %
\includegraphics{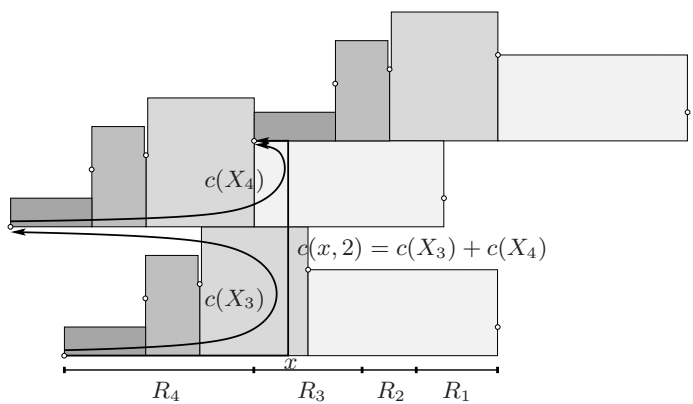}
\begin{picture}(0,0)(0,0)
\put(100,-160){
\begin{picture}(0,0)(0,0)
\put(-175,20){$R_4$}
\put(-120,20){$R_3$}
\put(-90,20){$R_2$}
\put(-65,20){$R_1$}
\put(-125,31){$x$}
\put(-155,55){$c(X_3)$}
\put(-155,100){$c(X_4)$}
\put(-120,74){$c(x,2)=c(X_3)+c(X_4)$}
\end{picture}}
\end{picture}
\vspace{140bp}
\caption{
\label{zorich:fig:first:return:cycles}
Decomposition of a long cycle into the sum of basic cycles. (We have
unfolded the flat surface along vertical trajectory emitted from the
point $x$)
}
\end{figure}

Actually, we prefer to close up a piece of trajectory  going from
$x\in  X$  to the first return  point  $T(x)\in X$ in a  slightly
different  way.  Instead of completing  the  path joining  the
endpoints it is more convenient to close this piece of trajectory
joining both points $x$ and  $T(x)$  to the left endpoint of  $X$
along $X$ (see Fig.~\ref{zorich:fig:first:return:cycles}). This modified
path defines the same homology cycle $c(x,1)$ as  the closed path
for which the points $x$ and $T(x)$ are joined directly.

Consider now the ``first return  cycle''  $c(x,1)$  as a function
$c(x)=c(x,1)$ of the starting point  $x\in  X$.  Let the interval
exchange transformation  $T:X\to  X$  decompose  $X$  into $\noi$
subintervals $X_1\sqcup \dots \sqcup  X_\noi$.  It is easy to see
that  the  function $c(x)$  is  piecewise  constant:  looking  at
Fig.~\ref{zorich:fig:first:return:cycles}  one  can  immediately  verify
that  if two points  $x_1$  and  $x_2$  are not  separated  by  a
discontinuity point (i.e.  if they belong to the same subinterval
$X_j$) they  determine  homologous  (and  even  homotopic) cycles
$c(x_1)=c(x_2)$. Each  subinterval $X_j$ determines its own cycle
$c(X_j)$.

The vertical trajectory emitted from a point $x\in  X$ returns to
$X$ at the points $T(x),  T^2(x),  \dots, T^N(x)$. It is easy  to
see that the cycle $c(x,2)$  corresponding  to  the second return
can  be   represented   as   a   sum  $c(x,2)=c(x)+c(T(x))$,  see
Fig.~\ref{zorich:fig:first:return:cycles}. Similarly  the cycle $c(x,N)$
obtained  by  closing  up  a  long  piece of vertical  trajectory
emitted from $x\in X$ and followed up to $N$-th return to $X$ can
be represented as a sum
\begin{equation}
\label{zorich:eq:c:N:equals:sum:c:i}
c(x,N)=c(x)+c(T(x))+\dots+c(T^{N-1}(x))
\end{equation}

According to  the  fundamental  Theorem of S.~Kerckhoff, H.~Masur
and   J.~Smillie~\cite{zorich:Kerckhoff:Masur:Smillie},   for  any   flat
surface the directional  flow in almost any direction is ergodic,
and  even  uniquely ergodic.  Hence, the  same  is  true  for the
corresponding  interval  exchange  transformation.  Applying  the
ergodic  theorem (see  Appendix~\ref{zorich:s:Ergodic:Theorem})  to  the
sum~\eqref{zorich:eq:c:N:equals:sum:c:i} and taking  into  consideration
that $c(x)$ is a piecewise-constant function we get
$$
c(x,N)\sim N\cdot\cfrac{1}{|X|}\int_X c(x)\, dx =
N \cdot \cfrac{1}{|X|}\, \Big( \lambda_1
c(X_1)+\dots+\lambda_\noi c(X_\noi) \Big)
$$
where $c(X_j)$ denotes  the ``first return cycle'' for the points
$x$       in        the       subinterval       $X_j$,        see
Fig.~\ref{zorich:fig:first:return:cycles}.   This   gives  an   explicit
formula for the
\emph{asymptotic cycle}\index{Asymptotic!cycle}
\begin{equation}
\label{zorich:eq:formula:for:asymptotic:cycle}
c=\lim_{N\to\infty}\cfrac{c(x,N)}{N}= \cfrac{1}{|X|}\, \Big(
\lambda_1 c(X_1)+\dots+\lambda_\noi c(X_\noi) \Big)
\end{equation}
Note that the asymptotic cycle does not depend on the starting
point $x\in X$.

\begin{Exercise}
Show  that  the  paths  ${\vec  v}_1,  \dots,  {\vec  v}_4$  (see
Fig.~\ref{zorich:fig:iet:octagon}) represent  a  basis  of cycles of the
corresponding  flat  surface  $S$  of  genus  $g=2$.  Show that the first
return cycles  $c(X_j)$  determine another basis of cycles and that
one can pass from one basis to the other using the following
relations (see Fig.~\ref{zorich:fig:iet:octagon}):
\begin{align*}
  c(X_1)&={\vec v}_1-{\vec v}_4 &
  c(X_2)&={\vec v}_1-{\vec v}_3-{\vec v}_4\\
  c(X_3)&=2{\vec v}_1+{\vec v}_2-{\vec v}_4 &
  c(X_4)&={\vec v}_1+{\vec v}_2+{\vec v}_3-{\vec v}_4
\end{align*}
Express  the  asymptotic cycle  in  the basis ${\vec v}_j$ in  terms of the
lengths $\lambda_j$ and then in terms of the angle by which the
regular octagon is turned with respect to the standard
presentation. Locally the vertical foliation goes to the North.
And globally?
\end{Exercise}

\subsection{Time Acceleration Machine (Renormalization):
Conceptual Description}
\label{zorich:ss:Time:acceleration:machine:Renormalization}

\index{Renormalization|(}

In the previous section we  have  seen why all trajectories of  a
typical  directional  flow wind around the surface following  the
same asymptotic cycle.  We have also  found an effective  way  to
evaluate  this  asymptotic   cycle:  we  have  seen  that  it  is
sufficient to find an interval exchange transformation $T:X\to X$
induced on any transverse segment $X$ as the first return  map of
the directional flow, and then  to  determine  the ``first return
cycles''  $c(X_j)$,  see Fig.~\ref{zorich:fig:first:return:cycles}.  The
linear combination  of  the  cycles  $c(X_j)$  taken with weights
proportional  to  the  lengths $\lambda_j=|X_j|$ of  subintervals
gives              the              asymptotic             cycle,
see~\eqref{zorich:eq:formula:for:asymptotic:cycle}.

Let us proceed now with a more delicate question of  deviation of
a  trajectory  of  the  directional  flow   from  the  asymptotic
direction.  Without  loss of generality we may  assume  that  the
directional flow under consideration is  the  vertical  flow.  We
know that a  very long cycle  $c(x,N)$ corresponding to  a  large
number $N$ of returns of the trajectory to the horizontal segment
$X$ stretches in the  direction  approaching the direction of the
asymptotic cycle $c$. We  want  to describe how $c(x,N)$ deviates
from               this               direction              (see
Sec.~\ref{zorich:ss:Deviation:from:Asymptotic:Cycle}                 and
Sec.~\ref{zorich:ss:Asymptotic:Flag:and:Dynamical:Hodge:Decomposition}).

We  have  already  seen that as  soon  as  we  have evaluated the
``first  return  cycles'' $c(X_j)$,  complete  information  about
cycles  representing   long   pieces   of   trajectories  of  the
directional flow is  encoded  in the corresponding trajectory $x,
T(x), \dots, T^{N-1}(x)$ of the interval exchange transformation;
see~\eqref{zorich:eq:c:N:equals:sum:c:i}.

An  interval  exchange transformations  gives  an  example  of  a
parabolic dynamical  system  which  is neither completely regular
(like rotation of a circle) nor completely chaotic (like geodesic
flow on a  compact manifold of constant negative curvature, which
is  a typical  example  of a {\it  hyperbolic}  system). In  some
aspects  interval   exchange  transformations  are  closer  to
rotations of a  circle:  say,  as it was proved  by  A.~Katok,  an
interval       exchange       transformation       is       never
mixing~\cite{zorich:Katok:not:mixing} (though, as it  was  very recently
proved   by   A.~Avila    and   G.~Forni   in~\cite{zorich:Avila:Forni},
generically  it  is weakly  mixing).  However,  the  behavior  of
deviation  from  the  ergodic  mean resembles the behavior  of  a
chaotic system.

Our  principal   tool   in   the   study   of  interval  exchange
transformations exploits certain self-similarity  of  these maps.
Choosing a shorter  horizontal interval $X'$ we make the vertical
flow  wind  for  a long time  before  the  first  return to $X'$.
However, the  new first  return map in a sense  would not be more
complicated than  the initial one:  it would be again an interval
exchange transformation $T':X'\to  X'$ of the same (or almost the
same) number of subintervals.

To  check  the latter statement let  us  study the nature of  the
points of  discontinuity of the first  return map $T:X\to  X$. An
interior point $x\in X$ is a point of discontinuity either if the
forward vertical  trajectory of $x$  hits one of the endpoints of
$X$  (as  on  Fig.~\ref{zorich:fig:iet})  or  if  the  forward  vertical
trajectory of  $x$ hits the conical  point before coming  back to
$X$ (see Fig.~\ref{zorich:fig:iet:octagon}). A conical  point having the
cone  angle   $2\pi(d+1)$     has   $d+1$   incoming  vertical
trajectories which land to this conical  point.  (Say,  the  flat
surface represented  on  Fig.~\ref{zorich:fig:iet:octagon}  has a single
conical  point  with the cone angle $6\pi$;  hence  this  conical
point has $3$ incoming vertical trajectories.)  Following them at
the backward direction till the  first  intersection  with $X$ we
find   $d+1$    points    of    discontinuity    on    $X$   (see
Fig.~\ref{zorich:fig:iet:octagon}).  Thus,  all   conical  points  taken
together produce $\sum_j (d_j+1)$ points of discontinuity on $X$.

Generically  two  more points  of  discontinuity  come  from  the
backward trajectories of  the endpoints of $X$. However, in order
to get as small number of discontinuity points as possible we can
choose $X$ in such way that either backward or forward trajectory
of  each of  the  two endpoints hits  some  conical point  before
coming  back  to  $X$.  This  eliminates   these  two  additional
discontinuity points.

\begin{Convention}
\label{zorich:conv:minimal:noi}
From now on we shall always choose any horizontal subinterval $X$
in such way  that the interval exchange transformation $T:X\to X$
induced by the first  return of the vertical flow to $X$  has the
minimal possible number
$$
\noi=\sum_j (d_j+1)+1=2g+(\text{number of conical points})-1
$$
of subintervals under exchange.
\end{Convention}

In the formula above we  used  the  Gauss--Bonnet formula telling
that $\sum_j d_j = 2g-2$, where $g$ is the genus of the surface.

Following  Convention~\ref{zorich:conv:minimal:noi}   we  shall  usually
place the left  endpoint of the  horizontal interval $X$  at  the
conical  singularity.  This leaves  a  discrete  choice  for  the
position of the right endpoint.

\paragraph{Renormalization}

We apply the following strategy in our study  of cycles $c(x,N)$.
Choose     some     horizontal     segment     $X$     satisfying
Convention~\ref{zorich:conv:minimal:noi}.       Consider        vertical
trajectories, which hit  conical  points. Follow them in backward
direction till  the  first  intersection  with  $X$. Consider the
resulting  decomposition $X=X_1\sqcup\dots  \sqcup  X_\noi$,  the
corresponding interval exchange transformation $T:X\to X$ and the
``first return cycles'' $c(X_j)$.

Consider   a   smaller   subinterval  $X'\subset  X$   satisfying
Convention~\ref{zorich:conv:minimal:noi}. Apply  the above procedure  to
$X'$;   let   $X'=X'_1\sqcup\dots   \sqcup   X'_\noi$   be    the
corresponding decomposition of $X'$. We  get  a  new partition of
our flat surface into a collection  of  $\noi$  rectangles  based
over subintervals $X'_1\sqcup \dots \sqcup X'_{\noi}$.

By  construction  the  vertical  trajectories of any  two  points
$x_0,x\in X'_k$ follow  the same high and narrow rectangle $R'_k$
of  the  new  building up to  their  first  return  to $X'$. This
implies  that  the  corresponding  new  ``first  return  cycles''
$c'(x_0)=c'(x)$ are the same and equal to $c'(X'_k)$.

Both vertical  trajectories  of  $x_0,  x\in  X'_k$ intersect the
initial  interval  $X$ many times before first  return  to  $X'$.
However, since these  trajectories  stay together, they visit the
same intervals $X_{j_k}$ in the same order $j_0, j_1, \dots, j_l$
(the   length   $l=l(k)$  of  this  trajectory  depends  on   the
subinterval $X'_k$).

This means  that  we  can  construct  an  $\noi\times\noi$-matrix
$\Cocycle_{jk}$ indicating  how  many times a vertical trajectory
emitted from a point  $x\in  X'_k$ have visited subinterval $X_j$
before  the  first return to $X'$. (By  convention  the  starting
point  counts,  while  the  first  return  point does not.)  Here
$X=X_1\sqcup\dots \sqcup X_\noi$ is the partition  of the initial
``long'' horizontal interval $X$  and  $X'=X'_1\sqcup\dots \sqcup
X'_\noi$ is the partition of the new ``short'' subinterval $X'$.

Having computed this  integer  matrix $\Cocycle$ we can represent
new ``first return cycles'' $c'(X'_k)$  in  terms  of the initial
``first return cycles'' $c(X_j)$ as
\begin{equation}
\label{zorich:eq:cprime:equals:A:c} c'(X'_k)=\Cocycle_{1 k} c(X_1)+
\dots + \Cocycle_{\noi k} c(X_\noi)
\end{equation}
Moreover, it is  easy to see that the lengths $\lambda'_k=|X'_k|$
of subintervals  of the new  partition are related to the lengths
$|X_j|$ of subintervals of the initial  partition  by  a  similar
relation
\begin{equation}
\label{zorich:eq:lambda:prime:equals:A:lambda} \lambda_j=\Cocycle_{j 1}
\lambda'_1 + \dots + \Cocycle_{j \noi} \lambda'_\noi.
\end{equation}
Note that to evaluate matrix $\Cocycle$ we, actually, do not need
to use the vertical flow: the  matrix  $\Cocycle$  is  completely
determined  by  the  initial  interval  exchange   transformation
$T:X\to X$ and by the position of the subinterval $X'\subset X$.

What we gain with this construction is the following. To consider
a  cycle  $c(x,N)$ representing  a  long  piece  of  leaf  of the
vertical  foliation  we followed the trajectory $x, T(x),  \dots,
T^N(x)$ of  the  initial interval exchange transformation $T:X\to
X$ and applied  formula~\eqref{zorich:eq:c:N:equals:sum:c:i}. Passing to
a shorter horizontal interval $X'\subset X$  we  can  follow  the
trajectory  $x,  T'(x), \dots, (T')^{N'}(x)$ of the new  interval
exchange transformation $T':X'\to X'$ (provided $x\in X'$). Since
the  subinterval  $X'$ is  much  shorter than  $X$  we cover  the
initial piece of  trajectory of the  vertical flow in  a  smaller
number $N'$ of steps. In other words, passing from $T$ to $T'$ we
accelerate the time: it is  easy  to see that the trajectory  $x,
T'(x), \dots,  (T')^{N'}(x)$  follows  the  trajectory  $x, T(x),
\dots, T^N(x)$ but jumps over many iterations of $T$ at a time.

Of course  this approach would be non efficient  if the new first
return map $T':X'\to X'$ would be much more  complicated than the
initial one.  But we  know that passing from $T$  to $T'$ we stay
within a family of interval exchange transformations of the fixed
number  $\noi$  of subintervals,  and,  moreover,  that  the  new
``first return cycles'' and  the  lengths of the new subintervals
are  expressed  in  terms of the  initial  ones  by  means of the
$\noi\times\noi$-matrix  $\Cocycle$,  which  depends only on  the
choice of $X'\subset X$ and which can be easily computed.

Our  strategy  can be  formalized  as follows.  In  the next  two
sections we  describe  a  simple explicit algorithm (generalizing
Euclidean   algorithm)  called   Rauzy--Veech   induction   which
canonically  associates  to  an interval exchange  transformation
$T:X\to X$ some  specific subinterval $X'\subset X$ and, hence, a
new  interval   exchange   transformation   $T':X'\to  X'$.  This
algorithm can be considered  as a map from the {\it space  of all
interval  exchange  transformations}\index{Interval exchange transformation!space of}
of a given number $\noi$  of
subintervals to  itself.  Applying  recursively this algorithm we
construct a  sequence  of  subintervals $X=X^{(0)}\supset X^{(1)}
\supset  X^{(2)}  \supset  \dots  $  and  a sequence of  matrices
$\Cocycle=\Cocycle(X^{(0)}), \Cocycle(X^{(1)}),  \dots $ describing
transitions     form     interval     exchange     transformation
$T^{(r)}:X^{(r)} \to X^{(r)}$ to interval exchange transformation
$T^{(r+1)}:X^{(r+1)}       \to       X^{(r+1)}$.        Rewriting
equations~\eqref{zorich:eq:cprime:equals:A:c}
and~\eqref{zorich:eq:lambda:prime:equals:A:lambda} in  a matrix form  we
get:
\begin{align}
\begin{pmatrix}
c(X^{(r+1)}_1)\\
\dots\\
c(X^{(r+1)}_\noi)
\end{pmatrix}
&=
\Bigg(
\Cocycle(X^{(r)})
\Bigg)^T
\cdot
\begin{pmatrix}
c(X^{(r)}_1)\\
\dots\\
c(X^{(r)}_\noi)
\end{pmatrix}
\notag\\
\label{zorich:eq:cocycle:matrix:form}
\\ 
\begin{pmatrix}
\lambda_1(X^{(r+1)})\\
\dots\\
\lambda_\noi(X^{(r+1)})
\end{pmatrix}
&=
\Bigg(
\Cocycle(X^{(r)})
\Bigg)^{-1}
\cdot
\begin{pmatrix}
\lambda_1(X^{(r)})\\
\dots\\
\lambda_\noi(X^{(r)})
\end{pmatrix}
\notag
\end{align}

Taking    a    product     $\Cocycle^{(s)}=\Cocycle(X^{(0)})\cdot
\Cocycle(X^{(1)})\cdot  \dots\cdot  \Cocycle(X^{(s-1)})$  we  can
immediately express the ``first return cycles''  to a microscopic
subinterval $X^{(s)}$  in  terms  of  the  initial ``first return
cycles''   to    $X$    by    a   linear   expression   analogous
to~\eqref{zorich:eq:cocycle:matrix:form}.  Note,  however,  that  before
coming  back   to  this  microscopic  subinterval  $X^{(s)}$  the
vertical flow  has to travel  for enormously long time. The first
return cycle to this very short  subinterval $X^{(s)}$ represents
the cycle  $c(x,N)$  corresponding  to  very  long trajectory $x,
T(x),   ...,   T^N(x)$   of   the   initial   interval   exchange
transformation with  $N\sim\exp(const\cdot  s)$.  In other words,
our renormalization procedure plays a role of a time acceleration
machine: instead of following patiently the  trajectory $x, T(x),
..., T^N(x)$ of the initial interval  exchange transformation for
the exponential  time  $N\sim\exp(const\cdot  s)$  we  obtain the
cycle $c(x,N)$ applying only $s$ steps of renormalization!

One can argue that in this  way we can describe only very special
parts of vertical trajectories: those which start and  end at the
same microscopically small subinterval  $X^{(s)}\subset  X$. This
can  be  overdone   by   the  following  technique.  Consider  an
enormously long trajectory $x, T(x), \dots,  T^N(x)$ which starts
and finishes at some generic points of $X$. One can choose $s(N)$
in such way that the trajectory would get to $X^{(s)}$ relatively
soon (in comparison with its length $N$); then  would return back
to $X^{(s)}$ many times; and would reach the  last point $T^N(x)$
relatively fast after the last visit  to  $X^{(s)}$.  That  means
that essentially (up to negligibly short starting part and ending
part)  one  can  assume  that  the  entire  trajectory  starts at
$X^{(s)}$  and  ends at $X^{(s)}$ (returning to this  subinterval
many times).

This  simple  idea  can  be  developed  and  rigorously  arranged
(see~\cite{zorich:Zorich:How:do} for details). To  avoid  overloading of
this survey with technicalities we  consider  only  a  simplified
problem  giving  a comprehensive description of the first  return
cycles  to  $X^{(s)}$.  The  nature  of  the
asymptotic  flag\index{Flag!asymptotic}\index{Asymptotic!flag}
is especially transparent in this case.

\subsection{Euclidean Algorithm as a Renormalization Procedure in Genus One}
\label{zorich:ss:Euclidean:algorithm}
To  illustrate  the  idea  of renormalization we start  with  the
``elementary'' case, when the Riemann surface  is  a  torus,  the
foliation is a standard irrational  foliation,  and  the  initial
transversal    $X$    is   a   meridian.   We   have   seen    at
Fig.~\ref{zorich:fig:rotation:of:the:circle} that in this case the first
return map $T:X\to X$ is just a rotation of a circle.

\begin{figure}[hbt]
\centering
\includegraphics{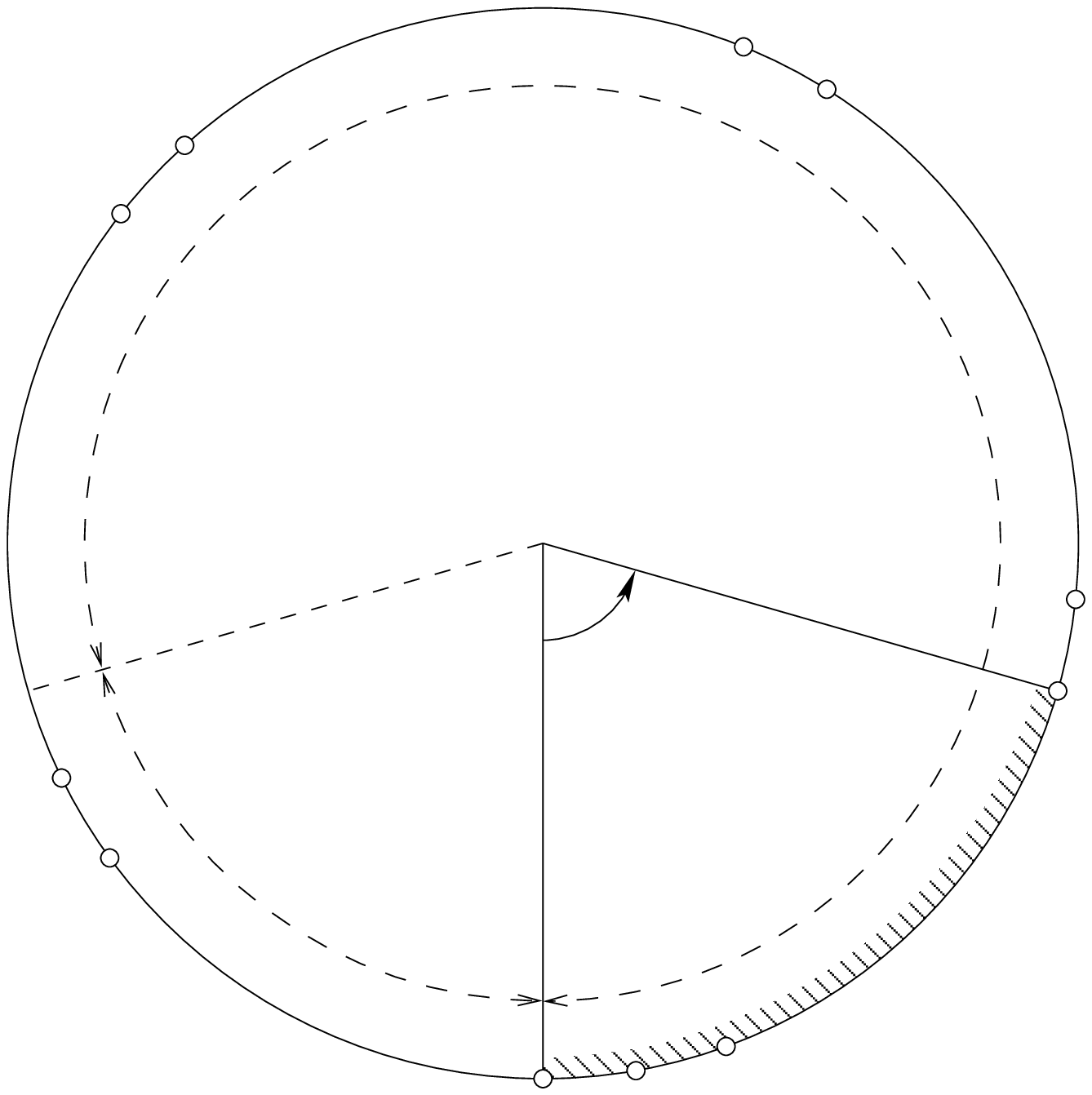}
\includegraphics{zorich_ciseaux6_1cm.eps}
\vspace{205bp} 
\begin{picture}(0,0)(0,0)
\put(94,104){$x$}
\put(92,120){$T^5 x$}
\put(88,138){$T^{10} x$}
\put(24,198){$Tx$}
\put(-2,202){$T^6 x$}
\put(-104,157){$T^2 x$}
\put(-113,140){$T^7 x$}
\put(-95,40){$T^3 x$}
\put(-78,22){$T^8 x$}
\put(35,10){$T^4 x$}
\put(53,20){$T^9 x$}
\put(0,110){$\alpha$}
\put(73,165){$X'$}
\put(-81,90){$X_1$}
\put(50,61){$X_2$}
\end{picture}
\caption{
\label{zorich:fig:circ}
Renormalization  for  rotation  of  a circle leads  to  Euclidean
algorithm and to Gauss measure
}
\end{figure}

Consider  rotation  of  a  circle  $T:S^1\to  S^1$  by  an  angle
$\alpha$. Let the  length  of the  circle  be normalized to  one.
Consider trajectory $x,  Tx,  T^2x, \dots$  of  a point $x$  (see
Fig.~\ref{zorich:fig:circ}).  Denote  the length of the arc $(x,Tx)$  by
$\lambda=\alpha/(2\pi)$.

Cutting the circle at  the point $x$ we get an interval  $X$; the
rotation of the  circle  generates a map of  the  interval $X$ to
itself which we denote by the same symbol $T:X\to X$.  Unbend $X$
isometrically to get a horizontal  interval  of  unit length in such
way that the counterclockwise orientation of the circle gives the
standard positive orientation of the horizontal interval. The
map $T$ acts  on $X$  as follows: it cuts the  unit interval $X$ into
two pieces  $X_1\sqcup  X_2$  of  lengths  $|X_1|=1-\lambda$  and
$|X_2|=\lambda$  and  interchanges  the  pieces  preserving   the
orientation, see Fig.~\ref{zorich:fig:circ}. In other words, the map $T$
is an interval exchange transformation of two subintervals. (To
avoid confusion we stress that $X_2=[T^{-1}x, x]$ and \emph{not}
$X=[x, Tx]$.)

Suppose now that we are  looking  at $X$ in the microscope  which
shows only  the  subinterval  $X'=[x,Tx[$  (corresponding  to the
sector of  angle  $\alpha$  at Fig.~\ref{zorich:fig:circ}). Consider the
trajectory $x, Tx,  T^2x,  \dots$ of the point  $x$  which is the
left extremity of $X_1$. For the  particular rotation represented
at Fig.~\ref{zorich:fig:circ}  the points $Tx, T^2 x, T^3  x, T^4 x$ are
outside  of  the  sector of our  vision;  the  next  point of the
trajectory which we see in $X'$ is the point $T^5 x$. This is the
first  return  point  $T'  x=T^5  x$  to  the  subinterval  $X'$.
Following the trajectory $T^5 x, T^6 x, ...$ further on  we would
not   see   several   more   points   and   then  we  shall   see
$T^{10}x=T'(T'x)$. This is  the  second return to the subinterval
$X'$.

Note that the distance between $x$ and $T'x=T^5 x$ is the same as
the distance between  $T' x=T^5 x$  and $T'(T' x)=T^{10}  x$;  it
equals  $(1-\{1/\lambda\})\cdot  \lambda$, where  $\{\ \}$ denotes
the fractional part of a real  number. It is easy to see that
$T':X'\to X'$ is again an interval exchange transformation of \emph{two}
subintervals $X'_1\sqcup  X'_2$.  The lengths of subintervals are
$|X'_1|=\{1/\lambda\}\cdot\lambda$  and   $(1-\{1/\lambda\})\cdot
\lambda$. After identification of the endpoints  the segment $X'$
becomes a circle  and  the  map $T'$ becomes a  rotation  of  the
circle  $T'$.  Having  started  with  a   rotation   $T$   in   a
\emph{counterclockwise} direction by  the angle $\alpha=2\pi\cdot
\lambda$ we get a  rotation  $T'$ in a \emph{clockwise} direction
by the angle $\alpha'=2\pi\cdot\left\{\cfrac{1}{\lambda}\right\}$
(please verify).

One should  not think that $T'=T^5$  identically. It is  true for
the points  of  first subinterval, $T'|_{X'_1}=T^5$. However, for
the   points   of   the   second  subinterval  $X'_2$   we   have
$T'|_{X'_2}=T^4$. In other words, for  the  points  of the sector
$\alpha$,  which  are  close  to the extremity $Tx$,  {\it  four}
iterations of $T$ bring them back to the sector. Thus, the matrix
$\Cocycle(X')$ of number of visits to subintervals has the form
$$
\begin{pmatrix}
\Cocycle_{11} & \Cocycle_{12}\\
\Cocycle_{21} & \Cocycle_{22}
\end{pmatrix}
=
\begin{pmatrix}
4 & 3\\
1 & 1
\end{pmatrix}
$$
(please draw $X_1\sqcup  X_2$  and $X'_1\sqcup X'_2$ and verify).
We  remind  that $\Cocycle_{jk}$  indicates   how  many  times  a
vertical trajectory emitted from a point $x\in X'_k$ have visited
subinterval  $X_j$  before  the  first return to $X'$,  where  by
convention the  starting  point  counts,  while  the first return
point does not.

Thus  we  get  a  renormalization procedure as described  in  the
previous section: confine the map  $T$  to  a smaller subinterval
$X'$; consider the  resulting first return map $T'$; rescale $X'$
to have unit length. Having started  with  an  interval  exchange
transformation    $T$    of    two    intervals    of     lengths
$(1-\lambda,\lambda)$,  where  $\lambda\in  (0,1)$ we get  (after
rescaling)  an  interval  exchange  transformation  $T'$  of  two
intervals  of  lengths  $\{1/\lambda\},1-\{1/\lambda\}$.  Or,  in
terms of rotations, having started with a \emph{counterclockwise}
rotation by the angle $\alpha=2\pi\lambda$  we  get  a  \emph{clockwise}
rotation                by                the               angle
$\alpha'=2\pi\cdot\left\{\cfrac{1}{\lambda}\right\}$.

One can  recognize  Euclidean  algorithm  in  our renormalization
procedure. Consider the ``space of rotations'', where rotations
are parametrized by the angle $2\pi\lambda$, $\lambda\in [0;1[$. The map
\begin{equation}
\label{zorich:eq:g}
g:\lambda\mapsto\left\{\cfrac{1}{\lambda}\right\}
\end{equation}
can be considered  as  a map from ``the  space  of rotations'' to
itself, or what is  the same, a map from ``the space  of interval
exchange transformations of two subintervals''\index{Interval exchange transformation!space of}
to itself. The map
$g$ is ergodic with respect to the invariant probability measure
\begin{equation}
\label{zorich:eq:dmu}
d\mu=\cfrac{1}{\log 2}\cdot\cfrac{d\lambda}{(\lambda+1)}
\end{equation}
on the parameter space $\lambda\in[0;1[$
which is  called the {\it  Gauss measure}. This map is intimately
related  with  the  development  of  $\lambda$   into a continued
fraction\index{Continued fraction}
$$
\lambda=\cfrac{1}{n_1+\cfrac{1}{n_2+\cfrac{1}{n_3 +
\cfrac{1}{\dots}}}}
$$

We  shall   see  another  renormalization  procedure  related  to
map~\eqref{zorich:eq:g}  in   the   next  sections,  in  particular,  in
Sec.~\ref{zorich:ss:Encoding:a:Continued:Fraction}.

\subsection{Rauzy--Veech Induction}
\label{zorich:ss:Rauzy:Veech:induction}

\index{Rauzy--Veech induction}

In  the   previous  section  we  have   seen  an  example   of  a
renormalization procedure  for interval exchange  transformations
of  two  intervals.   In  this  section  we  consider  a  similar
renormalization procedure  which now works for interval exchanges
of any number of subintervals.  As  we have seen in the  previous
section,  we  do not  need  to keep  information  about the  flat
surface to describe the renormalization algorithm.  Nevertheless,
we prefer to keep  track  of zippered rectangles decomposition of
the  surface  corresponding  to  the sequence of  the  horizontal
subintervals $X=X^{(0)}\supset X^{(1)} \supset \dots$ in order to
preserve geometric spirit of the algorithm.

\index{Renormalization|)}

Consider a flat surface $S$; choose  a  horizontal  interval  $X$
satisfying   Convention~\ref{zorich:conv:minimal:noi};   consider    the
corresponding  decomposition   of   the   surface  into  \emph{zippered
rectangles}\index{Zippered rectangle}
as      on     Fig.~\ref{zorich:fig:iet:octagon}.     Let
$X_1\sqcup\dots\sqcup X_\noi$ be  the corresponding decomposition
of the horizontal  segment $X$ in the base; let $\lambda_j=|X_j|$
denote the widths of subintervals.

\begin{Convention}
\label{zorich:cv:pi:and:pi:inv}
We associate to a decomposition of a flat surface into rectangles
a permutation $\pi$ in such way that the  top horizontal segments
of the rectangles are glued  to  the bottom side of the  interval
$X$ in the order $\pi^{-1}(1), \dots, \pi^{-1}(\noi)$.
\end{Convention}

\begin{Example}
In the  example  presented at Fig.~\ref{zorich:fig:iet:octagon} the four
rectangles $R_1, \dots, R_4$ appear at the bottom side of  $X$ in
the  order  $R_3,R_1,R_4,R_2$,  so  we associate to this  way  of
gluing a permutation
$$
\begin{pmatrix}
1 & 2 & 3 & 4 \\
3 & 1 & 4 & 2
\end{pmatrix}
   =
(3,1,4,2)^{-1}=(2,4,1,3)=\pi
$$
\end{Example}

\begin{Exercise}
Show that intersection indices $c(X_i)\circ c(X_j)$ of the
``first return cycles'' (see
Sec.~\ref{zorich:ss:Evaluation:of:the:asymptotic:cycle:using:an:iet})
are given by the following skew-symmetric matrix $\Omega(\pi)$
defined by the permutation $\pi$:
\begin{equation}
\label{zorich:eq:intersection:form}
\Omega_{ij}(\pi)=
\begin{cases}
1 &\text{if }i<j\text{ and }\pi^{-1}(i)>\pi^{-1}j\\
-1&\text{if }i>j\text{ and }\pi^{-1}(i)<\pi^{-1}j\\
0 &\text{otherwise}
\end{cases}
\end{equation}
Evaluate $\Omega(\pi)$ for the permutation in the Example above
and compare the result with a direct calculation for the cycles
$c(X_j)$, $j=1,\dots,4$, computed in the Exercise at the end of
Sec.~\ref{zorich:ss:Evaluation:of:the:asymptotic:cycle:using:an:iet}.
\end{Exercise}

Compare now the  width  $\lambda_\noi$ of the rightmost rectangle
$R_\noi$  with  the  width   $\lambda_{\pi^{-1}(\noi)}$   of  the
rectangle which is glued to the rightmost position  at the bottom
of  $X$.  As   a  new  subinterval  $X'\subset  X$  consider  the
subinterval $X'$,  which has the  same left extremity as $X$, but
which        is         shorter        than        $X$         by
$\min(\lambda_\noi,\lambda_{\pi^{-1}(\noi)})$.

\begin{figure}[htb]
%
%
\centering
\includegraphics{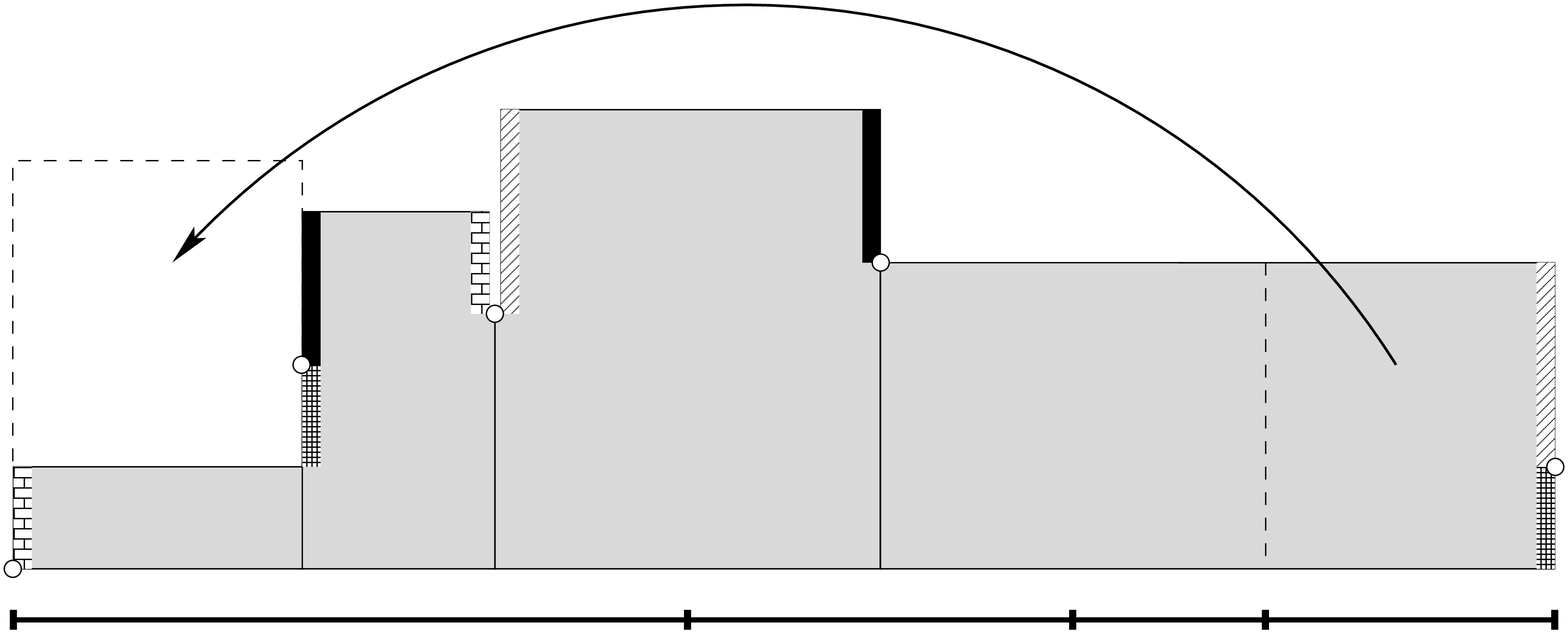}
\includegraphics{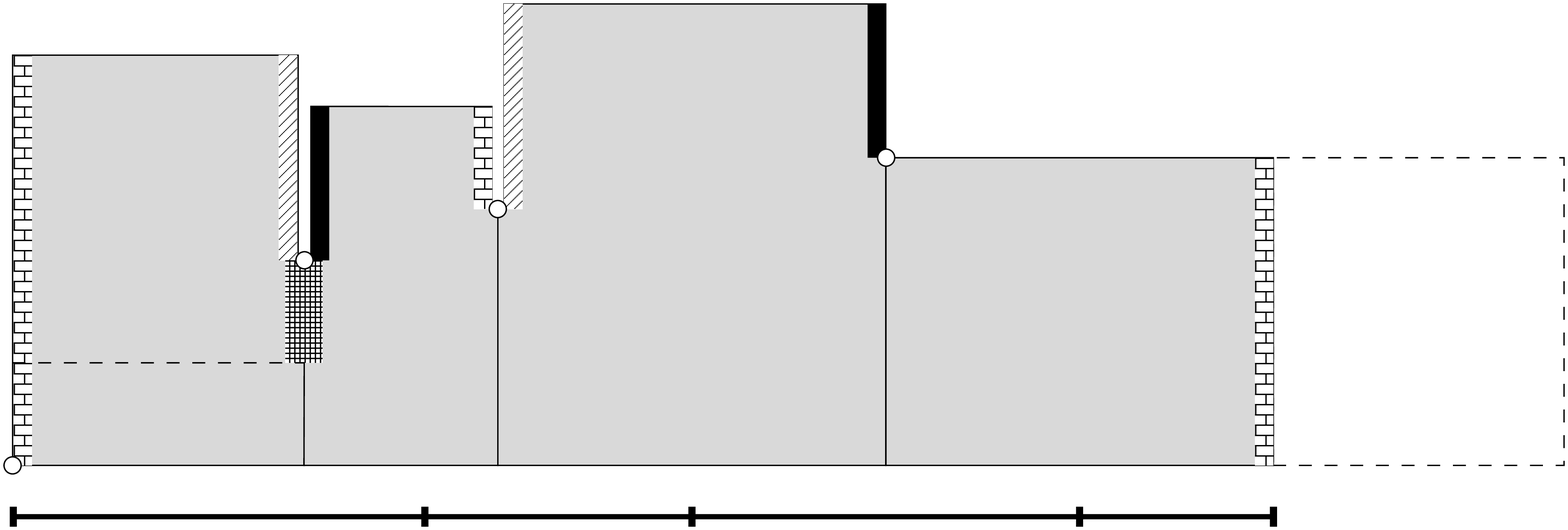}
\begin{picture}(0,0)(0,0)
\put(100,-195){
\begin{picture}(0,0)(0,0)
\put(-245,80){$R_1$}
\put(-195,80){$R_2$}
\put(-135,80){$R_3$}
\put(-10,80){$R_4$}
\end{picture}}
\put(100,-160){
\begin{picture}(0,0)(0,0)
\put(-210,20){$R_4$}
\put(-95,20){$R_3$}
\put(-32,20){$R_2$}
\put(20,20){$R_1$}
\end{picture}}
        %
   %
\put(100,-335){
\begin{picture}(0,0)(0,0)
\put(-245,80){$R'_1$}
\put(-195,80){$R'_2$}
\put(-135,80){$R'_3$}
\put(-55,80){$R'_4$}
\end{picture}}
\put(100,-302){
\begin{picture}(0,0)(0,0)
\put(-235,20){$R'_4$}
\put(-160,20){$R'_1$}
\put(-95,20){$R'_3$}
\put(-32,20){$R'_2$}
\end{picture}}
\end{picture}
\vspace{290bp}
\caption{
\label{zorich:fig:Rauzy:move:b}
Type I modification: the rightmost rectangle $R_4$ on  top of $X$
is  wider  than the rectangle $R_1=R_{\pi^{-1}(4)}$ glued to  the
rightmost position at the bottom of $X$.
}
\end{figure}

\begin{figure}[htb]
%
%
\centering
\includegraphics{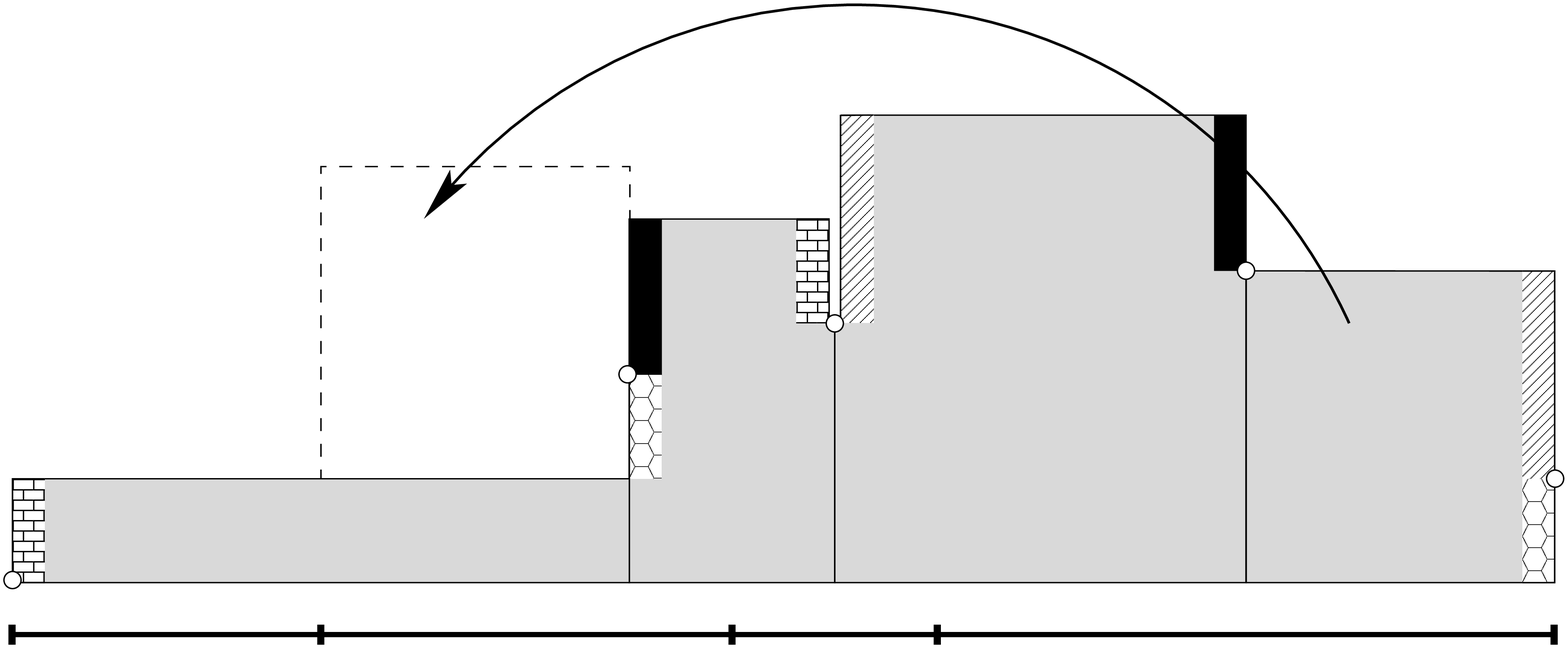}
\includegraphics{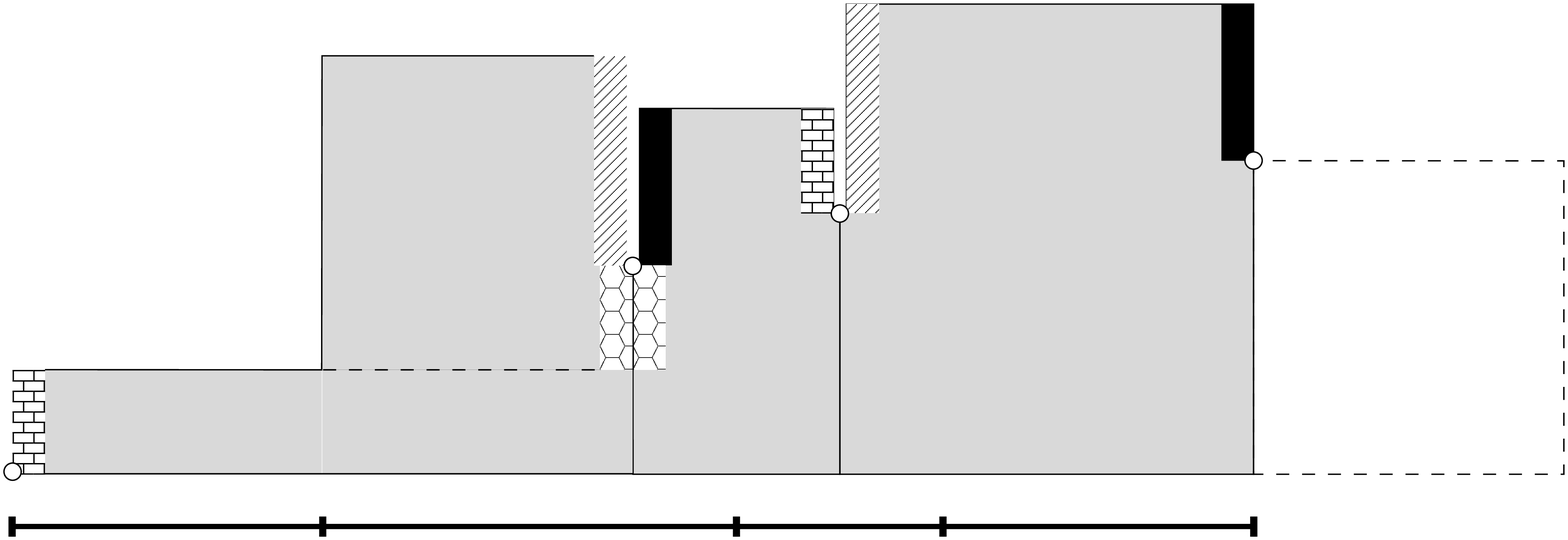}
\begin{picture}(0,0)(0,0)
\put(100,-200){
\begin{picture}(0,0)(0,0)
\put(-210,80){$R_1$}
\put(-125,80){$R_2$}
\put(-60,80){$R_3$}
\put(17,80){$R_4$}
\end{picture}}
\put(100,-165){
\begin{picture}(0,0)(0,0)
\put(-240,20){$R_4$}
\put(-165,20){$R_3$}
\put(-100,20){$R_2$}
\put(-20,20){$R_1$}
\end{picture}}
        %
   %
\put(100,-333){
\begin{picture}(0,0)(0,0)
\put(-245,80){$R'_1$}
\put(-180,80){$R'_2$}
\put(-125,80){$R'_3$}
\put(-60,80){$R'_4$}
\end{picture}}
\put(100,-301){
\begin{picture}(0,0)(0,0)
\put(-240,20){$R'_2$}
\put(-165,20){$R'_4$}
\put(-100,20){$R'_3$}
\put(-50,20){$R'_1$}
\end{picture}}
\end{picture}
\vspace{290bp}
\caption{
\label{zorich:fig:Rauzy:move:a}
Type II modification: the rightmost rectangle $R_4$ on top of
$X$ is narrower than the rectangle $R_1=R_{\pi^{-1}(4)}$ glued
to the rightmost position at the bottom of $X$.
}
\end{figure}

The  situation  when  $\lambda_\noi>\lambda_{\pi^{-1}(\noi)}$  is
represented  at  Fig.~\ref{zorich:fig:Rauzy:move:b}; the  situation when
$\lambda_\noi<\lambda_{\pi^{-1}(\noi)}$    is    represented   at
Fig.~\ref{zorich:fig:Rauzy:move:a}.

By construction  the first return map  $T':X\to X'$ has  the same
number $\noi$  of  subintervals  in  its decomposition. Observing
Fig.~\ref{zorich:fig:Rauzy:move:b} and~\ref{zorich:fig:Rauzy:move:b} one can
see      that       in       the       first      case,      when
$\lambda_\noi>\lambda_{\pi^{-1}(\noi)}$,  the new decomposition
\mbox{$X'_1\sqcup \dots\sqcup X'_\noi$} is obtained from the
original decomposition $X_1\sqcup \dots\sqcup X_\noi$
by  shortening  the  last interval by  $\lambda_{\pi^{-1}(\noi)}$
from    the     right.     In     the     second    case,    when
$\lambda_\noi<\lambda_{\pi^{-1}(\noi)}$,  the  new  decomposition
$X'_1\sqcup \dots \sqcup X'_\noi$ is obtained  from  the
original decomposition $X_1\sqcup \dots\sqcup X_\noi$ by
eliminating  the   last   subinterval $X_\noi$
and  by  partitioning  the
subinterval  $X_{\pi^{-1}(\noi)}$  into  two   ones   of  lengths
$\lambda_{\pi^{-1}(\noi)}-\lambda_\noi$    and     $\lambda_\noi$
correspondingly.

The order  in which the rectangles of the  new building are glued
to  the  bottom  of  the  interval  $X'$   changes.  The  new
permutation $\pi'$  can  be  described  as  follows. Consider the
initial permutation $\pi$ as a pair of orderings of a finite set:
a  ``top''  ordering  $1,2,\dots,\noi$   (corresponding   to  the
ordering  of  the rectangles  along  the  top  side  of  the base
interval $X$)  and  a  ``bottom''  ordering  $\pi^{-1}(1), \dots,
\pi^{-1}(m)$ (corresponding  to  the  ordering  of the rectangles
along the bottom  side  of the base interval  $X$).  In the first
case,  when  $\lambda_\noi>\lambda_{\pi^{-1}(\noi)}$,   the   new
permutation  $\pi'$  corresponds to the modification of the  {\it
bottom} ordering  by  cyclically  moving  one  step forward those
letters  occurring after  the  image of the  last  letter in  the
bottom line,  i.e., after the  letter $\noi$. In the second case,
when $\lambda_\noi<\lambda_{\pi^{-1}(\noi)}$, the new permutation
$\pi'$ corresponds to  the modification of the {\it top} ordering
by cyclically moving one step  forward  those  letters  occurring
after the image of  the last letter in the top line,  i.e., after
the letter $\pi^{-1}(\noi)$.

\begin{Example}
For the initial buildings at both  Figures~\ref{zorich:fig:Rauzy:move:b}
and~\ref{zorich:fig:Rauzy:move:a} the permutation $\pi$ corresponding to
the  initial  interval exchange transformation $T:X\to X$ is  the
same and equals
$$
\pi=
\begin{pmatrix}
1 & 2 &  3 & 4 \\
4 & 3 &  2 & 1
\end{pmatrix}
$$
Our modification produces permutation
$$
\begin{pmatrix}

1 & \ & 2 &     & 3 &     & 4 & \quad \\
4 & \ & 3 & \to & 2 & \to & 1 &
\end{pmatrix}
 =
\begin{pmatrix}
1& 2 &  3& 4\\
4& 1 &  3& 2
\end{pmatrix}=\pi'
$$
\includegraphics{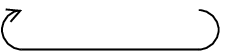}
in the first case (when $\lambda_\noi>\lambda_{\pi^{-1}(\noi)}$)
and permutation
$$
\begin{pmatrix}
1 & \phantom{\to} & 2 & \to & 3 & \to & 4 & \phantom{\to} \\
4 &               & 3 &     & 2 &     & 1 &
\end{pmatrix}
    =
\begin{pmatrix}
1& 4 & 2 & 3 \\
4& 3 & 2 & 1
\end{pmatrix}
    \sim
\begin{pmatrix}
1& 2 & 3 & 4 \\
2& 4 & 3 & 1
\end{pmatrix}
=\pi'
$$
\includegraphics{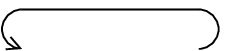}
in the second case (when
$\lambda_\noi<\lambda_{\pi^{-1}(\noi)}$).
\end{Example}

Note      that      in      the      second      case       (when
$\lambda_\noi<\lambda_{\pi^{-1}(\noi)}$)  passing  to   the   new
decomposition $X'_1\sqcup\dots\sqcup X'_\noi$ we  have  to change
the initial enumeration of the subintervals though physically all
subintervals but one stay  unchanged.  Another choice would be to
assign ``names''  to  subintervals  once  and  forever. Under the
first choice the permutations
$$
\begin{pmatrix}
1& 4 & 2 & 3 \\
4& 3 & 2 & 1
\end{pmatrix}
    \text{ and }
\begin{pmatrix}
1& 2 & 3 & 4 \\
2& 4 & 3 & 1
\end{pmatrix}
$$
coincide;  under   the   latter   choice  they  become  different
permutations.   The   article~\cite{zorich:Yoccoz:Les:Houches}  in   the
current volume adopts the second convention.

Similarly to the case of interval exchange transformations of two
intervals  the   induction   procedure  described  above  can  be
described entirely  in terms of interval exchange transformations
$T:X\to X$ and $T':X'\to X'$. Historically  it  was  proposed  by
G.~Razy~\cite{zorich:Rauzy}  in   these   latter   form   and  then  was
interpreted  by   W.~Veech~\cite{zorich:Veech:Annals:82}  in  terms   of
zippered   rectangles\index{Zippered rectangle}.
(Actually,   the   zippered   rectangles
decomposition    has    appeared   and    was    first    studied
in~\cite{zorich:Veech:Annals:82}.)


\subsection{Multiplicative Cocycle on the Space of Interval Exchanges}
\label{zorich:ss:space:of:iet}

\index{Interval exchange transformation!space of|(}

The      renormalization       procedure      constructed      in
Sec.~\ref{zorich:ss:Euclidean:algorithm} gives  a map $g$ from the space
of rotations of a circle to  itself, or, in other terms, from the
space of interval exchange transformations of two subintervals to
itself.  The  permutation  $\pi$  corresponding  to  an  interval
exchange transformations  of  two  intervals  $X_1\sqcup  X_2$ is
always  equal   to   $\pi=(2,1)$,   so   such  interval  exchange
transformation can be parametrized by  a  single  real  parameter
$\lambda\in(0,1)$, where $\lambda=|X_1|$. Here we assume that the
total length $|X_1|+|X_2|=|X|$ of the interval  $X$ is normalized
as $|X|=1$.

An  interval  exchange  transformation  of  $\noi$   subintervals
$X=X_1\sqcup\dots\sqcup X_\noi$ is parametrized  by  a collection
$\lambda_1, \dots, \lambda_\noi$ of positive numbers representing
the   lengths    of    subintervals    and   by   a   permutation
$\pi\in\mathfrak{S}_\noi$. Assuming  that the total length of the
interval $X$ is  normalized as $|X|=1$  we see that the {\it space
of interval exchange transformations}
is parametrized by a finite
collection       of       $(\noi-1)$-dimensional        simplices
$\Delta^{\noi-1}=\{(\lambda_1,    \dots,    \lambda_\noi)\     |\
\lambda_1+\dots+\lambda_\noi=1;   \lambda_j>0\}$,   where    each
simplex corresponds to some fixed permutation $\pi$.

\begin{figure}[htb]
\centering
\includegraphics{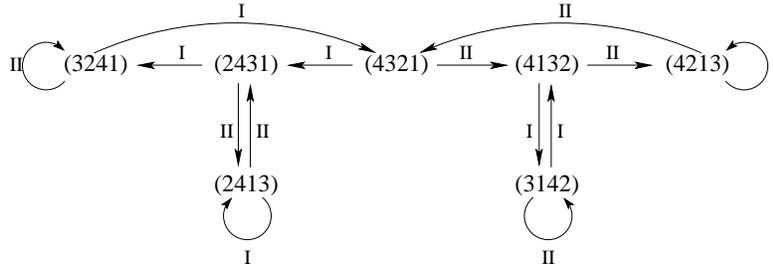}
\vspace{95bp}
\caption{
\label{zorich:fig:Rauzy:class}
Rauzy class of permutation $(4,3,2,1)$
}
\end{figure}

As a collection of permutations one can consider all permutations
obtained  from  a  given  one   by   applying   recursively   the
modifications described at  the end of the previous section. Such
collection  of  permutations  is  called  a
{\it  Rauzy  class}\index{Rauzy class}\index{0R10@$\mathfrak{R}$ -- (extended) Rauzy class}
$\mathfrak{R}\subset\mathfrak{S}_\noi$.
Figure~\ref{zorich:fig:Rauzy:class}  illustrates  the Rauzy  class  of the
permutation $(4,3,2,1)$, where the arrows indicated modifications
of the first and of the second type.

The renormalization  procedure  described in the previous section
(combined with rescaling of the  resulting  interval  $X'$ to the
unit length) defines a map
$$
\cT: \Delta^{\noi-1}\times\mathfrak{R} \to
\Delta^{\noi-1}\times\mathfrak{R}
$$
for  each  space of interval exchange transformations\index{Interval exchange transformation!space of}
to  itself.
The following Theorem of W.~A.~Veech (see~\cite{zorich:Veech:Annals:82})
is crucial in this story.

\begin{KeyTheorem}[W.~A.~Veech]
The map $\cT$ is ergodic  with  respect  to absolutely continuous
invariant measure.
\end{KeyTheorem}

\index{Interval exchange transformation|)}

\begin{NNRemark}
It is easy to see that the
Rauzy class\index{Rauzy class}\index{0R10@$\mathfrak{R}$ -- (extended) Rauzy class}
presented at
Fig.~\ref{zorich:fig:Rauzy:class} does not depend on the starting
permutation. Actually, the same is true for any Rauzy class.
Moreover, for almost any flat surface $S\in\cH^{comp}(d_1,
\dots, d_\noz)$ in any connected component\index{Moduli space!connected components of the strata}
of any stratum a
finite set of all permutations realizable by first return maps
to all possible horizontal segments $X$ satisfying
Convention~\ref{zorich:conv:minimal:noi} does not depend on the surface
$S$. This set is a disjoint union of a finite collection of
corresponding Rauzy classes $\mathfrak{R}_1, \dots ,
\mathfrak{R}_j$, where $j$ is the number of distinct entries
$d_{i_1}<d_{i_2}<\dots<d_{i_j}$. The union
$\mathfrak{R}_{ex}=\mathfrak{R}_1\sqcup\dots\sqcup\mathfrak{R}_j$
is called the
\emph{extended Rauzy class}\index{Rauzy class!extended Rauzy class}\index{0R10@$\mathfrak{R}$ -- (extended) Rauzy class};
it depends only on the connected component $\cH^{comp}(d_1,
\dots, d_\noz)$. In particular, connected components of the
strata are characterized by the extended Rauzy classes, where
the latter ones can be described in purely combinatorial terms.
\end{NNRemark}

With the  Theorem above we have almost  accomplished our scheme
for a renormalization procedure. There is only one trouble with
the map $\cT$:  the  measure mentioned in the Theorem  is
infinite (the total measure of  the  space of interval exchange
transformations is  infinite).  This  technical  problem  can be
fixed  by  the following trick.  We shall  modify  the renormalization
algorithm  described  in  the previous section by making several
modifications of the zippered rectangle\index{Zippered rectangle}
at a time. At  a  single
step of the  new  algorithm $\cG$ we apply  several steps of the
previous one $\cT$. Namely, we keep going as soon as we apply
consecutive transformations  $\cT$  of the same type $I$ or of
the  same  type $II$. Conceptually it  does  not change the
renormalization procedure, but now  the  renormalization
develops faster  than   before.  The  following  example
illustrates  the correspondence between the renormalization
procedures  $\cT$ and $\cG$:
$$
\begin{array}{cccccccccccccc}
(\lambda,\pi)&\xrightarrow{I}&
\cT(\lambda,\pi)&\xrightarrow{I}&
\cT^2(\lambda,\pi)&\xrightarrow{I}&
\cT^3(\lambda,\pi)&\xrightarrow{II}&
\cT^4(\lambda,\pi)&\xrightarrow{II}&
\cT^5(\lambda,\pi)&\xrightarrow{I}&\cdots
\\
||&&&&&&||&&&&||
\\
(\lambda,\pi)&
\multicolumn{5}{c}{\xrightarrow{\hspace*{3.5truecm}}}&
\cG(\lambda,\pi)&
\multicolumn{3}{c}{\xrightarrow{\hspace*{2truecm}}}&
\cG^2(\lambda,\pi)&\to&\cdots
\end{array}
$$

The accelerated procedure $\cG$ was introduced
in~\cite{zorich:Zorich:Gauss:map}, where the following Theorem was
proved.

\begin{NNTheorem}
The map $\cG$ is ergodic  with  respect  to absolutely continuous
invariant  probability   measure   on   each  space  of  zippered
rectangles\index{Zippered rectangle}.
\end{NNTheorem}

Now, when we  have elaborated almost  all necessary tools  we  are
ready  to  give  an  idea  of  the  proof  of  the  Theorem  from
Sec.~\ref{zorich:ss:Asymptotic:Flag:and:Dynamical:Hodge:Decomposition}
concerning
asymptotic flag\index{Flag!asymptotic}\index{Asymptotic!flag}.
The last element which is missing is
some   analysis    of   the   matrices    $\Cocycle(\lambda,\pi)$
in~\eqref{zorich:eq:cocycle:matrix:form}.

\paragraph{Multiplicative Cocycle}

In our interpretation of a renormalization procedure as  a map on
the  space  of  interval  exchange  transformations\index{Interval exchange transformation!space of|)}
we  have  to
consider   matrix   $\Cocycle$  as   a   matrix-valued   function
$\Cocycle(\lambda,\pi)$           on           the          space
$\Delta^{\noi-1}\times\mathfrak{R}$\index{0R10@$\mathfrak{R}$ -- (extended) Rauzy class}
of    interval    exchange
transformations.   Our    goal    (as    it   was   outlined   in
Sec.~\ref{zorich:ss:Time:acceleration:machine:Renormalization})  is   to
describe      the      properties       of      the      products
$\Cocycle^{(s)}=\Cocycle(X^{(0)})\cdot     \Cocycle(X^{(1)})\cdot
\dots\cdot  \Cocycle(X^{(s)})$   of  matrices  corresponding   to
successive  steps  of  renormalization.

We shall keep  the  same
notation  $\Cocycle(\lambda,\pi)$  for matrices  corresponding to
our fast renormalization procedure $\cG(\lambda,\pi)$. Consider the
product $\Cocycle^{(s)}(\lambda,\pi)$ of  values  of
$\Cocycle(\cdot)$   taken   along   the   orbit   $(\lambda,\pi),
\cG(\lambda,\pi), \dots, \cG^s(\lambda,\pi)$ of length $s$ of
the map $\cG$:
$$
\Cocycle^{(s)}(\lambda,\pi)= \Cocycle(\lambda,\pi)\cdot
\Cocycle(\cG(\lambda,\pi))\cdot \dots\cdot
\Cocycle(\cG^s(\lambda,\pi))
$$
Such $\Cocycle^{(s)}(\lambda,\pi)$ is called a {\it
multiplicative cocycle}\index{Cocycle!multiplicative}
over the map $\cG$:
$$
\Cocycle^{(p+q)}(\lambda,\pi)= \Cocycle^{(p)}(\lambda,\pi)\cdot
\Cocycle^{(q)}(\cG^p(\lambda,\pi))
$$
Using ergodicity of $\cG$ we  can  apply  multiplicative  ergodic
theorem\index{Multiplicative ergodic theorem}\index{Ergodic!multiplicative ergodic theorem}\index{Theorem!multiplicative ergodic}
(see Appendix~\ref{zorich:s:Multiplicative:ergodic:theorem})  to
describe   properties    of   $\Cocycle^{(s)}$.   Morally,    the
multiplicative ergodic theorem tells that for large values of $s$
the matrix $\Cocycle^{(s)}(\lambda,\pi)$ should be considered  as
a matrix conjugate  to the $s$-th  power of some  {\it  constant}
matrix. (See Appendix~\ref{zorich:s:Multiplicative:ergodic:theorem}  for
a  rigorous  formulation.) The logarithms of eigenvalues of  this
constant  matrix  are called  {\it  Lyapunov  exponents}  of  the
multiplicative cocycle.

Recall that  matrix  $\Cocycle^T(\lambda,\pi)$  was  defined as
a matrix representing the new  ``first  return cycles'' in terms
of the  old  ones, see~\eqref{zorich:eq:cocycle:matrix:form}. Actually,
it can be also interpreted as  a  matrix representing a change
of  a basis      in       the      first      relative
cohomology
$[\operatorname{Re}\omega]=(\lambda_1,\dots,\lambda_\noi)\in
H^1(S,\{\text{conical  singularities}\};\R{})$.  It is easy to
check that it respects the (degenerate)  symplectic form: the
intersection form~\eqref{zorich:eq:intersection:form}. Note that
symplectic matrices  have certain symmetry of eigenvalues. In
particular, it follows from the general theory that the
corresponding
Lyapunov exponents\index{0theta10@$\theta_1,\dots,\theta_g$ -- Lyapunov exponents of the Rauzy--Veech cocycle}\index{Exponent!Lyapunov exponent}\index{Lyapunov exponent}
have the following symmetry:
$$
\theta_1>\theta_2\ge \theta_3\ge \dots \ge \theta_g\ge
\underbrace{0=0=\dots=0}_{number\ of\ conical\ points\, -1} \ge
-\theta_g \ge \dots \ge -\theta_2
> -\theta_1
$$

Note that the first return cycles actually belong to the absolute
homology   group   $H_1(S;\Z{})\subset  H_1(S;\R{})\simeq\R{2g}$.
Passing to  this  $2g$-dimensional  space  we  get matrices which
already preserve a nondegenerate symplectic form. They define the
following subcollection
$$
\theta_1>\theta_2\ge \theta_3\ge \dots \ge \theta_g\ge
-\theta_g\ge \dots \ge -\theta_3>  -\theta_2 > -\theta_1
$$
of Lyapunov exponents.

The rest is an elementary linear algebra. We want to describe how
do the large powers $s$ of a symplectic matrix with eigenvalues
$$
\exp(\theta_1)>\exp(\theta_2)\ge\dots\ge\exp(\theta_g)
\ge\exp(-\theta_g)\ge\dots\ge\exp(-\theta_2)\ge\exp(-\theta_1)
$$
act on a $2g$-dimensional symplectic space.

We  know  that  the
Lyapunov   exponent\index{0theta10@$\theta_1,\dots,\theta_g$ -- Lyapunov exponents of the Rauzy--Veech cocycle}\index{Exponent!Lyapunov exponent}\index{Lyapunov exponent}
$\theta_g\ge   0$   is
nonnegative.  Assume\footnote{Actually,  this   assumption  is  a
highly nontrivial Theorem of G.~Forni~\cite{zorich:Forni:02}; see below,
see also Sec.~\ref{zorich:ss:Spectrum:of:Lyapunov:exponents}}
that it is actually strictly  positive:  $\theta_g>0$.  Then  for
half of dimensions our linear  map  is expanding and for half  of
dimensions it is contracting. In particular, under the assumption
that   $\theta_g>0$    we    conclude   that   the   linear   map
$\Cocycle^{(s)}(\lambda,\pi)$ projects  all  homology  space to a
Lagrangian  subspace  (spanned  by eigenvectors corresponding  to
positive Lyapunov exponents).

Assuming that the  spectrum of
Lyapunov exponents\index{0theta10@$\theta_1,\dots,\theta_g$ -- Lyapunov exponents of the Rauzy--Veech cocycle}\index{Exponent!Lyapunov exponent}\index{Lyapunov exponent}
is simple, that is assuming that
$$
\theta_1>\theta_2> \theta_3> \dots > \theta_g
$$
we get  the entire picture of deviation. A  generic vector in the
homology space stretches along the principal eigenvector (the one
corresponding  to   the  eigenvalue  $\theta_1$)  with  a  factor
$\exp(s\theta_1)$; it expands  along  the next eigenvector with a
factor $\exp(s\theta_2)$, etc, up to the order $g$; its deviation
from  the  Lagrangian  subspace  spanned   by   the   first   $g$
eigenvectors tends to zero. Hence, the norm $l$ of the image of a
generic vector under  $s$-th  power of our linear  map  is of the
order $l\sim  \exp(s\theta_1)$;  its  deviation from the subspace
$\cV_1$, which is spanned by the top eigenvector, is of the order
$\exp(s\theta_2)=l^\frac{\theta_2}{\theta_1}$; its deviation from
the subspace  $\cV_2$ spanned by the  two top eigenvectors  is of
the  order  $\exp(s\theta_3)=l^\frac{\theta_3}{\theta_1}$,   etc;
there is no deviation from the Lagrangian subspace spanned by the
top  $g$   eigenvectors.  In  particular,  the  exponent
$\nu_j$\index{0nu10@$\nu_1,\dots,\nu_g$ -- Lyapunov exponents related to the Teichm\"uller geodesic flow}\index{Exponent!Lyapunov exponent}\index{Lyapunov exponent}
responsible for deviation from the subspace  $\cV_{j-1}$ from the
Theorem                                                        in
Sec.~\ref{zorich:ss:Asymptotic:Flag:and:Dynamical:Hodge:Decomposition}
is obtained by normalization of the  Lyapunov exponent $\theta_j$
by the leading Lyapunov exponent\index{0theta10@$\theta_1,\dots,\theta_g$ -- Lyapunov exponents of the Rauzy--Veech cocycle}\index{Exponent!Lyapunov exponent}\index{Lyapunov exponent}
$\theta_1$:
\begin{equation}
\label{zorich:eq:nu:j:as:theta:j:over:theta:1}
\nu_j=\cfrac{\theta_j}{\theta_1}.
\end{equation}

This  completes  the proof of the  Theorem  in the case when  the
vertical  trajectory  starts  and  ends at the  same  microscopic
horizontal interval (in other words, when the piece of trajectory
is  ``almost  closed'').  Applying  some  additional  (relatively
involved) ergodic machinery one can complete  the  proof  of  the
Theorem for  arbitrary  long  pieces  of  vertical  trajectories;
see~\cite{zorich:Zorich:How:do} for a complete proof.

I  would  like   to  stress  that  the  original  Theorem  proved
in~\cite{zorich:Zorich:How:do}  is  conditional:   the  statement  about
Lagrangian   subspace   was   proved   modulo   conjecture   that
$\theta_g>0$; the statement about a complete
Lagrangian flag\index{Flag!Lagrangian}\index{Lagrangian!flag}
was proved modulo conjectural simplicity of the  spectrum of
Lyapunov\index{0theta10@$\theta_1,\dots,\theta_g$ -- Lyapunov exponents of the Rauzy--Veech cocycle}\index{Exponent!Lyapunov exponent}\index{Lyapunov exponent}
exponents.

Positivity     $\theta_g>0$    was     proved     by     G.~Forni
in~\cite{zorich:Forni:02}, and  simplicity  of the spectrum was recently
proved  by  A.~Avila  and M.~Viana~\cite{zorich:Avila:Viana}; see  more
details  in  Sec.~\ref{zorich:ss:Spectrum:of:Lyapunov:exponents}. As  it
was shown  in~\cite{zorich:Zorich:Gauss:map}  the  proof  of  the strict
inequality
$\theta_1>\theta_2$\index{0theta10@$\theta_1,\dots,\theta_g$ -- Lyapunov exponents of the Rauzy--Veech cocycle}\index{Exponent!Lyapunov exponent}\index{Lyapunov exponent}
immediately follows from  results
of W.~A.~Veech.

\begin{Exercise}
In      analogy       with       what       was      done      in
Sec.~\ref{zorich:ss:Euclidean:algorithm}   consider   the   Rauzy--Veech
induction  in  the torus case applying it  to  interval  exchange
transformations of two subintervals. We  have  seen  that in this
case the space of interval  exchange  transformations\index{Interval exchange transformation!space of}
is just an
interval $(0,1)$. Find an explicit formulae  for the Rauzy--Veech
renormalization map $\cT: (0,1)  \to  (0,1)$ and for the ``fast''
renormalization  map  $\cG:  (0,1)  \to (0,1)$. Explain  why  the
invariant measure is infinite for the map $\cT$.  Find a relation
between        $\cG$       and        the        Gauss        map
$g:x\mapsto\left\{\cfrac{1}{x}\right\}$. Let
$$
\cfrac{p_s}{q_s}=
\cfrac{1}{n_1+\cfrac{1}{n_2+\cfrac{1}{\dots +
\cfrac{1}{n_s}}}}
$$
be  the  $s$-th  best  rational  approximation  of a real  number
$x\in(0,1)$. In the torus case the spectrum of Lyapunov exponents
reduces to a single pair $\theta_1>-\theta_1$. Show that for
almost all  $x\in(0,1)$  the Lyapunov exponent
$\theta_1$\index{0theta10@$\theta_1,\dots,\theta_g$ -- Lyapunov exponents of the Rauzy--Veech cocycle}\index{Exponent!Lyapunov exponent}\index{Lyapunov exponent}
(called
in number theory  the  {\it L\'evy constant}, see~\cite{zorich:Levy}) is
responsible  for  the  growth  rate  of  the  denominator  of the
continued fraction\index{Continued fraction} expansion of $x$:
$$
\lim_{s\to\infty} \cfrac{\log q_s}{s} =
\theta_1=\cfrac{\pi^2}{12\log 2}
$$
\end{Exercise}

\subsection{Space  of   Zippered  Rectangles  and  Teichm\"uller geodesic flow}
\label{zorich:ss:Teichmuller:Geodesic:Flow}

We have proved the Theorem about an asymptotic Lagrangian flag of
subspaces responsible for  deviation  of the cycles $c(x,N)$ from
asymptotic direction.  We  have  also  proved  that the exponents
$\nu_j$ responsible  for  the  quantitative  description  of  the
deviation are expressed in terms of the Lyapunov exponents of the
multiplicative  cocycle  corresponding   to  our  renormalization
procedure: $\nu_j=\theta_j/\theta_1$.

There  remains  a  natural  question  why  should we choose  this
particular renormalization procedure and not a different one. One
more  natural  question  is  what   is   the   relation   between
renormalization procedure and the  flow  induced by the action of
the  diagonal   subgroup   $\begin{pmatrix}   \exp(t)  &  0\\0  &
\exp(-t)\end{pmatrix}$ on the space of  flat  surfaces  which  we
agreed to  call  the
\emph{Teichm\"uller  geodesic  flow}\index{Teichm\"uller!geodesic flow|(};
this
relation              was              announced               in
Sec.~\ref{zorich:ss:Asymptotic:Flag:and:Dynamical:Hodge:Decomposition}.
This section  answers  to  these  questions  which are, actually,
closely related.

In    our    presentation    we     follow     the    fundamental
paper~\cite{zorich:Veech:Annals:82} of W.~A.~Veech; the material at  the
end of the section is based  on the paper~\cite{zorich:Zorich:Gauss:map}
developing the initial paper~\cite{zorich:Veech:Annals:82}.

\paragraph{Space of Zippered Rectangles}

\index{Zippered rectangle|(}

We have seen that locally a flat surface $S$ can  be parametrized
by a collection of relative periods  of  the  holomorphic  1-form
$\omega$  representing  the flat surface, i.e. we  can  choose  a
small domain  containing  $[\omega]$  in  the relative cohomology
group $H^1(S,\{P_1,\dots,P_\noz\};\C{})$ as a coordinate chart in
the corresponding stratum $\cH(d_1,\dots,d_\noz)$\index{0H20@$\cH(d_1,\dots,d_\noz)$ -- stratum in the moduli space}\index{Stratum!in the moduli space}.

Decomposition  of  a flat surface into zippered rectangles gives
another  system  of  local  coordinates in the  stratum.  Namely,
choose      a      horizontal     segment      $X$     satisfying
Convention~\ref{zorich:conv:minimal:noi}                            from
Sec.~\ref{zorich:ss:Time:acceleration:machine:Renormalization}       and
consider the  corresponding  decomposition  of  $S$ into zippered
rectangles. Let $\lambda_1, \dots, \lambda_\noi$ be the widths of
the rectangles, $h_1, \dots, h_\noi$  their  heights,  and  $a_1,
\dots, a_\noz$ be the  altitudes  responsible for the position of
singularities  (we  zip  the  neighboring  rectangles  $R_j$  and
$R_{j+1}$ from  the bottom up to the altitude  $a_j$ and then the
rectangles      split      at      the      singularity,      see
Figures~\ref{zorich:fig:iet:octagon},            \ref{zorich:fig:Rauzy:move:b},
\ref{zorich:fig:Rauzy:move:a}). There is one more parameter describing a
decomposition  of  a  flat  surface into zippered  rectangles:  a
permutation $\pi\in\mathfrak{S}_\noi$.  This latter parameter  is
discrete.

The vertical parameters  $h_j, a_k$ and are not independent: they
satisfy some  linear  equations  and  inequalities.  Varying  the
continuous   parameters   $\lambda,h,a$  respecting   the  linear
relations  between  parameters $h$ and $a$  we  get a new set  of
coordinates in the stratum. These coordinates were introduced and
studied by W.~A.~Veech  in~\cite{zorich:Veech:Annals:82}. In particular,
it         was        proved         that         for         any
$(\lambda,\pi)\in\Delta^{\noi-1}\times\mathfrak{R}$\index{0R10@$\mathfrak{R}$ -- (extended) Rauzy class}
in the  space
of    interval    exchange   transformations\index{Interval exchange transformation!space of}
there    is    an
$\noi$-dimensional open  cone  of  solutions  $(h,a)$.  In  other
words, having any interval exchange transformation $T:X\to X$ one
can always construct a flat surface $S$ and  a horizontal segment
$X\subset S$ inside it such that the first return of the vertical
flow to $X$ gives the initial  interval exchange transformations.
Moreover,  there  is  a  $\noi$-dimensional family of  such  flat
surfaces  --   {\it   suspensions}  over  the  interval  exchange
transformation $T:X\to X$ (see~\cite{zorich:Masur:Annals:82},
\cite{zorich:Veech:Annals:82}).

Is  there  a  canonical  decomposition  of  a  flat  surface into
zippered rectangles? A choice of horizontal  segment $X\subset S$
completely  determines  a decomposition of a generic surface  $S$
into zippered rectangles, so our  question  is  equivalent to the
problem  of  a  canonical  choice  of  a  horizontal  segment $X$
satisfying Convention~\ref{zorich:conv:minimal:noi}. The choice which we
propose is {\it  almost} canonical; it leaves an arbitrariness of
finite order which is the same for almost all $S$ in the stratum.
Here  is  the  choice. Let us  place  the  left  extremity of the
horizontal segment $X$ at  one  of the conical singularities, and
let us  choose the  length $|X|$ of the segment  in such way that
$X$  would   be   the   shortest   possible  interval  satisfying
Convention~\ref{zorich:conv:minimal:noi} and condition $|X|\ge 1$.

In practice the interval $X$ can be constructed as follows: start
with a  sufficiently  long  horizontal  interval  having its left
extremity     at     a     conical    point    and     satisfying
Convention~\ref{zorich:conv:minimal:noi}.     Apply     the     ``slow''
Rauzy--Veech algorithm as  long  as the resulting subinterval has
length at  least $1$. For almost  all flat surfaces  after finite
number of steps we obtain the desired interval $X$.

The surface $S$  has finite number of conical singularities; each
conical singularity has finite number of horizontal prongs, so we
get arbitrariness of finite order. Thus, the resulting {\it space
of zippered rectangles} can be essentially viewed as a (ramified)
covering    over   the    corresponding    connected    component
$\cH^{comp}(d_1,\dots,d_\noz)$  of  the  stratum.  Passing  to  a
codimension  one  subspace  $\Omega$  defined  by  the  condition
$\lambda\cdot h=1$ we get a space of zippered  rectangles of area
one covering the space  $\cH_1^{comp}(d_1,\dots,d_\noz)$  of flat
surfaces  of   area  one.  Consider  a  codimension-one  subspace
$\Upsilon\subset\Omega$ of  zippered  rectangles  which have unit
area,   and   which   have   the   base   $X$   of   length  one,
$|X|=\lambda_1+\dots+\lambda_\noi=1$. The space  $\Upsilon$ has a
natural   structure   of  a   fiber   bundle   over   the   space
$\Delta^{\noi-1}\times\mathfrak{R}$\index{0R10@$\mathfrak{R}$ -- (extended) Rauzy class}
of    interval    exchange
transformations:   we   associate   to   a   zippered   rectangle
$(\lambda,h,a,\pi)\in\Upsilon$     the     interval      exchange
transformation $(\lambda,\pi)$.

We would like  to emphasize an  interpretation of $\Omega$  as  a
\emph{fundamental domain}  in  the  space  of \emph{all} zippered
rectangles  of  area  one.  As a fundamental domain  $\Omega$  is
defined  by the  additional  condition on the  base:  $X$ is  the
shortest          possible          interval           satisfying
Convention~\ref{zorich:conv:minimal:noi} such that $|X|\ge  1$.  In this
interpretation $\Upsilon$  is  the  boundary  of  the fundamental
domain.   Starting   with   an   arbitrary   zippered   rectangle
representation  satisfying  Convention~\ref{zorich:conv:minimal:noi}  we
can  apply   several   steps   of   Rauzy--Veech  algorithm  (see
Fig.~\ref{zorich:fig:Rauzy:move:b}   and   Fig.~\ref{zorich:fig:Rauzy:move:a}),
which does not change  the  surface $S$. After several iterations
we get to the fundamental domain $\Omega$.

We  would use  the  same notation for  $\Omega$  considered as  a
fundamental domain and for $\Omega$  considered  as  a  quotient,
when two boundary components of  $\Omega$  are  identified by the
modification      of     zippered      rectangles      as      on
Fig.~\ref{zorich:fig:Rauzy:move:b} and Fig.~\ref{zorich:fig:Rauzy:move:a}.

\paragraph{Teichm\"uller Geodesic  Flow  and its First Return Map
to a Cross-section}

Zippered  rectangles  coordinates  are extremely convenient  when
working with the  Teichm\"uller  geodesic flow,
which we identify
with the  action  of  the  diagonal subgroup $g_t=\begin{pmatrix}
\exp(t) &  0\\0  &  \exp(-t)\end{pmatrix}$. Namely, $g_t$ expands
the horizontal parameters  $\lambda$  by the factor $\exp(t)$ and
contracts the vertical parameters $h,a$ by the same factor.

Consider  a  zippered rectangle  $S=(\lambda,h,a,\pi)\in\Upsilon$
with the base $X$ of unit  length. Applying $g_t$ to $S$ with $t$
continuously increasing from  $t=0$  we shall eventually make the
length of the base of  the  deformed zippered rectangle $g_t S  =
(\exp(t)\lambda,\exp(-t)h,\exp(-t)a;\pi)$  too  long and  thus we
shall get  outside of the  fundamental domain $\Omega$. It is not
difficult to  determine an exact  time $t_0$ when it will happen.
We get to the boundary of the fundamental domain $\Omega$  at the
time
\begin{equation}
\label{zorich:eq:time:t0}
t_0(S)=-\log\big(1-\min(\lambda_\noi,\lambda_{\pi^{-1}(\noi)})\big).
\end{equation}
The  time $t_0$ is  chosen  in  such  way that  applying  to  the
zippered  rectangle  $g_{t_0}  S$  one step of  the  Rauzy--Veech
induction       (see       Fig.~\ref{zorich:fig:Rauzy:move:b}        and
Fig.~\ref{zorich:fig:Rauzy:move:a}) we get a new zippered rectangle with
the    base     $X'$     of     unit     length.     To    verify
formula~\eqref{zorich:eq:time:t0} for  $t_0$  it  is  sufficient to note
that expansion-contraction commutes  with Rauzy--Veech induction.
Thus, to evaluate $t_0$ we can \emph{first} apply one step of the
Rauzy--Veech      induction       and      \emph{then}      apply
expansion-contraction for  an appropriate time, which would bring
us back to $\Upsilon$, i.e.  which  would make the length of  the
base of the new building of zippered rectangles equal to one.

In other words,  starting at a point $S\in\Upsilon$ and following
the flow  for the time  $t_0(S)$ we  get to the  boundary of  the
fundamental domain in the space of  zippered  rectangles  and  we
have to instantly jump back to the point of $\Upsilon$ identified
with  $g_{t_0}  S$. One can recognize in  this  construction  the
first return map  $\cS: \Upsilon\to \Upsilon$ defined by the flow
$g_t$ on  the section $\Upsilon$: at  the time $t_0(S)$  the flow
$g_t$ emitted from a point $S\in\Upsilon$  returns  back  to  the
codimension-one subspace $\Upsilon$ transversal to the flow.

Morally  one  should  consider the map  $\cS$  as  a  map on some
subspace of flat surfaces. Note, that $\cS$ is  not applicable to
points of flat surfaces, it associates to a flat surface taken as
a whole another flat surface taken as a whole.

We  see  now  that  the  Rauzy--Veech  renormalization  procedure
$\cS:\Upsilon\to\Upsilon$  performed  on the  level  of  zippered
rectangles is nothing  but  discrete version of the Teichm\"uller
geodesic  flow.  Namely $\cS$  is  the first  return  map of  the
Teichm\"uller  geodesic   flow   to   a  section  $\Upsilon$.  By
construction          the          Rauzy--Veech         induction
$\cT:\Delta^{\noi-1}\times\mathfrak{R}\to
\Delta^{\noi-1}\times\mathfrak{R}$\index{0R10@$\mathfrak{R}$ -- (extended) Rauzy class}
on  the   space  of  interval
exchange transformations\index{Interval exchange transformation!space of}
is  just a projection of $\cS$. In other
words, the following diagram
$$
\begin{CD}
\Upsilon(\mathfrak{R})   @>\cS>>  \Upsilon(\mathfrak{R})\\
@VVV                             @VVV \\
\Delta^{\noi-1}\times\mathfrak{R} @>\cT>>
\Delta^{\noi-1}\times\mathfrak{R}
\end{CD}
$$
is   commutative,   and  the  invariant  measure  on  the   space
$\Delta^{\noi-1}\times\mathfrak{R}$    of    interval    exchange
transformations  is  a push  forward  of  the  natural  invariant
measure on the space $\Upsilon$ of zippered rectangles.

\paragraph{Choice of a Section}

Now we can return to the questions addressed at the  beginning of
this section. Ignoring an algorithmic aspect  of  the  choice  of
renormalization procedure we see that conceptually, it is defined
by a section of  the  Teichm\"uller geodesic flow. In particular,
the          ``fast''          renormalization          procedure
$\cG:\Delta^{\noi-1}\times\mathfrak{R}                        \to
\Delta^{\noi-1}\times\mathfrak{R}$\index{0R10@$\mathfrak{R}$ -- (extended) Rauzy class}
defined   in  the   previous
section   corresponds    to    a    choice    of   a   subsection
$\Upsilon'\subset\Upsilon$. Luckily it has  a  simple algorithmic
representation in terms of modification of  the interval exchange
transformation $T(\lambda,\pi)$,  and,  moreover, it has a simple
description in  terms  of  coordinates  $\lambda,h,a,\pi$  in the
space of zippered rectangles given by an extra  condition for the
parameter $a_\noi$.

Recall that parameters  $a_j$ are responsible for the position of
singularities:  we  zip  the  neighboring  rectangles  $R_j$  and
$R_{j+1}$  from  the   bottom  up  to  the  altitude  $a_j$,  see
Fig.~\ref{zorich:fig:iet:octagon}.  In  particular, by  construction all
$a_j$ for $j=1,\dots,\noi-1$ are positive. Parameter $a_\noi$ is,
however, different from the others: the rectangle $R_\noi$ is the
rightmost  rectangle  in the collection. If there  is  a  conical
singularity located at the right side of this rightmost rectangle
(see,  for  example, the zippered rectangle decomposition of  the
flat surface  on  the  top  part of Fig.~\ref{zorich:fig:Rauzy:move:b}),
then parameter $a_\noz$ is positive; it  indicates  as  usual  at
what height is located  the  singularity. However, the right side
of the rightmost rectangle might  contain  no  singularity.  This
means  that  the singularity  is  located  on  the  corresponding
vertical trajectory  {\it below} the  zero level of the base $X$.
The rectangle  which is glued to $X$ from  below at the rightmost
position is the rectangle  $R_{\pi^{-1}(\noi)}$;  the singularity
is located on the right side of this rectangle (see, for example,
the zippered rectangle decomposition  of  the flat surface on the
bottom part of Fig.~\ref{zorich:fig:Rauzy:move:b}). In this  case we let
$a_\noi$ be negative indicating how low we have  to descend along
downward vertical trajectory  emitted  from the right endpoint of
$X$ to hit the singularity.

The subsection  $\Upsilon'$  is  defined  by  the following extra
condition
\begin{align*}
\Upsilon' =\, &\{\,(\lambda,h,a,\pi)\in\Upsilon\, |\, a_\noi>0
\text{
when } \lambda_\noi>\lambda_{\pi^{-1}(\noi)}\}\, \sqcup\\
 \sqcup\, &\{\,(\lambda,h,a,\pi)\in\Upsilon\, |\, a_\noi<0 \text{ when
}\, \lambda_\noi<\lambda_{\pi^{-1}(\noi)}\}
\end{align*}

\begin{Exercise}
Check which zippered rectangles at Figures~\ref{zorich:fig:iet:octagon},
\ref{zorich:fig:Rauzy:move:b},   \ref{zorich:fig:Rauzy:move:a}   satisfy    the
condition  $a_\noi\cdot(\lambda_\noi-\lambda_{\pi^{-1}(\noi)})>0$
and which do not.
\end{Exercise}

It can be verified (see~\cite{zorich:Zorich:Gauss:map}) that the section
$\Upsilon'$ is  still a fiber bundle  over the space  of zippered
rectangles  and   that   the   corresponding   first  return  map
$\cS':\Upsilon'\to\Upsilon'$ projects to the map $\cG$:
\begin{equation}
\label{zorich:eq:suspension:CD}
\begin{CD}
\Upsilon'(\mathfrak{R})   @>\cS'>>  \Upsilon'(\mathfrak{R})\\
@VVV                             @VVV \\
\Delta^{\noi-1}\times\mathfrak{R} @>\cG>>
\Delta^{\noi-1}\times\mathfrak{R}
\end{CD}
\end{equation}

\begin{Exercise}
Verify that the definition of the renormalization procedure $\cG$
as a projection  of  the first  return  map of the  Teichm\"uller
geodesic flow to $\Upsilon'$ matches the  intrinsic definition of
$\cG$ given in Sec.~\ref{zorich:ss:space:of:iet}.
\end{Exercise}

\begin{figure}[htb]
\centering
\includegraphics{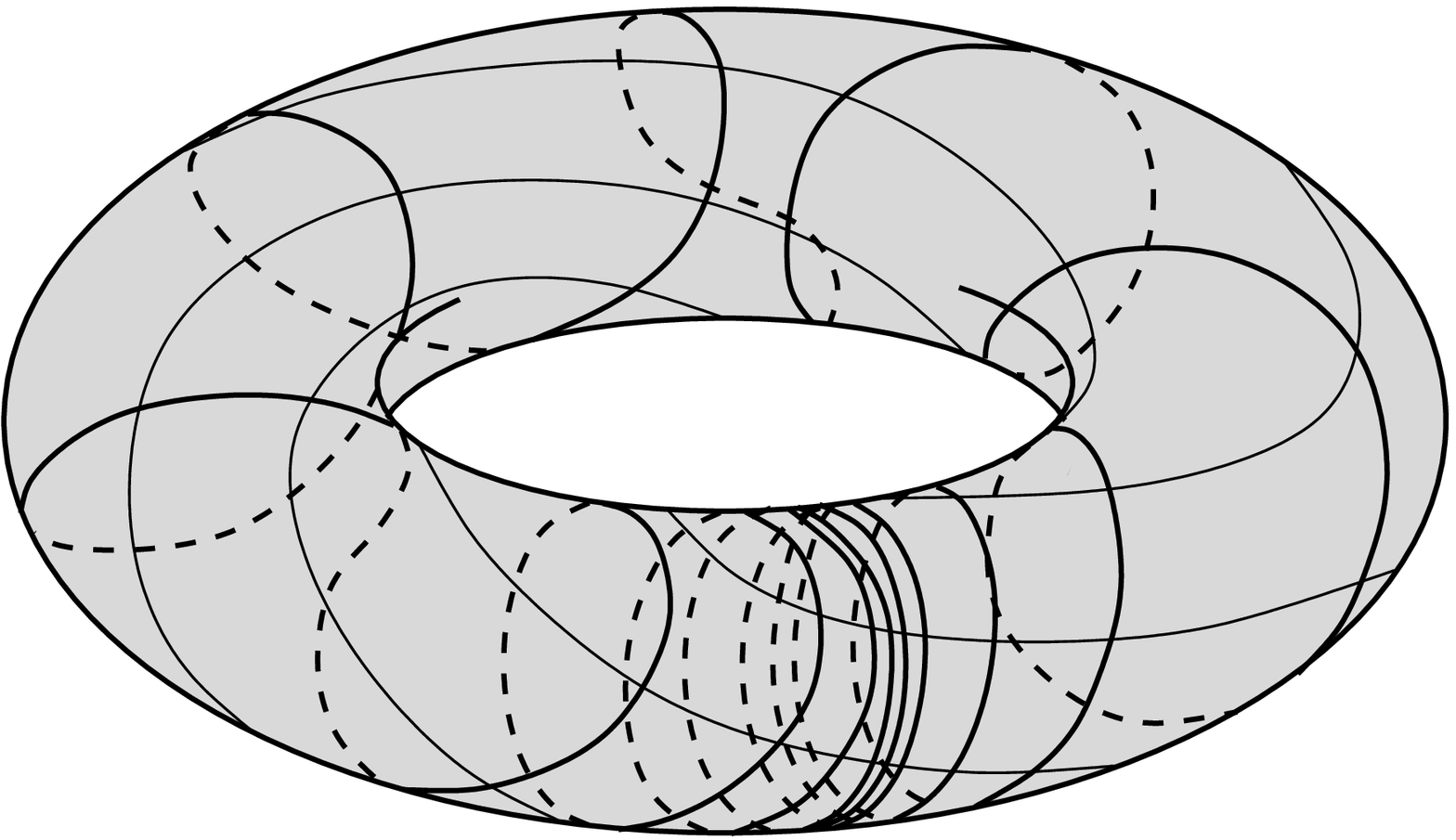}
\includegraphics{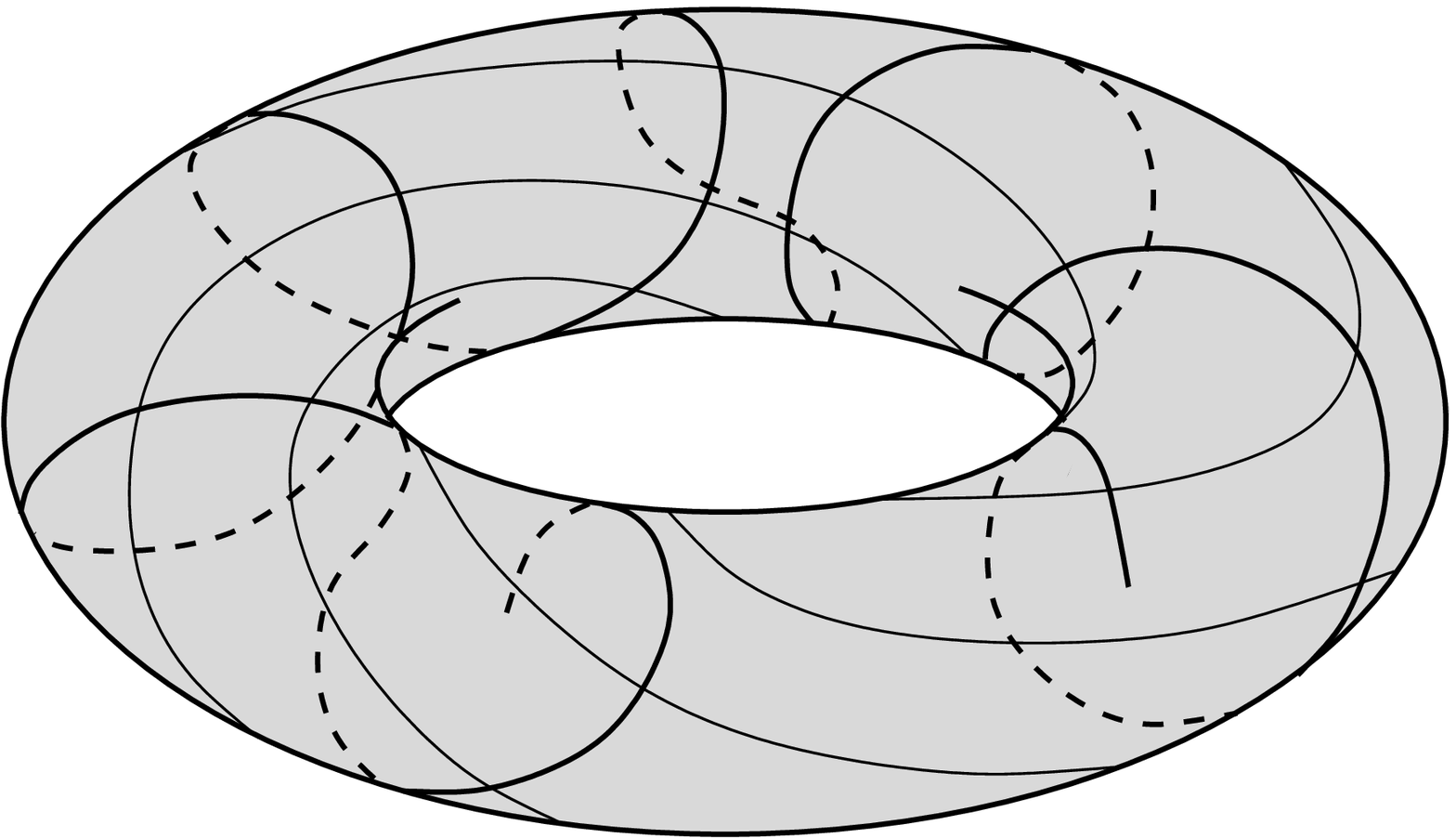}
\vspace{75bp}
\caption{
\label{zorich:fig:wrong:section}
The section  $Y$ on the  left picture has infinite measure though
the measure of the space is finite. The  subsection $Y'\subset Y$
on the right picture has finite measure. In both cases  the first
return map of the ergodic flow to the section is ergodic, but the
mean return time to the left subsection is zero
}
\end{figure}

Different choices of the  section  also explain why the invariant
measure  on   the  space  of  interval  exchange  transformations\index{Interval exchange transformation!space of}
$\Delta^{\noi-1}\times\mathfrak{R}$\index{0R10@$\mathfrak{R}$ -- (extended) Rauzy class}
was    infinite   for   the
Rauzy--Veech induction $\cT$  while is finite for the ``fast''
renormalization  procedure  $\cG$.  As  a model case  consider  a
directional flow  on a torus and  two different sections  to this
flow. Taking as a section  the  line $Y$ represented on the  left
picture of  Fig.~\ref{zorich:fig:wrong:section}  we  get  a  section  of
infinite measure  though the measure of  the torus is  finite and
the  flow is very  nice.  Taking  as  a section  a  finite  piece
$Y'\subset     Y$     as     on     the     right     side     of
Fig.~\ref{zorich:fig:wrong:section} we get a section of finite measure.

Similarly, the component  $\cH_1^{comp}(d_1,\dots,d_\noz)$ of the
stratum  has  finite  volume  and  hence  the  space  of zippered
rectangles   $\Omega$   which   is   a   finite   covering   over
$\cH_1^{comp}(d_1,\dots,d_\noz)$ also has finite volume. However,
the initial  section  $\Upsilon$ has infinite ``hyperarea'' while
section $\Upsilon'$ already has finite ``hyperarea''.

We  complete  this  section  with a several  comments  concerning
Lyapunov exponents. Though these comments are too brief to give a
comprehensive   proof   of   the   relation   between   exponents
$\nu_j$\index{0nu10@$\nu_1,\dots,\nu_g$ -- Lyapunov exponents related to the Teichm\"uller geodesic flow}\index{Exponent!Lyapunov exponent}\index{Lyapunov exponent}
responsible for  the  deviation and the Lyapunov exponents
of the Teichm\"uller geodesic flow,  they  present  the key idea,
which can be completed by an elementary calculation.

It  is  clear that the Lyapunov exponents  of  the  Teichm\"uller
geodesic       flow       $g_t$        on       the       stratum
$\cH_1^{comp}(d_1,\dots,d_\noz)$  coincide  with   the   Lyapunov
exponents on its finite covering $\Omega$. The Lyapunov exponents
of the flow $g_t$ differ from the Lyapunov exponents of the first
return map  $\cS'$ to the section  $\Upsilon'$ only by  a scaling
factor representing  the  ``hyperarea''  of  $\Upsilon'$. The map
$\cG$  on  the space of zippered rectangles  coincides  with  the
restriction  of  the  map  $\cS'$  on   the  zippered  rectangles
restricted  to  horizontal  parameters  $\lambda$  and   discrete
parameter $\pi$. Thus, the Lyapunov exponents  of  $\cG$  form  a
subcollection of  the  Lyapunov exponents of $\cS'$ corresponding
to the subspace of horizontal parameters $\lambda$. It remains to
note that the Lyapunov exponents of the map $\cG$ are  related to
the  Lyapunov  exponents  of the cocycle  $\Cocycle(\lambda,\pi)$
just by  the scaling factor, and  that we have  already expressed
the exponents $\nu_j$ responsible  for  the deviation in terms of
the Lyapunov  exponents of the multiplicative cocycle $\Cocycle$,
see~\eqref{zorich:eq:nu:j:as:theta:j:over:theta:1}.  Matching  all   the
elements of  this chain together we  get a representation  of the
exponents\index{0nu10@$\nu_1,\dots,\nu_g$ -- Lyapunov exponents related to the Teichm\"uller geodesic flow}\index{Exponent!Lyapunov exponent}\index{Lyapunov exponent}
$\nu_j$  in  terms  of the Lyapunov exponents  of  the
Teichm\"uller        geodesic         flow        given        in
Sec.~\ref{zorich:ss:Asymptotic:Flag:and:Dynamical:Hodge:Decomposition}.

\toread{Zippered rectangles and Lyapunov exponents}

More   details   on  Rauzy   classes\index{Rauzy class}\index{0R10@$\mathfrak{R}$ -- (extended) Rauzy class}
$\mathfrak{R}$,   zippered
rectangles, Lyapunov exponents of the Teichm\"uller geodesic flow
and  their  relation   might  be  found  in  original  papers  of
G.~Rauzy~\cite{zorich:Rauzy},           W.~Veech~\cite{zorich:Veech:Annals:82},
\cite{zorich:Veech:flow}    and    the   author~\cite{zorich:Zorich:Gauss:map},
\cite{zorich:Zorich:How:do}.

\index{Zippered rectangle|)}

\subsection{Spectrum of Lyapunov Exponents
(after M.~Kontsevich, G.~Forni, A.~Avila and M.~Viana)}
\label{zorich:ss:Spectrum:of:Lyapunov:exponents}

It  should  be  mentioned  that the statement that  the  subspace
$\cV_g$, such  that  $|dist(c(x,N),\cV_g)|\le const$ for any $N$,
has dimension \emph{exactly}  $g$  was formulated in the original
paper~\cite{zorich:Zorich:How:do}  as  a  conditional statement. It  was
based on the conjecture that the  Lyapunov  exponent
$\nu_g$\index{0nu10@$\nu_1,\dots,\nu_g$ -- Lyapunov exponents related to the Teichm\"uller geodesic flow}\index{Exponent!Lyapunov exponent}\index{Lyapunov exponent}
is strictly positive. This  conjecture  was later proved by G.~Forni
in~\cite{zorich:Forni:02}.

\begin{NNTheorem}[G.~Forni]
For any connected component $\cH^{comp}(d_1,\dots,d_\noz)$ of any
stratum of Abelian differentials the first $g$ Lyapunov exponents
of the Teichm\"uller geodesic flow are strictly greater than $1$:
$$
1+\nu_g>1
$$
\end{NNTheorem}

As an indication why this positivity is not something which
should be taken for granted we would like to give  some precision
about related results of G.~Forni.

The  Lyapunov  exponents of the Teichm\"uller geodesic flow  play
the role of (logarithms of)  eigenvalues  of  a virtual ``average
monodromy  of  the tangent bundle along the  flow''.  Instead  of
considering the tangent bundle to $\cH(d_1,\dots,d_\noz)$\index{0H20@$\cH(d_1,\dots,d_\noz)$ -- stratum in the moduli space}\index{Stratum!in the moduli space} one can
consider another vector bundle intimately related  to the tangent
bundle.  This  vector  bundle  has the space $H^1(S;\R{})$  as  a
fiber.   Since   we   know   how   to   identify   the   lattices
$H^1(S;\Z{})\subset    H^1(S;\R{})$    and   $H^1(S';\Z{})\subset
H^1(S';\R{})$ in  the fibers over  two flat surfaces $S$ and $S'$
which are close to each other in $\cH(d_1,\dots,d_\noz)$\index{0H20@$\cH(d_1,\dots,d_\noz)$ -- stratum in the moduli space}\index{Stratum!in the moduli space}, we know
how to transport the  fiber  $H^1(S;\R{})$ over the ``point'' $S$
in the  base  $\cH(d_1,\dots,d_\noz)$\index{0H20@$\cH(d_1,\dots,d_\noz)$ -- stratum in the moduli space}\index{Stratum!in the moduli space} to the fiber $H^1(S';\R{})$
over the ``point''  $S'$.  In other  words,  we have a  canonical
connection (called Gauss--Manin connection) in the vector bundle.
Hence we can  again  study the  ``average  monodromy of the  fiber
along  the  flow''.  It  is  not  difficult  to   show  that  the
corresponding Lyapunov  exponents  are  related  to  the Lyapunov
exponents of the  Teichm\"uller  flow. Namely, the new collection
of Lyapunov exponents\index{0nu10@$\nu_1,\dots,\nu_g$ -- Lyapunov exponents related to the Teichm\"uller geodesic flow}\index{Exponent!Lyapunov exponent}\index{Lyapunov exponent}
has the form:
$$
1\ge\nu_2\ge \dots \ge\nu_g\ge -\nu_g\ge \dots \ge -\nu_2\ge-1
$$
In  particular,  the collection  of  Lyapunov  exponents  of  the
Teichm\"uller geodesic flow  can be obtained as follows: take two
copies of the collection above;  add  $+1$ to all the entries  in
one copy; add $-1$ to all entries in another copy; take the union
of the resulting  collections. The theorem of G.~Forni tells that
for any connected component of any stratum we have $\nu_g>0$.

Consider    now    some     $SL(2,\R{})$-invariant     subvariety
$\cN\subset\cH(d_1,\dots,d_\noz)$.  Consider  the restriction  of
the vector  bundle with the  fiber $H^1(S;\R{})$ to $\cN$. We can
compute   the   ``average  monodromy  of  the  fiber  along   the
Teichm\"uller  flow''  restricted  to  $\cN$.  It   gives  a  new
collection  of
Lyapunov   exponents\index{0nu10@$\nu_1,\dots,\nu_g$ -- Lyapunov exponents related to the Teichm\"uller geodesic flow}\index{Exponent!Lyapunov exponent}\index{Lyapunov exponent}.
Since  the  ``holonomy''
preserves  the  natural   symplectic   form  in  the  fiber,  the
collection will be again symmetric:
$$
1\ge\nu'_2\ge \dots \ge\nu'_g\ge -\nu'_g\ge \dots \ge
-\nu'_2\ge-1
$$
G.~Forni has showed~\cite{zorich:Forni:Handbook} that there are examples
of  invariant   subvarieties   $\cN$   such  that  all  $\nu'_j$,
$j=2,\dots,g$, are  equal  to zero! Moreover, G.~Forni explicitly
describes the locus where the monodromy does not change the fiber
(or it exterior powers) too  much,  and where one may  get
multiplicities of Lyapunov exponents.

Another     conditional     statement     in     the     original
paper~\cite{zorich:Zorich:How:do}   concerns   \emph{strict}  inclusions
$\cV_1\subset \cV_2\subset\dots\subset \cV_g\subset H_1(S;\R{})$.
It was  based on the  other conjecture claiming that the Lyapunov
exponents have  simple  spectrum.  The  first  strict  inequality
$\nu_1>\nu_2$  is  an elementary corollary of general results  of
Veech; see~\cite{zorich:Zorich:Deviation}. The  other strict inequalities
are  much  more difficult to prove. Very  recently  A.~Avila  and
M.~Viana~\cite{zorich:Avila:Viana} have announced a proof of  simplicity
of  the  spectrum~\eqref{zorich:eq:simplicity:of:Lyapunov:exponents} for
any connected component of any stratum proving the conjecture which was
open for a decade.

\begin{NNTheorem}[A.~Avila, M.Viana]
For     any     connected    component     of     any     stratum
$\cH^{comp}(d_1,\dots,d_\noz)$ of Abelian differentials the first
$g$ Lyapunov exponents\index{0nu10@$\nu_1,\dots,\nu_g$ -- Lyapunov exponents related to the Teichm\"uller geodesic flow}\index{Exponent!Lyapunov exponent}\index{Lyapunov exponent}
are distinct:
\begin{equation}
\label{zorich:eq:simplicity:of:Lyapunov:exponents}
1+\nu_1>1+\nu_2>\dots>1+\nu_g
\end{equation}
\end{NNTheorem}

\paragraph{Sum of the Lyapunov exponents}

Currently there are no methods of calculation of Lyapunov
exponents for general dynamical systems. The Teichm\"uller
geodesic flow does not make an exception: there is some
knowledge of approximate values of the numbers
$\nu_j$\index{0nu10@$\nu_1,\dots,\nu_g$ -- Lyapunov exponents related to the Teichm\"uller geodesic flow}\index{Exponent!Lyapunov exponent}\index{Lyapunov exponent}
obtained
by computer simulations for numerous low-dimensional strata, but
there is no approach leading to explicit evaluation of these
numbers with exception for some very special cases.

Nevertheless, for any connected component of any stratum (and,
more generally, for any $GL^+(2;\R{})$-invariant suborbifold) it
is possible to evaluate the {\it sum} of the Lyapunov exponents
$\nu_1+\dots\nu_g$, where $g$ is the genus. The formula for this
sum was discovered by M.~Kontsevich in~\cite{zorich:Kontsevich};
it is given in terms of the following natural structures on the
strata $\cH(d_1,\dots,d_\noz)$\index{0H20@$\cH(d_1,\dots,d_\noz)$ -- stratum in the moduli space}\index{Stratum!in the moduli space}.

There is a natural action of $\C\ast$ on every stratum of the
moduli space of holomorphic 1-forms: we can multiply a
holomorphic form $\omega$ by a complex number. Let us denote by
$\cH_{(2)}(d_1,\dots,d_\noz)$ the quotient of
$\cH(d_1,\dots,d_\noz)$\index{0H20@$\cH(d_1,\dots,d_\noz)$ -- stratum in the moduli space}\index{Stratum!in the moduli space} over $\C\ast$. The space
$\cH_{(2)}(d_1,\dots,d_\noz)$ can be viewed as the space of flat
surfaces of unit area {\it without} choice of distinguished
direction.

There are two natural holomorphic vector bundles over
$\cH_{(2)}(d_1,\dots,d_\noz)$. The first one is the
$\C{\ast}$-bundle
$\cH(d_1,\dots,d_\noz)\to\cH_{(2)}(d_1,\dots,d_\noz)$. The
second one is the $\C{g}$-bundle, which fiber is composed of all
holomorphic 1-forms in the complex structure corresponding to a
flat surface $S\in\cH_{(2)}(d_1,\dots,d_\noz)$. Both bundles
have natural curvatures; we denote by $\gamma_1$ and $\gamma_2$
the corresponding closed curvature 2-forms.

Finally, there is a natural closed codimension two form $\beta$
on every stratum $\cH_{(2)}(d_1,\dots,d_\noz)$. To construct
$\beta$ consider the natural volume form $\Omega$ on
$\cH(d_1,\dots,d_\noz)$\index{0H20@$\cH(d_1,\dots,d_\noz)$ -- stratum in the moduli space}\index{Stratum!in the moduli space}. Four generators of the Lie algebra
$\mathfrak{gl}(2;\R{})$ define four distinguished vectors in the
tangent space $T_S \cH(d_1,\dots,d_\noz)$ at any ``point'' $S\in
\cH(d_1,\dots,d_\noz)$. Plugging these four vectors in the first
four arguments of the volume form $\Omega$ we get a closed
codimension four form on $\cH(d_1,\dots,d_\noz)$\index{0H20@$\cH(d_1,\dots,d_\noz)$ -- stratum in the moduli space}\index{Stratum!in the moduli space}. It is easy to
check that this form can be pushed forward along the
$\C{\ast}$-fibers of the bundle $\cH(d_1,\dots,d_\noz)\to
\cH_{(2)}(d_1,\dots,d_\noz)$ resulting in the closed codimension
two form on the base of this fiber bundle.

\begin{NNTheorem}[M.~Kontsevich]
For any connected component of any stratum the sum of the first
$g$ Lyapunov exponents\index{0nu10@$\nu_1,\dots,\nu_g$ -- Lyapunov exponents related to the Teichm\"uller geodesic flow}\index{Exponent!Lyapunov exponent}\index{Lyapunov exponent}
can be expressed as
$$
\nu_1+\dots+\nu_g=\cfrac{\int\beta\wedge \gamma_2}{\int \beta
\wedge \gamma_1}\ ,
$$
where the integration is performed over the corresponding
connected component of $\cH_{(2)}(d_1,\dots,d_\noz)$.
\end{NNTheorem}

As it was shown by G.~Forni, this formula can be generalized for
other $GL^+(2;\R{})$-invariant submanifolds.

The proof is based on two observations. The first one
generalizes the fact that dynamics of the geodesic flow on the
hyperbolic plane is in some sense equivalent to dynamics of
random walk. One can replace Teichm\"uller geodesics by geodesic
broken lines consisting of geodesic segments of unit length.
Having a broken line containing $n$ geodesic segments with the
endpoint at the point $S_n$ we emit from $S_n$ a new geodesic in
a random direction and stop at the distance one from $S_n$
at the new point $S_{n+1}$. This generalization suggested
by M.~Kontsevich was formalized and justified by G.~Forni.

Consider the vector bundle over the moduli space of holomorphic
1-forms with the fiber $H^1(S;\R{})$ over the ``point'' $S$. We
are interested in the sum $\nu_1+\dots+\nu_g$ of Lyapunov
exponents\index{0nu10@$\nu_1,\dots,\nu_g$ -- Lyapunov exponents related to the Teichm\"uller geodesic flow}\index{Exponent!Lyapunov exponent}\index{Lyapunov exponent}
representing mean monodromy of this vector bundle
along random walk. It follows from standard arguments concerning
Lyapunov exponents that this sum corresponds to the top Lyapunov
exponent of the exterior power of order $g$ of the initial
vector bundle. In other words, we want to measure the average growth
rate of the norm of a $g$-dimensional subspace in $H^1(S;\R{})$
when we transport it along trajectories of the random walk
using the Gauss--Manin connection.

Fix a Lagrangian subspace $L$ in the fiber $H^1(S_0;\R{})$
over a ``point'' $S_n$. Consider the set of points located at the
Teichm\"uller distance $1$ from $S_n$. Transport $L$ to each
point $S_{n+1}$ of this ``unit sphere'' along the corresponding
geodesic segment joining $S_n$ with $S_{n+1}$; measure the logarithm
of the change of the norm of $L$; take the average over the
``unit sphere''. The key observation of M.~Kontsevich
in~\cite{zorich:Kontsevich} is that for an
appropriate choice of the norm this average growth rate is the
same for all Lagrangian subspaces $L$ in $H^1(S_0;\R{})$ and
depends only on the point $S_n$. A calculation based on this
observation gives the formula above.

Actually, formula above can be rewritten in a much more explicit
form (which is a work in progress). The values of the sum given
by this more explicit formula perfectly match numerical
simulations. The table below gives the values of the sums of
Lyapunov exponents for some low-dimensional strata; this
computation uses the results of A.~Eskin and
A.~Okounkov~\cite{zorich:Eskin:Okounkov} for the volumes of the
strata.

$$
\begin{array}{|c|c|c|c|c|c|c|}
\multicolumn{7}{c}{\text{Conjectural values of }\nu_1+\dots+\nu_g \text{ for some strata}}\\
[-\halfbls]\multicolumn{7}{c}{\text{}}\\
\hline &&&&&&\\
\ \cH(2)\ &\ \cH(1,1)\ &\ \dots\ &\ \cH(4,1,1)\ &\ \dots\ &\
\cH(1,1,1,1,1,1)\ &\ \cH(1,1,1,1,1,1,1,1)
\\
[-\halfbls] &&&&&&\\ \hline &&&&&& \\ [-\halfbls] \cfrac{4}{3}&
\cfrac{3}{2} & \dots & \cfrac{1137}{550} & \dots &
\cfrac{839}{377} & \cfrac{235\,761}{93\,428}
\\ &&&&&&\\ \hline
\end{array}
$$

In particular, since $\nu_1=1$ this information gives the exact
value of the only nontrivial
Lyapunov exponent\index{0nu10@$\nu_1,\dots,\nu_g$ -- Lyapunov exponents related to the Teichm\"uller geodesic flow}\index{Exponent!Lyapunov exponent}\index{Lyapunov exponent}
$\nu_2$ for the
strata in genus two. Some extra arguments show that $\nu_2=1/3$
for the stratum $\cH(2)$ and for any $GL^+(2;R)$-invariant
submanifold in it; $\nu_2=1/2$ for the stratum $\cH(1,1)$ and
for any $GL^+(2;R)$-invariant submanifold in it.

\subsection{Encoding a Continued Fraction by a Cutting Sequence of a Geodesic}
\label{zorich:ss:Encoding:a:Continued:Fraction}

\index{Continued fraction|(}

We have seen that renormalization for a rotation of a  circle (or
equivalently for  an  interval  exchange  transformation  of  two
subintervals)  leads  to  the  Euclidean algorithm which  can  be
considered  in  this guise  as  a  particular  case  of  the fast
Rauzy--Veech induction.

The multiplicative cocycle
$$
\Cocycle^{(s)}=
\begin{pmatrix}
1 & n_1\\
0 & 1
\end{pmatrix}
\cdot
\begin{pmatrix}
1 & 0\\
n_2 & 1
\end{pmatrix}
\cdots
\begin{pmatrix}
1 &\ n_{2k-1}\\
0 &\ 1
\end{pmatrix}
\cdot
\begin{pmatrix}
1 &\ 0\\
n_{2k} &\ 1
\end{pmatrix}
\cdots
$$
considered  in  section~\ref{zorich:ss:space:of:iet} corresponds  to the
decomposition  of   a  real  number  $x\in(0,1)$  into  continued
fraction, $x=[0;n_1,n_2,\dots]$,
$$
x= \cfrac{1}{n_1+\cfrac{1}{n_2+\cfrac{1}{\dots
}}}
$$

A flat surface which realizes an interval exchange transformation
of two subintervals is a flat torus. The the moduli space of flat
tori  can  be  naturally identified with  $SL(2,\R{})/SL(2,\Z{})$
which  in its  turn  can be naturally  identified  with the  unit
tangent bundle  to  the  modular  surface  $\Hyp/SL(2,\Z{})$  see
Sec.~\ref{zorich:ss:Toy:Example:Family:of:Flat:Tori}.   Moreover,    the
Teichm\"uller  metric  on the space of tori  coincides  with  the
hyperbolic metric on  $\Hyp$,  and the
Teichm\"uller geodesic flow\index{Teichm\"uller!geodesic flow|)}
on
the moduli space of flat tori coincides with the geodesic flow on
the modular surface.

Hence, the construction  from  the previous section suggests that
the Euclidean  algorithm  corresponds  to the following geometric
procedure. There should  be a section $\Upsilon$ in the (covering
of)  the  unit tangent  bundle  to the  modular  surface and  its
subsection $\Upsilon'\subset\Upsilon$ such that the trajectory of
the geodesic flow emitted from a point of  $\Upsilon'$ returns to
$\Upsilon'$ after $n_1$ intersections with $\Upsilon$, then after
$n_2$ intersections with  $\Upsilon$, etc. In other words there is
a natural way  to  code  a continued fraction by  a  sequence  of
intersections   (so   called   ``cutting   sequence'')   of   the
corresponding       geodesic       with       some       sections
$\Upsilon'\subset\Upsilon$.

Actually, a geometric coding of a continued fraction by a cutting
sequence of a geodesics on a  surface is known since the works of
J.~Nielsen  and  E.~Artin  in  20s  and  30s. The  study  of  the
geometric coding  was developed in  the 80s and 90s by C.~Series,
R.~Adler, L.~Flatto and other authors. We refer to the expository
paper~\cite{zorich:Series} of  C.~Series for detailed description of the
following geometric coding algorithm.

Consider  a  tiling  of  the  upper  half  plane  with  isometric
hyperbolic  triangles as  at  Fig.~\ref{zorich:fig:cutting:sequence}.  A
fundamental domain  of the tiling  is a triangle with vertices at
$0,1$  and  $\infty$;  the  corresponding quotient surface  is  a
triple   cover   over   the   standard   modular   surface   (see
Fig.~\ref{zorich:fig:L:3}).  This  triangulation  of  $\Hyp$  by   ideal
triangles is also known as {\it Farey tessellation}.

\begin{figure}[htb]
\centering
\includegraphics{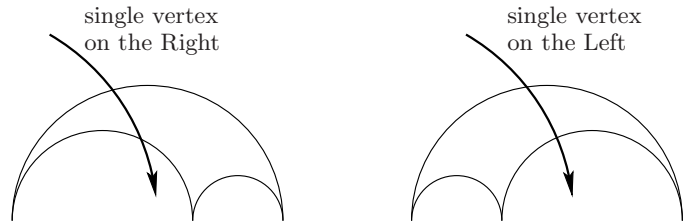}
\begin{picture}(0,0)(0,0)
\put(-100,0){\small single vertex}
\put(-100,-10){\small on the Right}
\put(60,0){\small single vertex}
\put(60,-10){\small on the Left}
\end{picture}

\vspace{85bp}
\caption{
\label{zorich:fig:coding:rule}
Coding rule:  when we cross a triangle leaving  one vertex on the
left and two on the right  we write symbol $L$; when there is one
vertex on the right and two on the left we write symbol $R$}
\end{figure}

\begin{figure}[htb]
\centering
\includegraphics{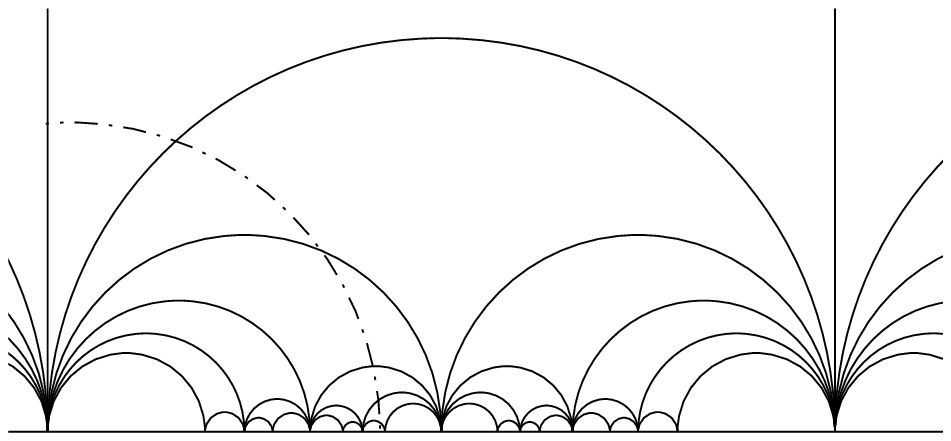}
\begin{picture}(0,0)(0,0)
\put(-137,-40){$iy$}
\put(-33,-135){$x$}
\put(-129,-135){$0$}
\put(99,-135){$1$}
\end{picture}
\vspace{140bp}
\caption{
\label{zorich:fig:cutting:sequence}
(After C.~Series.) The cutting sequence defined  by this geodesic
starts with  $R,R,L,L,R,L,L,R,R,L,\dots$  which  we abbreviate as
$R^2 L^2  R^1 L^2 R^2 L^1 \dots$. The  real number $x\in(0,1)$ at
which  lands   the  geodesic  has  continued  fraction  expansion
$x=[0;2,2,1,2,2,1,\dots]$
}
\end{figure}

Consider  a   real  number  $x\in(0,1)$.  Consider  any  geodesic
$\gamma$  landing to  the  real axis at  $x$  such that  $\gamma$
intersects with the  imaginary  axis; let  $iy$  be the point  of
intersection. Let us  follow  the geodesic $\gamma$ starting from
$iy$ in direction of  $x$. Each time when we cross a  triangle of
our tiling  let us  note by the symbol $L$  the situation when we
have a  single vertex  on the left and two  vertices on the right
(see  Fig~\ref{zorich:fig:coding:rule})   and  by  the  symbol  $R$  the
symmetric situation.

\begin{Example}
Following     the     geodesic      $\gamma$     presented     at
Fig.~\ref{zorich:fig:cutting:sequence}     from     some     $iy$     to
$x=(\sqrt{85}-5)/10)\approx   0.421954$   we   get   a   sequence
$R,R,L,L,R,L,L,R,R,L,\dots$ which  we  abbreviate as $R^2 L^2 R^1
L^2 R^2 L^1 \dots$.
\end{Example}

\begin{NNTheorem}[C.~Series]
Let $x\in(0,1)$ be irrational. Let $\gamma$ be a geodesic emitted
from    some     $iy$     and     landing     at     $x$;     let
$R^{n_1}L^{n_2}R^{n_3}L^{n_4}\dots$ be  the corresponding cutting
sequence. Then $x=[0;n_1,n_2,n_3,n_4,\dots]$.
\end{NNTheorem}

\toread{Geometric symbolic coding}

I can strongly recommend a paper of P.~Arnoux~\cite{zorich:Arnoux} which
clearly and rigorously explains the idea  of  suspension  in  the
spirit   of   diagram~\eqref{zorich:eq:suspension:CD}    handling    the
particular case of Euclidean algorithm  and  of  geodesic flow on
the Poincar\'e upper  half-plane. As a survey on geometric coding
we can recommend  the  survey of C.~Series~\cite{zorich:Series} (as well
as other surveys in this collection on related subjects).

\index{Continued fraction|)}

\section{Closed Geodesics and Saddle Connections on Flat Surfaces}
\label{zorich:s:Closed:Geodesics:and:Saddle:Connections:on:Flat:Surfaces}

\index{Geodesic!closed|(}
\index{Saddle!connection|(}
\index{Geodesic!counting of periodic geodesics|(}

Emitting a geodesic in an irrational direction on a flat torus we
get  an  irrational  winding  line;  emitting  it  in  a rational
direction we get a closed geodesic. Similarly, for a flat surface
of higher genus a  countable  dense set of directions corresponds
to closed geodesics.

In this section  we study how  many closed regular  geodesics  of
bounded length  live on a generic  flat surface $S$.  We consider
also \emph{saddle  connections}  (i.e.  geodesic segments joining
pairs of conical singularities) and count them.

We explain a curious phenomenon concerning saddle connections and
closed geodesics on flat surfaces:  they  often  appear in pairs,
triples,  etc  of parallel  saddle  connections  (correspondingly
closed geodesics) of equal length.

When  all   saddle   connections   (closed   geodesics)  in  such
\emph{configuration} become short the corresponding flat  surface
starts to degenerate and gets close to the boundary of the moduli
space. Thus, a description of possible configurations of parallel
saddle connections (closed  geodesics)  gives us a description of
the multidimensional ``cusps''\index{Cusp of the moduli space@Cusp of the moduli space, \emph{see also} Moduli space; principal boundary}
of the strata.

\subsection{Counting Closed Geodesics and Saddle Connections}
\label{zorich:ss:Counting:Closed:Geodesics:and:Saddle:Connections}

\paragraph{Closed Regular Geodesics Versus Irrational Winding
Lines}

Consider a flat torus  obtained  by identifying pairs of opposite
sides  of a  unit  square. A geodesic  emitted  in an  irrational
direction  (one  with irrational slope) is an irrational  winding
line; it is dense in the  torus. A geodesic emitted in a rational
direction is closed; all parallel  geodesic  are  also closed, so
directional flow in a rational direction  fills  the  torus  with
parallel periodic  trajectories.  The  set of rational directions
has measure  zero in the set of all  possible directions. In this
sense  directions  representing  irrational  winding  lines   are
typical   and   directions  representing   closed   geodesic   --
nontypical.

The situation with  flat surfaces of  higher genera $g\ge  2$  is
similar in many  aspects, though more complicated in details. For
example , for  any flat surfaces  $S$ almost all  directions  are
``irrational''; any  geodesic  emitted in an irrational direction
is dense in the surface. Actually,  even  stronger  statement  is
true:

\begin{NNTheorem}[S.~Kerckhoff, H.~Masur, J.~Smillie]
For any flat surface directional flow in almost any direction is
uniquely ergodic\index{Flow!uniquely ergodic}\index{Ergodic!uniquely ergodic}.
\end{NNTheorem}

For  the   torus   the   condition   that   directional  flow  is
\emph{minimal} (that is any trajectory going in this direction is
dense in the torus) is equivalent to the condition that  the flow
is \emph{uniquely  ergodic} (the natural Lebesgue measure induced
by the flat structure is the only finite  measure invariant under
directional   flow;   see  Appendix~\ref{zorich:s:Ergodic:Theorem}   for
details). Surprisingly a  directional flow on a surface of higher
genus  (already  for  $g=2$)  might  be  \emph{minimal}\index{Flow!minimal}
but  not
\emph{uniquely ergodic}\index{Flow!minimal}\index{Flow!uniquely ergodic}\index{Ergodic!uniquely ergodic}!
Namely,  for  some  directions  which  give  rise  to  a  minimal
directional flow it might be possible to divide  the surface into
two parts (of nonzero measure) in such way that some trajectories
would  mostly  stay in one part while  other  trajectories  would
mostly stay in the other.

Closed geodesics on flat surfaces of higher genera also have some
similarities  with ones  on  the torus. Suppose  that  we have  a
regular  closed  geodesic  passing  through a point  $x_0\in  S$.
Emitting a geodesic from a nearby point $x$ in the same direction
we obtain a parallel closed  geodesic  of the same length as  the
initial one. Thus, closed geodesics  also  appear  in families of
parallel closed geodesics.  However, in the torus case every such
family fills the entire torus while a family  of parallel regular
closed geodesics  on a flat surfaces  of higher genus  fills only
part of the surface. Namely,  it  fills a flat cylinder having  a
conical singularity on each of its boundaries.

\begin{Exercise}
Find  several  periodic  directions  on  the  flat  surface  from
Fig.~\ref{zorich:fig:suspension}.   Find   corresponding   families   of
parallel closed geodesics.  Verify that each of the surfaces from
Fig.~\ref{zorich:fig:two:different:parities}   decomposes   under    the
vertical flow  into  three cylinders (of different circumference)
filled with periodic trajectories. Find these cylinders.
\end{Exercise}

\paragraph{Counting Problem}

Take an arbitrary loop on a torus. Imagine that it is made from a
stretched  elastic  cord.  Letting  it contract we get  a  closed
regular geodesic (may be winding several times along itself). Now
repeat the experiment with a  more  complicated  flat surface. If
the initial  loop was very simple (or if  we are extremely lucky)
we again obtain a regular closed geodesic. However, in general we
obtain a closed broken line of geodesic segments with vertices at
a collection of conical points.

Similarly  letting  contract an elastic cord joining  a  pair  of
conical  singularities  we usually obtain a broken line  composed
from several  geodesic segments joining conical singularities. In
this sense torus is very different from a general flat surface.

A geodesic segment  joining  two conical singularities and having
no  conical  points  in  its  interior   is  called
\emph{saddle connection}\index{Saddle!connection}.
The case  when boundaries of  a saddle connection coincide is not
excluded:  a  saddle  connection  might join a conical  point  to
itself.

\begin{Convention}
In  this  paper we consider only saddle  connections  and  closed
regular  geodesics.  We  never  consider broken lines  formed  by
several geodesic segments.
\end{Convention}

Now we are ready to formulate the Counting Problem. Everywhere in
this section we normalize the area of flat surfaces to one.

\begin{CountingProblem}
Fix a flat surface $S$. Let $N_{sc}(S,L)$ be the number of saddle
connections on  $S$ of length at  most $L$. Let  $N_{cg}(S,L)$ be
the  number  of maximal  cylinders  filled  with  closed  regular
geodesics  of length  at  most $L$ on  $S$.  Find asymptotics  of
$N_{sc}(S,L)$ and $N_{cg}(S,L)$ as $L\to\infty$.
\end{CountingProblem}

It   was   proved   by   H.~Masur   (see~\cite{zorich:Masur:Lower:bound}
and~\cite{zorich:Masur:Growth:rate})  that  for  any  flat  surface  $S$
counting functions  $N(S,L)$  grow  quadratically in $L$. Namely,
there exist constants $0<const_1(S)<const_2(S)<\infty$ such that
$$
const_1(S) \le N(S,L)/L^2 \le const_2(S)
$$
for $L$  sufficiently  large.  Recently Ya.~Vorobets has obtained
in~\cite{zorich:Vorobets:uniform:bounds}  uniform  estimates   for   the
constants $const_1(S)$ and $const_2(S)$ which depend  only on the
genus of $S$.

Passing  from  \emph{all}  flat  surfaces  to  \emph{almost  all}
surfaces in  a given connected component  of a given  stratum one
gets a much more precise result; see~\cite{zorich:Eskin:Masur}.

\begin{NNTheorem}[A.~Eskin and H.~Masur]
For almost  all flat surfaces  $S$ in a given connected component
of  a  stratum  $\cH(d_1,\dots,d_\noz)$\index{0H20@$\cH(d_1,\dots,d_\noz)$ -- stratum in the moduli space}\index{Stratum!in the moduli space}  the  counting  functions
$N_{sc}(S,L)$ and $N_{cg}(S,L)$ have exact quadratic asymptotics
\begin{equation}
\label{zorich:eq:Siegel:Veech:const}
\lim_{L\to\infty}\cfrac{N_{sc}(S,L)}{\pi L^2}=const_{sc} \qquad
\lim_{L\to\infty}\cfrac{N_{cg}(S,L)}{\pi L^2}=const_{cg}
\end{equation}
where
Siegel--Veech constants
\index{Siegel--Veech!constant|(}
$const_{sc}$ and $const_{cg}$
are  the  same for  almost  all flat  surfaces  in the  component
$\cH^{comp}_1(d_1,\dots,d_\noz)$.
\end{NNTheorem}

We  multiply  denominator  by  $\pi$  to  follow  a  conventional
normalization.

\paragraph{Phenomenon of Higher Multiplicities}

Let us discuss now the following problem. Suppose that we  have a
regular  closed  geodesic  on  a flat surface $S$.  Memorize  its
direction, say, let  it be the North-West direction. (Recall that
by             Convention~\ref{zorich:conv:flat:surface}              in
Sec.~\ref{zorich:ss:Very:Flat:Surfaces} we  can  place  a compass at any
point of the surface and it will tell us what is the direction to
the North.)  Consider  the  maximal  cylinder  filled with closed
regular geodesics parallel to ours. Take a point $x$ outside this
cylinder  and  emit   a  geodesic  from  $x$  in  the  North-West
direction. There are two questions.\newline
\emph{-- How big is the chance to get a closed geodesic?}\newline
\emph{-- How big is the chance to get a closed geodesic of the
same length as the initial one?}

Intuitively it is clear that the answer to the first question is:
``the chances  are low'' and to the second  one ``the chances are
even    lower''.    This    makes    the    following     Theorem
(see~\cite{zorich:Eskin:Masur:Zorich}) somehow counterintuitive:

\begin{NNTheorem}[A.~Eskin, H.~Masur, A.~Zorich]
For almost all flat surfaces $S$ from any  stratum different from
$\cH_1(2g-2)$ or $\cH_1(d_1,d_2)$ the function $N_{two\underline\
cyl}(S,L)$  counting  the number of families of parallel  regular
closed geodesics filling two distinct maximal cylinders has exact
quadratic asymptotics
$$
\lim_{L\to\infty}\cfrac{N_{two\underline\ cyl}(S,L)}{\pi L^2}=
const_{two\underline\ cyl}
$$
where  Siegel--Veech   constants  $const_{two\underline\  cyl}>0$
depends only on the connected component of the stratum.

For almost all flat surface $S$  in any stratum one cannot find a
single  pair  of  parallel  regular  closed  geodesics on $S$  of
different length.
\end{NNTheorem}

There  is   general   formula   for  the  Siegel--Veech  constant
$const_{two\underline\  cyl}$  and for  similar  constants  which
gives explicit numerical answers for  all  strata  in low genera.
Recall that  the \emph{principal stratum} $\cH(1,\dots,1)$ is the
only stratum  of fill dimension in $\cH_g$\index{0H10@$\cH_g$ -- moduli space of holomorphic 1-forms}\index{Moduli space!of holomorphic 1-forms}; it  is the stratum of
holomorphic 1-forms  with simple zeros  (or, what is the same, of
flat surfaces with conical angles $4\pi$  at  all  cone  points).
Numerical values of the Siegel--Veech constants for the principal
stratum are presented in Table~\ref{zorich:tab:SV:constants}.

\begin{table}[hbt]
$$
\begin{array}{c|c|c|c|c}
& g=1 & g=2 & g=3 & g=4 \\
\hline &&&&\\
[-\halfbls]\text{single}&&&&\\
[-2\halfbls] &
\cfrac{1}{2}\cdot\cfrac{1}{\zeta(2)}\approx 0.304 &
\cfrac{5}{2}\cdot\cfrac{1}{\zeta(2)}\approx 1.52 &
\cfrac{36}{7}\cdot\cfrac{1}{\zeta(2)}\approx 3.13 &
\cfrac{3150}{377}\cdot\cfrac{1}{\zeta(2)}\approx 5.08 \\
[-1\halfbls]\text{cylinder}&&&&\\
\hline &&&&\\
[-\halfbls]\text{two}&&&&\\
[-2\halfbls] &
- & - &
\cfrac{3}{14}\cdot\cfrac{1}{\zeta(2)}\approx 0.13 &
\cfrac{90}{377}\cdot\cfrac{1}{\zeta(2)}\approx 0.145 \\
[-\halfbls]\text{cylinders}&&&&\\
\hline &&&&\\
[-\halfbls]\text{three}&&&&\\
[-2\halfbls] &
- & - & - &
\cfrac{5}{754}\cdot\cfrac{1}{\zeta(2)}\approx 0.00403 \\
[-\halfbls]\text{cylinders}&&&&\\
\hline
\end{array}
$$
\caption{
\label{zorich:tab:SV:constants} 
Siegel--Veech constants $const_{n\underline\ cyl}$ for
the principal stratum $\cH_1(1,\dots,1)$}
\end{table}

Comparing these  values we see,  that our intuition was not quite
misleading. Morally, in genus $g=4$  a  closed  regular  geodesic
belongs to  a  one-cylinder family with ``probability'' $97.1\%$,
to a two-cylinder family with  ``probability''  $2.8\%$  and to a
three-cylinder family  with  ``probability''  only $0.1\%$ (where
``probabilities'' are calculated  proportionally to Siegel--Veech
constants).

In theorem above we discussed closed regular geodesics. A similar
phenomenon is true for  saddle  connections. Recall that the cone
angle at a conical point on a flat surface is an integer multiple
of $2\pi$.  Thus, at a  point with  a cone angle\index{Cone angle}
$2\pi n$  every
direction is  presented $n$ times. Suppose  that we have  found a
saddle connection of length going from  conical  point  $P_1$  to
conical point $P_2$. Memorize its direction  (say, the North-West
direction)  and   its   length   $l$.   Then   with  a  ``nonzero
probability'' (understood in  the same sense as above) emitting a
geodesic  from  $P_1$ in one of the  remaining  $n-1$  North-West
directions  we  make it  hit  $P_2$  at  the  distance  $l$. More
rigorously, the Siegel--Veech constant counting configurations of
two parallel saddle  connections of equal length joining $P_1$ to
$P_2$ is nonzero.

The   explicit   formula   for    any    Siegel--Veech   constant
from~\cite{zorich:Eskin:Masur:Zorich} can  be  morally  described as the
follows.  Up   to   some  combinatorial  factor  responsible  for
dimensions, multiplicities  of zeroes and possible symmetries any
Siegel--Veech constant can be obtained as a limit
\begin{equation}
\label{zorich:eq:SV:formula}
\SVc=\lim_{\varepsilon\to 0} \cfrac{1}{\pi\varepsilon^2}\,
\cfrac{\Vol(\text{``$\varepsilon$-neighborhood of the cusp $\cC$
''})}{\Vol\cH^{comp}_1(d_1,\dots,d_\noz)}
\end{equation}
where  $\cC$  is  a  particular  \emph{configuration}  of  saddle
connections or closed geodesics.

Say, as a configuration $\cC$ one can consider a configuration of
two  maximal   cylinders  filled  with  parallel  closed  regular
geodesics of equal lengths. The $\varepsilon$-neighborhood of the
corresponding  cusp  is   the   subset  of  those  flat  surfaces
$S\in\cH_1^{comp}(d_1,\dots,d_\noz)$ which have at least one pair
of cylinders  filled  with  parallel  closed  geodesics of length
shorter than $\varepsilon$.

As  another  example one can consider a  configuration  of  three
parallel    saddle    connections    of    equal    lengths    on
$S\in\cH_1(1,1,4,8)$ joining  zero  $P_1$  of  degree $4$ (having
cone angle $10\pi$)  to  zero $P_2$  of  degree $8$ (having  cone
angle $18\pi$) separated by angles $2\pi,2\pi,6\pi$  at $P_1$ and
by      angles     $6\pi,10\pi,2\pi$      at      $P_2$.      The
$\varepsilon$-neighborhood
$\cH_1^{\varepsilon}(1,1,4,8)\subset\cH_1(1,1,4,8)$    of     the
corresponding  cusp  is  the  subset  of  those flat surfaces  in
$\cH_1(1,1,4,8)$  which  have  at  least  one  triple  of  saddle
connections  as   described   above   of   length   shorter  than
$\varepsilon$.

We explain the origin of the key
formula~\eqref{zorich:eq:SV:formula} in the next section. In
section~\ref{zorich:ss:Multiple:Isometric:Geodesics:and:Principal:Boundary}
we give an explanation of appearance of higher multiplicities.

\paragraph{Other counting problems (after Ya.~Vorobets)}

Having a flat surface $S$ of unit area we have studied above the
number of maximal cylinders $N_{cg}(S,L)$ filled with closed
regular geodesics  of length at most $L$ on  $S$. (In this
setting when we get in some direction several parallel maximal
cylinders of equal perimeter, we count each of them.) In the
paper~\cite{zorich:Vorobets:uniform:bounds} Ya.~Vorobets
considered other counting problems.

In particular, among all maximal periodic cylinders of length at
most $L$ (as above) he counted the number $N_{cg,\sigma}(S,L)$
of those ones, which have area greater than $\sigma$. He also
counted the total sum $N_{area}(S,L)$ of areas of all maximal
cylinders of perimeter at most $L$ and the number $N_x(S,L)$ of
regular periodic geodesics of length at most $L$ passing trough
a given point $x\in S$.

Ya.~Vorobets has also studied how the maximal cylinders filled
with closed geodesics are distributed with respect to their
direction and their area. He considered the induced families of
probability measures on the circle $S^1$, on the unit interval
$[0,1]$ and on their product $S^1\times[0;1]$. Given a subset
$U\subset S_1$, $V\subset [0;1]$, $W\subset S^1\times [0;1]$ the
corresponding measures $dir_L(U), ar_L(V), pair_L(W)$ tell the
proportion of cylinders of bounded perimeter having direction in
$U$, area in $V$, or the pair $(direction,\ area)$ in $W$
correspondingly.

Using the general approach of A.~Eskin and H.~Masur,
Ya.~Vorobets has proved in~\cite{zorich:Vorobets:uniform:bounds}
existence of exact quadratic asymptotics for the counting
functions introduced above. He has computed the corresponding
Siegel--Veech constants in terms of the Siegel--Veech constant
$const_{cg}$ in~\eqref{zorich:eq:Siegel:Veech:const} and found
the asymptotic distributions of directions and areas of the cylinders:

\begin{NNTheorem}[Ya.~Vorobets]
For almost any flat surface $S$ of unit area in any connected
component of any stratum $\cH^{comp}_1(d_1,\dots,d_\noz)$ and
for almost any point $x$ of $S$ one has
$$
\lim_{L\to\infty}\cfrac{N_{cg,\sigma}(S,L)}{\pi L^2}=
            c_{cg,\sigma} \quad
\lim_{L\to\infty}\cfrac{N_{area}(S,L)}{\pi L^2}=c_{area} \quad
\lim_{L\to\infty}\cfrac{N_x(S,L)}{\pi L^2}=
            c_x,
$$
where
$c_{cg,\sigma}=(1-\sigma)^{2g-3+\noz}\cdot const_{cg}\quad$ and
$\quad c_{area}=c_x =\cfrac{const_{cg}}{2g-2+\noz}$.

For almost any flat surface $S$ of unit area one has the following
week convergence of measures:
$$
dir_L\to \varphi \qquad ar_L\to \rho \qquad pair_L\to \varphi\times\rho,
$$
where $\varphi$ is the uniform probability measure on the circle and
$\rho$ is the probability measure on the unit interval $[0;1]$ with the density
$(2g-3+m)(1-x)^{2g-4+m}dx$.
\end{NNTheorem}

\toread{Directional flow}

Theorem  of  S.~Kerckhoff,  H.~Masur  and  J.~Smillie  is  proved
in~\cite{zorich:Kerckhoff:Masur:Smillie}. An example of minimal but  not
uniquely ergodic interval exchange transformations is constructed
by W.~Veech  in~\cite{zorich:Veech:1969} (using different  terminology);
an independent example (also using different  terminology) was constructed
at the same time by V.~I.~Oseledets.
For  flows  such  examples  are   constructed   in   the  paper  of
A.~Katok~\cite{zorich:Katok:1973}        and        developed         by
E.~Sataev in~\cite{zorich:Sataev:1975}. Another example was
discovered by M.~Keane~\cite{zorich:Keane:nonergodic}. For alternative
approach to the study of unique ergodicity of interval exchange transformations
see the paper of M.~Boshernitzan~\cite{zorich:Boshernitzan}.
A  very  nice construction of minimal but  not
uniquely ergodic interval exchange transformations
(in a  language  which is very close  to  the language of
this  paper)  can  be  found  in  the  survey   of  H.~Masur  and
S.~Tabachnikov~\cite{zorich:Masur:Tabachnikov}  or  in   the  survey  of
H.~Masur~\cite{zorich:Masur:Handbook:1B}.

\subsection{Siegel--Veech Formula}
\label{zorich:ss:Siegel:Veech:Formula}

We start from a slight formalization of our  counting problem. As
usual we start with a model case of the flat  torus.  As usual we
assume that our flat torus is  glued from a unit square. We count
closed  regular  geodesics on $\torus$ of a  bounded  length.  To
mimic count of saddle connections we  mark  two  points  $P_1\neq
P_2$ on $\T{2}$ and count geodesic  segments  of  bounded  length
joining $P_1$ and $P_2$.

Our  formalization  consists  in   the   following  construction.
Consider an  auxiliary  Euclidian  plane  $\R{2}$.  Having found a
regular  closed  geodesic  on  $\torus$  we  note  its  direction
$\alpha$ and length $l$ and draw a vector in $\R{2}$ in direction
$\alpha$ having length $l$. We  apply  a  similar construction to
``saddle connections''.  The  endpoints  of corresponding vectors
form two discrete subsets in $\R{2}$ which we  denote by $V_{cg}$
and $V_{sc}$.

\begin{figure}[htb]
\centering
\includegraphics{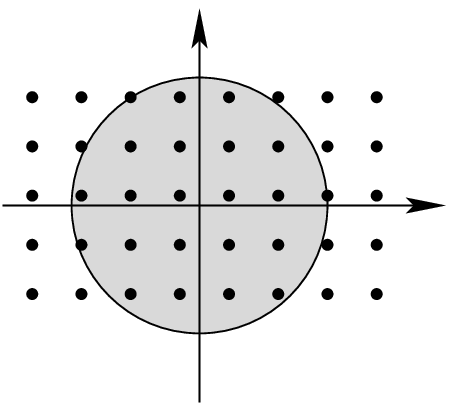}
\includegraphics{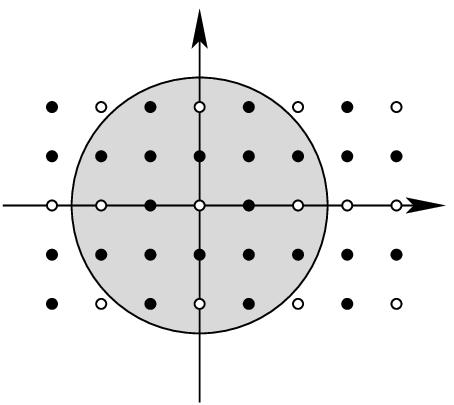}
\vspace{110bp}
\caption{
\label{zorich:fig:V:sc:and:V:cg}
Sets $V_{sc}$ and $V_{cg}$ for the flat torus}
\end{figure}

It is easy to see  that  for  the  torus case a generic choice of
$P_1$ and $P_2$ generates a set $V_{sc}$ which is just  a shifted
square   lattice,   see  Fig.~\ref{zorich:fig:V:sc:and:V:cg}.   The  set
$V_{cg}$ is a subset  of  \emph{primitive} elements of the square
lattice, see  Fig.~\ref{zorich:fig:V:sc:and:V:cg}.  Since  we count only
regular  closed  geodesics which do not turn  many  times  around
themselves we cannot obtain elements of the form  $(k n_1,k n_2)$
with $k, n_1, n_2 \in \Z{}$.

The  corresponding   counting  functions  $N_{sc}(\T{2},L)$   and
$N_{sc}(\T{2},L)$ correspond to the number of element of $V_{sc}$
and $V_{cg}$  correspondingly which get to  a disc of  radius $L$
centered  in  the  origin.  Both functions have  exact  quadratic
asymptotics.  Denoting  by  $\chi_L(v)$,  where  $v\in\R{2}$  the
indicator function of such disc we get
\begin{eqnarray}
\label{zorich:eq:asymptotics:for:T2}
N_{sc}(\T{2},L) &= \sum_{v\in V_{sc}} \chi_L(v) &\sim 1\cdot\pi L^2 \\
N_{cg}(\T{2},L) &= \sum_{v\in V_{cg}} \chi_L(v) &\sim
\cfrac{1}{\zeta(2)}\cdot \pi L^2 \notag
\end{eqnarray}

The   coefficients   in   quadratic   asymptotics   define    the
corresponding   Siegel--Veech    constants   $const_{sc}=1$   and
$const_{cg}=\cfrac{1}{\zeta(2)}=\cfrac{6}{\pi^2}$.   (Note   that
here we count every geodesic twice: once with one orientation and
the other one with the opposite orientation. This explains why in
this normalization we obtain  the  value of $const_{cg}$ twice as
much      as      $const_{cg}$      for     genus     one      in
Table~\ref{zorich:tab:SV:constants}.)

Consider now  a more general  flat surface $S$. Fix the geometric
type of  \emph{configuration}  $\cC$  of  saddle  connections  or
closed geodesics.  By  definition  all saddle connections (closed
geodesics) in $\cC$ are parallel and  have  equal  length.  Thus,
similar  to  the torus case, every  time  we see a collection  of
saddle connections (closed geodesics) of geometric  type $\cC$ we
can associate  to such collection a  vector in $\R{2}$.  We again
obtain a discrete set $V(S)\subset \R{2}$.

Now fix  $\cC$ and apply  this construction to every flat surface
$S$  in  the  stratum  $\cH_1(d_1,\dots,d_\noz)$\index{0H30@$\cH_1(d_1,\dots,d_\noz)$ -- ``unit hyperboloid''}\index{Stratum!in the moduli space}\index{Unit hyperboloid}.  Consider   the
following      operator      $       f      \mapsto      \hat{f}$
generalizing~\eqref{zorich:eq:asymptotics:for:T2}  from  functions  with
compact     support     on    $\R{2}$     to     functions     on
$\cH_1(d_1,\dots,d_\noz)$:
$$
\hat{f}(S)=\sum_{v\in V} f(v)
$$

\begin{NNLemma}[W.~Veech]
The functional\index{0d20@$d\nu_1$ -- volume element on the ``unit hyperboloid''}
$$
f\mapsto \int_{\cH^{comp}_1(d_1,\dots,d_\noz)} \hat{f}(S) d\nu_1
$$
\index{0SL@$SL(2,\R{})$-action on the moduli space}
\index{Action on the moduli space!ofSL@of $SL(2,\R{})$}
is $SL(2,\R{})$-invariant.
\end{NNLemma}

Having proved convergence of the integral above the Lemma follows
immediately from invariance of the  measure  $d\nu_1$\index{0d20@$d\nu_1$ -- volume element on the ``unit hyperboloid''}  under  the
action of $SL(2,\R{})$ and  from  the fact that $V(gS)=gV(s)$ for
any flat surface $S$ and any $g\in SL(2,\R{})$.

Now note  that  there very few $SL(2,\R{})$-invariant functionals
on functions with compact support in $\R{2}$. Actually, there are
two such functionals,  and the other ones are linear combinations
of these two. These two  functionals  are the value of $f(0)$  at
the origin and the integral $\int_{\R{2}} f(x,y) dx\,  dy$. It is
possible to see that the value $f(0)$ at the origin is irrelevant
for the functional from the Lemma above. Hence it is proportional
to $\int_{\R{2}} f(x,y) dx\, dy$.

\begin{NNTheorem}[W.~Veech]
For any function $f:\R{2}\to\R{}$ with compact support one has
\index{Formula!of Siegel--Veech}
\begin{equation}
\label{zorich:eq:Veech:counting:theorem}
\cfrac{1}{\Vol\cH^{comp}_1(d_1,\dots,d_\noz)}
\int_{\cH^{comp}_1(d_1,\dots,d_\noz)} \hat{f}(S) d\nu_1 =
C \int_{\R{2}} f(x,y) dx\, dy
\end{equation}
\end{NNTheorem}

Here   the   constant  $C$   in~\eqref{zorich:eq:Veech:counting:theorem}
\emph{does not  depend} on $f$;  it depends only on the connected
component $\cH^{comp}_1(d_1,\dots,d_\noz)$  and on the  geometric
type $\cC$ of the chosen configuration.

Note that it is an \emph{exact} equality. In particular, choosing
as $f=\chi_L$  the indicator function of a disc  of radius $L$ we
see that for  any given $L\in\R{}_+$ the \emph{average} number of
saddle  connections  not   longer   than  $L$  on  flat  surfaces
$S\in\cH^{comp}_1(d_1,\dots,d_\noz)$ is \emph{exactly} $C\cdot\pi
L^2$, where $C$ does not depend on $L$.

It would  be convenient to  introduce a special notation for such
$\hat f$. Let
$$
\hat{f}_L(S)=\sum_{v\in V} \chi_L(v)
$$
The Theorem  of  Eskin  and  Masur~\cite{zorich:Eskin:Masur} cited above
tells that for large values of  $L$  one  gets  $\hat{f}_L(S)\sim
\SVc  \pi  L^2$  for  almost   all   individual   flat   surfaces
$S\in\cH^{comp}_1(d_1,\dots,d_\noz)$ and that  the  corresponding
constant $\SVc$ coincides with the constant $C$ above.

Formula~\eqref{zorich:eq:Veech:counting:theorem}  can  be   applied   to
$\hat{f}_L$  for  any  value  of $L$. In particular,  instead  of
considering  large  $L$  we   can   choose  a  very  small  value
$L=\varepsilon$.        The        corresponding         function
$\hat{f}_\varepsilon(S)$ counts how many collections of  parallel
$\varepsilon$-short saddle connections (closed geodesics) of  the
type $\cC$ we can find on the flat surface $S$.

Usually there are no such saddle  connections (closed geodesics),
so for  most  flat  surfaces $\hat{f}_\varepsilon(S)=0$. For some
surfaces there is exactly one collection like this. We denote the
corresponding  subset  by  $\cH^{\varepsilon,thick}_1(\cC)\subset
\cH_1(d_1,\dots,d_\noz)$.  Finally,  for the  surfaces  from  the
remaining (very small) subset $\cH^{\varepsilon,thin}_1(\cC)$ one
has several  collections  of  short  saddle  connections  (closed
geodesics) of the type $\cC$. Thus,
$$
\hat{f}_\varepsilon(S)=\begin{cases}
0 & \text{for most of the surfaces } S\\
1 & \text{for } S\in \cH^{\varepsilon,thick}_1(\cC)\\
>1 & \text{for } S\in \cH^{\varepsilon,thin}_1(\cC)
\end{cases}
$$
and   we    can   rewrite~\eqref{zorich:eq:Veech:counting:theorem}   for
$\hat{f}_\varepsilon$ as
\begin{multline*}
\SVc\cdot\pi\varepsilon^2=
\SVc \int_{\R{2}} \chi_\varepsilon(x,y) dx\, dy=\\
=\cfrac{1}{\Vol\cH^{comp}_1(d_1,\dots,d_\noz)}
\int_{\cH^{comp}_1(d_1,\dots,d_\noz)} \hat{f}_\varepsilon(S)\ d\nu_1 =\\
=\cfrac{1}{\Vol\cH^{comp}_1(d_1,\dots,d_\noz)}
\int_{\cH^{\varepsilon, thick}_1(\cC)} 1\  d\nu_1 +\\
+\cfrac{1}{\Vol\cH^{comp}_1(d_1,\dots,d_\noz)}
\int_{\cH^{\varepsilon, thin}_1(\cC)} \hat{f}_\varepsilon(S)\
d\nu_1
\end{multline*}

It can be shown that  though  $\hat{f}_\varepsilon(S)$  might  be
large on  $\cH^{\varepsilon,  thin}_1(\cC)$  the  measure of this
subset is  so small (it is of the  order $\varepsilon^4$ that the
last integral above is negligible in comparison with the previous
one;  namely  it   is   $o(\varepsilon^2)$.  (This  is  a  highly
nontrivial result of Eskin and Masur~\cite{zorich:Eskin:Masur}.)  Taking
into consideration that
$$
\int_{\cH^{\varepsilon, thick}_1(\cC)} 1\  d\nu_1 =
\Vol\cH^{\varepsilon, thick}_1(\cC)=
\Vol\cH^{\varepsilon}_1(\cC)+o(\varepsilon^2)
$$
we can rewrite the chain of equalities above as
$$
\SVc\cdot\pi\varepsilon^2=
\frac{\Vol\cH^{\varepsilon}_1(\cC)}
{\Vol\cH^{comp}_1(d_1,\dots,d_\noz)}
+o(\varepsilon^2)
$$
   %
   %
which is equivalent to~\eqref{zorich:eq:SV:formula}.

\paragraph{Baby Case: Closed Geodesics on the Torus}

As  an  elementary application we can prove  that  proportion  of
primitive  lattice   points   among   all   lattice   points   is
$1/\zeta(2)$. In  other words, applying~\eqref{zorich:eq:SV:formula}  we
can  prove  asymptotic formula~\eqref{zorich:eq:asymptotics:for:T2}  for
the number of primitive lattice points in a disc of  large radius
$L$. As we have seen at  Fig.~\ref{zorich:fig:V:sc:and:V:cg} this number
equals  to  the number $N_{cg}(\T{2},L)$ of families of  oriented
closed geodesics  of length bounded  by $L$ on the standard torus
$\T{2}$.

We want  to  apply~\eqref{zorich:eq:SV:formula}  to  prove the following
formula for the corresponding  Siegel--Veech  constant $c^+_{cg}$
(where  superscript   $+$   indicates   that   we   are  counting
\emph{oriented} geodesics on $\T{2}$).
$$
c^+_{cg}= \lim_{L\to\infty} \cfrac{N_{cg}(\T{2},L)}{\pi L^2}=
\cfrac{1}{\zeta(2)} = \cfrac{6}{\pi^2}
$$

Note that the  moduli  space $\cH_1(0)$ of flat  tori  is a total
space  of  a unit  tangent  bundle to  the  modular surface  (see
Sec.~\ref{zorich:ss:Toy:Example:Family:of:Flat:Tori},
Fig.~\ref{zorich:fig:space:of:flat:tori};   see   also~\eqref{zorich:eq:cd}  in
Sec.~\ref{zorich:ss:SL2R:action:in:geometric:terms}    for     geometric
details). Modular  surface can be considered  as a space  of flat
tori of unit area without choice of direction to the North.

Measure  on  this circle  bundle  disintegrates  to  the  product
measure on  the fiber and the  hyperbolic measure on  the modular
curve.  In  particular,  $\Vol(\cH_1(0))=\pi\cdot  \pi/3$,  where
$\pi/3$ is the hyperbolic area of the modular surface.

Similarly,                 $\Vol(\cH^{\varepsilon}_1(0))=\pi\cdot
\Area(\text{Cusp}(\varepsilon))$,
where
\index{Cusp of the moduli space@Cusp of the moduli space, \emph{see also} Moduli space; principal boundary}
$\text{Cusp}(\varepsilon)$ is  a  subset  of  the modular surface
corresponding  to  those flat  tori  of unit  area  which have  a
geodesic       shorter        than       $\varepsilon$       (see
Fig.~\ref{zorich:fig:space:of:flat:tori}).

Showing that
$\Area(\text{Cusp}(\varepsilon))\approx\varepsilon^2$ we
apply~\eqref{zorich:eq:SV:formula} to get
$$
c_{cg}=\lim_{\varepsilon\to 0}= \cfrac{1}{\pi\varepsilon^2}\
\cfrac{\Area(\text{Cusp}(\varepsilon))}
      {\Area(\text{Modular surface})}=
\cfrac{1}{\pi\varepsilon^2}\
\cfrac{\varepsilon^2+o(\varepsilon^2)}
      {\pi/3} =
\cfrac{1}{2\zeta(2)}.
$$
Note that  the  Siegel--Veech  constant  $c_{cg}$  corresponds to
counting \emph{nonoriented}  closed  geodesics  on $\T{2}$. Thus,
finally we obtain the desired value $c^+_{cg}=2c_{cg}$.

In  the next  section we  give  an idea  of how  one can  compute
$\Vol\cH^{\varepsilon}_1(\cC)$   in   the   simplest   case,   In
Sec.~\ref{zorich:ss:Multiple:Isometric:Geodesics:and:Principal:Boundary}
we describe the  phenomenon  of higher multiplicities and discuss
the   structure   of   typical   cusps  of  the   moduli   spaces
$\cH(d_1,\dots,d_\noz)$\index{0H20@$\cH(d_1,\dots,d_\noz)$ -- stratum in the moduli space}\index{Stratum!in the moduli space}.

\index{Geodesic!closed|)}
\index{Saddle!connection|)}
\index{Geodesic!counting of periodic geodesics|)}
\index{Siegel--Veech!constant|)}

\subsection{Simplest Cusps of the Moduli Space}
\label{zorich:ss:Simplest:Cusps:of:the:Moduli:Space}

In this section  we  consider the  simplest  ``cusp'' $\cC$ on  a
stratum        $\cH(d_1,\dots,d_\noz)$\index{0H20@$\cH(d_1,\dots,d_\noz)$ -- stratum in the moduli space}\index{Stratum!in the moduli space}        and        evaluate
$\Vol\cH^{\varepsilon}_1(\cC)$ for this cusp.  Namely,  we assume
that the  flat surface has at  least two distinct  conical points
$P_1\neq P_2$;  let  $2\pi(d_1+1),  2\pi(d_1+1)$ be corresponding
cone  angles\index{Cone angle}.
As a  \emph{configuration}\index{Configuration of saddle connections or of closed geodesics}
$\cC$  we  consider  a
configuration when  we  have  a  single  saddle connection $\rho$
joining $P_1$ to  $P_2$ and no  other saddle connections  on  $S$
parallel to $\rho$. In our calculation we assume that the conical
points on every $S\in \cH(d_1,\dots,d_\noz)$ have names; we count
only saddle connections joining $P_1$ to $P_2$.

\index{Cusp of the moduli space@Cusp of the moduli space, \emph{see also} Moduli space; principal boundary}

Consider     some      $S\in\cH_1^{\varepsilon,thick}(\cC)\subset
\cH(d_1,\dots,d_\noz)$,  that  is  a  flat surface $S$  having  a
single  saddle  connection joining $P_1$ to $P_2$  which  is  not
longer  than  $\varepsilon$ and  having  no  other  short  saddle
connections or closed geodesics.

We are going to show that  there is a canonical way to shrink the
saddle    connection    on   $S\in\cH_1^{\varepsilon,thick}(\cC)$
coalescing  two  conical  points  into  one.  We shall see  that,
morally, this provides us with an (almost) fiber bundle
\begin{equation}
\label{zorich:cd:Thick:part:over:smaller:stratum}
\begin{CD}
\cH_1^{\varepsilon,thick}(d_1,d_2,d_3,\dots,d_\noz)\\
@VV{\tilde{D}^2_\varepsilon}V\\
\cH_1(d_1+d_2,d_3,\dots,d_\noz)
\end{CD}
\end{equation}
where $\tilde{D}^2_\varepsilon$  is  a  ramified  cover  of order
$(d_1+d_2+1)$   over   a   standard   metric   disc   of   radius
$\varepsilon$.  Moreover,  we  shall  see  that  the  measure  on
$\cH_1^{\varepsilon,thick}(d_1,d_2,d_3,\dots,d_\noz)$
disintegrates  into  a   product   of  the  standard  measure  on
$\tilde{D}^2_\varepsilon$    and    the   natural    measure   on
$\cH_1(d_1+d_2,d_3,\dots,d_\noz)$.  The latter  would  imply  the
following simple answer to our problem:
\begin{multline}
\label{zorich:eq:vol:of:simplest:cusp}
\Vol(\text{``$\varepsilon$-neighborhood of the cusp $\cC$''})=\\
=\Vol\cH^{\varepsilon}_1(\cC)=
\Vol\cH^{\varepsilon}_1(d_1,d_2,d_3,\dots,d_\noz) \approx \\
\approx (d_1+d_2+1)\cdot \pi\varepsilon^2\cdot
\Vol\cH_1(d_1+d_2,d_3,\dots,d_\noz)
\end{multline}

Instead of contracting an  isolated  short saddle connection to a
point we prefer to create it breaking a conical point  $P'\in S'$
of        degree        $d_1+d_2$       on        a       surface
$S'\in\cH_1(d_1+d_2,d_3,\dots,d_\noz)$ into two conical points of
degrees $d_1$ and $d_2$ joined by a short  saddle connections. We
shall see that this surgery is invertible, and thus we  shall get
a coalescing construction.  In the remaining part of this section
we describe this surgery following~\cite{zorich:Eskin:Masur:Zorich}.

\paragraph{Breaking up a Conical Point into Two}

\index{Conical!singularity|(}

We       assume       that        the       initial       surface
$S'\in\cH_1(d_1+d_2,d_3,\dots,d_\noz)$  does  not have  any short
saddle connections or short closed geodesics.

Consider a metric disc  of  a small radius $\varepsilon$ centered
at the  point $P'$, i.e. the set  of points  $Q'$ of the  surface
$S'$ such that Euclidean distance from $Q'$ to the point  $P'$ is
less  than   or   equal   to   $\varepsilon$.   We  suppose  that
$\varepsilon>0$   is   chosen   small   enough,   so   that   the
$\varepsilon$-disc does not  contain  any other conical points of
the metric; we assume also, that the disc which we defined in the
metric  sense  is  homeomorphic  to  a  topological  disc.  Then,
metrically our disc has a structure of a regular cone with a cone
angle\index{Cone angle}
$2\pi(d_1+d_2+1)$;  here  $d_1+d_2$  is the multiplicity of
the zero $P'$. Now cut the chosen disc (cone) out of the surface.
We shall  modify the flat  metric inside it preserving the metric
at the boundary,  and then paste  the modified disc  (cone)  back
into the surface.

\begin{figure}[ht]
  %
\includegraphics{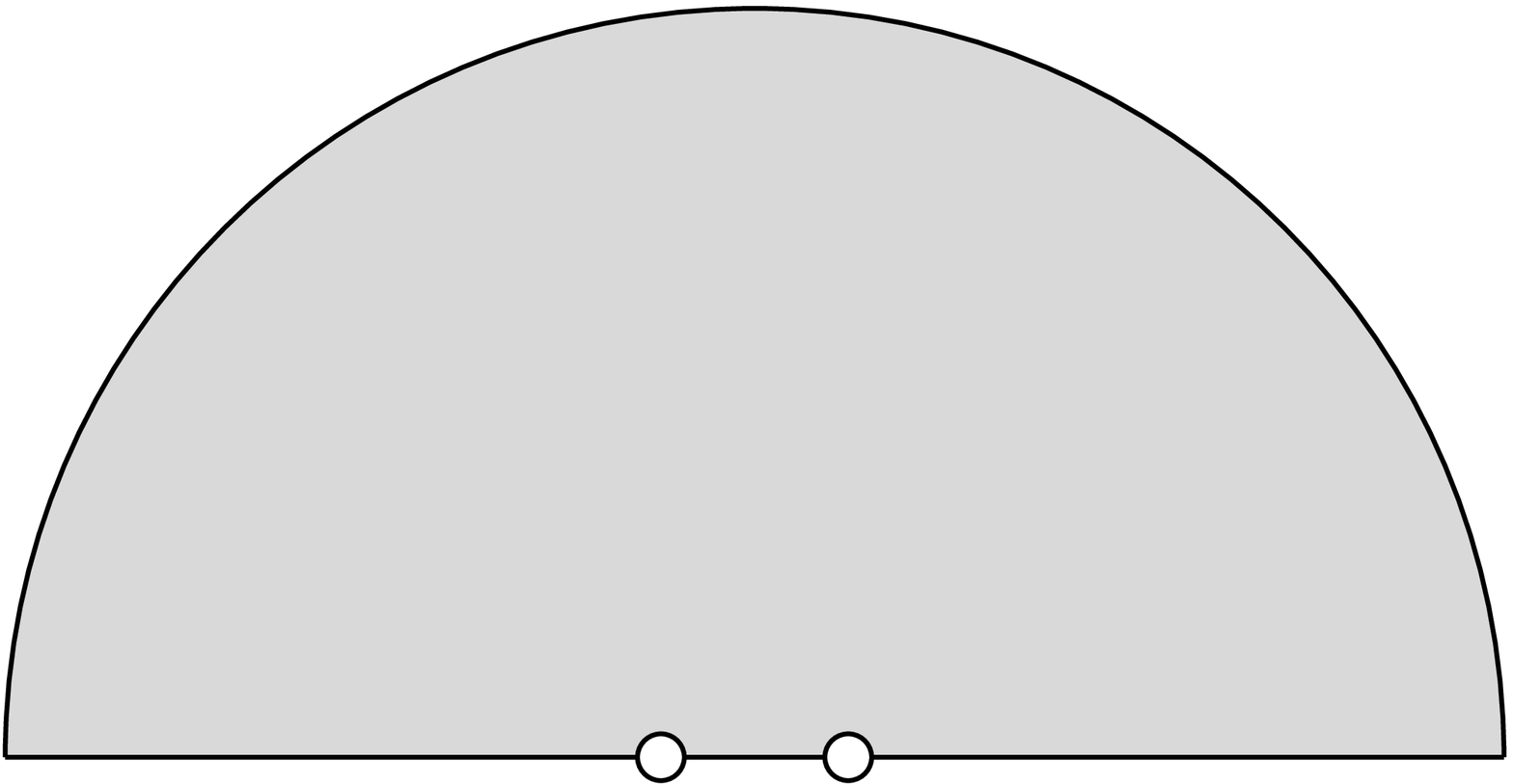}
\includegraphics{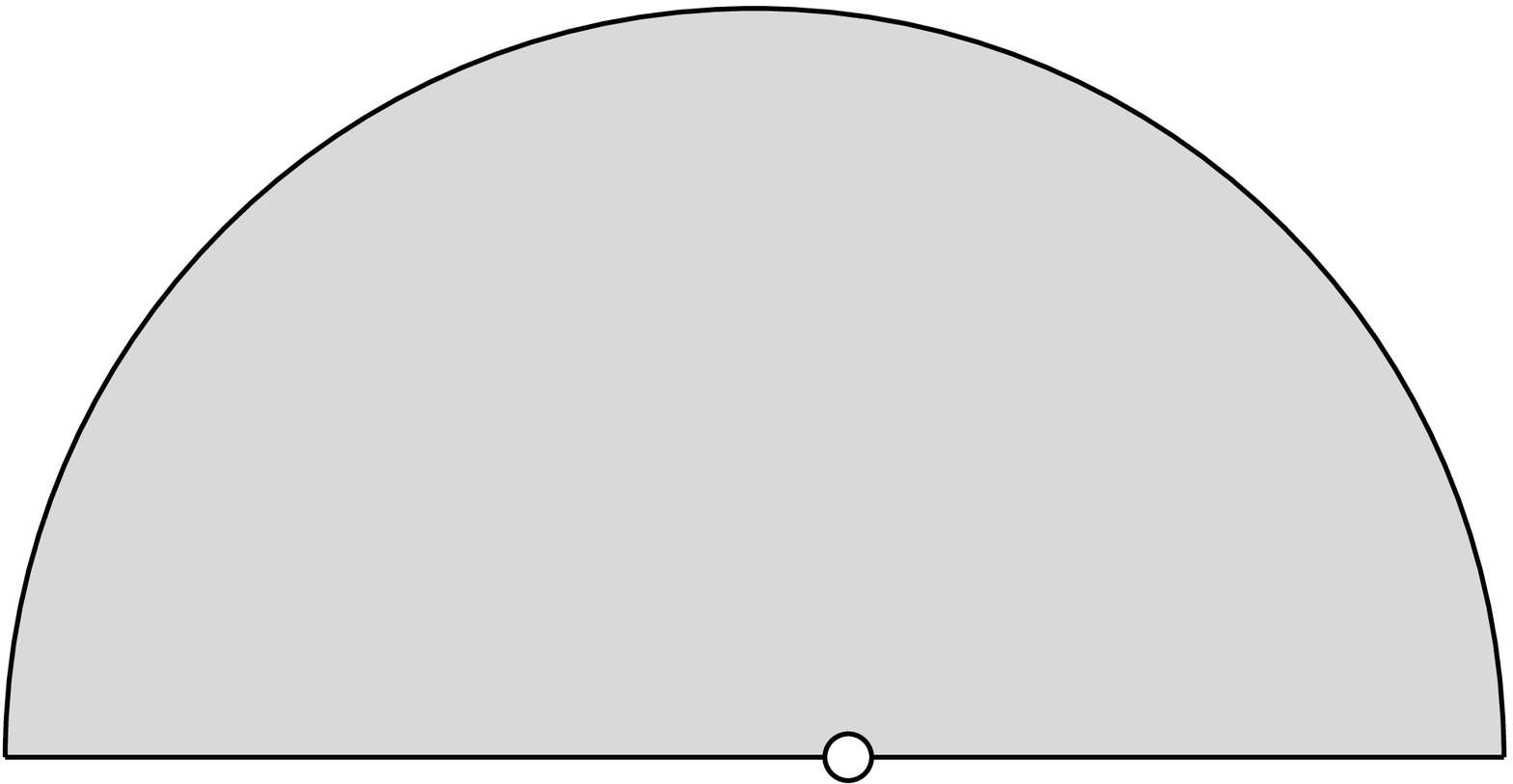}
\includegraphics{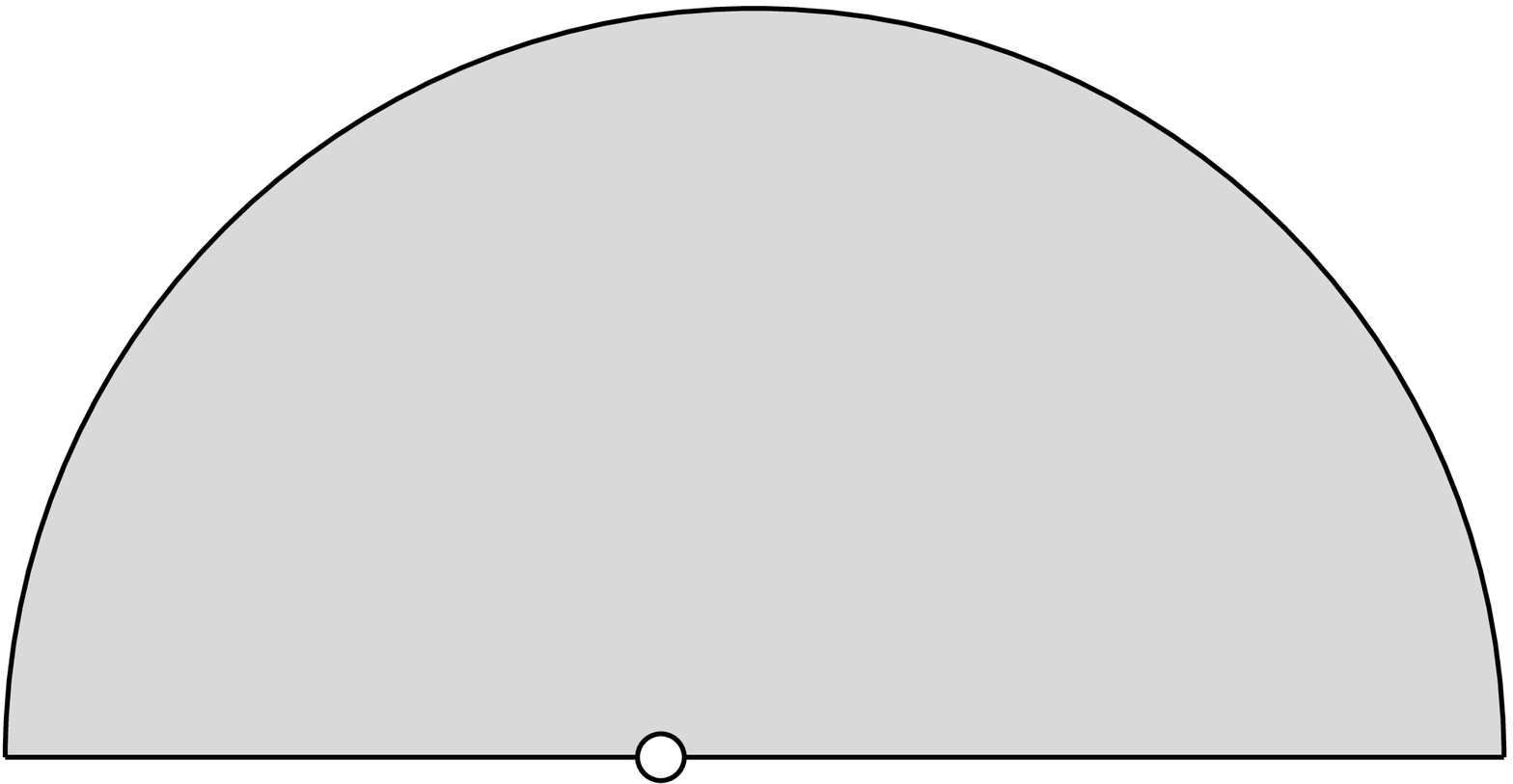}
\includegraphics{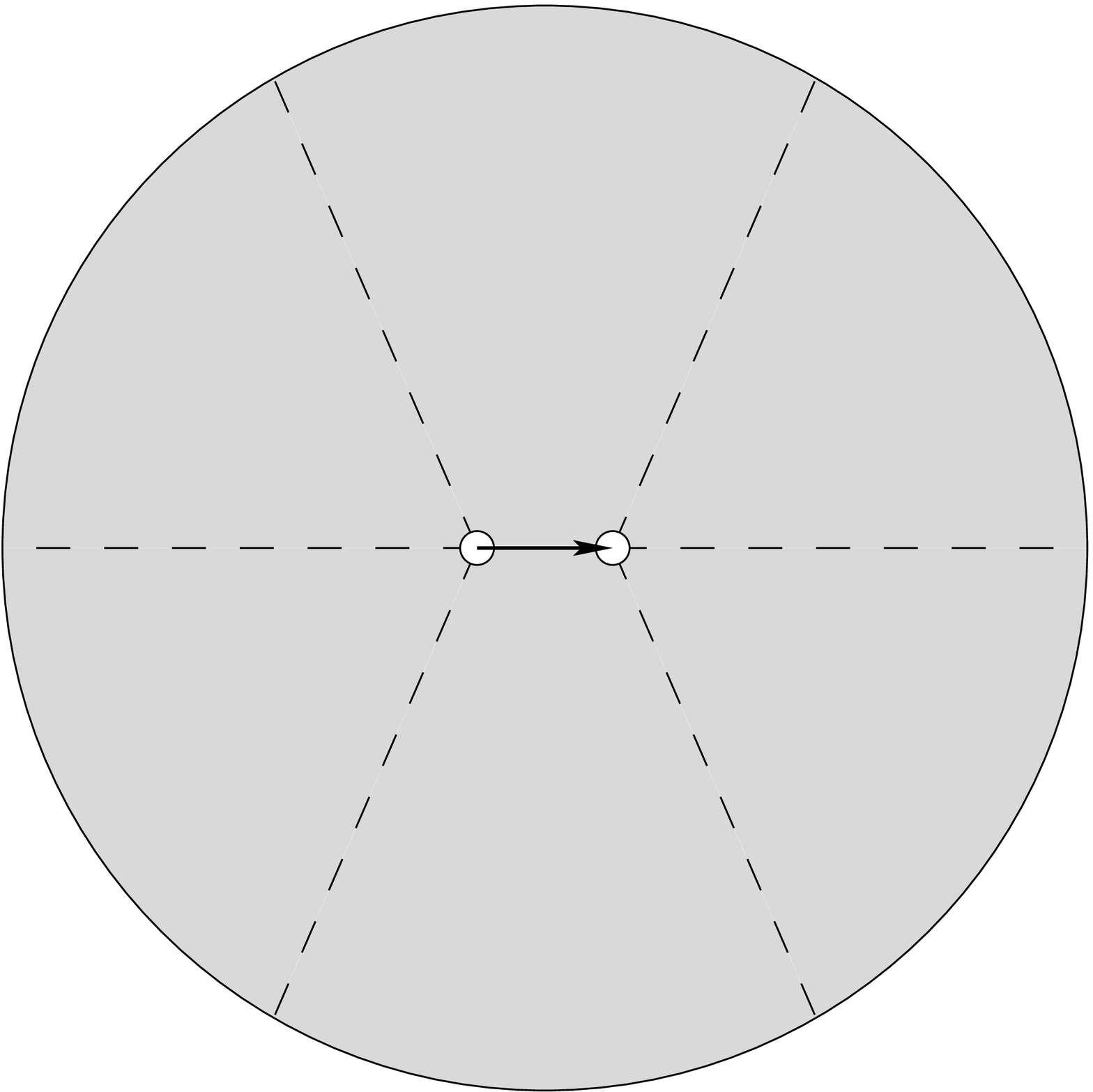}
\includegraphics{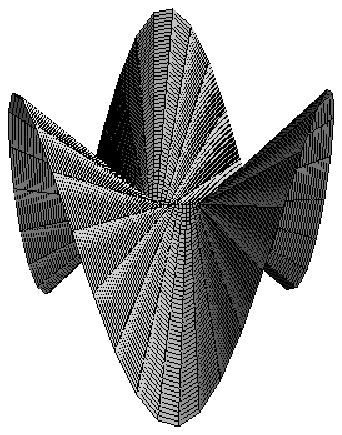}
    %
%
%
\begin{picture}(-169,-140)(-169,-140)
\put(10,-30) 
 {\begin{picture}(0,0)(0,0)
 \put(-145,-148){$\scriptstyle 2\delta$}
 \put(-160,-197){$\scriptstyle \varepsilon+\delta$}
 \put(-130,-197){$\scriptstyle \varepsilon-\delta$}
 \put(-160,-247){$\scriptstyle \varepsilon-\delta$}
 \put(-130,-247){$\scriptstyle \varepsilon+\delta$}
 \put(-39,-197){$\scriptstyle 2\delta$}
 \put(-65,-197){$\scriptstyle \varepsilon+\delta$}
 \put(-19,-197){$\scriptstyle \varepsilon+\delta$}
 \put(-64,-173){$\scriptstyle \varepsilon-\delta$}
 \put(-19,-173){$\scriptstyle \varepsilon-\delta$}
 \put(-65,-217){$\scriptstyle \varepsilon-\delta$}
 \put(-17,-217){$\scriptstyle \varepsilon-\delta$}
\end{picture}}
\end{picture}
\vspace{150bp} 
\caption{
\label{zorich:fig:breaking:up:a:zero}
Breaking up a conical point into two
(after~\cite{zorich:Eskin:Masur:Zorich}).
}
\end{figure}

Our cone  can be glued  from $2(k_i+1)$ copies of standard metric
half-discs of the radius $\varepsilon$,  see  the  picture at the
top  of  Fig.~\ref{zorich:fig:breaking:up:a:zero}.   Choose  some  small
$\delta$, where  $0<\delta<\varepsilon$  and  change  the  way of
gluing  the  half-discs  as  indicated on the bottom  picture  of
Fig.~\ref{zorich:fig:breaking:up:a:zero}. As patterns we  still  use the
standard  metric  half-discs,  but  we move slightly  the  marked
points on their diameters. Now we  use  two  special  half-discs;
they  have two  marked  points on the  diameter  at the  distance
$\delta$ from the center of the half disc. Each of  the remaining
$2k_i$  half-discs  has  a  single marked point at  the  distance
$\delta$ from the center of the half-disc. We are alternating the
half-discs with the marked point  moved  to the right and to  the
left from the  center.  The picture  shows  that all the  lengths
along identifications are matching; gluing the half-discs in this
latter way  we obtain a topological disc with  a flat metric; now
the flat metric has two cone-type  singularities  with  the  cone
angles $2\pi(d_1+1)$ and $2\pi(d_2+1)$.

Note that a  small tubular neighborhood  of the boundary  of  the
initial  cone   is   isometric   to   the  corresponding  tubular
neighborhood of the boundary of the resulting object. Thus we can
paste it back into the surface. Pasting it back we can turn it by
any angle $\varphi$, where $0\le\varphi< 2\pi(d_1+d_2+1)$.

We described how to break up a zero of multiplicity  $d_1+d_2$ of
an  Abelian   differential  into  two  zeroes  of  multiplicities
$d_1,d_2$. The construction  is local; it is parameterized by the
two free real parameters (actually, by one complex parameter): by
the small distance $2\delta$  between  the newborn zeroes, and by
the direction $\varphi$ of the short geodesic segment joining the
two newborn zeroes. In particular, as a parameter  space for this
construction  one  can  choose  a  ramified  covering  of  degree
$d_1+d_2+1$ over a standard metric disc of radius $\varepsilon$.

\index{Conical!singularity|)}

\subsection{Multiple Isometric Geodesics and Principal Boundary
of the Moduli Space}
\label{zorich:ss:Multiple:Isometric:Geodesics:and:Principal:Boundary}

In  this section  we  give an explanation  of  the phenomenon  of
higher multiplicities,  we consider typical degenerations of flat
surfaces  and  we  discuss  how  can  one  use  configurations of
parallel saddle connections or closed geodesics  to determine the
orbit of a flat surface $S$.

\index{Geodesic!closed|(}
\index{Saddle!connection|(}

\paragraph{Multiple Isometric Saddle Connections}

Consider a collection of saddle connections  and closed geodesics
representing      a     basis      of      relative      homology
$H_1(S,\{P_1,\dots,P_\noz\};\Z{})$ on a flat surface $S$. Recall,
that any geodesic on $S$ goes in some  constant direction. Recall
also that by Convention~\ref{zorich:conv:flat:surface} any flat  surface
is endowed with a distinguished direction to the North, so we can
place  a  compass  at any point  of  $S$  and  determine in which
direction  goes  our  geodesic.  Thus, every closed  geodesic  or
saddle        connection        determines        a        vector
$\vec{v}_j\in\R{2}\simeq\C{}$ which  goes  in  the same direction
and have the same  length  as the corresponding geodesic element.
Finally recall  that  collection  of planar vectors $\{\vec{v}_1,
\dots,  \vec{v}_{2g+\noz-1}\}$  considered   as  complex  numbers
provide    us    with    a    local    coordinate    system    in
$\cH(d_1,\dots,d_\noz)$\index{0H20@$\cH(d_1,\dots,d_\noz)$ -- stratum in the moduli space}\index{Stratum!in the moduli space}.  In  complex-analytic   language   these
coordinates are the
\emph{relative periods}\index{Period}
of holomorphic 1-form  representing  the flat surface $S$, namely
$\vec{v}_j=\int_{\rho_j}\omega$,    where    $\rho_j$   is    the
corresponding  geodesic  element  (saddle  connection  or  closed
geodesic).

\begin{figure}[htb!]
%
 %
 %
 \includegraphics{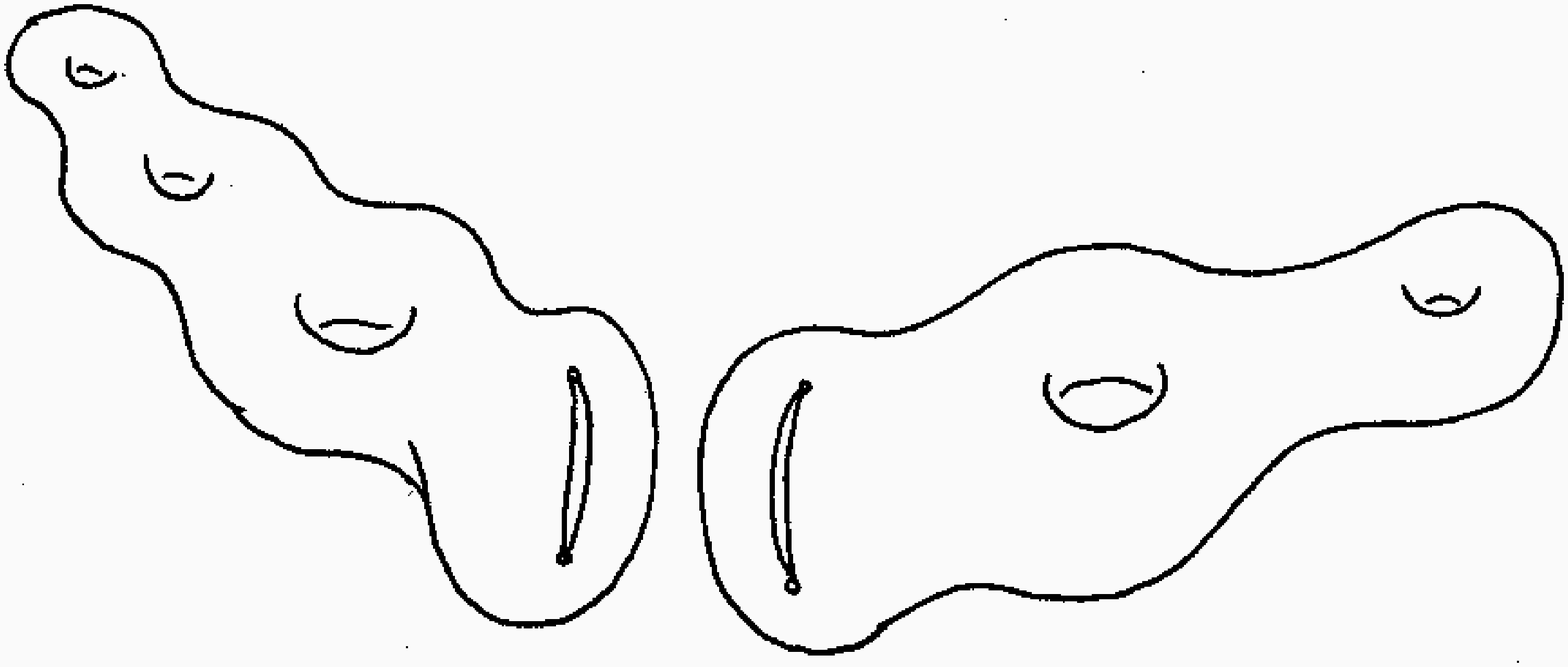}
 \includegraphics{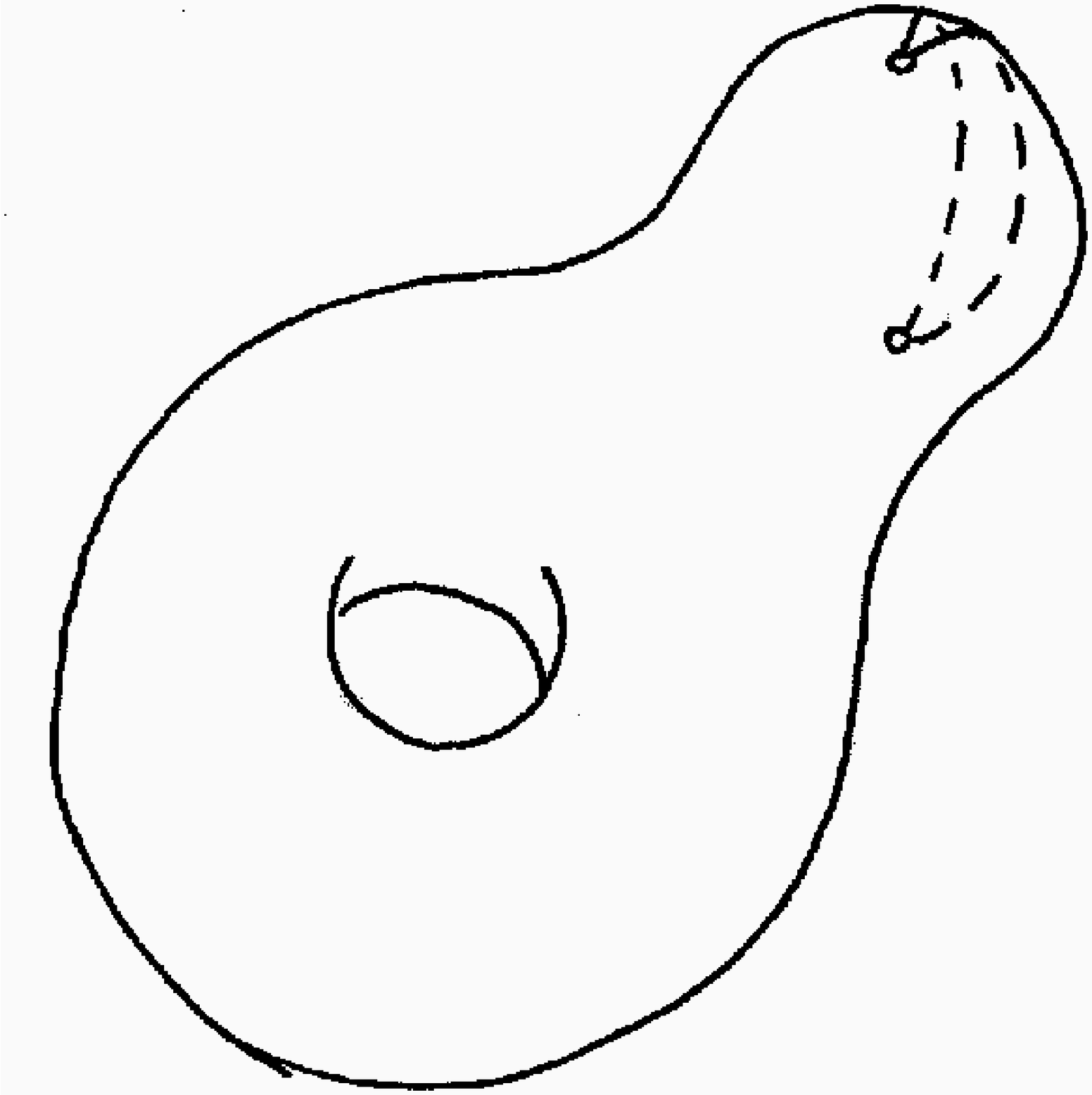}
 %
 %
 \includegraphics{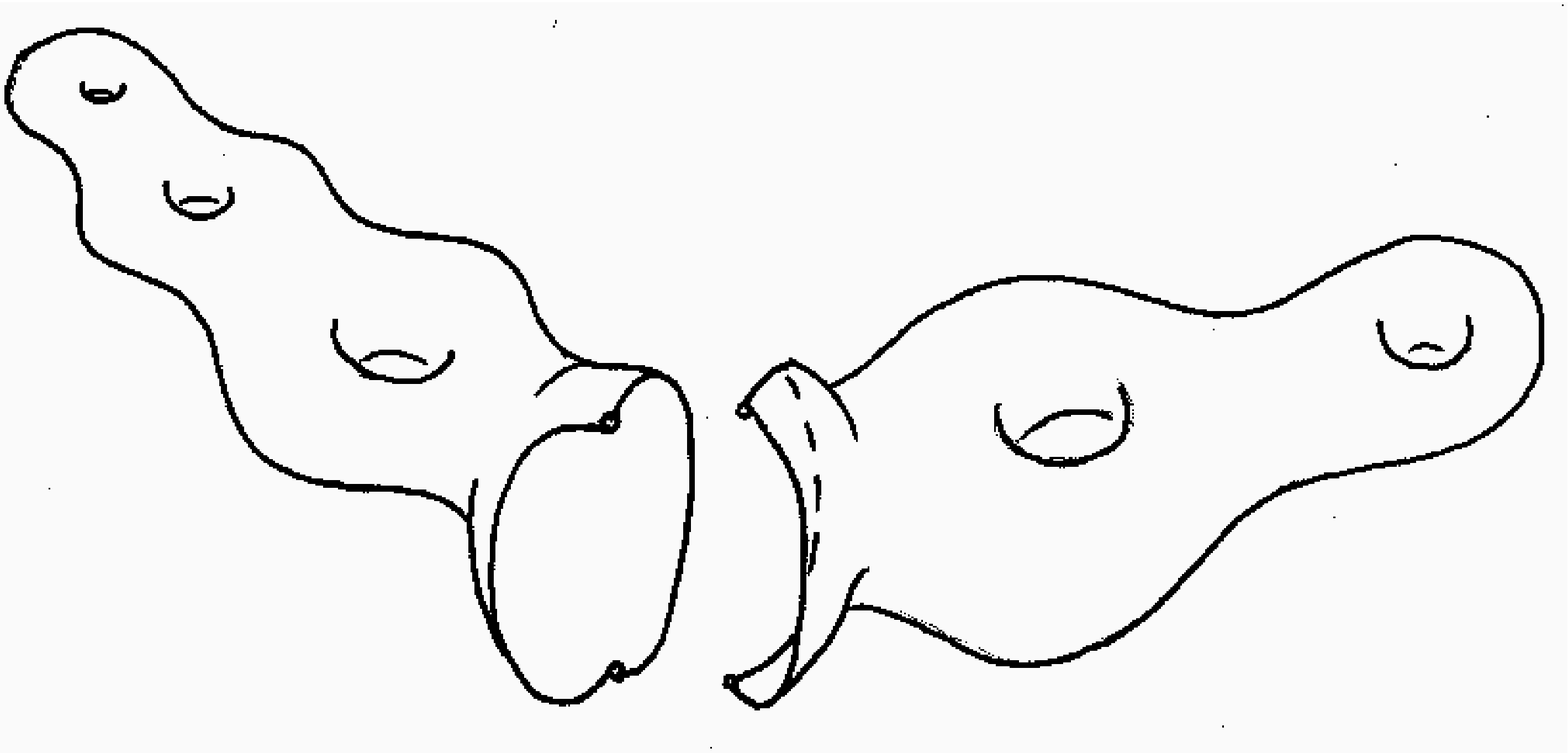}
 \includegraphics{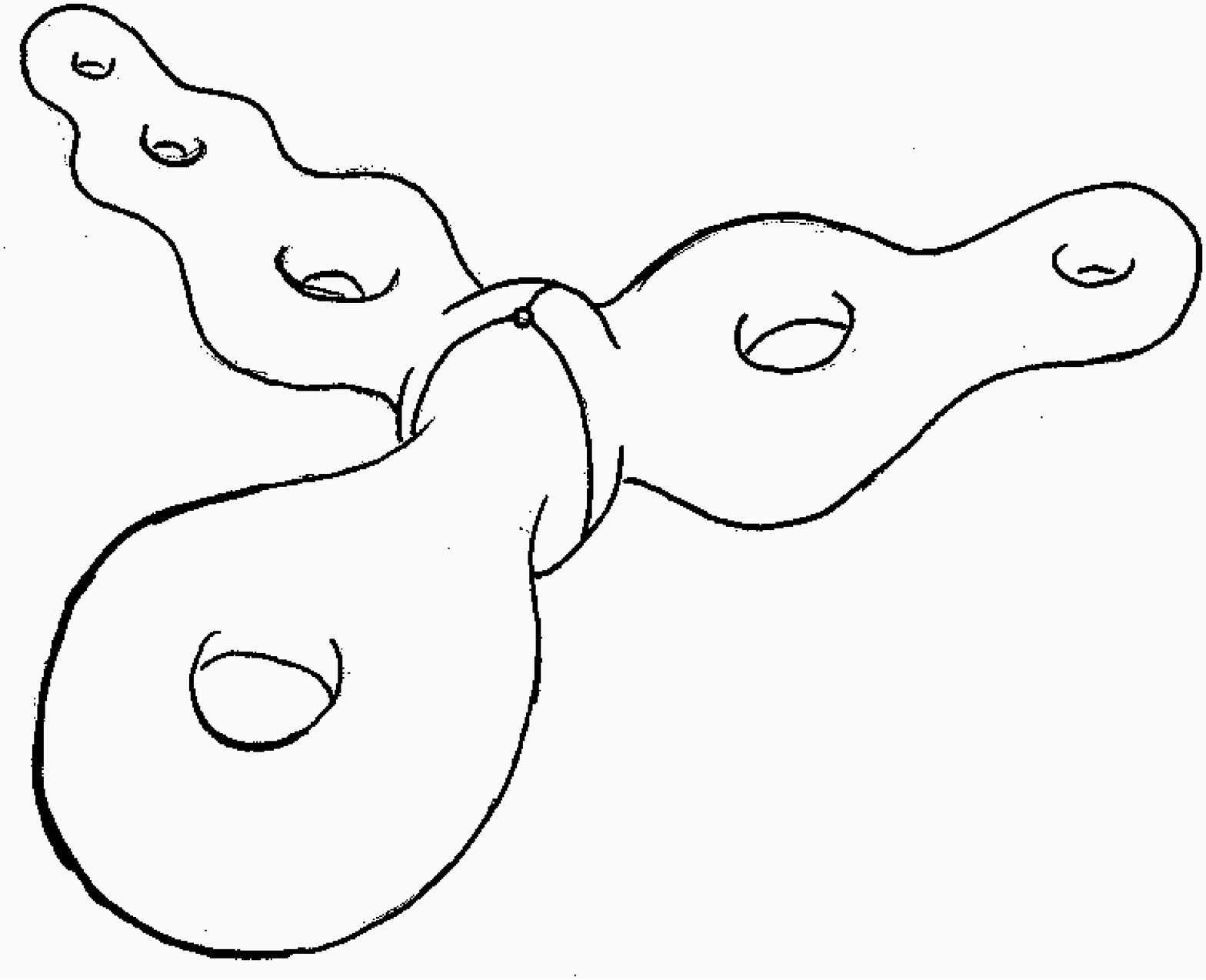}
%
%
\begin{picture}(-160,0)(-160,0)
\put(0,0){
 \begin{picture}(0,0)(0,0) 
 \put(-153,-36){$S'_3$}
 \put(-6,-46){$S'_2$} 
 \end{picture}
}
\put(5,0){
 \begin{picture}(0,0)(0,0) 
 \put(-94,-127){$S'_1$}
 \end{picture}
}
\put(0,0){
 \begin{picture}(0,0)(0,0) 
 \put(20,-12){$S'_3$}
 \put(155,-22){$S'_2$}
 \end{picture}
}
\put(0,0){
 \begin{picture}(0,0)(0,0) 
 \put(80,-180){$S'_1$}
 \put(20,-117){$S'_3$}
 \put(145,-126){$S'_2$}
 \put(84,-140){$\scriptstyle \gamma_1$}
 \put(70,-140){$\scriptstyle \gamma_2$}
 \put(79,-125){$\scriptstyle \gamma_3$}
 \end{picture}
}
\end{picture}
\vspace{195bp} 
\caption{
\label{zorich:fig:dist:sad:mult}
Multiple homologous saddle connections; topological picture}
\end{figure}

We  say  that two geodesic elements $\gamma_1, \gamma_2$  (saddle
connections or closed geodesics) are \emph{homologous}  on a flat
surface  $S$  if  they  determine  the  same homology classes  in
$H_1(S,\{P_1,\dots,P_\noz\};\Z{})$. In other words, $\gamma_1$ is
homologous  to  $\gamma_2$  if  cutting  $S$  by  $\gamma_1$  and
$\gamma_2$ we break the surface $S$ into two pieces. For example,
the saddle  connections  $\gamma_1,  \gamma_2,  \gamma_3$  on the
right surface  at  the bottom of Fig.~\ref{zorich:fig:dist:sad:mult} are
homologous.

The  following  elementary observation is very important for  the
sequel. Since the holomorphic 1-form $\omega$ representing $S$ is
closed, homologous geodesic elements $\gamma_1, \gamma_2$  define
the same period:
$$
\int_{\gamma_1}\omega = \vec{v} = \int_{\gamma_2}\omega
$$

We  intensively  used  the  following  result   of  H.~Masur  and
J.~Smillie~\cite{zorich:Masur:Smillie}  in  our  considerations  in  the
previous section.

\begin{NNTheorem}[H.~Masur, J.~Smillie]
There is a  constant $M$ such that for all $\varepsilon,\kappa>0$
the      subset      $\cH^\varepsilon_1(d_1,\dots,d_\noz)$     of
$\cH_1(d_1,\dots,d_\noz)$\index{0H30@$\cH_1(d_1,\dots,d_\noz)$ -- ``unit hyperboloid''}\index{Stratum!in the moduli space}\index{Unit hyperboloid}  consisting  of  those  flat  surfaces,
which have a saddle  connection  of length at most $\varepsilon$,
has volume at most $M\varepsilon^2$.

The volume  of the set of flat surfaces  with a saddle connection
of  length  at most  $\varepsilon$  and  a  nonhomologous  saddle
connection   with   length  at   most   $\kappa$   is   at   most
$M\varepsilon^2\kappa^2$.
\end{NNTheorem}

Morally, this Theorem says that a  subset  corresponding  to  one
complex parameter with  norm at most $\varepsilon$ has measure of
order $\varepsilon^2$, and a subset corresponding  to two complex
parameter with norm at most  $\varepsilon$  has  measure of order
$\varepsilon^4$. In  particular,  this theorem implies that $\Vol
\cH_1^{\varepsilon,thin}(d_1,\dots,d_\noz)$  is  of   the   order
$\varepsilon^4$.

In the previous section we considered the subset of flat surfaces
$S\in\cH_1^\varepsilon(d_1,d_2,d_3,\dots,d_\noz)$ having a single
short  saddle  connection  joining  zeroes of degrees  $d_1$  and
$d_2$.  We  associated  to   such   surface  $S$  a  new  surface
$S'\in\cH_1^\varepsilon(d_1+d_2,d_3,\dots,d_\noz)$ in the smaller
stratum. Note that, morally, surfaces $S$ and $S'$  have the same
periods with a  reservation  that $S$  has  one more period  than
$S'$:  the  extra small period represented by  our  short  saddle
connection.

Metrically surfaces  $S$ and $S'$ are  almost the same:  having a
surface $S$  we know how  to contract our short saddle connection
to a point; having a  surface  $S'$ and an abstract short  period
$\vec{v}\in\C{}\simeq\R{2}$   we   know   how   to   break    the
corresponding zero on  $S'$ into two  zeroes joined by  a  single
short saddle  connection  realizing  period  $\vec{v}$.  (In  the
latter construction we have some additional  discrete freedom: we
can  break  the zeroes  in  direction  $\vec{v}$  in  $d_1+d_2+1$
different ways.)

Our construction does not work when  we  have  two  nonhomologous
short geodesic elements on the  surface  $S$. But we do not  care
since according to the Theorem of  H,~Masur  and  J.~Smillie  the
subset   $\cH_1^{\varepsilon,thin}(d_1,\dots,d_\noz)$   of   such
surfaces has very small measure (of the order $\varepsilon^4$).

Now  consider  a  slightly  more general surgery  represented  by
Fig.~\ref{zorich:fig:dist:sad:mult}.   We   take  three   distinct  flat
surfaces, we  break  a zero on each of them as it was done in the
previous section. We do it  coherently  using  the same direction
and the same  distance  $\delta$ on  each  surface (left part  of
Fig.~\ref{zorich:fig:dist:sad:mult}). Then  we  slit  each surface along
the newborn  saddle connection and glue  the surfaces in  a close
compound   surface   as  indicated   on   the   right   part   of
Fig.~\ref{zorich:fig:dist:sad:mult}.

The resulting surface has three short parallel saddle connections
of equal  length.  By  construction  they  are \emph{homologous}:
cutting  along  any two of them  we  divide the surface into  two
parts.  Thus,  the  resulting  surface again has only  one  short
period!  Note  that  a  complete  collection  of  periods  of the
compound surface can be obtained as disjoint union
$$
\text{periods of }S'_1 \sqcup
\text{periods of }S'_2 \sqcup
\text{periods of }S'_3 \sqcup
\text{newborn short period}
$$
Hence, any flat surface $S$  with  three  short homologous saddle
connections and no other short geodesic elements can be viewed as
a nonconnected flat  surface  $S'_1\sqcup S'_2\sqcup S'_3$ plus a
memory about the  short period $\vec{v}\in\C{}$ which we use when
we break the zeroes (plus  some  combinatorial  arbitrariness  of
finite order).

The moduli  space of disconnected  flat surfaces of the same type
as $S'_1\sqcup S'_2\sqcup S'_3$  has  one dimension less than the
original  stratum   $\cH(d_1,\dots,d_\noz)$\index{0H20@$\cH(d_1,\dots,d_\noz)$ -- stratum in the moduli space}\index{Stratum!in the moduli space}.  Our  considerations
imply        the        following        generalization        of
formula~\eqref{zorich:eq:vol:of:simplest:cusp}   for   the   volume   of
``$\varepsilon$-neighborhood'' of the corresponding cusp:
\index{Cusp of the moduli space@Cusp of the moduli space, \emph{see also} Moduli space; principal boundary|(}
\begin{multline*}
   %
\Vol(\text{``$\varepsilon$-neighborhood of the cusp $\cC$''})=
\Vol\cH^{\varepsilon}_1(\cC)=\\
k \cdot\Vol\cH_1(\text{stratum of } S'_1)
\cdot\Vol\cH_1(\text{stratum of } S'_2)
\cdot\Vol\cH_1(\text{stratum of } S'_3) \cdot \pi\varepsilon^2
\end{multline*}
where  the  factor  $k$  is  an  explicit  constant  depending on
dimensions,   possible    symmetries,   and   combinatorics    of
multiplicities of zeroes. In particular, we get a subset of order
$\varepsilon^2$.

Now for  any  flat  surface  $S_0\in\cH(d_1,\dots,d_\noz)$ we can
state a counting problem, where we shall count  only those saddle
connections which appear in
configuration\index{Configuration of saddle connections or of closed geodesics}
$\cC$  of  triples  of
homologous saddle  connections breaking $S_0$ into three surfaces
of  the  same  topological  and geometric types as  $S'_1,  S'_2,
S'_3$.    Applying    literally    same    arguments    as     in
Sec.~\ref{zorich:ss:Siegel:Veech:Formula}
and~\ref{zorich:ss:Simplest:Cusps:of:the:Moduli:Space} we  can show that
such number of triples of homologous saddle connections of length
at most $L$ has quadratic growth rate and  that the corresponding
Siegel--Veech constant  $c(\cC)$  can  be  expressed  by the same
formula as above:
\begin{multline*}
   %
\SVc=\lim_{\varepsilon\to 0} \cfrac{1}{\pi\varepsilon^2}\,
\cfrac{\Vol(\text{``$\varepsilon$-neighborhood of the cusp $\cC$
''})}{\Vol\cH^{comp}_1(d_1,\dots,d_\noz)}=\\
=(\text{combinatorial factor})\
\cfrac{\prod_{k=1}^n\Vol\cH_1(\text{stratum of }S'_k)}
{\Vol\cH^{comp}_1(d_1,\dots,d_\noz)}
\end{multline*}
   %

\paragraph{Principal Boundary of the Strata}

\index{Moduli space!principal boundary of the moduli space@principal boundary of the moduli space, \emph{see also} Cusp of the moduli space|(}

The results above  give a description of typical degenerations of
flat surfaces.  A flat surface gets close to  the boundary of the
moduli  space  when  some  geodesic element (or a  collection  of
geodesic  elements)   get   short.  Morally,  we  have  described
something like ``faces'' of  the  boundary, while there are still
``edges'', etc, representing degenerations of higher codimension.
Flat surfaces which are close to this ``principal boundary'' of a
stratum $\cH(d_1,\dots,d_\noz)$\index{0H20@$\cH(d_1,\dots,d_\noz)$ -- stratum in the moduli space}\index{Stratum!in the moduli space} have the following structure.

If the short geodesic element is a saddle  connection joining two
distinct  zeroes,  then   the  surface  looks  like  the  one  at
Fig.~\ref{zorich:fig:dist:sad:mult}.  It  can  be decomposed to  several
flat surfaces with slits glued cyclically  one  to  another.  The
boundaries  of  the slits form short saddle  connections  on  the
compound surface. All these saddle connections join the same pair
of  points;  they  have  the  same  length  and  direction.  They
represent  homologous  cycles  in  the  relative  homology  group
$H_1(S,\{P_1,\dots,P_\noz\};\Z{})$.

\begin{figure}[ht]  
\centering
\includegraphics{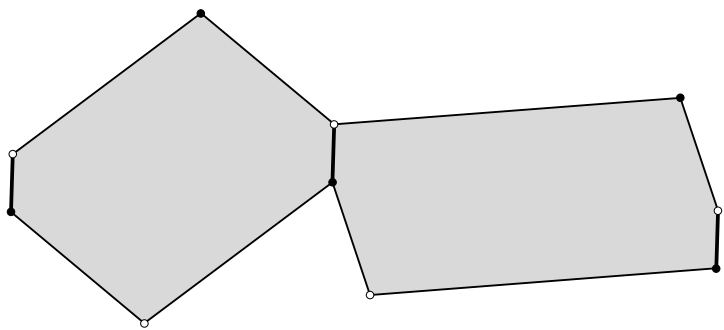}
\vspace{120bp}
\caption{
\label{zorich:fig:metric:surface:mult:2}
Flat  surface  with two short homologous saddle connections.  Any
small  deformation  of this  surface  also has  a  pair of  short
homologous saddle connections. }
\end{figure}

Figure~\ref{zorich:fig:metric:surface:mult:2} represents  a flat surface
$S_0\in\cH(1,1)$  unfolded  to a  polygon.  The  two  short  bold
segments represent  two homologous saddle connections. The reader
can easily check that on any surface $S\in\cH(1,1)$ obtained from
$S_0$  by a small  deformation  one  can  find a  pair  of  short
parallel saddle connections of equal length. Cutting $S$ by these
saddle connections we get a pair of tori with slits.

\begin{figure}[ht]
\centering
\includegraphics{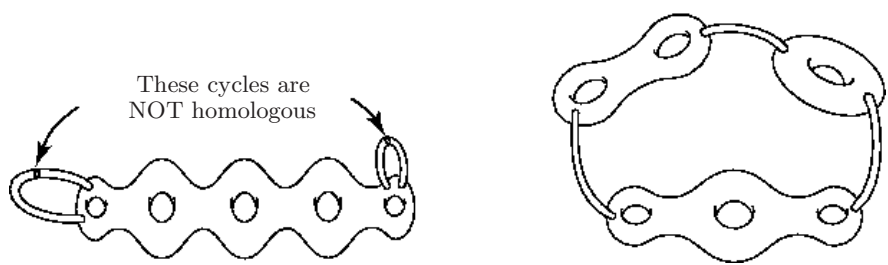}
\begin{picture}(0,0)(0,0)
\put(-120,-30){\ These cycles are}
\put(-120,-40){NOT homologous}
\end{picture}
\vspace{100bp}
\caption{
\label{zorich:fig:dance:ugly:drawing}
A typical (on the right) and nontypical (on the left)
degenerations of a flat surface. Topological picture
}
\end{figure}

We  did  not discuss in the  previous  section the case when  the
short geodesic element is a regular closed geodesic  (or a saddle
connection joining a conical point to  itself).  Morally,  it  is
similar to the case of saddle  connections,  but  technically  it
slightly  more  difficult. A  flat  surface  near  the  principal
boundary   of   this   type   is  presented  on  the   right   of
Fig.~\ref{zorich:fig:dance:ugly:drawing}.

\begin{figure}[htb]
\centering

\special{
psfile=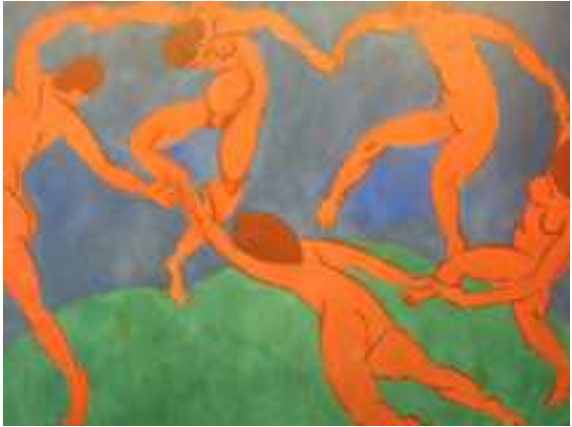    
hscale=100 
vscale=100 
voffset=-157 
hoffset=57   
}
\vspace{155bp}
\caption{
\label{zorich:fig:dance}
Flat surface near the principal boundary$^5$ 
}
\end{figure}

Similar to the  case of saddle  connections, the surface  can  be
also  decomposed   to  a  collection  of  well-proportional  flat
surfaces $S'_1, \dots, S'_n$ of lower genera. Each surface $S'_k$
has a pair of holes. Each  of these holes is realized by a saddle
connection joining a  zero to itself. The surfaces are cyclically
glued to a  ``necklace'', where two neighboring surfaces might be
glued directly or by a narrow cylinder. Since the waist curves of
all these  cylinders  and  all  saddle  connections  representing
boundaries of  surfaces  $S'_k$ are homologous, the corresponding
closed geodesics on  $S$ are parallel  and have equal  length.  A
more artistic\footnotemark\ image of a surface, which is located closed  to the
boundary of a stratum is represented on Fig.~\ref{zorich:fig:dance}.

The surface  on  the left of Fig.~\ref{zorich:fig:dance:ugly:drawing} is
close  to  an ``edge''of the moduli  space  in the sense that  it
represents  a  ``nontypical''  degeneration:  a  degeneration  of
codimension  two.  This  surface  has  two  nonhomologous  closed
geodesics  shorter  than  $\varepsilon$.  Due to the  Theorem  of
H.~Masur    and    J.~Smillie    cited    above,    the    subset
$\cH_1^{\varepsilon,  thin}(d_1,\dots,d_\noz)$  of such  surfaces
has measure of the order $\varepsilon^4$.

\paragraph{Configurations  of  Saddle  Connections and of  Closed
Geodesics as Invariants of Orbits}

\index{Configuration of saddle connections or of closed geodesics}

Consider a  flat surface $S_0$ and  consider its orbit  under the
action of  $SL(2,\R{})$. It is  very easy to construct this orbit
locally for those elements of  the  group  $GL^+(2,\R{})$ which are
close  to  identity.  It  is  a  fairly  complicated  problem  to
construct this  orbit  globally  in $\cH(d_1,\dots,d_\noz)$\index{0H20@$\cH(d_1,\dots,d_\noz)$ -- stratum in the moduli space}\index{Stratum!in the moduli space} and to
find its closure. Ergodic  theorem  of H.~Masur and W.~Veech (see
Sec.~\ref{zorich:ss:Action:of:SL2R:on:the:Moduli:Space}) tells that  for
almost every surface  $S_0$  the closure  of  the orbit of  $S_0$
coincides  with   the   embodying   connected  component  of  the
corresponding stratum. But  for some surfaces the closures of the
orbits are smaller.  Sometimes it is possible to distinguish such
surfaces  looking   at  the  configurations  of  parallel  closed
geodesics and  saddle connections. Say, for \emph{Veech surfaces}
which will be discussed in Sec.~\ref{zorich:ss:Veech:Surfaces} the orbit
of  $SL(2,\R{})$  is  already  closed  in  the stratum, so  Veech
surfaces  are  very  special.  This  property  has  an  immediate
reflection  in  behavior of parallel closed geodesics and  saddle
connections: as  soon as we have a saddle  connection or a closed
geodesic  in  some  direction  on  a
Veech  surface\index{Surface!Veech},
\emph{all}
geodesics in this direction  are  either closed or (finite number
of them) produce saddle connections.

Thus, it is useful to study  configurations  of  parallel  closed
geodesics on a surface  (which  includes the study of proportions
of corresponding  maximal  cylinders filled with parallel regular
closed  geodesics)  to  get  information  about  the  closure  of
corresponding orbit.

%
%
\footnotetext{
H. Matisse: La Danse.
The State Hermitage Museum, St. Petersburg.
(c) Succession H. Matisse/VG Bild-Kunst, Bonn, 2005
}

One can also use configurations of parallel closed geodesics on a
flat  surface  to determine  those  connected  component  of  the
stratum, to which belongs our surface  $S_0$. Some configurations
(say,  $g-1$  tori connected  in  a ``necklace''  by  a chain  of
cylinders, compare to Fig.~\ref{zorich:fig:dance}) are specific for some
connected components and never appear  for  surfaces  from  other
connected components. We return to  this  discussion  in the very
end of
Sec.~\ref{zorich:ss:Classification:of:Connected:Components:of:the:Strata}
where we use this idea to  distinguish
connected  components\index{Moduli space!connected components of the strata}
of the strata in the
moduli space of quadratic differentials\index{Moduli space!of quadratic differentials}.

\index{Geodesic!closed|)}
\index{Saddle!connection|)}
\index{Moduli space!principal boundary of the moduli space@principal boundary of the moduli space, \emph{see also} Cusp of the moduli space|)}
\index{Cusp of the moduli space@Cusp of the moduli space, \emph{see also} Moduli space; principal boundary|)}

\subsection{Application: Billiards in Rectangular Polygons}
\label{zorich:ss:Application:Billiards:in:Rectangular:Polygons}

\index{Billiard!counting of periodic trajectories|(}
\index{Trajectory!billiard trajectory|(}
\index{Billiard!in rectangular polygon|(}
\index{Billiard!trajectory!periodic|(}

Consider now  a  problem of counting
\emph{generalized diagonals}\index{Diagonal!generalized}
of  bounded  length  or  a  problem  of counting closed  billiard
trajectories  of  bounded  length  in  a  billiard in a  rational
polygon   $\Pi$.   Apply   Katok--Zemliakov   construction\index{Katok--Zemliakov construction}
(see
Sec.~\ref{zorich:ss:Billiards:in:Polygons}) and glue a very flat surface
$S$ from the billiard table  $\Pi$.  Every  generalized  diagonal
(trajectory joining two corners of the  billiard, possibly, after
reflections from the sides) unfolds  to  a  saddle connection and
every periodic trajectory unfolds to a closed regular geodesic.

It is very tempting to  use  the results described above for  the
counting problem  for  the billiard. Unfortunately, the technique
elaborated above is not applicable  to  billiards  directly.  The
problem  is  that  flat  surfaces  coming  from billiards form  a
subspace of large codimension in any stratum of flat surfaces; in
particular,  this  subspace has measure zero. Our ``almost  all''
technique does not see this subspace.

However,  the  problems  are  related  in   some  special  cases;
see~\cite{zorich:Eskin:Masur:Schmoll} treating  billiard in a  rectangle
with a barrier.  As another illustration we consider billiards in
``rectangular   polygons''.    These   results   represent    the
work~\cite{zorich:Athreya:Eskin:Zorich} which is in progress. We warn the reader
that we are  extremely  informal in  the  remaining part of  this
section.

\paragraph{Rectangular Polygons}

\index{Polygon!rectangular|(}

\begin{figure}[htb]
%
%
\includegraphics{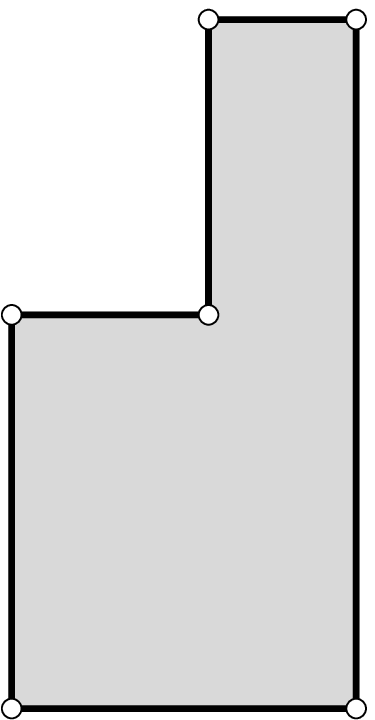}
\includegraphics{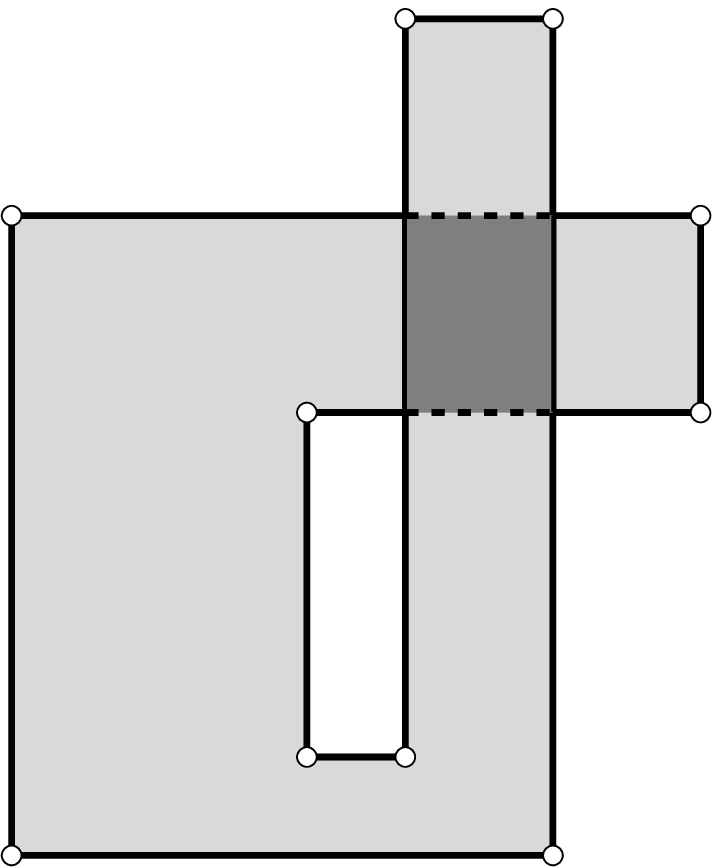}
\includegraphics{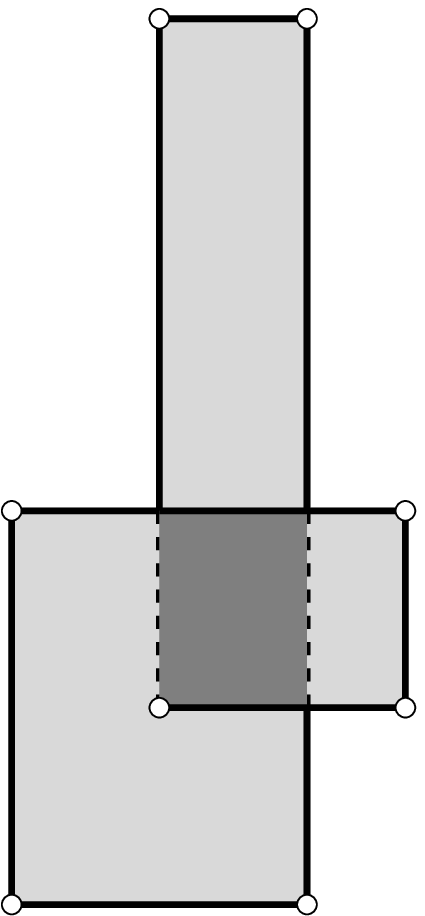}
\includegraphics{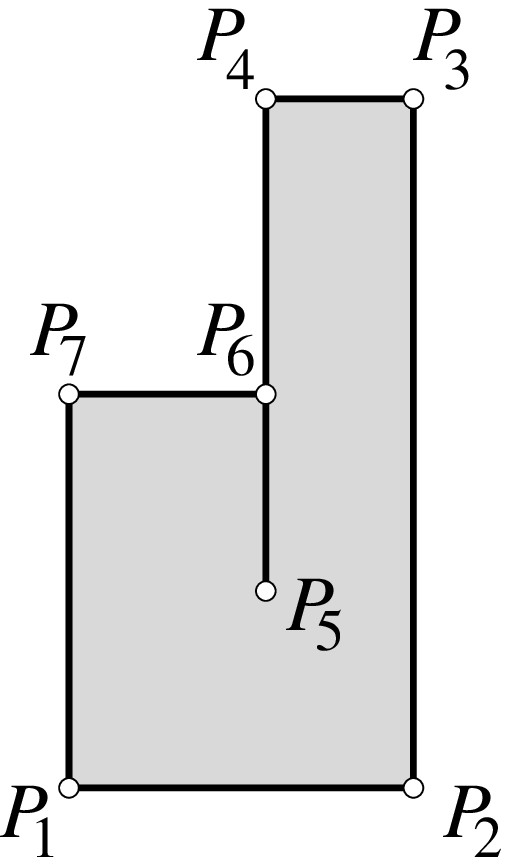}
\vspace{80bp}
\caption{
\label{zorich:fig:rectangular:polygons}
Rectangular  polygons
   }
\end{figure}

Figure~\ref{zorich:fig:rectangular:polygons}  suggests  several examples
of \emph{rectangular  polygons}\index{Polygon!rectangular}.
The ``polygons'' are allowed to
have ramification points  at  the boundary, with restriction that
the  angles  at ramification  points  are  integer  multiples  of
$\pi/2$. Note that we  do not identify the side $P_5 P_6$  with a
part of the  side  $P_4 P_5$ in the  right  polygon. This polygon
should be considered  with  a cut along the  side  $P_5 P_6$. The
corresponding billiard  has  a  ``barrier''  along  the side $P_5
P_6$.

Formally speaking, by a \emph{rectangular polygon}\index{Polygon!rectangular}
$\Pi$ we call
a topological disc  endowed  with a  flat  metric, such that  the
boundary $\partial\Pi$ is presented by  a  finite  broken line of
geodesic  segments  and  such  that  the  angle  between  any two
consecutive sides equals $k \pi/2$, where $k\in \N$.

Consider now  our  problem  for  a  standard rectangular billiard
table (the proportions  of  the sides do not  matter).  We emit a
trajectory from some corner of  the  table and want it arrive  to
another corner after several reflections from the sides.

When our trajectory  reflects  from a  side  it is convenient  to
prolong it as  a  straight  line by making a  reflection  of  the
rectangle  with   respect   to   the   corresponding   side  (see
Katok--Zemliakov                 construction                  in
Sec.~\ref{zorich:ss:Billiards:in:Polygons}).  Unfolding  our rectangular
table we  tile the plane  $\R{2}$ with a rectangular lattice. Our
problem can be reformulated as  a  problem  of counting primitive
lattice      points      (see      the     right     part      of
Fig.~\ref{zorich:fig:V:sc:and:V:cg}).

We are emitting our initial trajectory from some  fixed corner of
the billiard.  It means  that in the model with  a lattice in the
plane we are emitting a straight line from the origin  inside one
of the four quadrants. Thus, we are counting  the asymptotics for
the number  of primitive lattice  points in the intersection of a
coordinate quadrant  with a disc  of large radius $L$ centered at
the origin. This gives us  $1/4$  of the number of all  primitive
lattice points. Note that in our count we have fixed  the initial
corner, but we let our trajectory hit any of the  remaining three
corners. Thus, if we count only those generalized diagonals which
are launched  from some prescribed corner  $P_i$ and arrive  to a
prescribed corner $P_j$ (different from initial one) we get $1/3$
of  the  previous  number.  Hence,  the   number  $N_{ij}(L)$  of
generalized diagonals joining $P_i$  with  $P_j$ is $1/12$ of the
number   $N_{cg}(\T{2},L)\sim   (1/\zeta(2))\cdot  \pi   L^2$  of
primitive lattice points, see~\ref{zorich:eq:asymptotics:for:T2}.

In our  calculation we assumed that  the billiard table  has area
one. It  is clear that asymptotics  for our counting  function is
homogeneous with respect to the area of the  table. Adjusting our
formula for a rectangular billiard  table  of  the area different
from  $1$  we  get  the  following  answer  for   the  number  of
generalized  diagonals  of length at most $L$ joining  prescribed
corner  $P_i$  to a prescribed corner $P_j$  different  from  the
first one:
\begin{equation}
\label{zorich:eq:N:of:generalized:diagonals:for:T2}
N_{ij}(L)\approx \cfrac{1}{2\pi}\cdot
\cfrac{L^2}{\text{Area of the billiard table}}
\end{equation}

Now, having studied  a model case,  we announce two  examples  of
results from~\cite{zorich:Athreya:Eskin:Zorich} concerning rectangular polygons.

\begin{figure}[htb]
\centering
%
   %
   \includegraphics{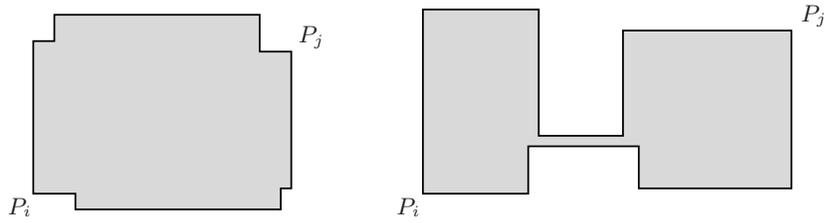}
\begin{picture}(0,0)(0,0)
\put(-155,-76){$P_i$}
\put(-45,-11){$P_j$}
\put(-8,-76){$P_i$}
\put(145,-3){$P_j$}
\end{picture}
\vspace{80bp}
\caption{
\label{zorich:fig:family:of:rectangular:polygons}
Family of rectangular polygons of the same geometry  and the same
area.  The  shape  of  these   polygons   is   quite   different.
Nevertheless for both billiard tables the  number of trajectories
of length at most $L$ joining the right-angle corner $P_i$ to the
right-angle corner $P_j$  is approximately the same as the number
of trajectories of length at most  $L$  joining  two  right-angle
corners of a rectangle of the same area.
   }
\end{figure}

Consider a family of rectangular polygons having exactly $k\ge 0$
angles   $3\pi/2$    and    all   other   angles   $\pi/2$   (see
Fig.~\ref{zorich:fig:family:of:rectangular:polygons}).    Consider     a
generic billiard table in this family (in the measure-theoretical
sense). Fix any  two corners $P_i\neq P_j$ having angles $\pi/2$.
The number $\tilde{N}_{ij}(L)$ of generalized diagonals of length
at most  $L$ joining $P_i$ to $P_j$ is  approximately the same as
for a rectangle:
$$
\tilde{N}_{ij}(L)\sim \cfrac{1}{2\pi}\cdot
\cfrac{L^2}{\text{Area of the billiard table}}
$$
We  have  to  admit  that  we  are slightly  cheating  here:  the
equivalence ``$\sim$'' which we can  currently  prove  is  weaker
than                                                ``$\approx$''
in~\eqref{zorich:eq:N:of:generalized:diagonals:for:T2}; nevertheless, we
do not want to go into technical details.

Note that the  shape  of the polygon within  the  family may vary
quite                      considerably,                      see
Fig.~\ref{zorich:fig:family:of:rectangular:polygons}, and  this does not
affect  the  asymptotic  formula.  However,  the  answer  changes
drastically when we  change  the family. For rectangular polygons
having  several  angles  of  the  form  $n\pi$  the  constant  in
quadratic asymptotics is more complicated. This is why we do
not  expect  any  elementary  proof  of  this formula (our  proof
involves evaluation of corresponding Siegel--Veech constant\index{Siegel--Veech!constant}).

\begin{figure}[htb]
\centering
\includegraphics{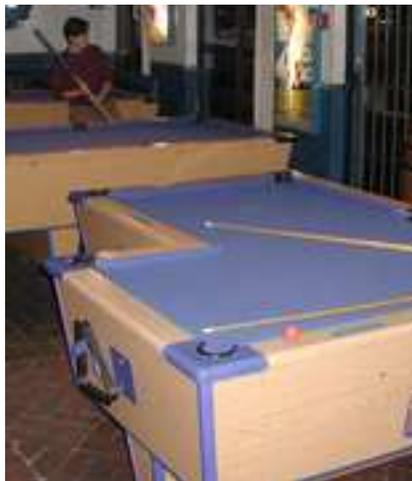}
   %
   %
\begin{picture}(0,0)(0,0)
\put(21,-187){\scriptsize $\bigcirc\hspace*{-6pt}c$ Moon Duchin}
\end{picture}
\vspace{185bp}
\caption{ \label{zorich:fig:L:shaped:billiard}
$L$-shaped billiard table
   }
\end{figure}

Actually, naive intuition does  not  help in counting problems
of this type. Consider, for example, an $L$-shaped billiard
table\index{Billiard!table!L-shaped} as on
Fig.~\ref{zorich:fig:L:shaped:billiard}.

The angle  at the vertex  $P_0$ is $3\pi/2$ which is \emph{three}
times larger  than the angle $\pi/2$  at the other  five vertices
$P_1, \dots, P_5$. However, the number
$$
\tilde{N}_{0j}(L)\sim \cfrac{2}{\pi}\cdot \cfrac{L^2}{\text{Area
of the billiard table}}
$$
of generalized  diagonals of length  at most $L$ joining $P_0$ to
$P_j$, where $1\le j \le 5$, is \emph{four} times bigger than the
number $\tilde{N}_{ij}(L)$ of generalized  diagonals  joining two
corners with the angles $\pi/2$.

\index{Billiard!counting of periodic trajectories|)}
\index{Trajectory!billiard trajectory|)}
\index{Polygon!rectangular|)}

\index{Billiard!in rectangular polygon|)}
\index{Billiard!trajectory!periodic|)}

\section{Volume of Moduli Space}
\label{zorich:s:Volume:of:the:Moduli:Space}

In Sec.~\ref{zorich:ss:Moduli:Space:of:Holomorphic:One:Forms} we defined
a volume element in the stratum  $\cH(d_1,\dots,d_\noz)$\index{0H20@$\cH(d_1,\dots,d_\noz)$ -- stratum in the moduli space}\index{Stratum!in the moduli space}. We used
linear    volume    element    in    cohomological    coordinates
$H^1(S,\{P_1,\dots,P_\noz\};\C{})$ normalized in such way that  a
fundamental domain of the lattice\index{Lattice!in the moduli space}
$$
H^1(S,\{P_1,\dots,P_\noz\};\,\Z{}\oplus\sqrt{-1}\,\Z{})\
\subset\ H^1(S,\{P_1,\dots,P_\noz\};\C{})
$$
has unit  volume. The unit lattice does not  depend on the choice
of  cohomological  coordinates,  its  vertices play the  role  of
\emph{integer      points}     in      the      moduli      space
$\cH(d_1,\dots,d_\noz)$\index{0H20@$\cH(d_1,\dots,d_\noz)$ -- stratum in the moduli space}\index{Stratum!in the moduli space}.
In  Sec.~\ref{zorich:ss:Square:tiled:surfaces}
we  suggest   a   geometric   interpretation   of  flat  surfaces
representing integer points of the strata.

Using this interpretation we give an idea for counting the
volume\index{Moduli space!volume of the moduli space}
(``hyperarea'')           of           the           hypersurface
$\cH_1(d_1,\dots,d_\noz)\subset\cH(d_1,\dots,d_\noz)$   of   flat
surfaces  of  area  one.  We  apply  the strategy  which  can  be
illustrated in a  model  example of evaluation of  the  area of a
unit sphere. We first count the asymptotics for the number $N(R)$
of  integer points  inside  a ball of  huge  radius $R$.  Clearly
$N(R)$ corresponds to the  volume of the ball, so if we  know the
asymptotics  for  $N(R)$  we  know  the  formula  for  the volume
$\Vol(R)$   of   the   ball   of  radius  $R$.   The   derivative
$\cfrac{d}{dR}\Vol(R)\big|_{R=1}$ gives  us  the area of the unit
sphere.

Similarly,  to  evaluate the ``hyperarea of a unit  hyperboloid''
$\cH_1(d_1,\dots,d_\noz)$\index{0H30@$\cH_1(d_1,\dots,d_\noz)$ -- ``unit hyperboloid''}\index{Stratum!in the moduli space}\index{Unit hyperboloid}   it   is   sufficient  to  count   the
asymptotics  for  the  number  of   integer   points   inside   a
``hyperboloid'' $\cH_R(d_1,\dots,d_\noz)$ of huge ``radius'' $R$.
The  role  of  the  ``radius''  $R$  is played  by  the  positive
homogeneous    real     function    $R=area(S)$    defined     on
$\cH(d_1,\dots,d_\noz)$\index{0H20@$\cH(d_1,\dots,d_\noz)$ -- stratum in the moduli space}\index{Stratum!in the moduli space}.

\index{0d10@$d\nu$ -- volume element in the moduli space}

Note that the volume $\nu\big(\cH_{\le R}(d_1,\dots,d_\noz)\big)$
of     a     domain    bounded     by     the     ``hyperboloid''
$\cH_R(d_1,\dots,d_\noz)$ is a homogeneous function of $R$ of the
weight $\dim_{\R{}}\cH(d_1,\dots,d_\noz)/2$ while the volume of a
ball of radius $R$ is a homogeneous function of $R$ of the weight
which  equals  the  dimension  of the space. This  difference  in
weights explains the factor $2$ in the formula below:
\begin{multline}
\label{zorich:eq:efinition:of:volume}
\Vol\big(\cH_1(d_1,\dots,d_\noz)\big)=
\left.2\cfrac{d}{dR}
\ \nu\big(\cH_{\le R}(d_1,\dots,d_\noz)\big)
\right|_{R=1}=\\ 
=\dim_{\R{}}(\cH_1(d_1,\dots,d_\noz))\cdot\nu\big(\cH_{\le 1}(d_1,\dots,d_\noz)\big)
\end{multline}

This  approach  to computation of  the volumes\index{Moduli space!volume of the moduli space}
was suggested  by
A.~Eskin  and  A.~Okounkov  and  by M.~Kontsevich and the author.
However,   the  straightforward  application  of  this  approach,
described   in   Sec.~\ref{zorich:ss:Square:tiled:surfaces},  gives  the
values  of  the  volumes only for several low-dimensional strata.
The  general  solution  of the  problem was found by A.~Eskin and
A.~Okounkov   who    used  in  addition   powerful   methods   of
representation    theory.    We  give an idea  of their method in
Sec.~\ref{zorich:ss:Approach:of:Eskin:and:Okounkov}.

\subsection{Square-tiled Surfaces}
\label{zorich:ss:Square:tiled:surfaces}

Let us  study the geometric properties  of the flat  surfaces $S$
represented       by
``integer        points''\index{Lattice!in the moduli space}
$S\in H^1(S,\{P_1,\dots,P_\noz\};\Z{}\oplus\sqrt{-1}\,\Z{})$         in
cohomological  coordinates.  Let  $\omega$  be  the   holomorphic
one-form representing  such  flat surface $S$. Since $[\omega]\in
H^1(S,\{P_1,\dots,P_\noz\};\Z{}\oplus\sqrt{-1}\,\Z{})$        all
periods  of  $\omega$  (including  relative  periods)  belong  to
$\Z{}\oplus\sqrt{-1}\,\Z{}$. Hence  the following map  $f_\omega$
from   the   flat    surface    $S$   to   the   standard   torus
$\torus=\C{}/(\Z{}\oplus\sqrt{-1}\,\Z{})$ is well-defined:
$$
f_{\omega} : P\mapsto \Big(\int_{P_1}^P \omega\Big)\ \text{mod }
\Z{}\oplus \sqrt{-1}\,\Z{}\,,
$$
where $P_1$  is one of the conical  points. It  is easy to  check
that  $f_{\omega}$  is a  ramified  covering,  moreover,  it  has
exactly $\noz$  ramification  points, and the ramification points
are exactly the  zeros $P_1, \dots, P_\noz$ of $\omega$. Consider
the flat torus  $\torus$  as a  unit  square with the  identified
opposite sides.  The  covering  $f_\omega:S\to  \torus$ induces a
tiling of the  flat  surface $S$ by unit  squares.  Note that all
unit squares are endowed with the following additional structure:
we  know exactly  which  edge is top,  bottom,  right, and  left;
adjacency of the squares respects this  structure  in  a  natural
way: we  glue vertices to  vertices and edges to edges, moreover,
the right edge of a square  is always identified to the left edge
of some square and top  edge  is always identified to the  bottom
edge of some  square.  We  shall call a flat  surface  with  such
tiling a
\emph{square-tiled surface}\index{Square-tiled surface}\index{Surface!square-tiled}.
We   see   that   the   problem   of   counting  the  volume   of
$\cH_1(d_1,\dots,d_\noz)$\index{0H30@$\cH_1(d_1,\dots,d_\noz)$ -- ``unit hyperboloid''}\index{Stratum!in the moduli space}\index{Unit hyperboloid} is equivalent to the following problem:
how  many  square-tiled  surfaces  of  a   given  geometric  type
(determined by number  and types of conical singularities) can we
construct   using    at    most    $N$    unit    squares.   Say,
Fig.~\ref{zorich:fig:L:3:surfaces} gives  the  list  of all square-tiled
surfaces of genus $g>1$ glued from at most $3$ squares; they all
belong to $\cH(2)$.

In  terms  of the  coverings  our Problem  can  be formulated  as
follows. Consider  the ramified coverings $p:S\to\torus$ over the
standard  torus  $\torus$.  Fix  the number $\noz$  of  branching
points,  and  ramification  degrees $d_1,\dots,d_\noz)$ at  these
points. Assume that all ramification points  $P_1, \dots, P_\noz$
project  to  the  same  point  of  the torus $\torus$.  Enumerate
ramified coverings of  any given ramification type having at most
$N\gg 1$  sheets.  Here  pairs  of  coverings forming commutative
diagrams as below are identified:
\begin{equation}
\label{zorich:eq:commutative:diag:for:coverings}
\begin{matrix}
S & \longleftrightarrow & S\\
&\searrow\quad  \swarrow\\
& \torus &
\end{matrix}
\end{equation}

\paragraph{Counting of Square-tiled Tori}

\begin{figure}[htb]
\centering
\includegraphics{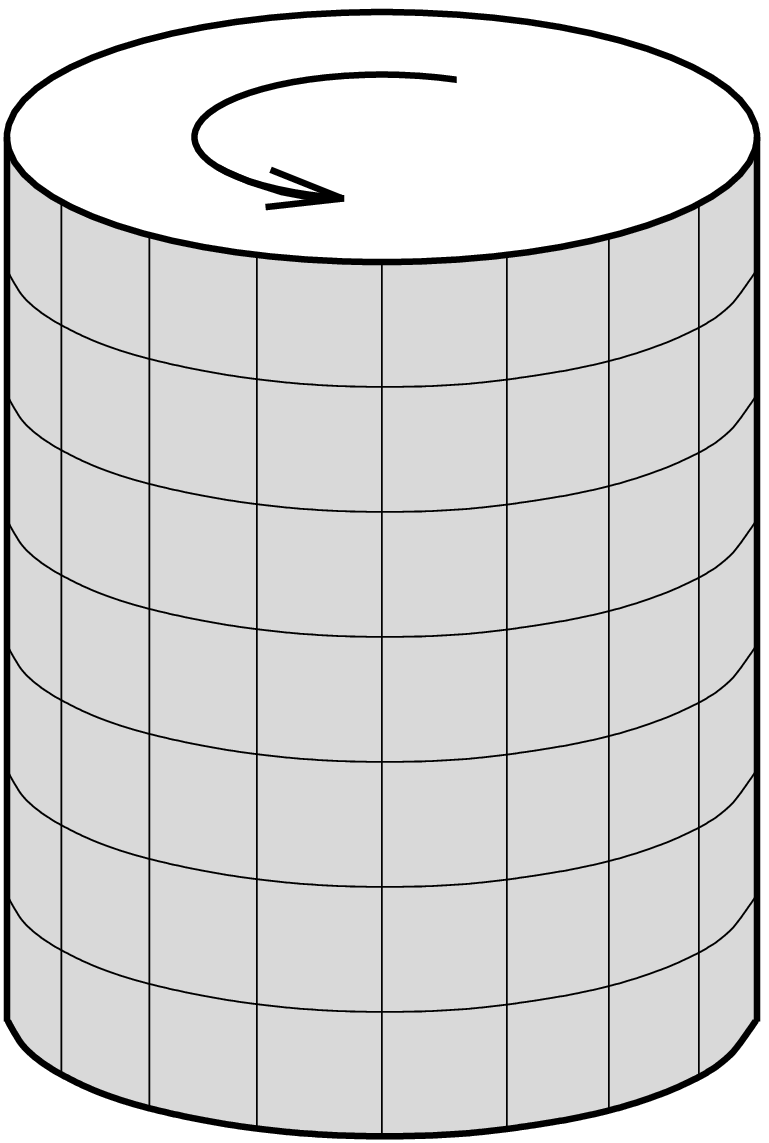}
\includegraphics{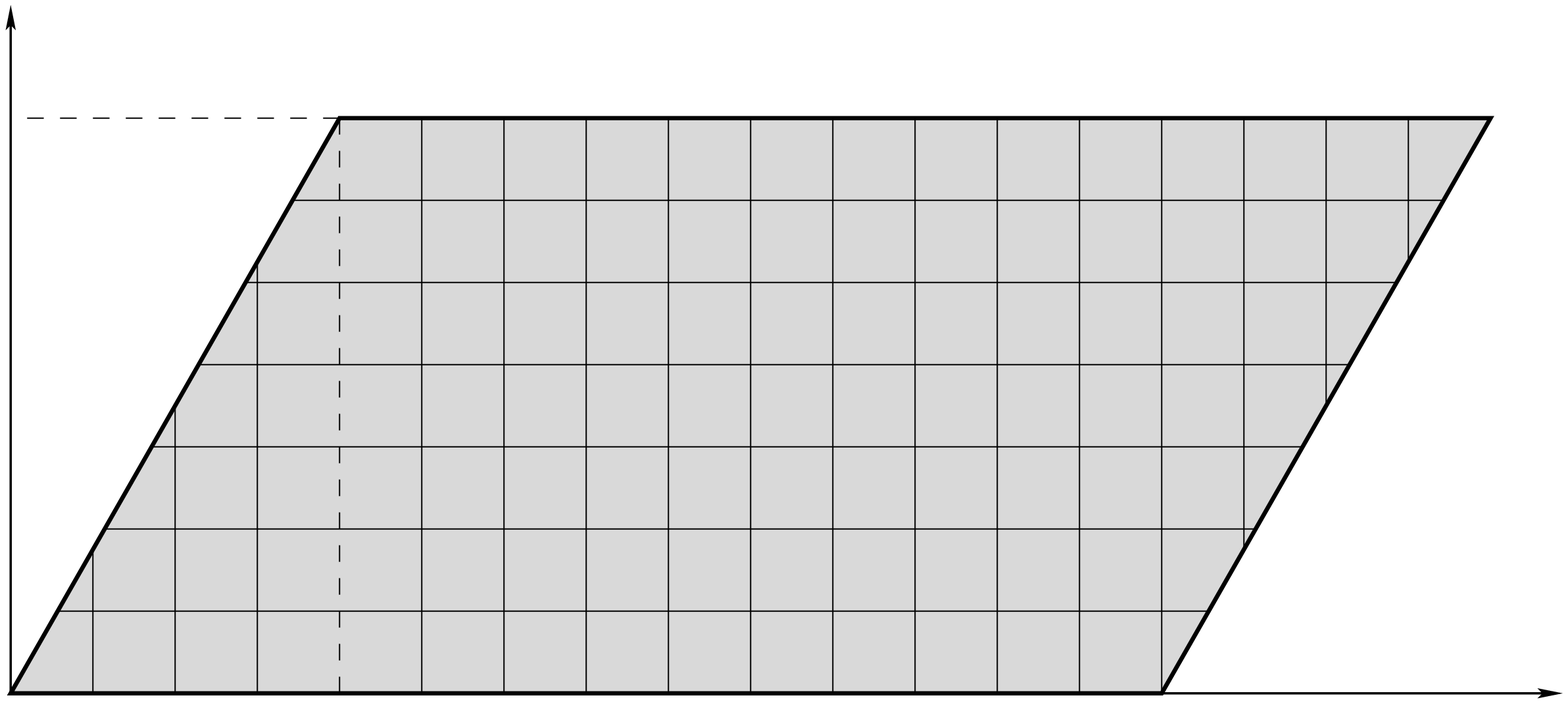}
\vspace{110bp}
\begin{picture}(120,0)(120,0) 
\put(180,0){\begin{picture}(0,0)(0,0)
\put(-111,87){$t$}
\put(-115,11){$w$}
\put(-80,55){$h$}
       %
       %
\put(-35,82){$h$}
\put(97,11){$w$}
\put(7,11){$t$}
\end{picture}}
\end{picture}
\caption{
\label{zorich:fig:1cyl}
A square-tiled surface is decomposed into several cylinders. Each
cylinder is parametrized  by its width (perimeter) $w$ and height
$h$. Gluing the cylinders together we get also  a twist parameter
$t$, where $0\le t<w$, for each cylinder}
\end{figure}

Let us count  the  number of square-tiled tori  tiled  by at most
$N\gg 1$  squares. In this case  our square-tiled surface  has no
singularities at all: we have  a  flat torus tiled with the  unit
squares in a regular way. Cutting our flat torus by  a horizontal
waist  curve  we  get a cylinder  with  a  waist  curve of length
$w\in\N$ and  a  height  $h\in\N$,  see  Fig.~\ref{zorich:fig:1cyl}. The
number of squares in the tiling equals $w\cdot h$. We  can reglue
the torus from  the  cylinder with  some  integer twist $t$,  see
Fig.~\ref{zorich:fig:1cyl}. Making  an  appropriate Dehn twist along the
waist  curve  we  can reduce the value of the twist $t$ to one of
the values $0, 1, \dots, w-1$. In other words, fixing the integer
perimeter $w$  and height $h$ of a cylinder  we get $w$ different
square-tiled tori.

Thus the number of square tiled tori constructed by using at most
$N$ squares is represented as
$$
\nu\big(\cH_{\le N}(0)\big)\sim
\sum_{\substack{w,h\in\N\\w\cdot h\le N}}
w = \sum_{\substack{w,h\in\N\\w\le \frac{N}{h}}} w \approx
\sum_{h\in\N} \cfrac{1}{2}\cdot\left(\cfrac{N}{h}\right)^2
= \cfrac{N^2}{2}
\cdot \zeta(2)  = \cfrac{N^2}{2} \cdot \cfrac{\pi^2}{6}
$$
Actually,  some  of the  tori  presented  by  the  first  sum are
equivalent by an affine  diffeomorphism,  so we are counting them
twice, or even several times. Say,  the  tori  $w=2;h=1;t=0$  and
$w=1;h=2;t=0$ are  equivalent.  However,  this happens relatively
rarely, and this correction does not affect the  leading term, so
we simply neglect it.

Applying   the   derivative   $\left.2\cfrac{d}{dN}\right|_{N=1}$
(see~\eqref{zorich:eq:efinition:of:volume}) we finally get the following
value for the volume\index{Moduli space!volume of the moduli space}
of the space of flat tori
$$
\Vol(\mathcal{H}_1(0)) = \cfrac{\pi^2}{3}
$$

\paragraph{Decomposition of a Square-Tiled Surface into Cylinders}

Let us study the geometry of square-tiles surfaces. Note that all
leaves  of  both  horizontal  and  vertical  foliation  on  every
square-tiled surface are closed.  In  particular the union of all
horizontal critical leaves (the ones adjacent  to conical points)
forms a finite graph $\Gamma$. The collection $P_1,\dots, P_\noz$
of conical points forms the  set  of the vertices of this  graph;
the  edges  of  the   graph   are  formed  by  horizontal  saddle
connections.  The  complement  $S-\Gamma$  is  a  union  of  flat
cylinders.

For    example,    for    the    square-tiled    surfaces    from
Fig.~\ref{zorich:fig:L:3:surfaces} we  get the following  decompositions
into horizontal cylinders. We have  one  surface  composed from a
single cylinder  filled with closed horizontal trajectories; this
cylinder has width (perimeter)  $w=3$  and hight $h=1$. Two other
surfaces  are  composed from two cylinders. The  heights  of  the
cylinders are $h_1=h_2=1$, the widths  are  $w_1=1$  and  $w_2=2$
correspondingly.   Observing   the   two-cylinder   surfaces   at
Fig.~\ref{zorich:fig:L:3:surfaces}  we  see  that  they  differ  by  the
\emph{twist}  parameter  $t_2$  (see Fig.~\ref{zorich:fig:1cyl}) of  the
wider  cylinder:  in one  case  $t_2=0$  and  in  the  other case
$t_2=1$. By construction the width $w_i$ and height  $h_i$ of any
cylinder are strictly positive  integers;  the value of the twist
$t_i$  is  a nonnegative  integer  bounded by  the  width of  the
cylinder: $0\le t_i < w_i$.

\paragraph{Separatrix Diagrams}

Let us  study in more details  the graphs $\Gamma$  of horizontal
saddle connections.

We start with an informal explanation. Consider the  union of all
saddle connections  for  the  horizontal  foliation,  and add all
critical points (zeroes of $\omega$). We obtain a finite oriented
graph $\Gamma$. Orientation on the edges comes from the canonical
orientation of the horizontal foliation. Moreover, graph $\Gamma$
is drawn  on an oriented  surface, therefore it carries so called
{\it ribbon structure} (even  if  we forget about the orientation
of edges), i.e. on the star  of each vertex $P$ a cyclic order is
given,  namely  the  counterclockwise  order in which  edges  are
attached  to  $P$.   The  direction  of  edges  attached  to  $P$
alternates (between directions toward $P$ and  from  $P$)  as  we
follow the counterclockwise order.

It is  well known that any  finite ribbon graph  $\Gamma$ defines
canonically  (up  to an isotopy) an oriented surface  $S(\Gamma)$
with boundary. To  obtain this surface  we replace each  edge  of
$\Gamma$  by  a thin oriented strip (rectangle)  and  glue  these
strips  together  using  the  cyclic  order  in  each  vertex  of
$\Gamma$. In  our case surface $S(\Gamma)$  can be realized  as a
tubular $\varepsilon$-neighborhood (in the  sense  of transversal
measure) of the  union of all saddle connections for sufficiently
small $\varepsilon>0$.

The  orientation  of  edges   of   $\Gamma$  gives  rise  to  the
orientation of  the  boundary  of  $S(\Gamma)$.  Notice that this
orientation is {\it not} the same as the canonical orientation of
the boundary of  an  oriented surface. Thus, connected components
of the boundary of  $S(\Gamma)$  are decomposed into two classes:
positively   and   negatively   oriented  (positively  when   two
orientations of  the boundary components coincide and negatively,
when  they  are   different).   The  complement  to  the  tubular
$\varepsilon$-neighborhood of $\Gamma$ is a finite disjoint union
of  open  cylinders  foliated  by  oriented  circles. It gives  a
decomposition    of     the     set     of    boundary    circles
$\pi_0(\partial(S(\Gamma)))$  into  pairs  of  components  having
opposite signs of the orientation.

Now we are ready to give a formal definition:

A
{\it separatrix  diagram}\index{Separatrix!diagram}\index{Diagram!separatrix diagram}
is a finite oriented  ribbon  graph $\Gamma$, and a decomposition
of the set of boundary components of $S(\Gamma)$ into pairs, such
that
\begin{itemize}
\item[--]
the orientation of edges at any vertex is alternated with respect
to the cyclic order of edges at this vertex;
\item[--] there is one positively  oriented  and  one  negatively
oriented boundary component in each pair.
\end{itemize}

Notice that ribbon graphs which appear as a part of the structure
of a separatrix diagram are  very  special. Any vertex of such  a
graph has  even degree, and  the number of boundary components of
the associated surface with  boundary  is even. Notice also, that
in general the graph of a separatrix diagram is {\it not} planar.

Any separatrix  diagram  $(\Gamma,  pairing)$  defines  a  closed
oriented surface together with an embedding of $\Gamma$  (up to a
homeomorphism) into this  surface. Namely, we glue to the surface
with boundary  $S(\Gamma)$  standard oriented cylinders using the
given pairing.

In pictures representing diagrams we  encode  the  pairing on the
set of boundary components painting corresponding  domains in the
picture by some  colors (textures in the black-and-white text) in
such a way that every  color  appears exactly twice. We will  say
also that paired components have the {\it same color}.

\begin{figure}[htb]
%
\includegraphics{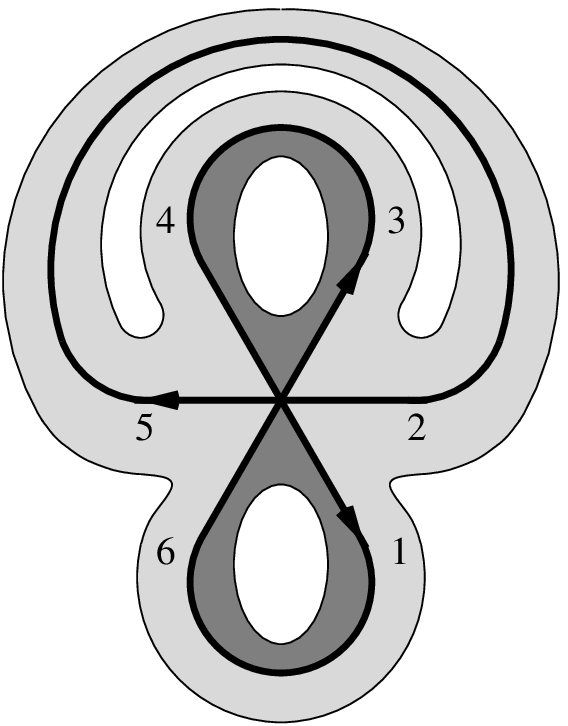}
\includegraphics{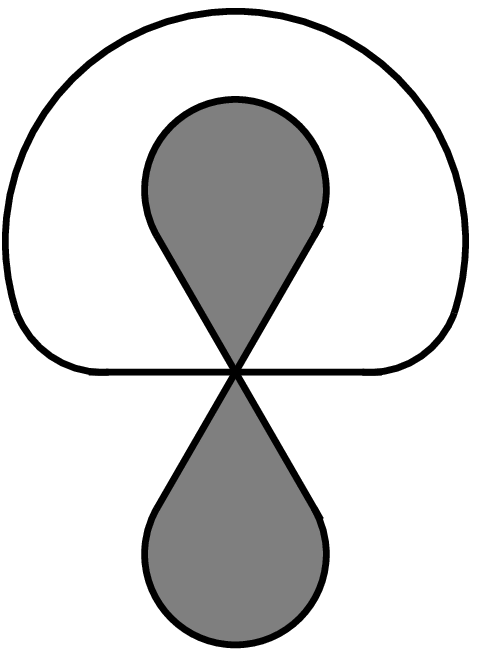}
\vspace{160bp}
\caption{
\label{zorich:fig:g2}
An example of  a separatrix diagram.  A detailed picture  on  the
left can be encoded by a schematic picture on the right. }
\end{figure}

\begin{Example}
\label{zorich:ex:g2}
The ribbon  graph presented at Figure~\ref{zorich:fig:g2} corresponds
to the horizontal foliation  of an Abelian differential on a
surface of genus $g=2$.  The Abelian differential  has a single
zero  of degree $2$. The ribbon graph has two pairs of boundary
components.
\end{Example}

Any  separatrix   diagram   represents   an  orientable  measured
foliation with only closed  leaves  on a compact oriented surface
without boundary.  We say that  a diagram is {\it realizable} if,
moreover, this measured foliation can be chosen as the horizontal
foliation  of  some  Abelian  differential. Lemma below  gives  a
criterion of realizability of a diagram.

Assign to each saddle connection a real variable standing for its
``length''. Now any boundary component is  also  endowed  with  a
``length'' obtained as sum of the ``lengths'' of all those saddle
connections which  belong to this component.  If we want  to glue
flat cylinders to the boundary components,  the  lengths  of  the
components in every pair should match each other.  Thus for every
two boundary components paired together  (i.e.  having  the  same
color) we  get a linear  equation: ``the length of the positively
oriented component equals  the  length of the negatively oriented
one''.

\begin{NNLemma}
        %
A diagram  is realizable if  and only if the corresponding system
of linear equations  on  ``lengths'' of saddle connections admits
strictly positive solution.
\end{NNLemma}

The proof is obvious.

\begin{Example}
The diagram  presented  at  Fig.~\ref{zorich:fig:g2}  has  three  saddle
connections, all of  them are loops. Let $p_{16}, p_{52}, p_{34}$
be their ``lengths''. There are two pairs of boundary components.
The corresponding system of linear equations is as follows:
$$
\left\{\begin{array}{l}
p_{34}=p_{16}\\
p_{16}+p_{52}=p_{34}+p_{52}
\end{array}\right.
$$
\end{Example}

\begin{figure}[htb]
\centering
\includegraphics{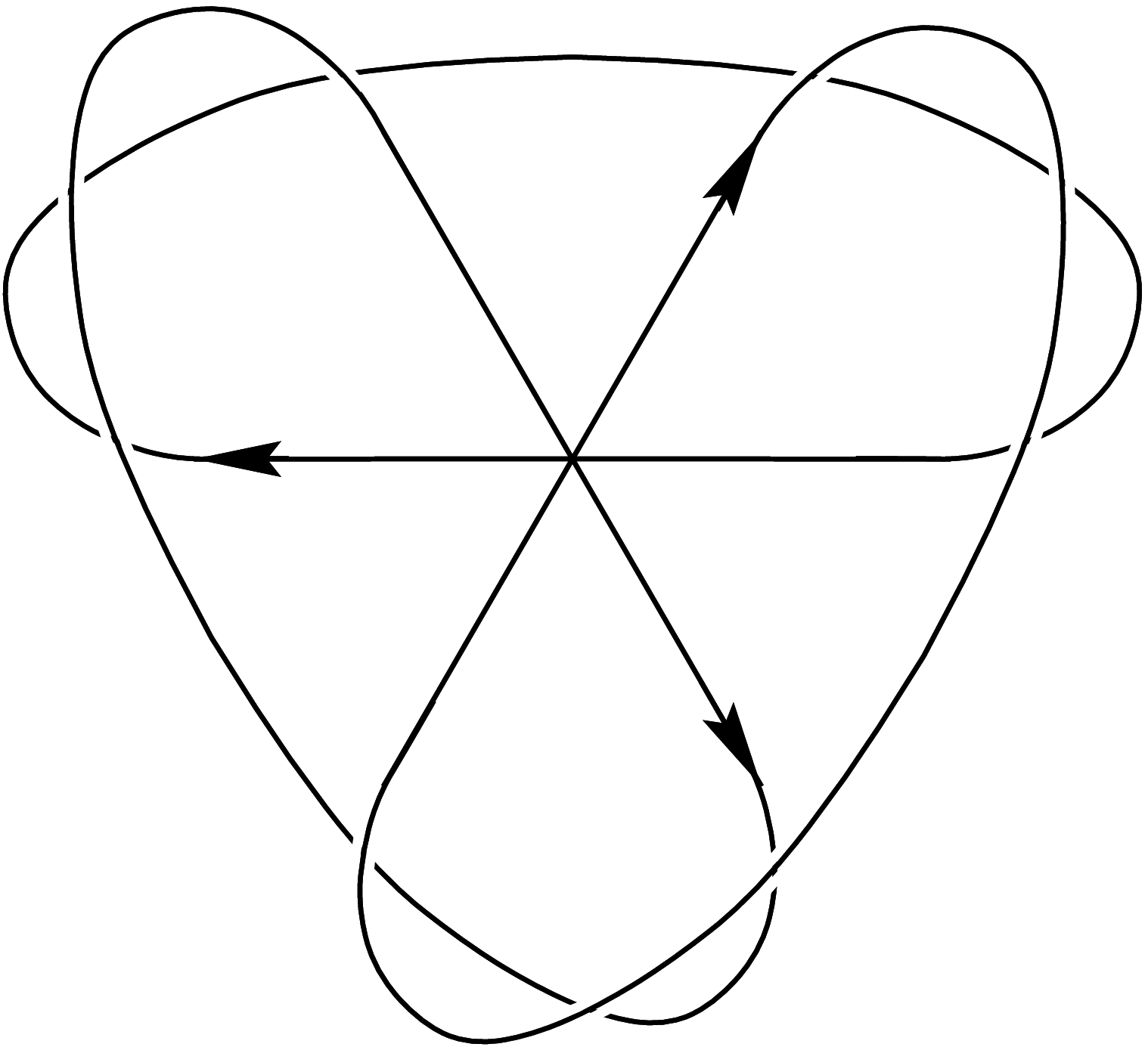}
\includegraphics{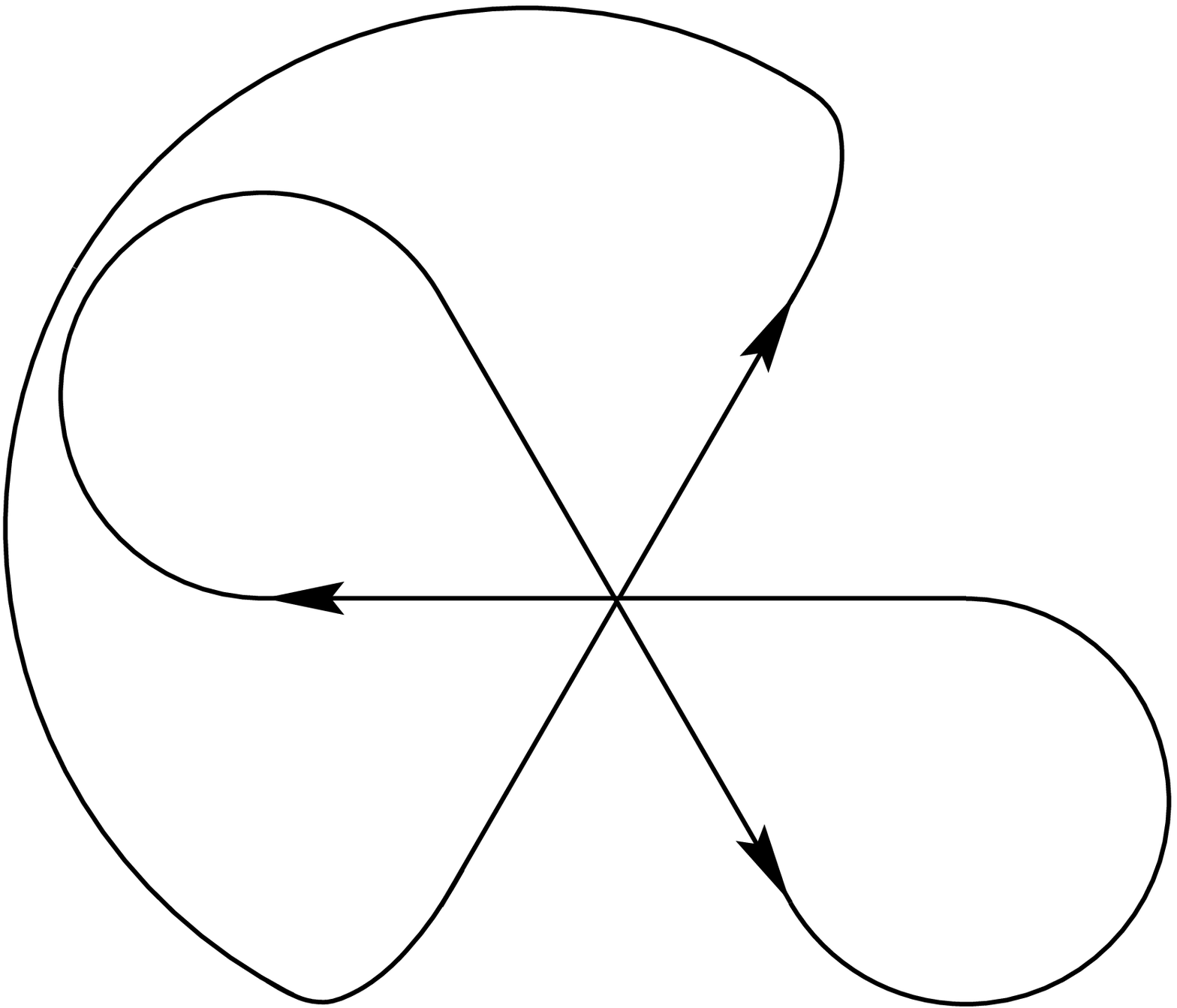}
\includegraphics{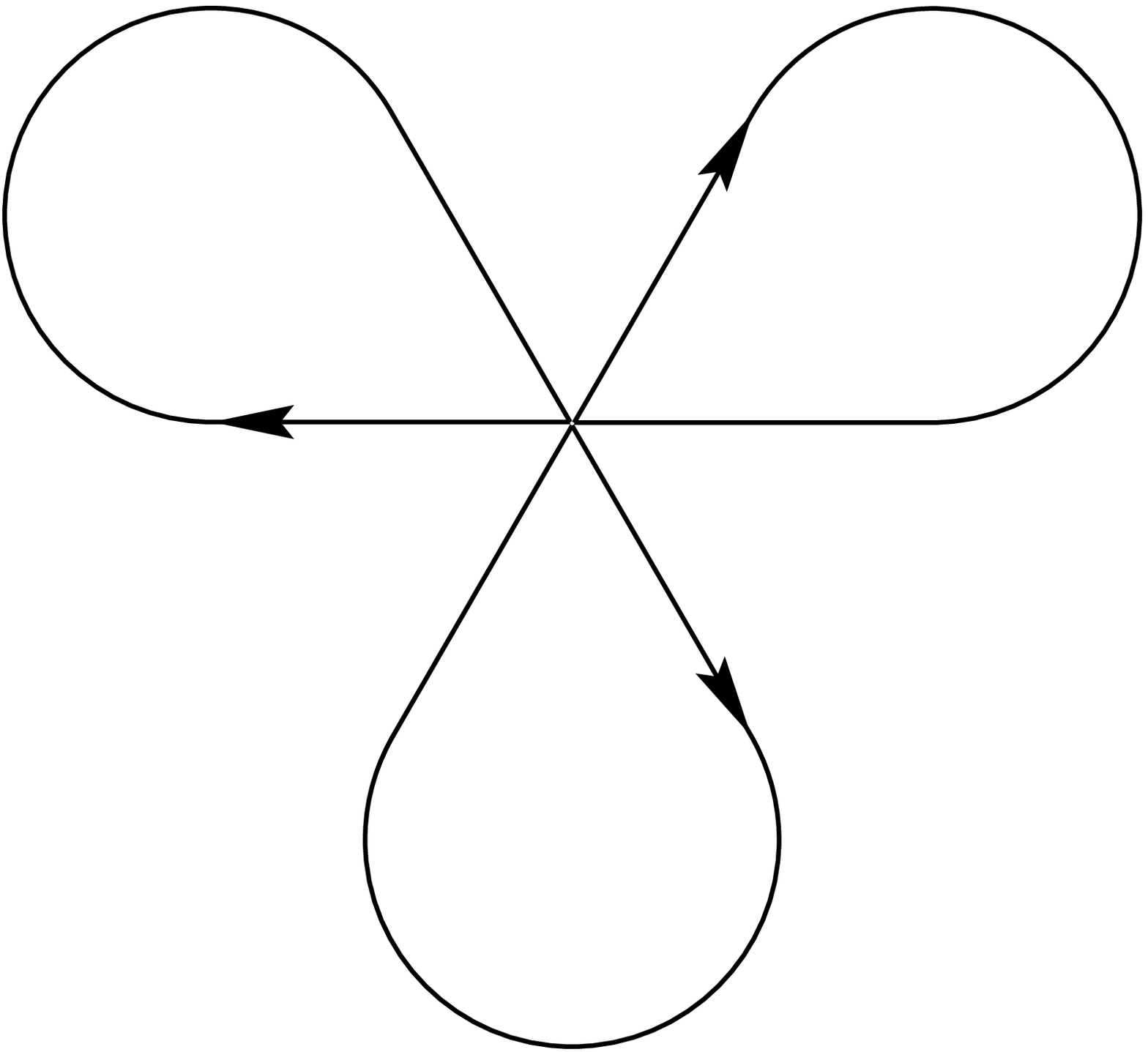}
\vspace{115bp} 
\begin{picture}(112,0)(112,0)
\put(180,0){\begin{picture}(0,0)(0,0)
\put(-129,20){$p_1$}
\put(-95,75){$p_2$}
\put(-162,75){$p_3$}
\put(-12,100){$p_1$}
\put(35,85){$p_2$}
\put(-12,15){$p_1$}
\end{picture}}
\end{picture}
\caption{
\label{zorich:fig:separatrix:diagrams}
The  separatrix  diagrams  from  left  to  right  a)  represent a
square-tiled  surface   glued   from   one   cylinder   of  width
$w=p_1+p_2+p_3$; b)  represent  a square-tiled surface glued from
two cylinder of widths $w_1=p_1$ and  $w_2=p_1+p_2$;  c)  do  not
represent any square-tiled surface} \end{figure}

\begin{Exercise}
Check       that       two      separatrix       diagrams      at
Fig.~\ref{zorich:fig:separatrix:diagrams}  are realizable,  and  one  --
not. Check that there are no other realizable separatrix diagrams
for the surfaces from  the  stratum $\cH(2)$. Find all realizable
separatrix diagrams for the stratum $\cH(1,1)$.
\end{Exercise}

\paragraph{Counting of Square-tiled Surfaces in $\cH(2)$}

To consider one more example  we  count  square-tiled surfaces in
the stratum $\cH(2)$. We have seen that in this stratum there are
only  two  realizable separatrix diagrams; they are presented  on
the left and in the center of Fig.~\ref{zorich:fig:separatrix:diagrams}.

Consider those  square tiled surfaces from $\mathcal{H}(2)$ which
correspond       to       the       left       diagram       from
Fig.~\ref{zorich:fig:separatrix:diagrams}. In this case the ribbon graph
corresponding  to  the separatrix diagram has single ``top''  and
single ``bottom'' boundary component, so  our  surface  is  glued
from a single cylinder. The waist curve of the cylinder is of the
length $w=p_1+p_2+p_3$, where  $p_1, p_2, p_3$ are the lengths of
the separatrix loops. As usual, denote the height of the cylinder
by $h$. The twist $t$ of the cylinder has an integer value in the
interval $[0,w-1]$.  Thus the number  of surfaces of this type of
the area bounded by $N$ is asymptotically equivalent to the sum
$$
\cfrac{1}{3}\sum_{\substack{p_1,p_2,p_3,h\in\N\\(p_1+p_2+p_3)h\le N}} (p_1+p_2+p_3)
\sim \cfrac{N^4}{24}\cdot \zeta(4) = \cfrac{N^4}{24}\cdot \cfrac{\pi^4}{90}
$$
(see  more  detailed  computation in~\cite{zorich:Zorich:square:tiled}).
The coefficient $1/3$ compensates the arbitrariness of the choice
of enumeration  of  $p_1,p_2,p_3$ preserving the cyclic ordering.
Similar to  the torus case we  have neglected a  small correction
coming from counting equivalent surfaces several times.

\begin{Exercise}
Check that for $p_1=p_2=p_3=1$  all  possible values of the twist
$t=0,1,2$    give    equivalent    flat   surfaces;   see    also
Fig.~\ref{zorich:fig:L:3:surfaces}
\end{Exercise}

Consider now a ribbon  graph  corresponding to the middle diagram
from Fig.~\ref{zorich:fig:separatrix:diagrams}. It has  two  ``top'' and
two ``bottom''  boundary  components. Thus, topologically, we can
glue in a pair of  cylinders  in two different ways. However,  to
have a flat structure on  the  resulting surface we need to  have
equal lengths  of  ``top''  and  ``bottom''  boundary components.
These lengths are  determined by the lengths of the corresponding
separatrix  loops. It is  easy  to  check  that one  of  the  two
possible gluings of cylinders  is  forbidden: it implies that one
of the separatrix loops has zero length, and hence the surface is
degenerate.

The other gluing is  realizable. In this case there is a  pair of
separatrix    loops     of     equal     lengths    $p_1$,    see
Fig.~\ref{zorich:fig:separatrix:diagrams}. The square-tiled  surface  is
glued from two cylinders: one having a waist curve $w_1=p_1$, and
the  other  one having  waist  curve  $w_2=p_1+p_2$.  Denote  the
heights and twists of  the  corresponding cylinders by $h_1, h_2$
and $t_1,t_2$.  The twist of the first cylinder  has the value in
the interval $[0,w_1-1]$; the twist  of  the  second cylinder has
the  value  in  the  interval  $[0,w_2-1]$.  Thus  the  number of
surfaces  of  2-cylinder type  of  the  area  bounded  by  $N$ is
asymptotically equivalent to the sum
$$
\sum_{\substack{p_1,p_2,h_1,h_2\\p_1 h_1+(p_1+p_2)h_2\le N}} p_1(p_1+p_2) =
\sum_{\substack{p_1,p_2,h_1,h_2\\p_1(h_1+h_2)+p_2 h_2\le N}} p_1^2+p_1 p_2
$$
It is not difficult  to represent these two sums in terms  of the
\emph{multiple  zeta   values}  $\zeta(1,3)$  and    $\zeta(2,2)$
(see    detailed    computation   in~\cite{zorich:Zorich:square:tiled}).
Applying     the     relations    $\zeta(1,3)=\zeta(4)/4$     and
$\zeta(2,2)=3\zeta(4)/4$ we get the following asymptotic  formula
for our sum:
$$
\sum_{\substack{p_1,p_2,h_1,h_2\\p_1(h_1+h_2)+p_2 h_2\le N}} p_1^2+p_1 p_2
\sim\cfrac{N^4}{24} \big( 2\cdot\zeta(1,3)+\zeta(2,2)\big)
=\cfrac{N^4}{24} \cdot \cfrac{5}{4} \cdot \zeta(4)
=\cfrac{N^4}{24} \cdot \cfrac{5}{4} \cdot \cfrac{\pi^4}{90}
$$
Joining  the   impacts   of   the   two   diagrams  and  applying
$\left.2\cdot                         \cfrac{d}{dN}\right|_{N=1}$
(see~\eqref{zorich:eq:efinition:of:volume}) we  get the following  value
for the volume\index{Moduli space!volume of the moduli space}
of the stratum $\cH(2)$:
$$
\Vol(\mathcal{H}_1(2))= \cfrac{\pi^4}{120}
$$

\toread{Separatrix diagrams}

Technique  of   separatrix   diagrams   is  extensively  used  by
K.~Strebel    in~\cite{zorich:Strebel}    and     in    some    articles
like~\cite{zorich:Kontsevich:Zorich}.

\subsection{Approach of A.~Eskin and A.~Okounkov}
\label{zorich:ss:Approach:of:Eskin:and:Okounkov}

It is time  to  confess  that evaluation of the  volumes  of  the
strata by means of naive counting square-tiled surfaces suggested
in the previous  section is not  efficient in general  case.  The
number  of  realizable  separatrix  diagrams  grows   and  it  is
difficult to express the  asymptotics  of the sums for individual
separatrix  diagrams  in  reasonable  terms  (say,  in  terms  of
multiple zeta values). In general case  the  problem  was  solved
using  the   following   approach   suggested   by  A.~Eskin  and
A.~Okounkov in~\cite{zorich:Eskin:Okounkov}.

Consider a  square-tiled  surface  $S\in\cH(d_1, \dots, d_\noz)$.
Enumerate the squares in some way. For the square number  $j$ let
$\pi_{r}(j)$ be  the number of its neighbor to  the right and let
$\pi_{u}(j)$ be  the number of  the square atop the square number
$j$.             Consider             the              commutator
$\pi'=\pi_{r}\pi_{u}\pi^{-1}_{r}\pi^{-1}_{u}$  of  the  resulting
permutations. When the total number of squares is big enough, for
most of the squares

Geometrically the resulting permutation $\pi'$ corresponds to the
following path:  we start from a square number  $j$, then we move
one step right, one step up, one step left, one step down, and we
arrive to  $\pi'(j)$. When the  total number of squares is large,
then for majority of the squares  such path brings us back to the
initial square; for these values of index $j$ we get $\pi'(j)=j$.
However, we may have more than $4$ squares adjacent to  a conical
singularity $P_k\in S$. For the squares adjacent to a singularity
the path right-up-left-down does not bring us back to the initial
square.   It   is    easy    to   check   that   the   commutator
$\pi'=\pi_{r}\pi_{u}\pi^{-1}_{r}\pi^{-1}_{u}$ is represented as a
product of $\noz$ cycles of lengths  $(d_1+1), \dots, (d_\noz+1)$
correspondingly.

We conclude  that  a  square-tiled  surface  $S\in\cH(d_1, \dots,
d_\noz)$   can   be   defined   by   a   pair   of   permutations
$\pi_{r},\pi_{u}$,      such      that       the       commutator
$\pi_{r}\pi_{u}\pi^{-1}_{r}\pi^{-1}_{u}$  decomposes  into  given
number  $\noz$  of  cycles  of  given  lengths  $(d_1+1),  \dots,
(d_\noz+1)$. Choosing another  enumeration  of the squares of the
same  square-tiled  surface  $S$  we  get  two  new  permutations
$\tilde\pi_{r},\tilde\pi_{u}$. Clearly  the permutations in  this
new pair are conjugate  to the initial ones by means of  the same
permutation $\rho$ responsible  for  the change of enumeration of
the         squares:         $\tilde\pi_{r}=\rho\pi_{r}\rho^{-1},
\tilde\pi_{u}=\rho\pi_{u}\rho^{-1}$.

We see that the  problem  of enumeration of square-tiled surfaces
can  be reformulated  as  a problem of  enumeration  of pairs  of
permutations of at most  $N$  elements such that their commutator
decomposes into  a given number  of cycles of given lengths. Here
the pairs of permutations are considered  up  to  a  simultaneous
conjugation. This problem was solved by  S.~Bloch and A.~Okounkov
by using methods of representation theory.  However,  it  is  not
directly applicable  to  our problem. Describing the square-tiled
surfaces  in  terms of pairs of  permutations  one has to add  an
additional explicit  constraint  that  the resulting square-tiled
surface is \emph{connected}! Taking a random pair of permutations
of very  large number $N\gg 1$  of elements realizing  some fixed
combinatorics of the commutator we  usually  get  a  disconnected
surface!

The necessary correction is quite nontrivial. It was performed by
A.~Eskin and A.~Okounkov in~\cite{zorich:Eskin:Okounkov}. In the further
paper            A.~Eskin,            A.~Okounkov             and
R.~Pandharipande~\cite{zorich:Eskin:Okunkov:Pandharipande}   give    the
volumes\index{Moduli space!volume of the moduli space}
of all individual connected components of the strata; see
also very nice and accessible survey~\cite{zorich:Eskin:Handbook}.

For  a  given  square-tiled  surface $S$ denote by  $Aut(S)$  its
automorphism group. Here we count only  those automorphisms which
isometrically send each square  of  the tiling to another square.
For most of the  square-tiled  surfaces $Aut(S)$ is trivial; even
for those rare square-tiled surfaces, which have nontrivial inner
symmetries  the  group $Aut(S)$ is obviously finite. We  complete
this section with the following arithmetic Theorem which confirms
two conjectures of M.~Kontsevich.

\begin{NNTheorem}[A.~Eskin, A.~Okounkov, R.~Pandharipande]
For every connected component of every stratum the generating
function
$$
\sum_{N=1}^{\infty} q^N
\sum_{\substack{N\!\text{-square-tiled}\\ \text{surfaces }S}}
\cfrac{1}{|Aut(S)|}
$$
is a quasimodular form, i.e.  a  polynomial  in Eisenstein series
$G_2(q)$, $G_4(q)$, $G_6(q)$.

Volume\index{Moduli space!volume of the moduli space}
$\Vol(\cH_1^{comp}(d_1,\dots,d_\noz)$ of  every  connected
component of every  stratum is a rational multiple of $\pi^{2g}$,
where $g$ is the genus, $d_1+\dots+d_\noz=2g-2$.
\end{NNTheorem}

\section{Crash Course in Teichm\"uller Theory}
\label{zorich:s:Crash:Course:in:Teichmuller:Theory}

In  this  section we  present  the  Teichm\"uller  theorem  about
extremal quasiconformal map and define Teichm\"uller metric. This
enables us to explain  finally  in what sense the ``Teichm\"uller
geodesic  flow''\index{Teichm\"uller!geodesic flow}
(which  we  initially defined in terms  of  the
action of the subgroup of  diagonal  matrices  in $SL(2,\R{})$
on flat surfaces)
\index{Action on the moduli space!ofSL@of $SL(2,\R{})$}\index{0SL@$SL(2,\R{})$-action on the moduli space}
is a \emph{geodesic} flow.

\index{Quasiconformal!coefficient of quasiconformality|(}
\index{Quasiconformal!extremal quasiconformal map|(}

\subsection{Extremal Quasiconformal Map}
\label{zorich:ss:Extremal:Quasiconformal:Map}

\paragraph{Coefficient of Quasiconformality}

\begin{figure}[htb]
\centering
\includegraphics{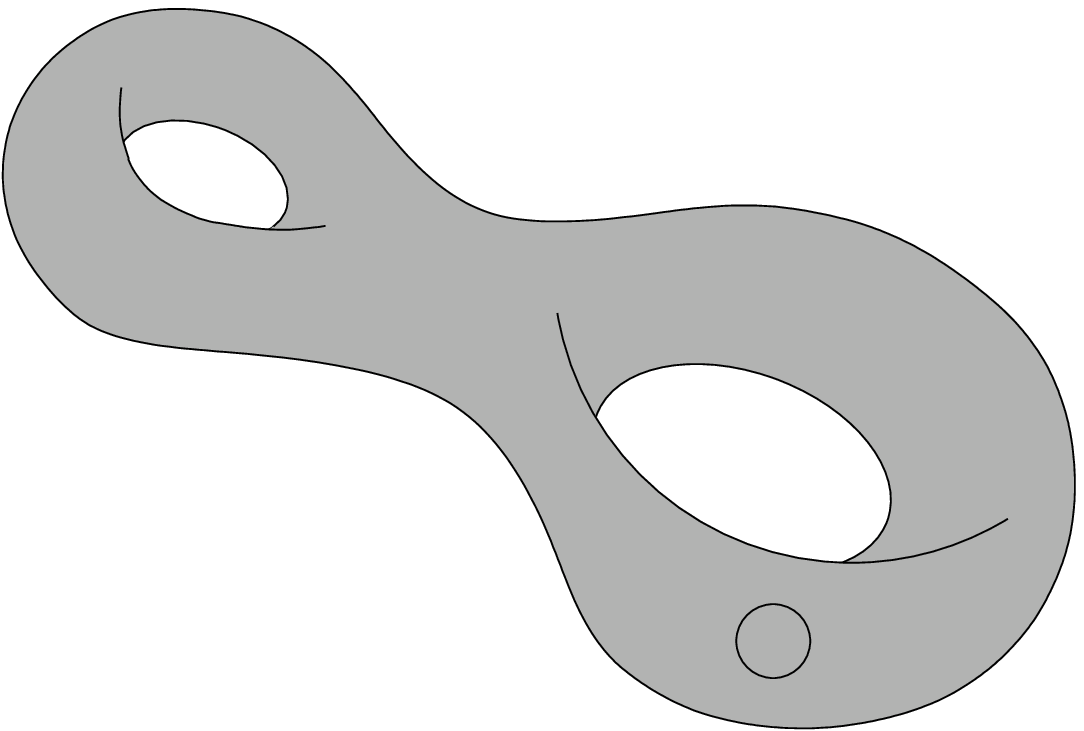}
\includegraphics{zorich_krendel_new.eps}
\includegraphics{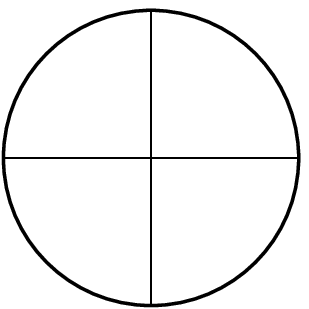}
\includegraphics{zorich_circle.eps}
\begin{picture}(0,0)(0,0)
\put(-12,-79){$\to$}
\end{picture}
\vspace{100bp}
\caption{
\label{zorich:fig:quasiconformal:map}
Quasiconformal map}
\end{figure}

Consider  a  closed  topological  surface  of  genus $g$ and  two
complex  structures  on  it.  Let  $S_0$   and   $S_1$   be   the
corresponding Riemann  surfaces.  When the complex structures are
different there are  no  conformal maps  from  $S_0$ to $S_1$.  A
smooth  map  $f:   S_0\to  S_1$  (or,  being  more  precise,  its
derivative $Df$) sends an  infinitesimal  circle at $x\in S_0$ to
an infinitesimal ellipse, see Fig.~\ref{zorich:fig:quasiconformal:map}.

\emph{Coefficient of quasiconformality}\index{Quasiconformal!coefficient of quasiconformality}
of $f$ at $x\in S_0$ is the ratio
$$
K_x(f)=\frac{a}{b}
$$
of demi-axis of this ellipse. \emph{Coefficient of
quasiconformality} of $f$ is
$$
K(f)=\sup_{x\in S_0} K_x(f)
$$
Though $S_0$  is a compact Riemann surface we  use $\sup$ and not
$\max$  since  the smooth  map  $f$ is  allowed  to have  several
isolated  critical   points   where   $Df   =0$   (and  not  only
$\det(Df)=0$) and where $K_x(f)$ is not defined.

\index{Quasiconformal!coefficient of quasiconformality|)}
\index{Quasiconformal!extremal quasiconformal map|)}

\paragraph{Half-translation Structure}

There is  a class of flat metrics which  is slightly more general
than  the  \emph{very  flat}  metrics which we consider  in  this
paper. Namely,  we can allow to a  flat metric  to have the  most
simple nontrivial linear  holonomy which is only possible: we can
allow to a  tangent  vector to change its  sign  after a parallel
transport along some closed  loops  (see the discussion of linear
holonomy\index{Holonomy}
in Sec.~\ref{zorich:ss:Very:Flat:Surfaces}).

Surfaces endowed with such flat structures are called
\emph{half-translation surfaces}\index{Surface!half-translation|(}\index{Half-translation surface|(}.
A holomorphic one-form (also called a holomorphic differential or
an Abelian differential) is an  analytic  object  representing  a
\emph{translation}  surface  (in our  terminology,  a  \emph{very
flat} surface). A holomorphic
\emph{quadratic} differential\index{Holomorphic!quadratic differential}
is an analytic object representing a half-translation surface.

In local coordinate $w$ a quadratic  differential  has  the  form
$q=q(w) (dw)^2$. In other words, the tensor rule for $q$ is
\begin{equation}
\label{zorich:eq:quadratic:differential}
q=q(w)(dw)^2=q\left(w(w')\right)\cdot
\left(\frac{dw}{dw'}\right)^2 \cdot\, (dw')^2
\end{equation}
under a change of coordinate $w=w(w')$.

One should not confuse $(dw)^2$  with  a  wedge product $dw\wedge
dw$  which equals  to  zero!  It  is just  a tensor  of  the type
described by the  tensor  rule~\eqref{zorich:eq:quadratic:differential}.
In  particular,   any   holomorphic  one-form  defined  in  local
coordinates  as  $\omega=\omega(w)  dw$  canonically  defines   a
quadratic differential\index{Holomorphic!quadratic differential}
$\omega^2 = \omega^2(w) (dw)^2$.

Reciprocally, a holomorphic quadratic differential $q=q(w)(dw)^2$
locally defines a pair of holomorphic one forms $\pm\sqrt{q(w)}\,
dw$ in any  simply-connected  domain where $q(w)\neq 0$. However,
for a generic holomorphic quadratic differential neither of these
1-forms  is  globally defined: trying to extend  the  local  form
$\omega_+=\sqrt{q(w)}\, dw$ along a  closed  path we may get back
with the form $\omega_- =-\sqrt{q(w)}\, dw$.

Recall  that  there  is  a  bijection  between  \emph{very  flat}
(=translation) surfaces  and  holomorphic  1-forms.  There  is  a
similar   bijection   between   half-translation   surfaces   and
holomorphic quadratic  differentials, where similar to the ``very
flat''  case  a   flat   surface  corresponding  to  a  quadratic
differential is polarized: it is endowed  with canonical vertical
and horizontal directions. (They can be defined locally using the
holomorphic  one-forms  $\omega_\pm=\pm\sqrt{q(w)}\, dw$.)  Note,
however, that we cannot distinguish anymore  between direction to
the North and to the South,  or between direction to the East and
to the West  unless the quadratic  differential $q$ is  a  global
square  of  a  holomorphic  1-form $\omega$. In  particular,  the
vertical and horizontal foliations  are  \emph{nonorientable} for
generic quadratic differentials\index{Holomorphic!quadratic differential}.

\paragraph{Teichm\"uller Theorem}

Choose any two complex structures on  a  topological  surface  of
genus $g\ge 1$; let $S_0$ and $S_1$ be  the corresponding Riemann
surfaces. Developing ideas of Gr\"otzsch Teichm\"uller has proved
a Theorem which we adapt to our language.

Note that \emph{flat structure} used  in  the  formulation of the
Theorem below is slightly more general  than  one  considered  in
Sec.~\ref{zorich:ss:Very:Flat:Surfaces}              and              in
Convention~\ref{zorich:conv:flat:surface}:   it   corresponds    to    a
half-translation structure and to a holomorphic  \emph{quadratic}\index{Holomorphic!quadratic differential}
differential (see above in this section). In particular, speaking
about a flat metric compatible with a given  complex structure we
mean a  flat  metric  corresponding  to  a quadratic differential
holomorphic in the given complex structure.

\index{Quasiconformal!coefficient of quasiconformality|(}
\index{Quasiconformal!extremal quasiconformal map|(}

\begin{NNTheorem}[Teichm\"uller]
For any pair  $S_0,  S_1$ of Riemann surfaces  of  genus $g\ge 1$
there exist an  extremal map $f_0:S_0\to S_1$ which minimizes the
coefficient of quasiconformality $K(f)$.

For this extremal map $f_0$ the  coefficient of quasiconformality
is constant on $S_0$ (outside of a finite  collection of singular
points)
$$
K_x(f_0)=K(f_0) \quad \forall x \in S_0 -\{P_1, \dots, P_\noz\}
$$

One  can   choose  a  flat  metric
(half-translation  structure\index{Surface!half-translation|)}\index{Half-translation surface|)})
compatible  with  the complex structure in which foliation  along
big  (correspondingly   small)   demi-axis  of  ellipses  is  the
horizontal  (correspondingly  vertical)  foliation  in  the  flat
metric.

In   flat   coordinates  the   extremal   map   $f_0$   is   just
expansion-contraction with coefficient $\sqrt{K}$.
\end{NNTheorem}

\begin{figure}[htb]
\centering
\includegraphics{zorich_octagon_houches.eps}
\includegraphics{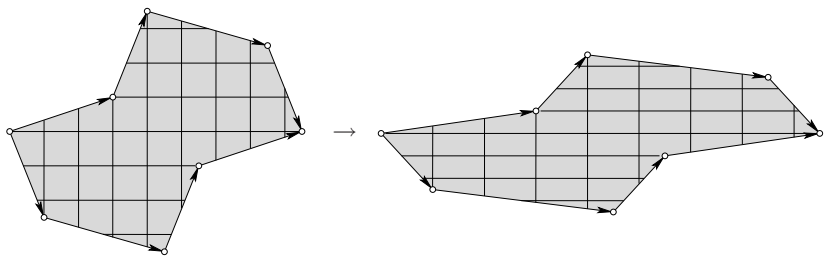}
\begin{picture}(0,0)(0,0)
\put(-33,-41){$\to$}
\end{picture}
\vspace{90bp}
\caption{
\label{zorich:fig:extremal:map}
In  flat  coordinates  the  extremal   map   $f_0$   is  just  an
expansion-contraction linear map
}
\end{figure}

\subsection{Teichm\"uller Metric and Teichm\"uller Geodesic Flow}
\label{zorich:ss:Teichmuller:Metric:and:Teichmuller:Geodesic:Flow}

Now  everything   is  ready  to  define  the
\emph{Teichm\"uller metric}\index{Teichm\"uller!metric}\index{Metric!Teichm\"uller metric}.
In  this  metric  we  measure  the  distance between two  complex
structures as
$$
dist(S_0,S_1) =  \cfrac{1}{2}\, \log K(f_0),
$$
where $f_0: S_0\to S_1$ is the extremal map.

\index{Quasiconformal!coefficient of quasiconformality|)}
\index{Quasiconformal!extremal quasiconformal map|)}

It means  that  a  holomorphic  quadratic  differential defines a
direction of deformation  of the complex structure and a geodesic
in the  Teichm\"uller  metric.  Namely,  a  holomorphic quadratic
differential defines a flat metric.  A  one-parameter  family  of
maps  which  in  the  flat  coordinates  are defined by  diagonal
matrices
$$
g_t=\begin{pmatrix} e^t & 0 \\ 0 & e^{-t} \end{pmatrix}
$$
is  a  one-parameter family  of  extremal  maps,  so  it  forms a
Teichm\"uller geodesic\index{Teichm\"uller!geodesic flow}.
According to the definition above we have
$$
dist(S_0, g_t S_0) = t.
$$

\begin{NNRemark}
Note that the Teichm\"uller metric is not a Riemannian metric but
a Finsler  metric: it does not correspond to  a quadratic form in
the tangent space, but just to a norm  which depends continuously
on the point of the space of complex structures.
\end{NNRemark}

It is known that the space  of complex structures on a surface of
genus $g\ge  2$ has complex  dimension $3g-3$. We have seen, that
the space  of  pairs  (complex  structure,  holomorphic quadratic
differential) can be identified with a tangent space to the space
of  complex  structures, in particular, it has complex  dimension
$6g-6$. Taking into consideration the functorial  behavior of the
space  of   pairs   (complex   structure,  holomorphic  quadratic
differential)\index{Moduli space!of quadratic differentials}
one can check, that it should be  identified with a
total space of a \emph{cotangent} bundle.

\section{Hope for a Magic Wand and Recent Results}
\label{zorich:s:Hope:for:a:Magic:Stick:and:Recent:Results}

This section  is devoted to  one of the most challenging problems
in  the  theory of  flat  surfaces: to  the  problem of  complete
classification  of   the   closures   of   \emph{all}  orbits  of
$GL^+(2,\R{})$\index{Action on the moduli space!ofGL@of $GL^+(2,\R{})$}\index{0GL@$GL^+(2,\R{})$-action on the moduli space}
on  the  moduli  spaces  of  Abelian (and quadratic)
differentials.  This  problem was very recently solved for  genus
two in the works of K.~Calta  and of C.~McMullen; we give a short
survey           of           their          results           in
Sec.~\ref{zorich:ss:Revolution:in:Genus:Two}
and~\ref{zorich:ss:Teichmuller:Discs}.

%
\subsection{Complex Geodesics}
\label{zorich:ss:SL2R:action:in:geometric:terms}

In  this  section  we  are  following  the geometric approach  of
C.~McMullen         developed          in~\cite{zorich:McMullen:Hilbert}
and~\cite{zorich:McMullen:genus:2}.

Fix the genus $g$  of the surfaces. We have seen in  the previous
section that  we can identify the  space $\cQ$\index{0M10@$\cQ$ -- moduli space of quadratic differentials}\index{Moduli space!of quadratic differentials}
of  pairs (complex
structure,  holomorphic  quadratic  differential) with the  total
space of the  (co)tangent  bundle to  the  moduli space
$\cM$\index{0M10@$\cM_g$ -- moduli space of complex structures}\index{Moduli space!of complex structures} of
complex structures.  Space  $\cH$  of  pairs  (complex structure,
holomorphic quadratic  differential)  can  be  identified  with a
subspace in $\cQ$  of those quadratic differentials, which can be
represented  as  global  squares  of  holomorphic  1-forms.  This
subspace forms a vector subbundle  of  special  directions in the
(co)tangent space  which we denote  by the same symbol $\cH$. The
``unit hyperboloid''\index{Unit hyperboloid}
$\cH_1\subset \cH$  of  holomorphic 1-forms corresponding to flat
surfaces  of  unit area  can  be  considered  as  a  subbundle of
\emph{unit  vectors}  in   $\cH$.   It  is  invariant  under  the
Teichm\"uller  geodesic  flow  --  the  geodesic   flow  for the
Teichm\"uller metric\index{Teichm\"uller!geodesic flow}.

One can check that an $SL(2;\R{})$-orbit
\index{0SL@$SL(2,\R{})$-action on the moduli space|(}
\index{Action on the moduli space!ofSL@of $SL(2,\R{})$|(}
in $\cH_1$ descends to a commutative diagram
\begin{equation}
\label{zorich:eq:cd}
\begin{CD}
SL(2;\R{})  @>>>  \cH_1\\
@VVV         @VVV    \\
SL(2;\R{})/SO(2;\R{}) \simeq\Hyp @>>> \cM,
\end{CD}
\end{equation}
which we interpret as
\begin{equation*}
\begin{CD}
\begin{pmatrix}
\text{Unit tangent}\\
\text{bundle to}\\
\text{hyperbolic plane}
\end{pmatrix}
  @>>>
\begin{pmatrix}
\text{Unit tangent}\\
\text{subbundle to}\\
\text{moduli space}
\end{pmatrix}
\\
@VVV         @VVV    \\
\begin{pmatrix}
\text{Hyperbolic}\\
\text{plane}
\end{pmatrix}
@>>>
\begin{pmatrix}
\text{Moduli}\\
\text{space}
\end{pmatrix}
\end{CD}
\end{equation*}

Moreover, it can be checked  that  the map $\Hyp\to \cM$\index{0M10@$\cM_g$ -- moduli space of complex structures}\index{Moduli space!of complex structures}
in  this
diagram is an isometry with  respect  to  the standard hyperbolic
metric  on  $\Hyp$  and  Teichm\"uller  metric  on  $\cM$.  Thus,
following C.~McMullen it is natural  to  call  the projections of
$SL(2;\R{})$-orbits to $\cM$  (which  coincide with the images of
$\Hyp$)  \emph{complex   geodesics}\index{Geodesic!complex geodesic}.   Another   name  for  these
projections is
\emph{Teichm\"uller discs}\index{Teichm\"uller!disc}.

\subsection{Geometric Counterparts of Ratner's Theorem}
\label{zorich:ss:Geometric:Counterparts:of:Ratners:Theorem}

Though it  is proved that  the moduli space of complex structures\index{Moduli space!of complex structures}
is  not  a hyperbolic  manifold (see~\cite{zorich:Masur:not:hyperbolic})
there is a strong  hope  that with respect to $SL(2,\R{})$-action
on $\cH$\index{Moduli space!of holomorphic 1-forms}
and on
$\cQ$\index{0M10@$\cQ$ -- moduli space of quadratic differentials}\index{Moduli space!of quadratic differentials}
the moduli space behaves as if it is.

In this  section we present  several facts about group actions on
homogeneous spaces  and  several  related  facts  about  geodesic
submanifolds.  We  warn  the  reader  that  our selection is  not
representative; it opens only a  narrow  slit  to the fascinating
world of  interactions  of  group  actions,  rigidity, hyperbolic
geometry, dynamics and number theory.

We  start  with  an  informal  formulation  of  part  of Ratner's
Theorem (see much better exposition adopted to our subject in the
survey of A.~Eskin~\cite{zorich:Eskin:Handbook}).

A discrete subgroup  $\Gamma$  of  a Lie group $G$  is  called  a
\emph{lattice}\index{Lattice!subgroup}
if a homogeneous space $G/\Gamma$ has finite volume.

\begin{NNTheorem}[M.~Ratner]
Let $G$ be  a  connected Lie group and  $U$  a connected subgroup
generated  by   unipotent   elements.   Then,   for  any  lattice
$\Gamma\subset  G$  and any $x\in G/\Gamma$ the  closure  of  the
orbit $U x$ in $G/\Gamma$  is  an orbit of some closed  algebraic
subgroup of $G$.
\end{NNTheorem}

We  would  like to point out  why  this theorem is so  remarkably
powerful. Considering a dynamical system, even an ergodic one, it
is possible to  get  a  lot of information about  a  generic  (in
measure-theoretical sense) trajectory. However, usually there are
plenty of trajectories  having  rather particular behavior. It is
sufficient to  consider geodesic flow  on a surface with cusps to
find trajectories  with  closures  producing  fairly  wild  sets.
Ratner's theorem proves,  that the closure of \emph{any} orbit of
a  unipotent  group acting  on  a  homogeneous  space  is  a nice
homogeneous space.

Ratner's theorem  has  numerous important relations with geometry
of homogeneous spaces. As an illustration we have chosen a result
of N.~Shah~\cite{zorich:Shah} and a generalization  of  his  result  for
noncompact      hyperbolic      manifolds       obtained       by
T.~Payne~\cite{zorich:Payne}.

\begin{NNTheorem}[N.~Shah]
In a compact manifold of constant negative curvature, the closure
of  a   totally  geodesic,  complete  (immersed)  submanifold  of
dimension  at  least  $2$  is   a   totally   geodesic   immersed
submanifold.
\end{NNTheorem}

\begin{NNTheorem}[T.~Payne]
Let  $f:M_1\to  M_2$ be  a  totally  geodesic  immersion  between
locally symmetric spaces of noncompact type, with $M_2$ of finite
volume. Then the  closure  of  the image of $f$  is  an  immersed
submanifold. Moreover, when $M_1$  and  $M_2$ have the same rank,
the closure of the image is a totally geodesic submanifold.
\end{NNTheorem}

\subsection{Main Hope}
\label{zorich:ss:Main:Hope}

\paragraph{Main Conjecture and its Possible Applications}

If only  the moduli space of complex structures  $\cM$\index{0M10@$\cM_g$ -- moduli space of complex structures}\index{Moduli space!of complex structures}
would be a
homogeneous space we would immediately apply the Theorem above to
diagram~\eqref{zorich:eq:cd} and  would  solve  considerable part of our
problems. But it is not. Nevertheless, there is a strong hope for
an analogous Theorem.

\index{Action on the moduli space!ofGL@of $GL^+(2,\R{})$|(}
\index{0GL@$GL^+(2,\R{})$-action on the moduli space|(}

\begin{Problem}
   %
Classify    the     closures     of
$GL^+(2,\R{})$-orbits in $\cH_g$\index{0H10@$\cH_g$ -- moduli space of holomorphic 1-forms}\index{Moduli space!of holomorphic 1-forms}
and in $\cQ_g$\index{0M10@$\cQ$ -- moduli space of quadratic differentials}\index{Moduli space!of quadratic differentials}.
Classify  the orbit  closures    of    the    unipotent    subgroup
$\begin{pmatrix}1&t\\0&1\end{pmatrix}_{t\in\R{}}$ on $\cH_g$\index{0H10@$\cH_g$ -- moduli space of holomorphic 1-forms}\index{Moduli space!of holomorphic 1-forms} and on
$\cQ_g$\index{0M10@$\cQ$ -- moduli space of quadratic differentials}\index{Moduli space!of quadratic differentials}.
\end{Problem}

The following  Conjecture is one of  the key conjectures  in this
area.

\begin{NNConjecture}
The closure $\cC(S)$  of a $GL^+(2,\R{})$-orbit of any flat surface
$S\in\cH$ (or $S\in\cQ$\index{0M10@$\cQ$ -- moduli space of quadratic differentials}\index{Moduli space!of quadratic differentials})
is a complex-algebraic suborbifold.
\end{NNConjecture}

\begin{NNRemark}
We do not discuss here the problems related with possible \emph{compactifications}
of the moduli spaces $\cH_g$\index{0H10@$\cH_g$ -- moduli space of holomorphic 1-forms}\index{Moduli space!of holomorphic 1-forms} and
$\cQ_g$\index{0M10@$\cQ$ -- moduli space of quadratic differentials}\index{Moduli space!of quadratic differentials}.
A complex-analytic description of
a compactification of $\cQ_g$\index{0M10@$\cQ$ -- moduli space of quadratic differentials}\index{Moduli space!of quadratic differentials}
can be found in the papers of J.~Fay~\cite{zorich:Fay}
and H.~Masur~\cite{zorich:Masur:compactification}.
\end{NNRemark}

Recall  that   according   to   Theorem   of  M.~Kontsevich  (see
Sec.~\ref{zorich:ss:General:Philosophy})   any    $GL^+(2,\R{})$-invariant
complex suborbifold in $\cH$ is represented by an affine subspace
in  cohomological  coordinates. Thus, if the Conjecture above  is
true,  the  structure  of   orbit   closures  of  the  action  of
$GL^+(2,\R{})$  on  $\cH$  and  on  $\cQ$\index{0M10@$\cQ$ -- moduli space of quadratic differentials}\index{Moduli space!of quadratic differentials}
(and of $SL(2,\R{})$  on
$\cH_1$ and  on $\cQ_1$) would  be as  simple as in  the case  of
homogeneous spaces.

We have not  discussed the aspects of Ratner's Theorem concerning
the measures;  it states more  than we cited above. Actually, not
only orbit closures have a nice form, but  also invariant ergodic
measures; namely, all such measures are just the natural measures
supported  on  orbits  of  closed  subgroups.  Trying  to  make a
parallel with Ratner's  Theorem  one should extend the Conjecture
above to invariant measures.

In the most  optimistic  hopes the  study  of an individual  flat
surface  $S\in\cH(d_1,\dots,d_\noz)$   would  look  as   follows.
(Frankly speaking, here  we  slightly exaggerate in our scenario,
but after all  we  are describing  the  dreams.) To describe  all
geometric properties of  a  flat surface  $S$  we first find  the
orbit closure $\cC(S)=\overline{GL^+(2,\R{})\,  S}\subset  \cH(d_1,
\dots, d_\noz)$ (our optimistic hope assumes  that  there  is  an
efficient  way  to  do  this).  Then  we consult a  (conjectural)
classification list and find  $\cC(S)$  in some magic table which
gives all information about $\cC(S)$ (like volume, description of
cusps, Siegel--Veech constants, Lyapunov exponents, adjacency  to
other  invariant  subspaces, etc). Using this information we  get
answers to all  possible questions which  one can ask  about  the
initial flat surface $S$.

Billiards  in  rational polygons  give  an  example  of  possible
implementation of this optimistic scenario. Fixing  angles of the
polygon which defines a billiard table we can  change the lengths
of its sides. We  get a family $\cB$ of polygons which  induces a
family $\tilde\cB$  of flat surfaces obtained by Katok--Zemlyakov
construction   (see   Sec.~\ref{zorich:ss:Billiards:in:Polygons}).  This
family    $\tilde\cB$    belongs   to    some    fixed    stratum
$\tilde\cB\subset\cH(d_1,  \dots,  d_\noz)$. However,  it  has  a
nontrivial codimension in the stratum, so $\tilde\cB$ has measure
zero  and  one cannot  use  ergodic theorem  naively  to get  any
information about  billiards  in corresponding polygons. Having a
version of ergodic  theorem  which treats \emph{all} orbits (like
in Ratner's Theorem)  presumably it would  be possible to  get  a
powerful tool for the study of rational billiards.

\begin{Exercise}
Consider   the   family  $\cB$   of   billiard   tables   as   on
Fig.~\ref{zorich:fig:family:of:rectangular:polygons}.   Determine    the
stratum  $\cH(d_1,   \dots,   d_\noz)$   to   which   belong  the
corresponding  flat  surfaces and compute the codimension of  the
resulting family $\tilde\cB\subset\cH(d_1, \dots, d_\noz)$.
\end{Exercise}

\paragraph{Content of Remaining Sections}

The Conjecture above is trivial for genus one, since in this case
the ``Teichm\"uller  space  of  Riemann  surfaces  of genus one''
coincides  with  an  upper  half-plane,  and   the  entire  space
coincides with a single
Teichm\"uller  disc\index{Teichm\"uller!disc}
(image  of $\Hyp$ in
diagram~\eqref{zorich:eq:cd}).

Very  recently   C.~McMullen   proved  the  Conjecture  in  genus
two~\cite{zorich:McMullen:Hilbert}, \cite{zorich:McMullen:genus:2}, and this is
a  highly nontrivial  result.  We give a  short  report of
revolutionary   results    of   K.~Calta~\cite{zorich:Calta}   and    of
C.~McMullen~\cite{zorich:McMullen:Hilbert}--\cite{zorich:McMullen:decagon:proof}
in  Sec.~\ref{zorich:ss:Revolution:in:Genus:Two}  below. However,  their
techniques, do  not  allow any straightforward generalizations to
higher genera: Riemann  surfaces of genus two are rather special,
in particular, every such surface is hyperelliptic.

In the next two
sections~\ref{zorich:ss:Classification:of:Connected:Components:of:the:Strata}
and~\ref{zorich:ss:Veech:Surfaces} we try to give an  idea  of  what  is
known about invariant  submanifolds in higher genera (which is an
easy task since, unfortunately, little is known).

Having an  invariant  submanifold  $\cK\subset\cH_g$ (or $\cK\subset\cQ_g$\index{0M10@$\cQ$ -- moduli space of quadratic differentials}\index{Moduli space!of quadratic differentials})
in genus $g$
one    can    construct    a     new     invariant    submanifold
$\tilde\cK\subset\cH_{\tilde g}$ (correspondingly $\tilde\cK\subset\cQ_{\tilde g}$)
in  higher  genus $\tilde g >g$
replacing every  $S\in\cK$  by  an  appropriate ramified covering
$\tilde S$  over $S$ of some fixed combinatorial  type. We do not
want to specify here what does  a  ``fixed  combinatorial  type''
mean; what  we claim is that  having an invariant  manifold $\cK$
there is  some procedure which  allows to construct a whole bunch
of new  invariant  submanifolds  $\tilde\cK$  for  higher  genera
$\tilde g >g$.

In some cases all quadratic differentials in the
invariant submanifold
$\tilde\cK$ obtained by a ramified covering construction
from some $\cK\subset\cQ_g$\index{0M10@$\cQ$ -- moduli space of quadratic differentials}\index{Moduli space!of quadratic differentials}
might become
global squares of Abelian differentials. Hence, using special ramified
coverings one can construct $GL^+(2,\R{})$-invariant   submanifolds
$\tilde\cK\subset\cH_{\tilde g}$ from invariant submanifolds
$\cK\subset\cQ_g$.

What is  really interesting to understand is what invariant
manifolds  form   the   ``roots''  of  such  constructions.  Such
invariant manifolds are often called the \emph{primitive} ones.

In the following two
sections~\ref{zorich:ss:Classification:of:Connected:Components:of:the:Strata}
and~\ref{zorich:ss:Veech:Surfaces} we consider the two extremal classes
of primitive invariant submanifolds: the largest ones and the
smallest ones. Namely, in
Sec.~\ref{zorich:ss:Classification:of:Connected:Components:of:the:Strata}
we present a classification of connected components of the strata
$\cH(d_1,  \dots,   d_\noz)$.   It  follows  from  ergodicity  of
$SL(2,\R{})$-action on $\cH_1(d_1, \dots, d_\noz)$ that the orbit
closure of almost any  flat  surface in $\cH(d_1, \dots, d_\noz)$
coincides with  the  embodying  connected  component of $\cH(d_1,
\dots, d_\noz)$.

In  section~\ref{zorich:ss:Veech:Surfaces}  we   consider  the  smallest
possible   $GL^+(2,\R{})$-invariant   submanifolds:   those   which
correspond to  closed  orbits.  Teichm\"uller  discs  obtained as
projections of  such orbits to  the moduli space $\cM$\index{0M10@$\cM_g$ -- moduli space of complex structures}\index{Moduli space!of complex structures}
of complex
structures form  the  ``closed  complex  geodesics''\index{Geodesic!complex geodesic}  --- Riemann
surfaces with cusps.

\index{Action on the moduli space!ofGL@of $GL^+(2,\R{})$|)}
\index{0GL@$GL^+(2,\R{})$-action on the moduli space|)}

\index{Action on the moduli space!ofSL@of $SL(2,\R{})$|)}
\index{0SL@$SL(2,\R{})$-action on the moduli space|)}

\subsection{Classification of Connected Components of the Strata}
\label{zorich:ss:Classification:of:Connected:Components:of:the:Strata}

In  order  to formulate the classification theorem for  connected
components of  the  strata  $\cH(d_1,\dots,d_\noz)$\index{0H20@$\cH(d_1,\dots,d_\noz)$ -- stratum in the moduli space}\index{Stratum!in the moduli space}  we  need  to
describe  the  classifying  invariants.  There are two  of  them:
\emph{spin structure}  and \emph{hyperellipticity}. Both  notions
are applicable only to part of the strata: flat surfaces from the
strata  $\cH(2d_1,\dots,2d_\noz)$  have \emph{even}  or \emph{odd
spin structure}.  The  strata $\cH(2g-2)$ and $\cH(g-1,g-1)$ have
special \emph{hyperelliptic connected component}.

\paragraph{Spin Structure}

Consider    a     flat     surface    $S$    from    a    stratum
$\cH(2d_1,\dots,2d_\noz)$. Let  $\rho:  S^1  \to  S$  be a smooth
closed path on $S$; here $S^1$ is a standard circle. Note that at
any point of the surfaces $S$ we know where is the ``direction to
the North''. Hence, at any point $\rho(t)=x\in S$ we can  apply a
compass and measure the direction of the tangent vector $\dot x$.
Moving along  our path $\rho(t)$  we make the tangent vector turn
in the  compass. Thus we get a map  $G(\rho):S^1\to S^1$ from the
parameter circle to the circumference of the compass. This map is
called  the   \emph{Gauss   map}.   We  define  the  \emph{index}
$\ind(\rho)$ of the path $\rho$ as a degree  of the corresponding
Gauss map (or, in other words as the algebraic number of turns of
the tangent vector around the compass) taken modulo $2$.
$$
\ind(\rho)= \deg G(\rho) \mod 2
$$

It  is  easy  to  see  that  $\ind(\rho)$  does   not  depend  on
parameterization.  Moreover,  it  does  not  change  under  small
deformations of the  path. Deforming the path more drastically we
may change its  position with respect to conical singularities of
the flat  metric. Say,  the initial path might go  on the left of
$P_k$ and its deformation might pass on the right of  $P_k$. This
deformation changes  the  $\deg  G(\rho)$.  However,  if the cone
angle  at  $P_k$  is  of  the  type  $2\pi(2d_k+1)$,  then  $\deg
G(\rho)\mod 2$ does not change!  This  observation  explains  why
$\ind(\rho)$ is well-defined  for  a free homotopy class $[\rho]$
when  $S\in\cH(2d_1,\dots,2d_\noz)$  (and hence,  when  all  cone
angles are odd multiples of $2\pi$).

Consider a collection of closed  smooth  paths  $a_1, b_1, \dots,
a_g,  b_g$   representing   a   symplectic   basis   of  homology
$H_1(S,\Z{}/2\Z{})$.   We   define  the   \emph{parity   of   the
spin-structure} of  a flat surface  $S\in\cH(2d_1,\dots,2d_\noz)$
as
\begin{equation}
\label{zorich:eq:parity:of:the:spin:structure}
\phi(S)=\sum_{i=1}^g
\left(\ind(a_i)+1\right)\left(\ind(b_i)+1\right) \mod 2
\end{equation}

\begin{NNLemma}
The value  $\phi(S)$ does not  depend on the choice of symplectic
basis  of  cycles   $\{a_i,b_i\}$.   It  does  not  change  under
continuous deformations of $S$ in $\cH(2d_1,\dots,2d_\noz)$.
\end{NNLemma}

Lemma above shows  that  the parity of the  spin  structure is an
invariant of  connected  components  of  strata  of those Abelian
differentials, which have zeroes of even degrees.

\begin{figure}[htb]
\centering
\includegraphics{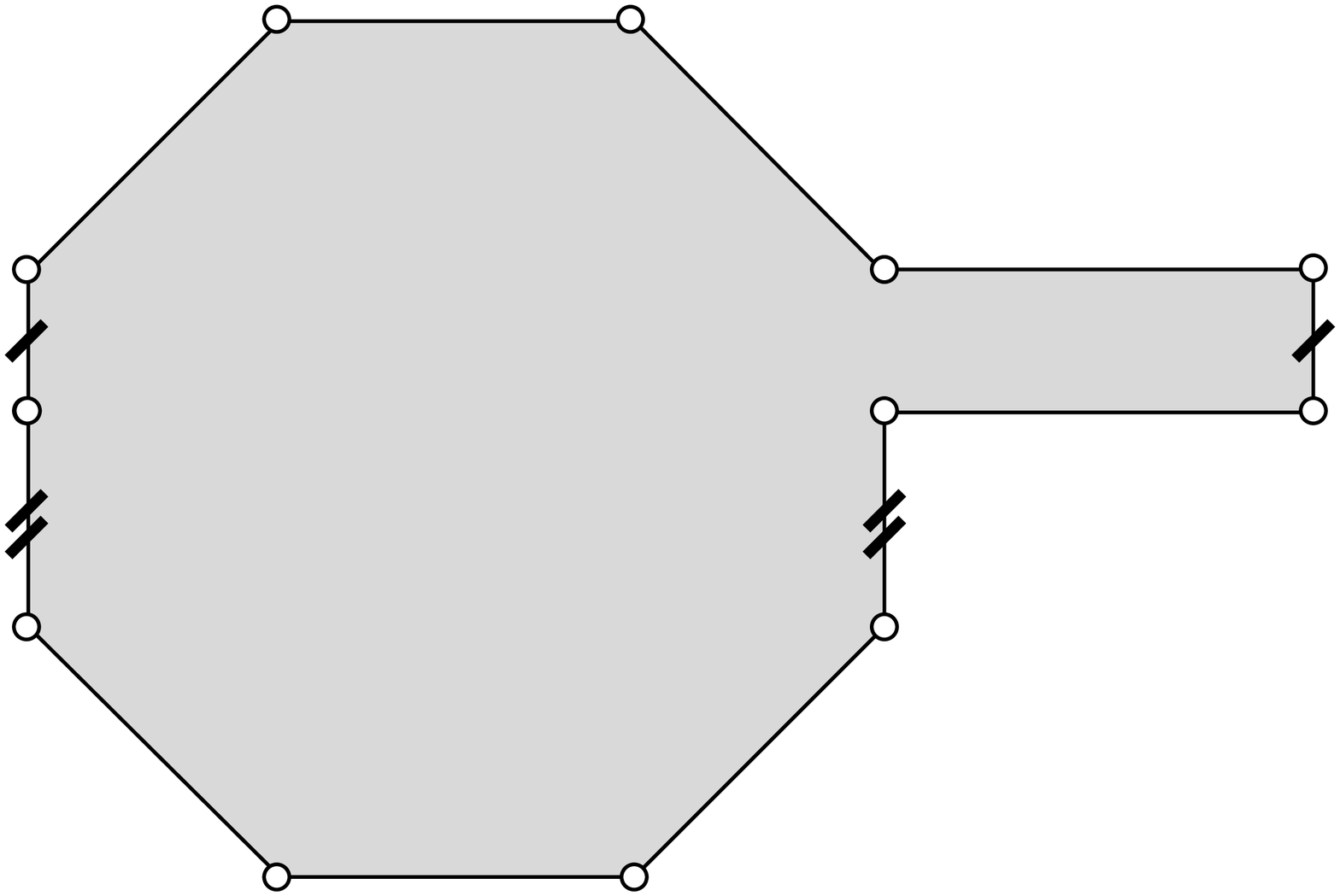}
\includegraphics{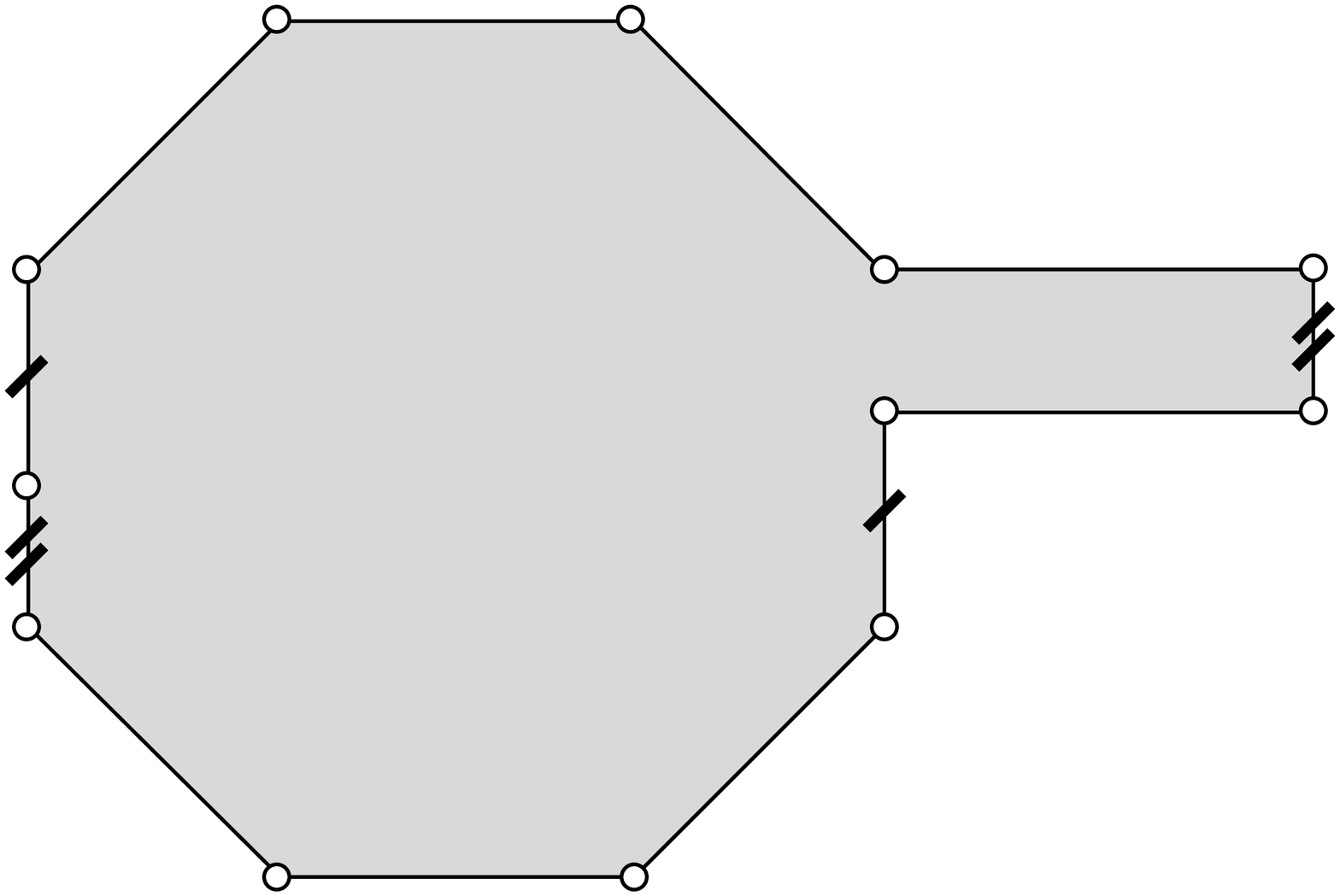}
\vspace{90bp}
\caption{
\label{zorich:fig:two:different:parities}
Attaching  a  handle to the flat surface  $S_0\in\cH(2)$  in  two
different  ways  we  get  two  flat  surfaces  in  $\cH(4)$  with
different  parities  of spin structure. Hence the resulting  flat
surfaces live  in  different  connected  components  of $\cH(4)$}
\end{figure}

\begin{Exercise}
Consider      two      flat       surfaces      presented      at
Fig.~\ref{zorich:fig:two:different:parities}.  They are  obtained  by  a
surgery which attaches a handle to a flat surface obtained from a
regular octagon. Note,  however, that the handles are attached in
two different ways  (see  the identifications of vertical sides).
Check that both surfaces belong to the same stratum $\cH(4)$.

Consider  a  symplectic basis of cycles of  the  initial  surface
(corresponding to the regular octagon) realized by paths which do
not pass through  the conical singularity. Show that a symplectic
basis of cycles for  each of two new surfaces can be  obtained by
completion  of  the initial basis with  a  pair of cycles $a,  b$
representing  the  attached handle, where the cycle  $a$  is  the
waist curve of  the handle. Calculate $\ind(a)$ and $\ind(b)$ for
each of two surfaces. Check that $ind(b)$ are different, and thus
our two flat surfaces  have  different parities of spin structure
and hence belong  to  different connected components of $\cH(4)$.
\end{Exercise}

\toread{Spin structure}

We  have  hidden  under   the   carpet  geometry  of  the  ``spin
structure''  defining  the ``parity-of-the-spin-structure''.  The
reader  can  find  details  in~\cite{zorich:Kontsevich:Zorich}  and   in
original         papers         of        M.Atiyah~\cite{zorich:Atiyah},
J.~Milnor~\cite{zorich:Milnor},       D.Mumford~\cite{zorich:Mumford}       and
D.~Johnson~\cite{zorich:Johnson}.         Recent         paper        of
C.~McMullen~\cite{zorich:McMullen:spin} contains further applications of
spin structures to flat surfaces.

\paragraph{Hyperellipticity}

A flat surface $S$ may  have  a symmetry; one specific family  of
such flat surfaces, which  are  ``more symmetric than others'' is
of a  special interest for us. Recall that  there is a one-to-one
correspondence between  flat  surfaces and pairs (Riemann surface
$M$      ,      holomorphic       1-form      $\omega$),      see
Sec.~\ref{zorich:ss:Dictionary:of:Complex:Analytic:Language}.  When  the
corresponding   Riemann   surface  is   \emph{hyperelliptic}  the
\emph{hyperelliptic  involution}  $\tau:M\to  M$   acts   on  any
holomorphic 1-form $\omega$ as $\tau^\ast\omega=-\omega$.

We say that  a  flat surface  $S$  is a
\emph{hyperelliptic  flat surface}\index{Hyperelliptic!involution}\index{Hyperelliptic!surface}
if  there  is an isometry $\tau:S\to  S$  such that $\tau$ is  an
involution,   $\tau\circ\tau=\id$,   and  the   quotient  surface
$S/\tau$  is   a   topological   sphere.   In   flat  coordinates
differential of such involution obviously satisfies $D\tau=-\Id$.

\begin{Exercise}
Check that the flat surface $S$ from Fig.~\ref{zorich:fig:suspension} is
hyperelliptic,  and  that  the  central symmetry of  the  polygon
induces the hyperelliptic involution of $S$.
\end{Exercise}

In  a   general  stratum  $\cH(d_1,\dots,d_\noz)$\index{0H20@$\cH(d_1,\dots,d_\noz)$ -- stratum in the moduli space}\index{Stratum!in the moduli space}   hyperelliptic
surfaces  form   a  small  subspace  of  nontrivial  codimension.
However, there are two special  strata,  namely  $\cH(2g-2)$  and
$\cH(g-1,g-1)$,  for  which hyperelliptic  surfaces  form  entire
\emph{hyperelliptic connected components}\index{Hyperelliptic!connected component}
$\cH^{hyp}(2g-2)$ and $\cH^{hyp}(g-1,g-1)$ correspondingly.

Note that in  the  stratum $\cH(g-1,g-1)$ there are hyperelliptic
flat surfaces of two different types.  A hyperelliptic involution
$\tau S\to S$  may fix the  conical points or  might  interchange
them.  It  is not difficult to  show  that for surfaces from  the
\emph{connected component} $\cH^{hyp}(g-1,g-1)$ the hyperelliptic
involution \emph{interchanges} the conical singularities.

The remaining  family  of  those  hyperelliptic  flat surfaces in
$\cH(g-1,g-1)$, for  which the hyperelliptic involution keeps the
saddle points fixed,  forms  a subspace of nontrivial codimension
in  the  complement $\cH(g-1,g-1)-\cH^{hyp}(g-1,g-1)$.  Thus, the
hyperelliptic connected  component $\cH^{hyp}(g-1,g-1)$ does  not
coincide with the space of all hyperelliptic flat surfaces.

\paragraph{Classification Theorem for Abelian Differentials}

Now, having introduced the classifying invariants  we can present
the classification of  connected  components of strata of Abelian
differentials.

\begin{NNTheorem}[M.~Kontsevich and A.~Zorich]
All connected components of any stratum  of Abelian differentials
on a curve of genus $g\ge 4$ are described by the following list:

The stratum $\mathcal{H}(2g-2)$ has  three  connected components:
the   hyperelliptic   one,  $\mathcal{H}^{hyp}(2g-2)$,   and  two
nonhyperelliptic   components:   $\mathcal{H}^{even}(2g-2)$   and
$\mathcal{H}^{odd}(2g-2)$  corresponding  to even  and  odd  spin
structures.

The stratum  $\mathcal{H}(2d,2d)$,  $d\ge  2$ has three connected
components:  the  hyperelliptic one,  $\mathcal{H}^{hyp}(2d,2d)$,
and two  nonhyperelliptic components: $\mathcal{H}^{even}(2d,2d)$
and $\mathcal{H}^{odd}(2d,2d)$.

All      the       other       strata       of      the      form
$\mathcal{H}(2d_1,\dots,2d_\noz)$ have  two connected components:
$\mathcal{H}^{even}(2d_1,\dots,2d_\noz)$                      and
$\mathcal{H}^{odd}(2d_1,\dots,2d_n)$, corresponding to  even  and
odd spin structures.

The stratum $\mathcal{H}(2d-1,2d-1)$, $d\ge 2$, has two connected
components;  one   of  them:  $\mathcal{H}^{hyp}(2d-1,2d-1)$   is
hyperelliptic;  the  other  $\mathcal{H}^{nonhyp}(2d-1,2d-1)$  is
not.

All the  other strata of Abelian  differentials on the  curves
of genera $g\ge 4$ are nonempty and connected.
\end{NNTheorem}

In the case of small genera $1\le g\le 3$ some components are
missing in comparison with the general case.

\begin{NNTheorem}
The moduli space  of Abelian differentials  on a curve  of  genus
$g=2$    contains    two     strata:    $\mathcal{H}(1,1)$    and
$\mathcal{H}(2)$.  Each  of them is connected and coincides  with
its hyperelliptic component.

Each of  the  strata  $\mathcal{H}(2,2)$, $\mathcal{H}(4)$ of the
moduli space  of Abelian differentials  on a curve of genus $g=3$
has two  connected  components:  the  hyperelliptic  one, and one
having odd spin structure. The other  strata  are  connected  for
genus $g=3$.
\end{NNTheorem}

Since there  is  a  one-to-one  correspondence  between connected
components of the  strata  and
\emph{extended Rauzy classes}\index{Rauzy class}\index{0R10@$\mathfrak{R}$ -- (extended) Rauzy class}
(see
Sec.~\ref{zorich:ss:space:of:iet} and paper~\cite{zorich:Yoccoz:Les:Houches} in
this collection) the Classification Theorem above classifies also
the extended Rauzy classes.

\paragraph{Classification Theorem for Quadratic Differentials}

Note that for any partition $d_1+\dots+d_\noz=2g-2$ of a positive
even  integer  $2g-2$  the  stratum  $\cH(d_1,\dots,d_\noz)$\index{0H20@$\cH(d_1,\dots,d_\noz)$ -- stratum in the moduli space}\index{Stratum!in the moduli space}   of
Abelian  differentials  is nonempty.  For  meromorphic  quadratic
differentials   with   at  most  simple  poles  there  are   four
\emph{empty} strata! Namely,

\begin{NNTheorem}[H.~Masur and J.~Smillie]
Consider a  partition of the number  $4g-4$, where $g\ge  0$ into
integers $d_1+\dots+d_\noz$  with  all  $d_j  \in \mathbb{N} \cup
\{-1 \}$.  The  corresponding  stratum $\cQ(d_1,\dots,d_\noz)$\index{0M10@$\cQ$ -- moduli space of quadratic differentials}\index{Moduli space!of quadratic differentials}
is
non-empty with the following four exceptions:
$$
 \cQ(\emptyset), \cQ(1,-1)\ (\text{in genus } g=1)
 \quad \text{and}\quad
 \cQ(4), \cQ(1,3)\ (\text{in genus } g=2)
$$
\end{NNTheorem}

Classification  of   connected   components   of  the  strata  of
meromorphic quadratic differentials with at most simple poles was
recently obtained by E.~Lanneau~\cite{zorich:Lanneau}.

\begin{NNTheorem}[E.~Lanneau]
Four    exceptional    strata     $\cQ(-1,9)$,     $\cQ(-1,3,6)$,
$\cQ(-1,3,3,3)$   and   $\cQ(12)$    of   meromorphic   quadratic
differentials\index{0M10@$\cQ$ -- moduli space of quadratic differentials}\index{Moduli space!of quadratic differentials}
contain  exactly  two connected components; none of
them hyperelliptic.

Three series of strata
$$
\begin{array}{ll}
\cQ(2(g-k)-3,2(g-k)-3,2k+1,2k+1)\quad       &k \geq -1,\ g\geq 1,\  g-k \geq 2\notag\\
\cQ(2(g-k)-3,2(g-k)-3,4k+2)            &k \geq 0, \ g\geq 1\  \text{ and } g-k \geq 1\notag\\
\cQ(4(g-k)-6,4k+2)                     &k \geq 0, \ g \geq 2\ \text{ and } g-k \geq 2\notag
\end{array}
$$
contain hyperelliptic connected  components. The strata from these
series in genera $g\ge 3$ and the strata $\cQ(-1,-1,3,3)$, $\cQ(-1,-1,6)$
in genus $g=2$ contain  exactly  two  connected  components;  one   of  them  --
hyperelliptic, the other one -- not.

The remaining strata from these series,
namely, $\cQ(1,1,1,1)$, $\cQ(1,1,2)$, $\cQ(2,2)$ in genus $g=2$ and
$\cQ(1,1,-1,-1)$, $\cQ(-1,-1,2)$ in genus $g=1$ coincide with their
hyperelliptic connected component.

All other strata  of  meromorphic quadratic differentials with at
most simple poles are connected.
\end{NNTheorem}

Recall that having  a  meromorphic quadratic differential with at
most simple poles one  can associate to it a surface with  a flat
metric which is  slightly more general that our usual \emph{very}
flat metric (see  Sec.~\ref{zorich:ss:Extremal:Quasiconformal:Map} for a
discussion of
half-translation structures\index{Surface!half-translation}\index{Half-translation surface}).

It is easy to  verify  whether a half-translation surface belongs
to a hyperelliptic  component or not. However, currently there is
no simple  and  efficient  way  to  distinguish  half-translation
surfaces from the four exceptional components.

\begin{Problem}
   %
Find an invariant  of  the half-translation structure which would
be easy to evaluate and which  would distinguish half-translation
surfaces  from   different
connected  components\index{Moduli space!connected components of the strata}
of  the  four
exceptional  strata  $\cQ(-1,9)$, $\cQ(-1,3,6)$,  $\cQ(-1,3,3,3)$
and $\cQ(12)$\index{0M10@$\cQ$ -- moduli space of quadratic differentials}\index{Moduli space!of quadratic differentials}.
\end{Problem}

Currently there are  two  ways to determine to  which  of the two
connected components of an exceptional  stratum  belongs  a  flat
surface $S$.

The first approach suggests to find a ``generalized permutation''
for an analog of the first return map of the  vertical  flow to a
horizontal segment and then to find it in one of the two extended
Rauzy classes\index{Rauzy class}
(see Sec.~\ref{zorich:ss:space:of:iet})
corresponding to  two  connected  components. Note, however, that
already for  the  stratum  $\cQ(-1,9)$\index{0M10@$\cQ$ -- moduli space of quadratic differentials}\index{Moduli space!of quadratic differentials}
the  corresponding  Rauzy
classes    contain    $97\,544$,   and    $12\,978$   generalized
permutations; the Rauzy  classes  for the components of $\cQ(12)$
contain already $894\,117$ and $150\,457$ elements.

In the second approach one studies configurations of saddle
connections (see
Sec.~\ref{zorich:ss:Multiple:Isometric:Geodesics:and:Principal:Boundary})
on the surface  $S$  and tries to find  a  configuration which is
forbidden  for  one  of  the  two  connected  components  of  the
corresponding stratum.

For  example,  for   surfaces  from  one  of  the  two  connected
components of $\cQ(-1,9)$ as soon as we have  a saddle connection
joining  the simple  pole  with the zero  we  necessarily have  a
closed geodesic  going in the  same direction. Thus, if we manage
to find on a surface $S\in\cQ(-1,9)$ a saddle connection which is
not accompanied by a parallel closed geodesic, $S$  belong to the
other connected component of $\cQ(-1,9)$.

\subsection{Veech Surfaces}
\label{zorich:ss:Veech:Surfaces}

\index{Action on the moduli space!ofSL@of $SL(2,\R{})$|(}
\index{0SL@$SL(2,\R{})$-action on the moduli space|(}

For   almost   every   flat   surface   $S$    in   any   stratum
$\cH_1(d_1,\dots,d_\noz)$\index{0H30@$\cH_1(d_1,\dots,d_\noz)$ --
``unit hyperboloid''}\index{Stratum!in the moduli space}\index{Unit hyperboloid}
the orbit $SL(2,\R{})\cdot S$ is  dense
in the  stratum and for any $g_1\neq g_2  \in SL(2,\R{})$ we have
$g_1  S  \neq g_2  S$.  However, some  flat  surfaces have  extra
symmetries. When a flat surface $S_0$ has an affine automorphism,
i.e. when  for some $g_0\in SL(2,\R{})$ we  get $g_0  S = S$  the
orbit of $S_0$ is smaller than usual.

The \emph{stabilizer} $Stab(S)\in SL(2;\R{})$, that is a subgroup
of  those $g\in  SL(2,\R{})$  for which $g  S=S$,  is called  the
\emph{Veech group}\index{Veech!group}
of the flat surface $S$ and is denoted $SL(S)$. In representation
of the flat surface $S$ in  terms of a pair (Riemann surface $X$,
holomorphic 1-form $\omega$ on it) the Veech group  is denoted as
$SL(X,\omega)$       following        the       notation       of
C.~McMullen~\cite{zorich:McMullen:Hilbert}.

Some exceptional flat surfaces $S$  possess  very  large group of
symmetry  and  their  orbits  are very small. The  flat  surfaces
having the largest possible symmetry group are called \emph{Veech
surfaces}. More precisely, a flat surface is called a
\emph{Veech surface}\index{Veech!surface|(}\index{Surface!Veech|(}
if its Veech group $SL(S)$ is a lattice in $SL(2,\R{})$  (that is
the quotient $SL(2,\R{})/SL(S)$ has finite volume).

\begin{NNTheorem}[J.~Smillie]
An  $SL(2,\R{})$-orbit  of  a  flat  surface  $S$  is  closed  in
$\cH_1(d_1,\dots,d_\noz)$\index{0H30@$\cH_1(d_1,\dots,d_\noz)$ -- ``unit hyperboloid''}\index{Stratum!in the moduli space}\index{Unit hyperboloid} if and only if $S$ is a Veech surface.
\end{NNTheorem}

Forgetting  polarization  (direction to the ``North'') of a  flat
surface we get a
\emph{Teichm\"uller disc}\index{Teichm\"uller!disc}
of $S$ (see~\eqref{zorich:eq:cd} and the comments below it)
$$
\begin{array}{rcl}
&\backslash\, SL(2,\R{}) /&\\
[-\halfbls] SO(2,\R{})&& SL(S)
\end{array}
\quad = \quad
\begin{array}{rl}
{\mathbb H} /&\\
[-\halfbls]  & SL(S)
\end{array}
$$

A flat  surface $S$ is a Veech surface  if its Teichm\"uller disc
$\Hyp/SL(S)$ has finite volume. However, even for a Veech surface
the Teichm\"uller disc is never compact:  it necessarily contains
at least one cusp. The Teichm\"uller discs of  Veech surfaces can
be considered  as \emph{closed complex geodesics}\index{Geodesic!complex geodesic} (see discussion
at the end of Sec.~\ref{zorich:ss:SL2R:action:in:geometric:terms}).

\begin{figure}[htb]
\centering
\includegraphics{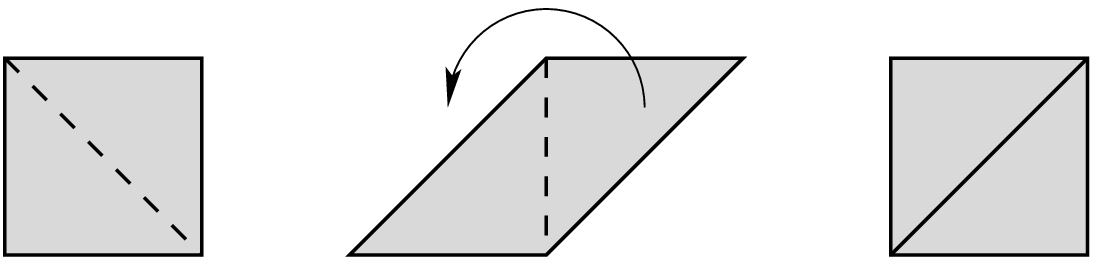}
\includegraphics{zorich_ciseaux6_1cm.eps}
\includegraphics{zorich_arrow.eps}
\begin{picture}(0,0)(0,0)
\put(-130,-72){$a$}
\put(-130,-4){$a$}
\put(-162,-35){$b$}
\put(-95,-35){$b$}
\put(-130,-30){$c$}
\put(-29,-72){$a$}
\put(30,-4){$a$}
\put(-32,-35){$b$}
\put(39,-35){$b$}
\put(5,-30){$c$}
\put(125,-72){$a$}
\put(125,-4){$a$}
\put(93,-30){$c$}
\put(160,-30){$c$}
\put(125,-35){$b$}
\put(60,-45){\Huge$=$}
\end{picture}
\vspace{90bp} 
\caption{
\label{zorich:fig:torus:stabilizer}
This linear transformation belongs to the Veech group\index{Veech!group}
of $\T{2}$
}
\end{figure}

Consider an elementary  example. As a  flat surface take  a  flat
torus       obtained       from        a       unit       square.
Fig.~\ref{zorich:fig:torus:stabilizer}    shows    why    the    element
$g_+=\begin{pmatrix}1   &   1\\0&1\end{pmatrix}\in    SL(2,\Z{})$
belongs to a stabilizer of $S$.

Similarly      the      element      $g_-=\begin{pmatrix}1      &
0\\1&1\end{pmatrix}\in  SL(2,\R{})$  also belongs  to  the  Veech
group $SL(\T{2})$ of $\T{2}$. Since  the  group  $SL(2,\Z{})$  is
generated by $g_+$ and $g_-$ we  conclude that $SL(2,\Z{})\subset
SL(\T{2})$. It  is  easy  to  check  that, actually, $SL(2,\Z{})=
SL(\T{2})$.  As  the  Teichm\"uller  disc of $\T{2}$ we  get  the
modular           curve           $\Hyp/SL(2,\Z{})$          (see
Fig.~\ref{zorich:fig:space:of:flat:tori})  which,  actually,   coincides
with the moduli space of complex structures on the torus.

\begin{figure}[htb]
\centering
\includegraphics{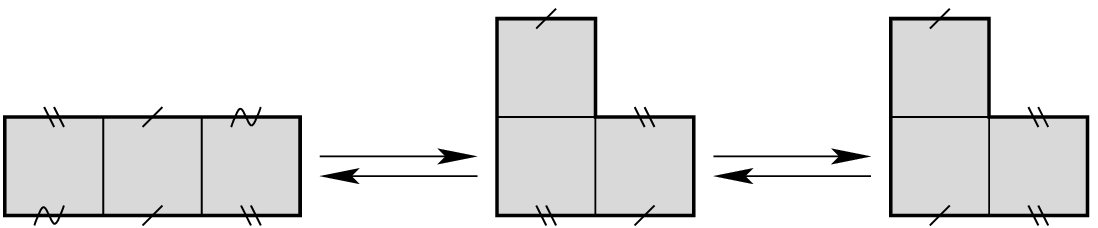}
\begin{picture}(0,0)(0,0)
\put(10,-10)
 {\begin{picture}(0,0)(0,0)
 \put(-70,-25){$\begin{pmatrix}0&1\\-1&0\end{pmatrix}$}
 \put(50,-25){$\begin{pmatrix}1&1\\0&1\end{pmatrix}$}
 \end{picture}}
\end{picture}
\vspace{80bp}
\caption{
\label{zorich:fig:L:3:surfaces}
There are three $3$-square-tiled surfaces\index{Surface!square-tiled}\index{Square-tiled surface}
in $\cH(2)$. Our
picture shows that they all belong to the same
$SL(2;\Z{})$-orbit }
\end{figure}

Consider a slightly more complicated example.

\begin{figure}[htb]
\includegraphics{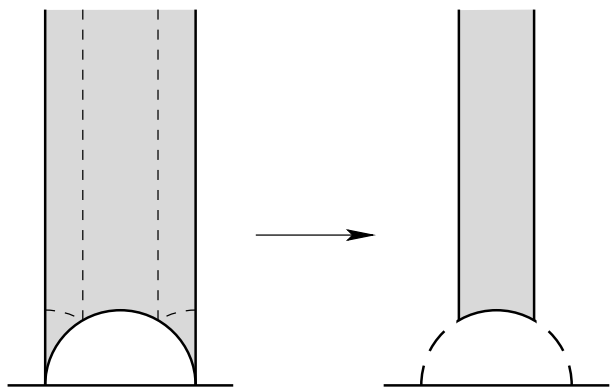}
\vspace{140bp}
\caption{
\label{zorich:fig:L:3}
Teichm\"uller discs of a $3$-square-tiled  surface  is  a  triple
cover over the modular curve
}
\end{figure}

\begin{Exercise}
Verify     that    square-tiled     surfaces     presented     at
Fig.~\ref{zorich:fig:L:3:surfaces} belong to the stratum $\cH(2)$.  Show
that there are  no other 3-square-tiled surfaces. Verify that the
linear transformations indicated  at  Fig.~\ref{zorich:fig:L:3:surfaces}
act as it  is  described on the Figure;  check  that the surfaces
belong to the same $SL(2,\R{})$-orbit. Find Veech groups of these
three surfaces. Show that these flat surfaces are Veech surfaces.
Verify that  the  corresponding  Teichm\"uller  disc  is a triple
cover over the modular curve (see Fig.~\ref{zorich:fig:L:3}).
\end{Exercise}

\index{Action on the moduli space!ofSL@of $SL(2,\R{})$|)}
\index{0SL@$SL(2,\R{})$-action on the moduli space|)}

\paragraph{Primitive Veech Surfaces}

It is not difficult to generalize the Exercise above and  to show
that       any
\emph{square-tiled       surface}\index{Surface!square-tiled}\index{Square-tiled surface}
(see Sec.~\ref{zorich:ss:Square:tiled:surfaces})  is  necessarily   a   Veech
surface.

A square-tiled surface is a ramified covering over  a flat torus,
such that  all ramification points project  to the same  point on
the flat torus, which is a Veech surface. One can generalize this
observation.  Having  a  Veech  surface $S$ one can  construct  a
ramified covering  $\tilde  S\to  S$  such  that all ramification
points on $\tilde S$ project to conical singularities on $S$. One
can  check  that any such $\tilde  S$  is a Veech surface.  Thus,
having a Veech surface we  can  construct a whole bunch of  Veech
surfaces in higher genera.

Veech  surfaces  which cannot  be  obtained  from  simpler Veech
surfaces    by    the   covering    construction    are called
\emph{primitive}. For a  long time (and till recent revolution
in genus two, see Sec.~\ref{zorich:ss:Revolution:in:Genus:Two},
the list of known primitive  Veech surfaces was  very short.
Very recently C.~McMullen has found infinitely many Veech
surfaces in genera $3$ and $4$ as well,
see~\cite{zorich:McMullen:prym}. All other known primitive Veech
surfaces  of genus $g>2$ can be obtained by Katok--Zemlyakov
construction (see Sec.~\ref{zorich:ss:Billiards:in:Polygons})
from triangular billiards\index{Billiard!triangular}
of the first three types in the list below:
$$
\begin{array}{lll}
\left(\cfrac{\pi}{n}\ ,\ \cfrac{n-1}{2n}\ \pi\ ,\ \cfrac{n-1}{2n}\
\pi\right), &\text{ for } n\ge 6 &
\text{(discovered by W.~Veech)}\\
\left(\cfrac{\pi}{n}\ ,\ \cfrac{\pi}{n}\ ,\ \cfrac{n-2}{n}\
\pi\right), &\text{ for } n\ge 7  &
\text{(discovered by W.~Veech)}\\
\left(\cfrac{\pi}{n}\ ,\ \cfrac{\pi}{2n}\ ,\ \cfrac{2n-3}{2n}\
\pi\right), &\text{ for }n\ge 4\quad &
\text{(discovered by Ya.~Vorobets})\\
\left(\cfrac{\pi}{3}\ ,\ \cfrac{\pi}{4}\ ,\
\cfrac{5\pi}{12}\right) &&
\text{(discovered by W.~Veech)})\\
\left(\cfrac{\pi}{3}\ ,\ \cfrac{\pi}{5}\ ,\
\cfrac{7\pi}{15}\right) &&
\text{(discovered by Ya.~Vorobets)}\\
\left(\cfrac{2\pi}{9}\ ,\ \cfrac{3\pi}{9}\ ,\
\cfrac{4\pi}{9}\right) && \text{(discovered by R.~Kenyon}\\
[-\halfbls]&&\hfill\text{and J.~Smillie)}\\
\left(\cfrac{\pi}{3}\ ,\ \cfrac{\pi}{12}\ ,\
\cfrac{7\pi}{12}\right) &&
\text{(discovered by W.~P.~Hooper)}\\
\end{array}
$$

--- The flat surface corresponding to the isosceles triangle
with the angles $\pi/n,(n-1)\pi/(2n), (n-1)\pi/(2n)$ belongs to
the hyperelliptic component $\cH^{hyp}(2g-2)$ when $n=2g$ and to
the hyperelliptic component $\cH^{hyp}(g-1,g-1)$ when $n=2g+1$.
The surface can be unwrapped to the regular $2n$-gon with
opposite sides identified by parallel translations, see
Fig.~\ref{zorich:fig:iet:octagon}.

--- The flat surface corresponding to the isosceles triangle
with the angles $\pi/n,\pi/n, (n-2)\pi/(2n)$ belongs to the
hyperelliptic component $\cH^{hyp}(2g-2)$ when $n=2g+1$ and to
the hyperelliptic component $\cH^{hyp}(g-1,g-1)$ when $n=2g+2$.
The surface can be unwrapped to a pair of regular $n$-gons glued
by one side. Each side of one polygon is identified by a
parallel translation with the corresponding side of the other
polygon, see Fig.~\ref{zorich:fig:different:unfoldings}.

--- The flat surface corresponding to the obtuse triangle with the angles
$\pi/n,\pi/(2n), (2n-3)\pi/(2n)$ belongs to one of two
nonhyperelliptic components of the stratum $\cH(2g-2)$ where
$n=g+1$.

--- The flat surface corresponding to the acute triangle
$\pi/3,\pi/4, 5\pi/12$ belongs to the nonhyperelliptic component
$\cH^{odd}(4)$; here $g=3$.

--- The flat surface corresponding to the acute triangle
$\pi/3,\pi/5, 7\pi/15$ belongs to the nonhyperelliptic component
$\cH^{even}(6)$; here $g=4$.

--- The flat surface corresponding to the acute triangle
$2\pi/9,3\pi/9, 4\pi/9$ belongs to the stratum $\cH(3,1)$; here
$g=3$.

--- The flat surface corresponding to the obtuse triangle
$\pi/3,\pi/12,7\pi/12$ belongs to the stratum $\cH(6)$; here
$g=4$ (the information that this is a Veech surface
is taken from~\cite{zorich:McMullen:prym}).

The details on unwrapping of these surfaces and on cylinder
decompositions of some of them can be found in the paper of
Ya.~Vorobets~\cite{zorich:Vorobets:billiards:in:triangles}.

It is proved that unwrapping triangular billiards in other
acute, rectangular or isosceles triangles does not give new
Veech surfaces in genera $g>2$
(see~\cite{zorich:Kenyon:Smillie}, \cite{zorich:Puchta},
\cite{zorich:Vorobets:billiards:in:triangles} and further
references in these papers). For obtuse triangles the question
is open.

We  discuss  genus  $g=2$  separately in the next  section:  very
recently                K.~Calta~\cite{zorich:Calta}                 and
C.~McMullen~\cite{zorich:McMullen:Hilbert} have found a countable family
of primitive Veech surfaces in  the  stratum  $\cH(2)$ and proved
that  the list  is  complete. However, even  in  genus $g=2$  the
situation with  the  stratum $\cH(1,1)$ is drastically different:
using the results of
M.~Moeller~\cite{zorich:Moeller:Galois}--\cite{zorich:Moeller:Periodic:points}
very             recently             C.~McMullen             has
proved~\cite{zorich:McMullen:decagon:proof} the following result.

\begin{NNTheorem}[C.~McMullen]
The only primitive Veech surface in the stratum $\cH(1,1)$ is the
surface  represented  by  the  regular  decagon  with  identified
opposite sides.
\end{NNTheorem}

Thus, it  is not  clear, what one should expect  as a solution of
the following general problem.

\begin{Problem}
Find all primitive Veech surfaces.
\end{Problem}

An algebro-geometric  approach  to  Veech  surfaces  suggested by
M.~M\"oller   in~\cite{zorich:Moeller:Galois}   and~\cite{zorich:Moeller:Hodge}
might help to shed some light on this Problem.

\toread{Veech surfaces}

I          recommend         the survey
paper~\cite{zorich:Hubert:Schmidt:Handbook} of  P.~Hubert and
T.~Schmidt as an introduction to  Veech  surfaces. A canonical
reference for square-tiled  surfaces\index{Surface!square-tiled}\index{Square-tiled surface}
(also   called {\it
arithmetic  Veech surfaces})     is     the     paper of
E.~Gutkin      and C.~Judge~\cite{zorich:Gutkin:Judge}. More
information  about   Veech surfaces   can   be   found in   the
pioneering   paper   of W.~Veech~\cite{zorich:Veech:Eisenstein}
and    in   the paper    of
Ya.~Vorobets~\cite{zorich:Vorobets:billiards:in:triangles}. For
the most recent  results  concerning Veech  groups  and geometry
of  the Teichm\"uller  discs  see  the  original papers of
P.~Hubert  and
T.~Schmidt~\cite{zorich:Hubert:Schmidt:Veech:groups},
\cite{zorich:Hubert:Schmidt:Invariants},
\cite{zorich:Hubert:Schmidt:Infinitely:generated},
\cite{zorich:Hubert:Schmidt:Infinitely:generated:2}, of
C.~McMullen~\cite{zorich:McMullen:inf:generated}, \cite{zorich:McMullen:prym}
and  of P.~Hubert  and
S.~Leli\`evre~\cite{zorich:Hubert:Lelievre:teich;discs},
\cite{zorich:Hubert:Lelievre:Noncongruence}.

\index{Veech!surface|)}
\index{Surface!Veech|)}

\subsection{Kernel Foliation}
\label{zorich:ss:Kernel:Foliation}

In this section we describe some natural holomorphic foliation
on the moduli space of Abelian differentials. In higher genera little is known about this
foliation (though it seems to be worth of study). We use this
foliation in the next section to describe
$GL(2,\R{})$-invariant submanifolds of ``intermediate type''
discovered by K.~Calta and by C.~McMullen in genus two.

We have  seen that any stratum  $\cH(d_1, \dots, d_\noz)$  can be
locally parameterized  by  a  collection  of  basic {\it relative
periods}\index{Period!relative}
of the holomorphic one-form  $\omega$,  or,  in other words,
that  a neighborhood  $\cU([\omega])\subset  H^1(S;\{P_1,  \dots,
P_\noz\};\C{})$ gives a local chart in $\cH(d_1, \dots, d_\noz)$.

Let  $S\in\cH(1,1)$.  Let  closed  paths  $a_1,  a_2,  b_1,  b_2$
represent  a  basis  of  cycles  in  $H_1(S;\Z{})$. Any path  $c$
joining conical singularities  $P_1$  and $P_2$ represents a {\it
relative cycle}  in  $H_1(S,\{P_1,  P_2\};\Z{})$.  Let $A_1, A_2,
B_1,  B_2,  C\in \C{}$  be  the {\it  periods}  of $\omega$:  the
integrals   of    $\omega$    over   $a_1,a_2,   b_1,   b_2,   c$
correspondingly.

\begin{figure}
%
%
\centering
    \includegraphics{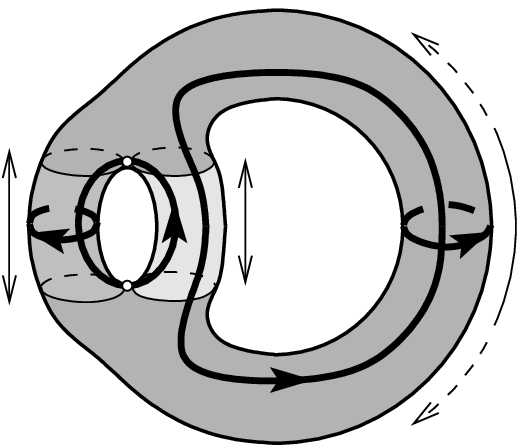}
    \includegraphics{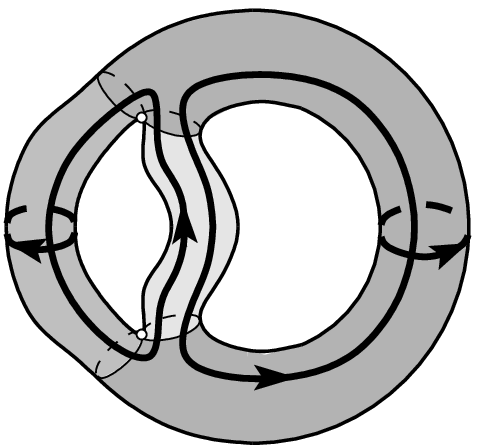}
\begin{picture}(56,0)(56,0) 
 \put(-79,-64){$h_1$}
 \put(0,-64){$h_2$}
 \put(76,-40){$h_3$}
 \put(-52,-89){$w_1$}
 \put(-31,-89){$w_2$}
 \put(-59,-72){$\mathbf{a_1}$}
 \put(-33,-64){$\mathbf{b_1}$}
 \put(32,-62){$\mathbf{a_2}$}
 \put(2,-12){$\mathbf{b_2}$}
 \end{picture}
\vspace{130bp}
\caption{
\label{zorich:fig:periods}
A  deformation of  a flat surface  inside the  kernel foliation
keeps the absolute periods unchanged}
\end{figure}

\begin{Example}
The  collection of  cycles  $a_i, b_i$, $i=1,2$,  on  the
surfaces  from Fig.~\ref{zorich:fig:periods}   represent   a
basis   of   cycles   in $H_1(S;\Z{})$.  All  horizontal
geodesics on  these surfaces are closed; each surface  can  be
decomposed into three cylinders  filled  with  horizontal
geodesics.  Let  $w_1,  w_2, w_3=w_1+w_2$ be the widths
(perimeters) of these cylinders; $h_1, h_2, h_3$  be their
heights; $t_1,  t_2, t_3$ their  twists, see
Fig.~\ref{zorich:fig:1cyl}. It is easy to check that
\begin{align*}
A_1&=\int_{a_1}\omega = -w_1             &B_1 =\int_{b_1}\omega = (t_1-t_2)+\sqrt{-1}(h_2-h_1)\notag\\
A_2&=\int_{a_2}\omega = -(w_1+w_2)\quad &B_2 =\int_{b_2}\omega =
(t_2+t_3)-\sqrt{-1}(h_2+h_3)
\end{align*}
   %
\end{Example}
The
\emph{kernel foliation}\index{Kernel foliation}\index{Foliation!kernel foliation}
in $\cH(1,1)$ is the foliation defined in local coordinates by
equations
$$
\begin{cases}A_1&=const_{11}\\A_2&=const_{21}\end{cases} \qquad
\begin{cases}B_1 =const_{12}\\B_2 =const_{22}\end{cases}
$$
In other  words, this is a foliation which  is obtained by fixing
\emph{all} absolute  periods  and  changing  the  relative period
$C=\int_{P_1}^{P_2}\omega$. Similarly,  the \emph{kernel foliation}
in arbitrary stratum $\cH(d_1, \dots,  d_\noz)$  is  a  foliation
which in  cohomological  coordinates  is  represented by parallel
complex  $(\noz-1)$-dimensional   affine  subspaces obtained   by
changing all relative periods while fixing the absolute ones.

Passing to  a finite cover  over $\cH(d_1, \dots, d_\noz)$ we can
assume that  all  zeroes  $P_1,  \dots,  P_\noz$ are \emph{named}
(i.e.  having two  zeroes  $P_j, P_k$ of  same  degrees, we  know
exactly which of the two is $P_j$ and which is $P_k$). Now we can
fix  an  arbitrary subcollection of zeroes and  define  a  kernel
``subfoliation'' along relative periods  corresponding  to chosen
subcollection.

Recall that the area of a flat surface is expressed in terms of
the absolute periods (see Riemann bilinear relation in
Table~\ref{zorich:tab:dictionary:geometry:to:complex:analysis} in
Sec.~\ref{zorich:ss:Dictionary:of:Complex:Analytic:Language}).
Thus, moving along leaves of kernel foliation we do not change the
area of the surface. In particular, we can consider the kernel
foliation as a foliation of the ``unit hyperboloid''
$\cH_1(d_1, \dots, d_\noz)$.

\paragraph{Exercises on Kernel Foliation}

\begin{Exercise}
To deform the  flat surface on the left of Fig.~\ref{zorich:fig:periods}
along the kernel  foliation we have to keep all $A_1,A_2,B_1,B_2$
unchanged. Hence, we cannot change the widths (perimeters) of the
cylinders, since they are expressed in terms of  $A_1$ and $A_2$.
Increasing the height  of the second cylinder by $\varepsilon$ we
have to  \emph{increase} the height  of the first cylinder by the
same         amount         $\varepsilon$         to         keep
$B_1=(t_1-t_2)+\sqrt{-1}\left((h_2+\varepsilon)-(h_1+\varepsilon)
\right)$  unchanged;  we also have to \emph{decrease} the  height
$h_3$  of  the third cylinder by $\varepsilon$  to  preserve  the
value of $B_2$.  Similarly,  \emph{increasing} the twist $t_1$ by
$\delta$ we  have to \emph{increase}  the twist $t_2$ by the same
amount $\delta$ and to \emph{decrease} $t_3$ by $\delta$.
\end{Exercise}

\begin{Exercise}
It is  convenient to consider the  kernel foliation in  the total
moduli  space   $\cH_g$\index{0H10@$\cH_g$ -- moduli space of holomorphic 1-forms}\index{Moduli space!of holomorphic 1-forms}   of   all  holomorphic  1-forms  without
subdivision of  $\cH_g$\index{0H10@$\cH_g$ -- moduli space of holomorphic 1-forms}\index{Moduli space!of holomorphic 1-forms}  into  strata  $\cH(d_1, \dots, d_\noz)$,
where $\sum_j d_j = 2g-2$. In particular, to deform a  surface
$S\in\cH(2)\subset \cH_2$ along  the kernel foliation  we have to break the double
zero into two  simple zeroes preserving the absolute periods. The
corresponding        surgery        is        presented        at
Fig.~\ref{zorich:fig:breaking:up:a:zero}.

The leaves  of the kernel  foliation are naturally endowed with a
flat structure, which  has conical singularities at the points of
intersection  of  the  leaf   with  the  smaller strata and  with
degenerate strata.

Assuming that the zeroes  $P_1,  P_2$ of a surface $S\in\cH(1,1)$
are  \emph{named}  show  that  the  intersection  of  the  kernel
foliation with  the  stratum  $\cH(2)$  corresponds  to a conical
point with the cone angle  $6\pi$,  while  the intersections with
the two strata of degenerate flat surfaces (determine which ones)
are just the regular points of the flat structure.
\end{Exercise}

In the exercise below we use a polygonal representation of a
flat  surface  (compare to Fig. 8 in  the  paper~\cite{zorich:Calta}  of
K.~Calta).

\begin{Exercise}
Consider a regular decagon. Imagine that there are springs inside
its  sides  so that we can  shrink  or  expand  the  sides keeping them
straight segments. Imagine that we hammer a nail in the center of
each  side.  Though the centers of  the  sides are now fixed  our
decagon is still  flexible:  we can pull a  vertex  and the whole frame
will follow, see Fig.~\ref{zorich:fig:kernel:foliation}b. We assume that
under any such deformation each nail stays exactly  in the middle
of the corresponding side.

\begin{figure}
%
%
\centering
\includegraphics{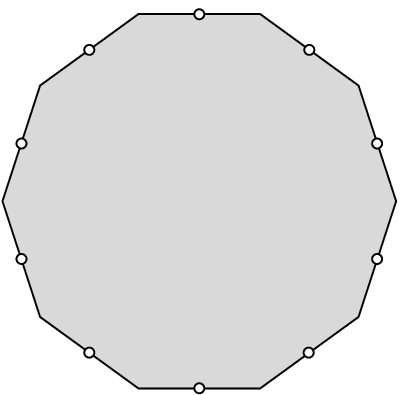}
\includegraphics{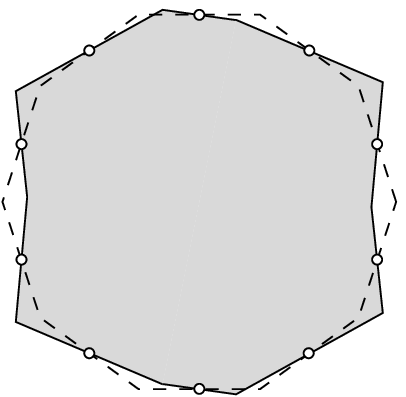}
\includegraphics{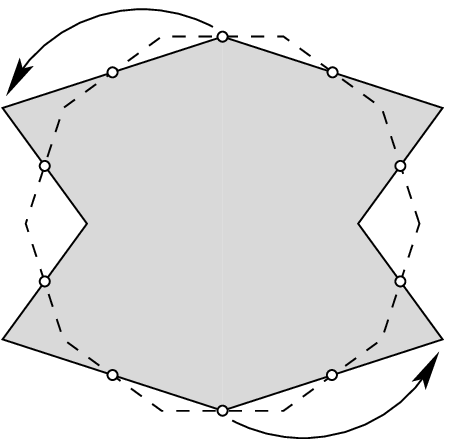}
\includegraphics{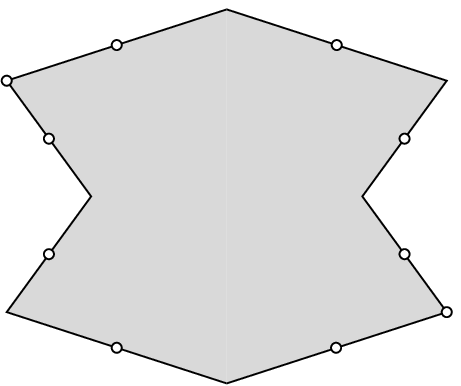}
\includegraphics{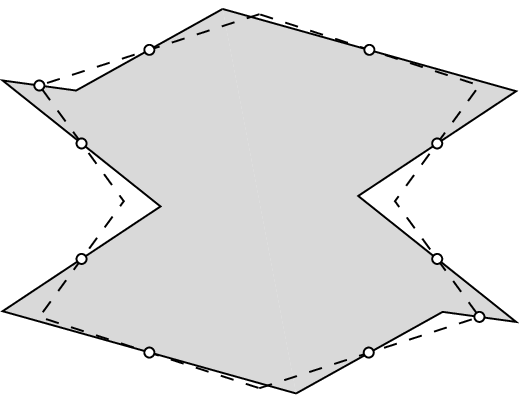}
\includegraphics{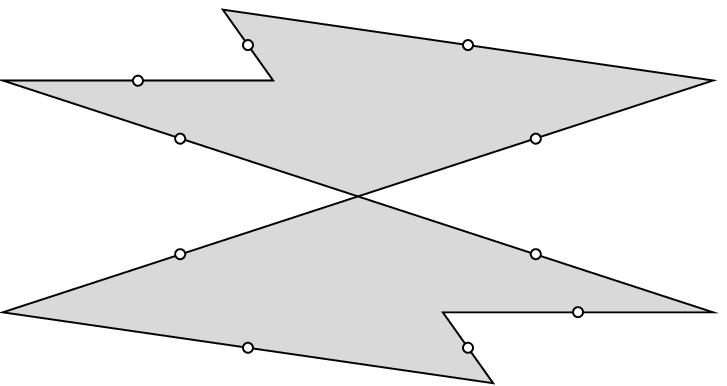}
\begin{picture}(25,0)(25,0)
  \put(-125,-47){a}
  \put(165,-47){b}
  \put(-125,-151){c}
  \put(165,-151){d}
  \put(-125,-264){e}
  \put(165,-264){f}
\end{picture}
\vspace{310bp} 
\caption{
\label{zorich:fig:kernel:foliation}
This cartoon movie represents a path living inside a leaf of the kernel
foliation. When  at the stage c it  lands to  a surface from  the
stratum $\cH(2)$  we remove the nails  from one pair  of vertices
and hammer  them in the other  symmetric pair of  vertices (stage
d). Then we continue the trip inside the kernel foliation}
\end{figure}

Prove that the deformed polygon is again centrally symmetric with
the same center of symmetry. Prove that the opposite sides of the
deformed polygon are parallel  and  have equal length. Prove that
the resulting  flat surface lives in the same  leaf of the kernel
foliation. Show that the ``nails'' (the centers of the sides) and
the  center  of  symmetry  are  the  Weierstrass  points  of  the
corresponding Riemann  surface (fixed points of the hyperelliptic
involution).

Move a vertex towards the  center  of the side which this  vertex
bounds.  The  side  becomes  short   (see   the   upper  side  on
Fig.~\ref{zorich:fig:kernel:foliation}b)  and  finally  contracts  to  a
point (Fig.~\ref{zorich:fig:kernel:foliation}c). What flat surface do we
get?

Show      that       the      deformation      presented at
Fig.~\ref{zorich:fig:kernel:foliation}c brings  us to a surface $S_0$
in the stratum $\cH(2)$. Remove the pair of clues, which is
hammered at the vertices of the resulting octagon. Hammer them
to another pair         of         symmetric vertices (see
Fig.~\ref{zorich:fig:kernel:foliation}d). We can declare that we have a
new  centrally-symmetric  decagon with a pair of  sides of  zero
lengths. Stretching  this  pair of  sides and making them have
positive length (see Fig.~\ref{zorich:fig:kernel:foliation}e) we
continue our  trip inside the kernel  foliation. Show that for a
given small value $\vec  v\in\C{}$ of the saddle connection
joining two conical singularities, there are exactly three
different surfaces obtained as a  small deformation of the
surface $S_0\in\cH(2)$ as on Fig.~\ref{zorich:fig:kernel:foliation}c
along the leaf of the kernel foliation  and  having  a  saddle
connection $\vec v$. (Use previous Exercise and
Fig.~\ref{zorich:fig:breaking:up:a:zero}.)

Following our path in the kernel foliation we get to the surface
$S_1$ as on Fig.~\ref{zorich:fig:kernel:foliation}f. Is this surface
singular? To what stratum (singular stratum) it belongs?
\end{Exercise}

\paragraph{Compact Leaves of the Kernel Foliation}

Let us show that  the leaf $\cK$ of  the kernel foliation
passing through a square-tiled surface\index{Surface!square-tiled}\index{Square-tiled surface}
$S(\omega_0)$ of genus
two is a compact square-tiled surface. To simplify notations
suppose that the absolute periods of the holomorphic 1-form
$\omega_0$, representing the flat surface $S(\omega_0)$,
generate the entire integer lattice $\Z{} + \sqrt{-1}\,\Z{}$.

Consider the ``relative period''\index{Period!relative} map $p: \cK \to \torus$ from
the corresponding leaf $\cK$ of the kernel foliation containing
$S(\omega_0)$ to the torus. The map $p$ associates to a flat
surface $S(\omega)\in\cK$ the relative period $C$ taken modulo
integers,
$$
\cK\ni S(\omega) \stackrel{p}{\mapsto} C =
\int_{P_1}^{P_2}\omega \mod \Z{} + \sqrt{-1}\,\Z{}\quad \in\quad
\C{}/(\Z{} + \sqrt{-1}\,\Z{})=\torus
$$
Since the flat surface $S(\omega)$ belongs to the same leaf
$\cK$ as the square-tiled surface $S(\omega_0)\in\cK$, the
absolute periods of $\omega$ are the same as the ones of
$\omega_0$, and hence the integral above taken modulo integers
does not depend on the path on $S(\omega)$ joining $P_1$ and
$P_2$. It is easy to check that the map $p$ is a finite ramified
covering over the torus $\torus$, and thus the leaf $\cK$ is a
square-tiled surface.

Those flat surfaces $S(\omega)\in\cK$, which have integer
relative period $C\in \Z{} + \sqrt{-1}\,\Z{}$, have \emph{all}
periods in $\Z{} + \sqrt{-1}\,\Z{}$. Hence, these flat surfaces
are square-tiled. Since $S(\omega)$ and $S(\omega_0)$ have the
same area, the number $N$ of squares tiling $S(\omega)$ and
$S(\omega_0)$ is the same. Thus, $\cK$ has a structure of a
square-tiled surface such that the vertices of the tiling are
represented by $N$-square-tiled surfaces $S(\omega)\in \cK$.

To discuss the geometry of $\cK$ we need to agree about
enumeration of zeroes $P_1, P_2$ of a surface $S\in\cH(1,1)$. We
choose the convention where the zeroes are \emph{named}. That
is, given two zeroes of order $1$ we know which of them is $P_1$
and which of them is $P_2$. Under this convention the
square-tiled surface $\cK$ is a translation surface; it is
represented by a \emph{holomorphic one-form}. (Accepting the
other convention we would obtain the quotient of $\cK$ over the
natural involution exchanging the names of the zeroes. In this
latter case the zeroes of $S\in\cH(1,1)$ are not
distinguishable; the leaf of the kernel foliation gives a flat
surface represented by a \emph{quadratic differential}.)

The lattice points of the square-tiled surface $\cK$ are
represented by $N$-square-tiled surfaces $S\in\cK$ of several
types. We have the lattice points represented by $N$-square
tiled surfaces from $\cH(1,1)$. These points are the regular
points of the flat metric on $\cK$.

There are points of intersection of $\cK$ with $\cH(2)$. Such
point $S(\omega)\in\cK\cap\cH(2)$ is always represented by a
square-tiled surface, and hence gives a vertex of the tiling of
$\cK$. We have seen (see
Fig.~\ref{zorich:fig:breaking:up:a:zero}) that given a surface
$S(\omega)\in\cH(2)$ and a small complex period $C$ one can
construct three different flat surfaces $S(\omega_1),
S(\omega_2), S(\omega_3)\in\cH(1,1)$ with the same absolute
periods as $\omega$ and with the relative period $C$. Thus, the
points $S(\omega)\in\cK\cap\cH(2)$ correspond to conical points
of $\cK$ with the cone angles $3\cdot 2\pi$ when the zeroes are
named (and with the cone angle $3\pi$, when they are not named).
The total number of such points was computed in the paper of
A.~Eskin, H.~Masur and
M.~Schmol~\cite{zorich:Eskin:Masur:Schmoll}; it equals\index{Conical!singularity}
\begin{equation}
\label{zorich:number:of:zeroes:on:K} \text{number of conical
points on } \cK =
\cfrac{3}{8}\,(N-2)N^2\prod_{p\,|N}\left(1-\cfrac{1}{p^2}\right)
\end{equation}

There remain vertices of the tiling of $\cK$ represented by
degenerate square-tiled surfaces $S(\omega)$. It is not
difficult to show that these points are regular for the flat
metric on $\cK(\omega)$ when the zeroes are named. (They
correspond to conical singularities with the cone angle $\pi$,
when the zeroes are not named). The degenerate $N$-square-tiled
surfaces $S(\omega)$ are of two types. It might be an
$N$-square-tiled torus with two points of the tiling identified.
It might be a pair of square-tiled tori with a vertex of the
tiling on one torus identified with a vertex of the tiling on
the other torus. Here the total number of squares used to tile
these two tori is $N$. The total number of the vertices of the
tiling of  $\cK$ of this type is computed in the
paper~\cite{zorich:Schmoll:kernel:foliation}; it equals
$$
\text{number of special points on } \cK =
\cfrac{1}{24}\,(5N+6)N^2\prod_{p\,|N}\left(1-\cfrac{1}{p^2}\right)
$$

Summarizing we conclude that the translation surface $\cK$ lives
in the stratum $\cH(\underbrace{2,\dots,2}_{k})$, where the
number $k$ of conical points is given by
formula~\eqref{zorich:number:of:zeroes:on:K}.

We complete this section with an interpretation of a compact
leaf $\cK$ as a space of torus coverings; this interpretation
was introduced by A.~Eskin, H.~Masur and M.~Schmol in
~\cite{zorich:Eskin:Masur:Schmoll} and developed by M.~Schmoll
in~\cite{zorich:Schmoll:kernel:foliation},
\cite{zorich:Schmoll:kernel:foliation:2}.

We have seen that a nondegenerate flat surface
$S(\omega_0)\in\cH(1,1)$ representing a vertex of the square
tiling of $\cK$ is an $N$-square-tiled surface. Hence,
$S(\omega)$ is a ramified covering over the standard torus
$\torus$ of degree $N$ having two simple ramification points,
which project to the same point of the torus. A non vertex point
$S(\omega)\in\cK$ is also a ramified covering over the standard
torus $\torus$. To see this consider once more the period map,
but this time applied to $S(\omega)$:
$$
S(\omega)\ni P\stackrel{proj}{\mapsto} \int_{P_1}^{P}\omega \mod
\Z{} + \sqrt{-1}\,\Z{}\quad \in\quad \C{}/(\Z{} +
\sqrt{-1}\,\Z{})=\torus.
$$
Since all absolute periods of $\omega$ live in $\Z{} +
\sqrt{-1}\,\Z{}$ the integral taken modulo integers does not
depend on the path joining the marked point (conical
singularity) $P_1$ with a point $P$ of the flat surface
$S(\omega)$. The map $proj$ is a ramified covering.

The degree of the covering can be computed as the ratio of areas
of $S(\omega)$ and of the torus $\torus$, which gives $N$. The
covering has precisely two simple ramification points, which are
the conical points $P_1, P_2$ of $S(\omega)$. This time they
project to two different points of the torus.

Recall that by convention we assume that the absolute periods of
$\omega$ generate the entire lattice $\Z{} + \sqrt{-1}\,\Z{}$.
They corresponds to primitive covers: the ones which do not
quotient through a larger torus.

\begin{NNProposition}[M.~Schmoll]
Consider primitive branched covers over the standard torus
$\torus$. Fix the degree $N$ of the cover. Let the cover have
exactly two simple branch points.

\index{Action on the moduli space!ofSL@of $SL(2,\R{})$}
\index{0SL@$SL(2,\R{})$-action on the moduli space}

The space of such covers is connected; its natural
compactification coincides with the corresponding leaf $\cK$ of
the kernel foliation. The Veech group of the square-tiled\index{Surface!square-tiled}\index{Square-tiled surface}
surface $\cK$ coincides with $SL(2,\Z{})$.
\end{NNProposition}

Connectedness of the space of covers is not quite obvious
(actually, it was proved earlier in other terms by
W.~Fulton~\cite{zorich:Fulton}). An observation above shows,
that the leaf $\cK$ coincides with a connected component of the
space of covers. Thus, connectedness of the space of covers
implies that this space coincides with $\cK$. The group
$SL(2,\Z{})$ acts naturally on the space of covers; in
particular it maps the space of covers to itself. This implies
that $SL(2,\Z{})$ belongs to the Veech group of the square-tiled
surface $\cK$. It is easy to show, that it actually coincides
with $SL(2,\Z{})$.

\begin{NNCorollary}[M.~Schmoll]
Consider square-tiled surfaces $S(\omega)$ of genus two such
that the absolute periods of $\omega$ span the entire integer
lattice $\Z{} + \sqrt{-1}\,\Z{}$. For any given $N>3$ all such
$N$-square tiled surfaces belong to the same compact connected
leaf $\cK(N)$ of the kernel foliation.
\end{NNCorollary}

For more information on kernel foliation of  square-tiled
surfaces in genus two see  the papers of A.~Eskin, H.~Masur and
M.~Schmol~\cite{zorich:Eskin:Masur:Schmoll} and of
M.~Schmoll~\cite{zorich:Schmoll:kernel:foliation},
\cite{zorich:Schmoll:kernel:foliation:2}. In particular, the
latter papers propose a beautiful formula for Siegel--Veech
constants of any flat surface $S\in\cK(N)$ in terms of geometry
of the cylinder decomposition of the square-tiled surface
$\cK(N)$.

\subsection{Revolution in Genus Two (after K.~Calta and C.~McMullen)}
\label{zorich:ss:Revolution:in:Genus:Two}

In  this  section  we  give  an   informal   survey   of   recent
revolutionary results in genus $g=2$ due to K.~Calta~\cite{zorich:Calta}
and to C.~McMullen~\cite{zorich:McMullen:Hilbert}.

\index{Veech!surface|(}
\index{Surface!Veech|(}

Using different  methods  they  found  a  countable collection of
primitive Veech surfaces in the  stratum  $\cH(2)$,  proved  that
this collection describes all Veech surfaces,  and gave efficient
algorithms  which   recognize  and  classify  Veech  surfaces  in
$\cH(2)$.

This  result  is  in  a  sharp  contrast  with   the  Theorem  of
C.~McMullen~\cite{zorich:McMullen:decagon:proof}   cited   above,  which
tells that  in the other  stratum $\cH(1,1)$ in genus $g=2$ there
is {\it only one} primitive Veech surface.

This discovery of  an infinite family of primitive Veech surfaces
in the stratum $\cH(2)$ is also in a sharp contrast with our poor
knowledge of primitive Veech surfaces in  higher  genera:  as  we
have seen in the previous  section,  primitive  Veech surfaces in
higher genera  $g\ge 3$ are  currently known only in some special
strata (mostly hyperelliptic),  and  even in these special strata
we   know   only  finite  number  of  primitive  Veech   surfaces
(basically, only one).

Another  remarkable  result  is  a discovery by K.~Calta  and  by
C.~McMullen of  nontrivial  examples of invariant submanifolds of
intermediate dimension:  larger  than  closed  orbits and smaller
than the entire stratum.

One more revolutionary result in  genus  two  is a Classification
Theorem due  to C.~McMullen~\cite{zorich:McMullen:genus:2} which  proves
that  a  closure  of  {\it  any}  $GL^+(2,\R{})$-orbit  is  a  nice
complex-analytic  variety  which is either an entire stratum,  or
which has one of the types mentioned above.

\subparagraph{Algebro-geometric Approach}
To avoid overloading of this survey I had  to sacrifice beautiful
algebro-geometric part  of  this  story developed by C.~McMullen;
the       reader       is       addressed       to       original
papers~\cite{zorich:McMullen:Hilbert}--\cite{zorich:McMullen:decagon:proof} and
to a short overview presented in~\cite{zorich:Hubert:Schmidt:Handbook}.

\paragraph{Periods of Veech Surfaces in Genus $g=2$}

\index{Action on the moduli space!ofGL@of $GL^+(2,\R{})$|(}
\index{0GL@$GL^+(2,\R{})$-action on the moduli space|(}

If  $S$  is  a Veech surface then the flat surface $gS$ is also a
Veech surface for any $g\in GL^+(2,\R{})$. Thus, speaking about a
finite  or  about a countable collection of  Veech  surfaces  we,
actually, choose some family of representatives  $\{S_k\}$ of the
orbits $GL^+(2,\R{})\cdot S$ of Veech surfaces.

The question, which  elements  of our collection $\{S_k\}$ belong
to  the  same  $GL^+(2,\R{})$-orbit  and  which  ones  belong  to
different orbits  is a matter of  a separate nontrivial  study. A
solution was  found by C.~McMullen in~\cite{zorich:McMullen:spin}; it is
briefly presented in the next Sec.~\ref{zorich:ss:Teichmuller:Discs}. In
this section we present  effective  algorithm due to K.~Calta and
to C.~McMullen which enables to  determine  whether  a given flat
surface in $\cH(2)$ is a Veech surface.

Following K.~Calta we  say  that a flat surface  $S$  can be {\it
rescaled} to  a flat  surface $S'$ if $S$ and  $S'$ belong to the
same $GL^+(2,\R{})$-orbit. We say that a flat surface $S$ is {\it
quadratic} if for any homology cycle $c\in H_1(S;\Z{})$ we have
$$
\int_c\omega = (p+q\sqrt{d})+i(r+s\sqrt{d}), \quad\text{where }
d\in\N,\quad p,q,r,s\in\Q
$$
In other words,  we say  that a  flat  surface $S$  defined by  a
holomorphic 1-form $\omega$  is {\it quadratic} if all periods of
$\omega$ live in $\Q(\sqrt{d})+i\Q(\sqrt{d})$.

We can considerably  restrict the area  of our search  using  the
following Lemma of W.~Thurston.

\begin{NNLemma}[W.~Thurston]
Any Veech surface in genus $g=2$ (no matter primitive or  not, in
the  stratum  $\cH(2)$  or  $\cH(1,1)$)  can  be  rescaled  to  a
quadratic surface.
\end{NNLemma}

Using  this   Lemma  K.~Calta  suggest  the  following  algorithm
deciding whether a given flat surface  $S\in\cH(2)$  is  a  Veech
surface or not.

\paragraph{Algorithm of Calta}

Recall, that if $S$ is a Veech surface (of arbitrary genus), then
by              Veech              alternative               (see
Sec.~\ref{zorich:ss:Implementation:of:General:Philosophy}) a directional
flow in any direction  is  either minimal or completely periodic, that is
a presence of a closed geodesic going in some direction implies that all geodesics
going in this direction are periodic.
Moreover, it was proved by Veech that as soon as there a saddle connection
going in some direction, this direction is also completely periodic.
In both cases the surface
decomposes into a  finite  collection of cylinders; each boundary
component  of  each cylinder contains a conical singularity  (see
Sec.~\ref{zorich:ss:Square:tiled:surfaces}).

The algorithm works as follows. Having a flat surface $S\in\cH(2)$
it is  easy to  find {\it some} closed geodesic  on $S$ (which is
allowed  to  be  a  closed  geodesic  saddle  connection).  Since
``rescaling''  the   surface   $S$   (i.e.   applying   a  linear
transformation from $GL^+(2,\R{})$) we replace a Veech surface by
a Veech surface, we  can turn $S$ in such way that  the direction
of  the  closed geodesic will become horizontal.  We  denote  the
resulting surface by the same symbol $S$.

Since $S$ lives  in  $\cH(2)$ it has a  single  conical point $P$
with  the  cone angle $6\pi$. In particular,  there  are  exactly
three  geodesics  leaving  the  conical  point  in  the  positive
horizontal direction (to  the  East). If  at  least one of  these
three horizontal geodesics does not come back to  $P$ the surface
$S$ is {\it not} a Veech surface. Otherwise our test continues.

As we  have seen in Sec.~\ref{zorich:ss:Square:tiled:surfaces} there are
two possible ways in which  three  horizontal  geodesics  emitted
from  $P$  to  the  East  can return  to  $P$.  Either  all three
geodesics return  at the angle $3\pi$, or one  of them returns at
the angle $3\pi$ and two  others  return at the angle $\pi$,  see
Fig.~\ref{zorich:fig:separatrix:diagrams}. In both  cases all horizontal
geodesics are closed. In the  first  case  the surface decomposes
into a  single cylinder; in the second case  the surface is glued
from two cylinders, see Sec.~\ref{zorich:ss:Square:tiled:surfaces}.

If  the  surface is  decomposed  into a  single  cylinder, it  is
sufficient  to   compare   the  lengths  $p_1,p_2,p_3$  of  three
horizontal           saddle           connections,            see
Fig.~\ref{zorich:fig:separatrix:diagrams}.  The  flat  surface $S$ is  a
Veech surface if and only  if  $p_1,p_2,p_3$  are  commensurable.
Moreover, if $p_1,p_2,p_3$ are commensurable, we  can rescale $S$
to a  square-tiled\index{Surface!square-tiled}\index{Square-tiled surface}
surface. It  can be done in several elementary
steps. First we rescale $S$  in  the  horizontal direction making
$p_1,p_2,p_3$ rational and then  integer.  Then we rescale $S$ in
the  vertical  direction making the height $h$  of  the  cylinder
integer.  Finally,  we  apply  an  appropriate  parabolic  linear
transformation $\begin{pmatrix} 1  &  s\\ 0 & 1\end{pmatrix}$. It
does  not  change   neither   $p_1,p_2,p_3$  nor  $h$,  but  when
$s\in\R{}$ varies the twist $t$  also  varies  continuously  (see
Fig.~\ref{zorich:fig:1cyl}) taking  all values in $\R{}$; in particular,
we can achieve $t=0$. We get a square-tiled surface.

Consider the second case, when  the  surface  decomposes into two
cylinders.   Let   $w_1,w_2,h_1,h_2,t_1,t_2$    be   the   widths
(perimeters),   heights    and    twists    of   this   cylinders
correspondingly (compare to Sec.~\ref{zorich:ss:Square:tiled:surfaces}).
In a  complete analogy with  the one-cylinder case we can rescale
the surface horizontally,  then  vertically, and finally apply an
appropriate parabolic linear transformation in order  to make the
width (perimeter) $w_1$ and height  $h_1$  of  the first cylinder
equal to  one, $h_1=w_1=1$, and  the twist $t_1=0$ equal to zero.
Applying an appropriate Dehn twist to the second  cylinder we can
assure $0\le t_2 <w_2$.

If  after  our  rescaling  all   parameters   $w_2,   h_2,   t_2$
characterizing  the  second  cylinder  do  not  get  to  the same
quadratic field $\Q(\sqrt{d})$ for some $d\in\N$, the surface $S$
is not a Veech surface.

If $w_2,h_2,t_2$ are rational (i.e. if $d$ is a complete square),
the surface $S$ can be rescaled to a square-tiled surface.

The remaining case is  treated by one of the key Theorems  in the
paper of K.~Calta~\cite{zorich:Calta}.

\begin{NNTheorem}[K.~Calta]
Let all  parameters  $w_j,h_j,t_j$,  $j=1,2$  of  a  two-cylinder
decomposition of a flat  surface  $S\in\cH(2)$ belong to the same
quadratic field  $\Q(\sqrt{d})$  with  $d\in\N$  not  a  complete
square. Then $S$ is a Veech surface if and only if the parameters
satisfy the following system of equations:
\begin{equation}
\label{zorich:eq:Calta:equations}
\begin{cases}
w_1 \bar h_1 &=-w_2 \bar h_2, \\
\bar w_1 t_1 + \bar w_2 t_2  &=
w_1\bar t_1 + w_2\bar t_2,
\end{cases}
\end{equation}
(where         the         bar        denotes         conjugation
$\overline{p+q\sqrt{d}}=p-q\sqrt{d}$   in   $\Q(\sqrt{d})$   with
$p,q\in\Q$).
\end{NNTheorem}

Actually, we  kept the system of equation above  as it is written
in  the  original  paper~\cite{zorich:Calta}.  In  this  form it can  be
adopted  to  a more general normalization of  parameters:  it  is
sufficient to  rescale surface $S$  to bring all $w_j, h_j, t_j,\
j=1,2$ to a quadratic field.

\begin{NNRemark}
Similar {\it necessary} conditions for Veech surfaces in $\cH(2)$
were  obtained  by D.~Panov  independently  of  K.~Calta  and  of
C.~McMullen.
\end{NNRemark}

Since by  Lemma of Thurston any Veech surface  in $\cH(2)$ can be
rescaled  to  a  quadratic  surface, taking a collection  of  all
quadratic  surfaces  decomposed  into  two  horizontal  cylinders
satisfying the condition above, we  get  representatives  of  the
$GL^+(2,\R{})$-orbits of {\it all} flat surfaces in $\cH(2)$.

\begin{figure}
%
%
\centering
    \includegraphics{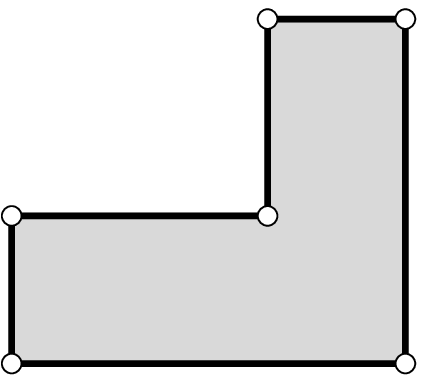}
    \includegraphics{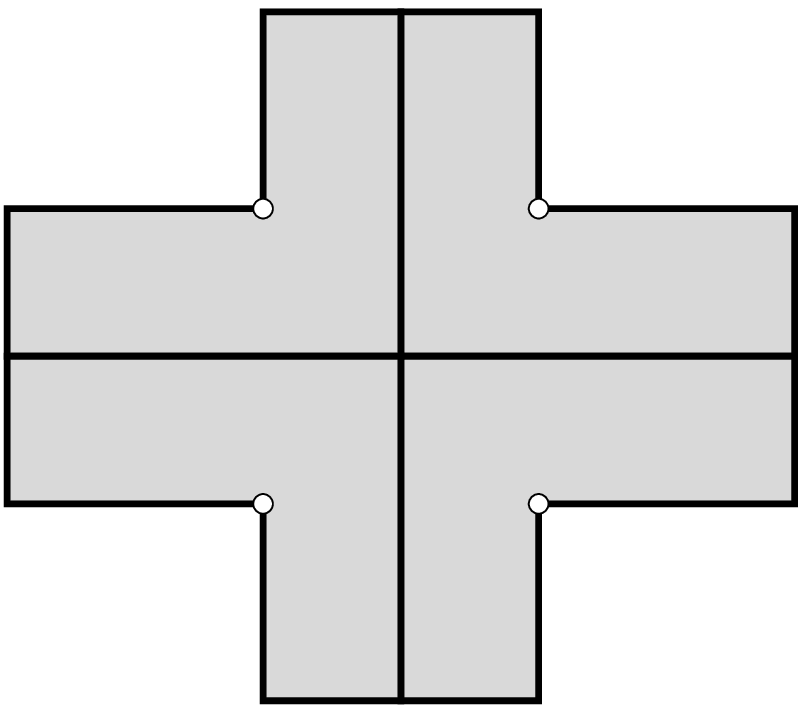}
\begin{picture}(0,3)(0,3)
 \put(-62,-28){$a$}
 \put(-150,-50){$1$}
 \put(-80,10){$1$}
 \put(-110,-71){$b$}
\end{picture}
\vspace{140bp}
\caption{
\label{zorich:fig:McM:billiard}
L-shaped billiard\index{Billiard!table!L-shaped}
table $P(a,b)$ and its unfolding into a flat surface $S\in\cH(2)$
(after C.~McMullen~\cite{zorich:McMullen:Hilbert})
}
\end{figure}

\begin{Exercise}
Show  that   the  Katok--Zemlyakov  construction  applied  to  an
$L$-shaped billiard\index{Billiard!table!L-shaped}
as on  Fig.~\ref{zorich:fig:McM:billiard}  (see also
Fig.~\ref{zorich:fig:L:shaped:billiard})     generates     a     surface
$S\in\cH(2)$.  Show  that  this  surface is decomposed  into  two
cylinders  filled  by closed horizontal geodesics and that  these
cylinders have parameters $w_1=2, h_1=2a-2, t_1=0$  for the first
cylinder and  $w_2=2b,  h_2=2,  t_2=0$  for  the second cylinder.
Using  the  condition  above  prove  the   following  Theorem  of
C.~McMullen~\cite{zorich:McMullen:Hilbert}:
\end{Exercise}

\begin{NNTheorem}[C.~McMullen]
The    $L$-shaped     billiard     table     $P(a,b)$    as    on
Fig.~\ref{zorich:fig:McM:billiard} generates a Veech surface if and only
if $a$ and $b$ are rational or
$$
a=x+z\sqrt{d}\ \text{ and }\ b=y+z\sqrt{d}
$$
for some $x,y,z\in\Q$ with $x+y=1$ and $d\ge 0$ in $\Z{}$.
\end{NNTheorem}

\paragraph{Kernel Foliation in Genus 2}

The following  elementary  observation  explains  our interest to
kernel   foliation   in    the    content   of   our   study   of
$GL^+(2,\R{})$-invariant subvarieties. Let  $\cN\subset  \cH(d_1,
\dots,  d_n)$  be  a  $GL^+(2,\R{})$-invariant  submanifold.  The
``germ'' of the kernel foliation at  $\cN$  is  equivariant  with
respect to the action of $GL^+(2,\R{})$.

In other words this statement can be described as follows. Denote
by $t(\delta)$ a  translation  by $\delta$ along kernel foliation
defined in a neighborhood of $S\in\cH(1,1)$. Here $\delta\in\C{}$
is  a  small  parameter.  Let  $g\in  GL^+(2,\R{})$  be  close to
identity. Then
$$
g\circ t(\delta)\cdot S= t(g\delta)\circ g S
$$
where $g$  acts on a complex number  $\delta$ as  on a vector  in
$\R{2}$. Moving along  the kernel foliation, and then applying an
element $g$ of the  group is the same as applying first  the same
element  of  the  group  and  then  moving  along  an appropriate
translation along  the  kernel  foliation. A similar construction
works in a general stratum.

We get the following very tempting picture. Suppose, that we have
a closed $GL^+(2,\R{})$-orbit $\cN$ in $\cH(1,1)$ or in $\cH(2)$.
For example,  suppose that $\cN$ is an orbit  of a Veech surface.
Consider  a  union  of  leaves  of  the kernel foliation  passing
through   $\cN$.   Due   to   the   remark   above   this   is  a
$GL^+(2,\R{})$-invariant subset!

The  weakness of  this  optimistic picture is  that  there is  no
\emph{a  priori}  reason  to  hope that the  resulting  invariant
subset in  $\cH(1,1)$ would be  closed. And here the magic comes.
The   following   statement   was    proved    independently   by
K.~Calta~\cite{zorich:Calta} and by C.~McMullen~\cite{zorich:McMullen:Hilbert}.

\begin{NNTheorem}[K.~Calta; C.~McMullen]
For {\it any} Veech surface $S_0\in\cH(2)$ the union of leaves
of the kernel  foliation passing through the
$GL^+(2,\R{})$-orbit of $S_0$ is a closed
$GL^+(2,\R{})$-invariant complex orbifold  $\cN$  of  complex
dimension $3$.
\end{NNTheorem}

In particular, the complex dimension of the resulting orbifold is
the sum of the  complex  dimension of the $GL^+(2,\R{})$-orbit of
$S_0$  and  of  the  complex  dimension  of the kernel  foliation
$3=2+1$.

K.~Calta   has   found   the   following   beautiful    geometric
characterization  of   flat   surfaces  living  in  an  invariant
subvariety  $\cN$  as  above.  Surfaces  in  any  such  $\cN$ are
\emph{completely periodic}: as soon  as  there is a single closed
trajectory in some direction, \emph{all} geodesics  going in this
direction are closed. This condition is  necessary and sufficient
condition for a surface to live in an  invariant subvariety $\cN$
as above. In particular, this shows that not  only Veech surfaces
have this property.

An  algorithm   analogous  to  the  algorithm  determining  Veech
surfaces in $\cH(2)$  (see above) allows to K.~Calta to determine
whether a given surface $S\in\cH(1,1)$ is  completely periodic or
not (and hence, whether it belongs  to  an  invariant  subvariety
$\cN$  as  above or  not).  As  before  one  starts  with finding
\emph{some}  closed   geodesic or \emph{some} saddle connection.
If  the  surface  is  completely
periodic in the corresponding direction,
then all other geodesics  going  in  this direction are
periodic and  the  surface  decomposes  into  cylinders. After an
appropriate rotation this periodic direction becomes  horizontal.
Without loss of generality, we may  assume  that  $S$  decomposes
into three cylinders. By $w_i$, $h_i$ and $t_i$ with $1 \le i \le
3$, we denote  the widths, heights and twists. After renumbering,
we may  assume  that $w_3 = w_1 + w_2$. Define $s_1 = h_1 + h_3$,
$s_2 = h_2 + h_3$, $\tau_1 = t_1 + t_3$, $\tau_2 = t_2 + t_3$. If
the surface is completely periodic  its  absolute  periods can be
rescaled to get to $\Q(\sqrt{d})+i \Q(\sqrt{d})$  (compare to the
algorithm for Veech  surfaces).  Leaving the elementary case when
$d\in\N$ is  a  complete  square,  the  following  characteristic
equations     obtained     by      K.~Calta     (analogous     to
equations~\eqref{zorich:eq:Calta:equations} above) tell whether our flat
surface is completely periodic or not:
\begin{align*}
w_1 \bar s_1 &=-w_2 \bar s_2, \\
\bar w_1 \tau_1 + \bar w_2 \tau_2 &= w_1 \bar \tau_1 + w_2 \bar
\tau_2 , \qquad {0 \le \tau_i < w_i+w_3}.
\end{align*}

Note that we consider any invariant subvariety $\cN$  as above as
a  subvariety   in  $\cH_2=\cH(1,1)\sqcup\cH(2)$.  Clearly,   the
intersection  $\cN\cap\cH(1,1)$ and  $\cN\cap\cH(2)$  results  in
close $GL^+(2,\R{})$-invariant subvarieties  in the corresponding
strata. Note that $dim_{\C{}}\cN=3$. Since  $\cN$  is  a union of
leaves   of    the    kernel    foliation   this   implies   that
$dim_{\C{}}\cN\cap\cH(2)=2$. This  means that \emph{any}  surface
$S\in\cN\cap\cH(2)$ is a Veech surface!

\paragraph{Classification Theorem of McMullen}

We  complete  this section with a description  of  the  wonderful
result of  C.~McMullen~\cite{zorich:McMullen:genus:2} realizing a  dream
of a complete classification of closures of $GL^+(2,\R{})$-orbits
in  genus  $g=2$. The classification is astonishingly simple.  We
slightly  reformulate  the original Theorem using the notions  of
kernel foliation and of completely periodic  surface (and, hence,
using implicitly results of K.~Calta~\cite{zorich:Calta}).

\begin{NNTheorem}[C.~McMullen]
\begin{itemize}
\item
If   a   surface   $S\in\cH(2)$   is   a   Veech   surface,   its
$GL^+(2,\R{})$-orbit   is  a   closed   complex   $2$-dimensional
subvariety;
\item
Closure of $GL^+(2,\R{})$-orbit of any surface $S\in\cH(2)$ which
is not a Veech surface is the entire stratum $\cH(2)$;
\item
If   a   surface   $S\in\cH(1,1)$   is  a  Veech   surface,   its
$GL^+(2,\R{})$-orbit   is  a   closed   complex   $2$-dimensional
subvariety;
\item
If a surface  $S\in\cH(1,1)$  is  not a Veech surface  but  is  a
completely  periodic   surface,   then   the   closure   of   its
$GL^+(2,\R{})$-orbit   is  a   closed   complex   $3$-dimensional
subvariety $\cN$ foliated by leaves  of  the  kernel foliation as
described above;
\item
If a surface  $S\in\cH(1,1)$ is not completely periodic, then the
closure  of   its  $GL^+(2,\R{})$-orbit  is  the  entire  stratum
$\cH(1,1)$.
\end{itemize}
\end{NNTheorem}

Actually,  the   Theorem   above   is  even  stronger:  connected
components of  these  invariant  submanifolds  are basically also
classified. We have seen that  the  $GL^+(2,\R{})$-orbit  of  any
Veech surfaces in  $\cH(2)$ has a representative with all periods
in a quadratic  field.  The \emph{discriminant} $D=b^2-4c>0$ is a
positive integer:  a  discriminant of the corresponding quadratic
equation $x^2+bx+c=0$ with integer coefficients. The discriminant
is an invariant of an $GL^+(2,\R{})$-orbit. Since for any integer
$b$  the  number  $b^2\mod  4$  can be  either  $0$  or  $1$, the
discriminant  $D\mod  4  =  0,1$.  The  values  $D=1,4$  are  not
realizable,   so    the    possible    values    of    $D$    are
$5,8,9,12,13,\dots$.

We  postpone  the   description   of  results  of  P.~Hubert  and
S.~Leli\`evre~\cite{zorich:Hubert:Lelievre:teich;discs}       and       of
C.~McMullen~\cite{zorich:McMullen:spin} on classification  of the orbits
of Veech surfaces in $\cH(2)$ to the next section. Here  we state
the following  result of C.~McMullen~\cite{zorich:McMullen:genus:2}.  By
$\cN(D)$ denote the 3-dimensional invariant submanifold  obtained
as a union of leaves of the kernel foliation passing  through the
orbits  of  all   Veech   surfaces  corresponding  to  the  given
discriminant $D$.

\begin{NNTheorem}[C.~McMullen]
The invariant subvariety  $\cN(D)$  is nonempty and connected for
any $D=0,1\mod 4$, $D\in\N$, $D\ge 5$.
\end{NNTheorem}

\paragraph{Ergodic Measures}

\index{Action on the moduli space!ofSL@of $SL(2,\R{})$|(}
\index{0SL@$SL(2,\R{})$-action on the moduli space|(}

Actually, the  Classification  Theorem  of  C.~McMullen  is  even
stronger: it also classifies the invariant measures. Consider the
``unit   hyperboloids''   $\cH_1(2)$    and   $\cH_1(1,1)$:   the
subvarieties of  real  codimension one representing flat surfaces
of  area  $1$.  The  group  $SL(2,\R{})$  acts on $\cH_1(2)$  and
$\cH_1(1,1)$ preserving  the  measure  induced  on  these  ``unit
hyperboloids'',                                               see
Sec.~\ref{zorich:ss:Moduli:Space:of:Holomorphic:One:Forms}.    Let    us
discuss what other $SL(2,\R{})$-invariant measures do we know.

When we have  an $SL(2,\R{})$-invariant subvariety, we can get an
invariant measure  concentrated  on this subvariety. For example,
as  we  know an  $SL(2,\R{})$-orbit  of a  Veech  surface $S$  is
closed; it  is isomorphic to the quotient $SL(2,\R{})/\Gamma(S)$,
where, by definition of a Veech surface, this quotient has finite
volume.  Thus,  Haar  measure  on $SL(2,\R{})$ induces  a  finite
invariant  measure  on the $SL(2,\R{})$-orbit of a Veech  surface
$S$.

Consider now a ``unit hyperboloid''  $\cN_1\subset  \cN$  in  the
manifold  $\cN$  obtained as  a  union  of  leafs  of  the kernel
foliation passing  through  $SL(2,\R{})$-orbit of a Veech surface
$S$. Note  that by Riemann bilinear relations the  area of a flat
surface  can  be expressed in terms of  absolute  periods.  Thus,
moving along the kernel foliation  we  do not change the area  of
the surface. We have seen that every leaf of the kernel foliation
is flat. Consider  the  corresponding Euclidean volume element in
each  leaf.  The  group  $SL(2,\R{})$ maps leaves of  the  kernel
foliation to leaves and respects this volume element. Thus we get
an   invariant   measure   on   $\cN_1\subset   \cN$;   near   an
$SL(2,\R{})$-orbit  of  a  Veech  surface it disintegrates  to  a
product measure.

We have associated  to any connected invariant subvariety of each
of four types  as  above a natural $SL(2,\R{})$-measure supported
on it.  One more result of C.~McMullen in~\cite{zorich:McMullen:genus:2}
tells that there are no other ergodic measures.

\paragraph{Other Properties}

The  invariant  subvarieties have  numerous  wonderful  geometric
properties. In particular, their projections to  the moduli space
$\cM$\index{0M10@$\cM_g$ -- moduli space of complex structures}\index{Moduli space!of complex structures}
of complex structures on  a  surface of genus two are  also
nice   subvarieties.    C.~McMullen    has    showed   that   the
$GL^+(2,\R{})$-orbit  of   a   Veech   surface   projects  to  an
isometrically immersed algebraic curve and $\cN(D)$ projects to a
complex surface.  Such surfaces (of \emph{complex} dimension two)
are called \emph{Hilbert modular surfaces}.

One   more   surprising   phenomenon   proved   by    C.~McMullen
in~\cite{zorich:McMullen:genus:2} concerns Veech groups of flat surfaces
in genus $g=2$.

\begin{NNTheorem}[C.~McMullen]
If the Veech group $\Gamma(S)$ of  a  flat  surface  $S\in\cH(2)$
contains a  hyperbolic element, the flat  surface $S$ is  a Veech
surface; in particular,  its Veech group $\Gamma(S)$ is a lattice
in $SL(2,\R{})$.

If the Veech group  $\Gamma(S)$  of a flat surface $S\in\cH(1,1)$
contains a hyperbolic element, and the flat surface $S$ is  not a
Veech surface, then $S$ is completely periodic. In  this case the
Veech group $\Gamma(S)$ is infinitely generated.
\end{NNTheorem}

The  Veech group  of  a flat surface  $S$  contains a  hyperbolic
element  if  and  only  if  $S$  admits  an  affine  pseudoanosov
diffeomorphism.

We complete this section  with  the following natural problem for
higher genera $g\ge 3$.

\begin{Problem}   Let   $\cK$   be   a   $GL^+(2,\R{})$-invariant
subvariety   in   some   stratum    of    holomorphic   one-forms
$\cH(d_1,\dots,d_\noz)\subset\cH_g$. Consider the  union $\cU$ of
leaves of the kernel foliation passing through $\cK$.  Is $\cU$ a
closed subvariety in $\cH(1,\dots,1)$? Similar question for other
strata.
\end{Problem}

\subsection{Classification of Teichm\"uller Discs of Veech Surfaces in $\cH(2)$}
\label{zorich:ss:Teichmuller:Discs}

\index{Surface!square-tiled|(}
\index{Square-tiled surface|(}

It is easy to check that  any  \emph{square-tiled  surface}  (see
Sec.~\ref{zorich:ss:Square:tiled:surfaces}) is a Veech surface, and thus
an $SL(2,\R{})$-orbit of any square-tiled surface is closed. Such
orbit   contains   other   square-tiled   surfaces.   Since   the
$SL(2,\R{})$-action does not change  the  area of a surface these
other square-tiled surfaces are tiled with  the  same  number  of
squares,        see        Fig.~\ref{zorich:fig:L:3:surfaces}         in
Sec.~\ref{zorich:ss:Veech:Surfaces}  and  the  Exercise related to  this
Figure.

For a fixed integer  $n$  the number of $n$-square-tiled surfaces
is finite. It would be interesting to know (and this is certainly
a part of  general  Problem from Sec.~\ref{zorich:ss:Main:Hope}) how the
square-tiled surfaces  are  arranged into orbits of $SL(2,\R{})$.
Say, we have seen in the previous section that there  are exactly
three    $3$-square-tiled     surfaces    in    $\cH(2)$     (see
Fig.~\ref{zorich:fig:L:3:surfaces}) and  that  they  belong  to the same
orbit. For $n=4$ there are already nine $4$-square-tiled surfaces
in $\cH(2)$ and they still belong to the same $SL(2,\R{})$-orbit.
The corresponding Teichm\"uller disc is a $9$-fold cover over the
modular curve. For $n=5$ there are  twenty seven $5$-square-tiled
surfaces in $\cH(2)$ and they split into two  different orbits of
$SL(2,\R{})$.

Generalizing      a       result      of      P.~Hubert       and
S.~Leli\`evre~\cite{zorich:Hubert:Lelievre:teich;discs}   obtained   for
prime  number  $n$ C.~McMullen has recently proved the  following
conjecture of P.~Hubert and S.~Leli\`evre.

\begin{NNTheorem}[C.~McMullen]
All $n$-square-tiled surfaces in $\cH(2)$, which  cannot be tiled
with $p\times q$-rectangles with $p$ or $q$ greater than $1$, get
to the  same $SL(2,\R{})$-orbit when $n\ge 4$ is  even and get to
exactly two distinct orbits when $n\ge 5$ is odd.
\end{NNTheorem}

Actually, C.~McMullen  has   classified in~\cite{zorich:McMullen:spin}
the orbits of all Veech surfaces  in  $\cH(2)$. As we have seen
in   the previous  section  Veech surfaces in $\cH(2)$ are
characterized by an integer parameter, called the
\emph{discriminant}\index{Discriminant}
$D$. For $n$-square-tiled surfaces the discriminant equals $D=n^2$.

Veech surfaces which cannot be rescaled to a square-tiled
surface are called
\emph{nonarithmetic}\index{Veech!surface!nonarithmetic}
Veech surfaces. Any nonarithmetic Veech surface in $\cH(2)$ can be rescaled to a
flat surface having all periods in a quadratic field (see the
Lemma of W.~Thurston in the previous section). The discriminant
corresponding to a nonarithmetic Veech surface is the
discriminant of this quadratic field. Of course a
$GL^+(2,\R{})$-orbit of a nonarithmetic Veech surface might have
different representatives $S$, such that all periods of $S$
belong a quadratic field. Nevertheless, for Veech surfaces in
genus $g=2$ the discriminant is well-defined: it is   an
invariant of a $GL^+(2,\R{})$-orbit. The discriminant is a
positive integer $D=0,1\mod 4$, $D\ge   5$.

C.~McMullen has proved  the following classification
Theorem~\cite{zorich:McMullen:spin}:

\begin{NNTheorem}[C.~McMullen]
For   $D=1\mod   8$,  $D>9$,  all  Veech  surfaces  in   $\cH(2)$
corresponding to discriminant $D$ get  to  exactly  two  distinct
$GL^+(2,\R{})$-orbits. For  other values $D=0,1\mod 4$, $D\ge 5$,
they belong to the same $GL^+(2,\R{})$-orbit.
\end{NNTheorem}

\begin{figure}[htb]
\centering
\includegraphics{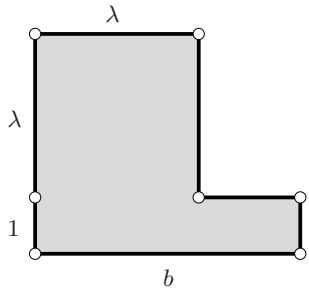}
\begin{picture}(0,5)(0,5)
\put(-64,-76){$1$}
\put(-64,-35){$\lambda$}
\put(-27,5){$\lambda$}
\put(-5,-95){$b$}
\end{picture}
\vspace{100bp}
\caption{
\label{zorich:fig:L:shaped:represent}
The $L$-shaped\index{Billiard!table!L-shaped}
billiard table $L(b,e)$ generating
a ``canonical'' Veech surface
}
\end{figure}

Moreover,   C.~McMullen   proposed    the   following   canonical
representative     for     any     such     $GL^+(2,\R{})$-orbit,
see~\cite{zorich:McMullen:spin}.   Consider   an   $L$-shaped\index{Billiard!table!L-shaped}
billiard
$L(b,e)$ as on Fig.~\ref{zorich:fig:L:shaped:represent}, where
\begin{equation}
\label{zorich:eq:b:e:lambda}
\begin{cases}
b,e\in\Z{}; \quad \lambda=(e+\sqrt{e^2+4b})/2\\
e=-1,0\text{ or } 1\\
e+1<b\\
\text{if }e=1\text{ then } b\text{ is even}
\end{cases}
\end{equation}

The billiard table  $L(b,e)$ generates a flat surface $S(b,e)$ in
$\cH(2)$, see Fig.~\ref{zorich:fig:McM:billiard}.

\begin{NNTheorem}[C.~McMullen]
The  flat  surface $S(b,e)$ generated by the $L$-shaped\index{Billiard!table!L-shaped}
billiard
table        $L(b,e)$         with        parameters        $b,e$
satisfying~\eqref{zorich:eq:b:e:lambda} is a Veech  surface.  Any closed
$GL^+(2,\R{})$-orbit  in  $\cH(2)$ is represented by one of  such
$S(b,e)$ and this representation is unique.  The discriminant $D$
of $S(b,e)$ equals $D=e^2+4b$.
\end{NNTheorem}

\begin{Exercise}
Using  linear  transformations  and  scissors  rescale  the flat
surface   obtained    from    the    ``double   pentagon'' (see
Fig.~\ref{zorich:fig:different:unfoldings}) to the flat surface
obtained from      a       ``golden      cross''      $P(a,b)$,
with $a=b=\cfrac{1+\sqrt{5}}{2}$; see
Fig.~\ref{zorich:fig:McM:billiard}  for the definition of $P(a,b)$.
Using  linear  transformations  and  scissors  rescale  any of
these flat surfaces to a surface obtained from the billiard
table $L(1,-1)$, see Fig.~\ref{zorich:fig:L:shaped:represent}, proving
that these Veech surfaces have discriminant $D=5$. (For
solutions see Fig.~4 in~\cite{zorich:McMullen:Hilbert}; see also
and~\cite{zorich:McMullen:decagon}).
\end{Exercise}

\begin{Exercise}[T.~Schmidt]
Using  linear  transformations  and  scissors  rescale  the flat
surface   obtained    from    the  regular octagon to a surface
obtained from the billiard table $L(2,0)$, see
Fig.~\ref{zorich:fig:L:shaped:represent}, proving that these Veech
surfaces have discriminant $D=8$. (See~\cite{zorich:McMullen:decagon}).
\end{Exercise}

\paragraph{Stratum $\cH(1,1)$}

We have  discussed in details  Veech surfaces in $\cH(2)$. We did
not  discussed  the Veech  surface  for  the  stratum  $\cH(1,1)$
because for  square-tiled  surfaces  in  $\cH(1,1)$  (also called
\emph{arithmetic  Veech  surfaces})  the  classification  of  the
$GL^+(2,\R{})$-orbits is knot known yet...

The  classification  of \emph{nonarithmetic}  Veech  surfaces  in
$\cH(1,1)$ is, however, known (since very recently), and is quite
surprising. Using the results of
M.~Moeller~\cite{zorich:Moeller:Galois}--\cite{zorich:Moeller:Periodic:points}
C.~McMullen has proved in~\cite{zorich:McMullen:decagon:proof} the following Theorem.

\begin{NNTheorem}[C.~McMullen]
Up  to  a  rescaling  the  only  primitive nonarithmetic  Veech  surface in
$\cH(1,1)$  is  the  one  obtained  from  the regular decagon  by
identification  of   opposite   sides.   In   other   words,  any
primitive nonarithmetic  Veech   surface   in  $\cH(1,1)$  belongs  to  the
$GL^+(2,\R{})$-orbit of  the  surface  obtained  from the regular
decagon.
\end{NNTheorem}

For higher  genera nothing is known  neither about the  number of
$SL(2,\R{})$-orbits of $n$-square-tiled surfaces, nor about their
geometry.

\begin{Problem}
Classify  orbits  of  square-tiled  surfaces in any  stratum,  in
particular in $\cH(1,1)$.
\end{Problem}

\toread{Square-tiled Surfaces}
\index{Surface!square-tiled}
An elementary introduction can  be
found                              in~\cite{zorich:Zorich:square:tiled}.
Paper~\cite{zorich:Hubert:Lelievre:teich;discs}    of    P.~Hubert   and
S.Leli\`evre  and~\cite{zorich:McMullen:spin}  of  C.~McMullen  classify
orbits of square-tiled surfaces  in  $\cH(2)$. See also the paper
of  G.~Schmith\"usen~\cite{zorich:Schmithusen}   for  an  algorithm   of
evaluation of  the Veech group  of a square-tiled surface and for
examples of square-tiled surfaces having $SL(2,\R{})$  as a Veech
group. Another such example due  to  M.~M\"oller  is presented in
the   survey~\cite{zorich:Hubert:Schmidt:Handbook}   of  P.~Hubert   and
T.~Schmidt.

\index{Veech!surface|)}
\index{Surface!Veech|)}

\index{Surface!square-tiled|)}
\index{Square-tiled surface|)}

\section{Open Problems}
\label{zorich:s:Open:Problems}

\paragraph{\textbf{Flat  surfaces with  nontrivial  holonomy  and
billiards in general polygons}}

\begin{FProblem}[Geodesics on general flat surfaces;
Sec.~\ref{zorich:ss:Flat:Surfaces}]

Describe behavior  of  geodesics  on  general  flat surfaces with
nontrivial  holonomy.  Prove (or disprove) that geodesic flow  is
ergodic on a typical (in any reasonable sense) flat surface.

Does  any  (almost any)  flat  surface has  at  least one  closed
geodesic which does not pass through singular points?

If yes, are there many regular closed geodesics? Namely, find the
asymptotics for the  number of closed geodesics of bounded length
as a function of the bound.
\end{FProblem}

\begin{FProblem}[Billiards in general polygons;
Sec.~\ref{zorich:ss:Billiards:in:Polygons}]

\index{Billiard!in polygon}
\index{Billiard!trajectory!generic}
\index{Billiard!trajectory!periodic}
\index{Billiard!triangular}

Describe the behavior of a generic regular billiard trajectory in
a generic triangle, in  particular,  prove (or disprove) that the
billiard flow is ergodic\index{Ergodic}.

Does any (almost  any) billiard table  has at least  one  regular
periodic trajectory?  If  the  answer  is  affirmative, does this
trajectory survive under deformations of the billiard table\index{Billiard!table}?

\index{Trajectory!billiard trajectory}
\index{Billiard!counting of periodic trajectories}

If  a  periodic  trajectory  exists,  are   there  many  periodic
trajectories like  that?  Namely,  find  the  asymptotics for the
number of periodic trajectories  of  bounded length as a function
of the bound.
\end{FProblem}

More problems on billiards can be found in the survey of
E.~Gutkin~\cite{zorich:Gutkin:survey:on:billiards:2}.

\begin{FProblem}[Renormalization of billiards in polygons;
Sec.~\ref{zorich:ss:Billiards:in:Polygons}
and Sec.~\ref{zorich:s:Renormalization:Rauzy:Veech:Induction}]

Is there  a  natural dynamical system (renormalization\index{Renormalization}
procedure) acting on the space of billiards in polygons?
\end{FProblem}

\paragraph{\textbf{Classification    of    orbit   closures
in $\cH_g$ and $\cQ_g$}}

\begin{FProblem}[Orbit closures for moduli spaces;
Sec.~\ref{zorich:ss:Main:Hope}]

Is it true that the  closures  of  $GL^+(2,\R{})$-orbits in $\cH_g$\index{0H10@$\cH_g$ -- moduli space of holomorphic 1-forms}\index{Moduli space!of holomorphic 1-forms}
and  $\cQ_g$\index{0M10@$\cQ$ -- moduli space of quadratic differentials}\index{Moduli space!of quadratic differentials}
are  always  complex-analytic  (complex-algebraic?)
orbifolds? Classify these closures. Classify ergodic measures for
the action of $SL(2,\R{})$ on ``unit hyperboloids''.

\index{Action on the moduli space!ofGL@of $GL^+(2,\R{})$|)}
\index{0GL@$GL^+(2,\R{})$-action on the moduli space|)}

Suppose that these orbit closures  are  described  by an explicit
list. Find natural intrinsic invariants of  a  flat  surface  $S$
which would allow to determine the closure of the orbit of $S$ in
the list.
\end{FProblem}

To be honest, even having obtained a conjectural classification
above, one would need to develop a serious further machinery to
get full variety of interesting applications. The situation with the problem below
is quite different: a reasonable solution of this problem
would immediately give a burst of applications since such a machinery
already exists. For the experts interested in ergodic aspects, counting problems,
etc, the measure-theoretic analogue of Ratner's Theorem discussed
below is the biggest open problem in the area.

\begin{FProblem}[\textbf{A.~Eskin}: ``Ratner's theorem'' for moduli spaces;
Sec.~\ref{zorich:ss:Main:Hope}]

Does the unipotent subgroup of $SL(2,\R{})$ act nicely on
$\cH_g$\index{0H10@$\cH_g$ -- moduli space of holomorphic 1-forms}\index{Moduli space!of holomorphic 1-forms} and  $\cQ_g$ or not? Are the orbit closures always
nice (for example, real-analytic) orbifolds or one can get
complicated closures (say, Kantor sets)?

If the action is ``nice'', classify the closures of orbits of
the unipotent subgroup $\begin{pmatrix}1 & t \\ 0 &
1\end{pmatrix}_{t\in\R{}}$ in
$\cH_g$ and  $\cQ_g$\index{0M10@$\cQ$ -- moduli space of quadratic differentials}\index{Moduli space!of quadratic differentials}.
Classify
ergodic measures for the action of the unipotent subgroup on
``unit hyperboloids''.

Solve the problem for genus $g=2$.
\end{FProblem}

The problem above is solved in the particular case when the unipotent flow
acts on a $SL(2,\R{})$-invariant submanifold in $\cH_g$\index{0H10@$\cH_g$ -- moduli space of holomorphic 1-forms}\index{Moduli space!of holomorphic 1-forms} obtained
by a ramified covering construction from a Veech surface; see the papers
of A.~Eskin, H.~Masur and M.~Schmoll~\cite{zorich:Eskin:Masur:Schmoll} and of
A.~Eskin, J.~Marklof, D.~Witte~Morris~\cite{zorich:Eskin:Marklof:Witte}.

The next problem concerns possibility of construction of $GL^+(2,\R{})$-invariant
submanifolds in higher genera using kernel foliation.

It follows from the results of K.~Calta and C.~McMullen that for
{\it any} Veech surface $S_0\in\cH(2)$ the union of complex
one-dimensional leaves of the kernel  foliation passing through
the complex two-dimensional $GL^+(2,\R{})$-orbit $\cO(S_0)$ of $S_0$
is a closed complex orbifold $\cN$  of complex dimension $3$, see
Sec.~\ref{zorich:ss:Revolution:in:Genus:Two}. By construction it
is $GL^+(2,\R{})$-invariant.

\begin{FProblem}[Kernel foliation;
Sec.~\ref{zorich:ss:Revolution:in:Genus:Two}]

Let $\cO\subset\cH(d_1,\dots,d_\noz)\subset\cH_g$ be a
$GL(2,\R{})$-invariant submanifold (suborbifold). Let $\cH(d'_1,
\dots, d'_n)\subset \cH_g$ be a bigger stratum adjacent to the
first one, $d_j=d'_{k_1}+\dots+d'_{k_j}$, $i=1, \dots, \noz$.

Consider the closure of the union of leaves of the kernel
foliation in $\cH(d'_1, \dots, d'_n)$ (or in $\cH_g$\index{0H10@$\cH_g$ -- moduli space of holomorphic 1-forms}\index{Moduli space!of holomorphic 1-forms}) passing
through $\cO$. We get a closed $GL^+(2,\R{})$-invariant subset
$\cN\subset\cH(d'_1, \dots, d'_n)$ (correspondingly
$\cN\subset\cH_g$).

Is $\cN$ a complex-analytic (complex-algebraic) orbifold? When
$\cN$ does not coincide with the entire connected component of
the stratum $\cH(d'_1, \dots, d'_n)$ (correspondingly $\cH_g$\index{0H10@$\cH_g$ -- moduli space of holomorphic 1-forms}\index{Moduli space!of holomorphic 1-forms})?
When $\dim_{\C{}}\cN=\dim_{\C{}}\cO+(n-\noz)$ (correspondingly
$\dim_{\C{}}\cN=\dim_{\C{}}\cO+(2g-2-\noz)$? Here $(n-\noz)$ and
$(2g-2-\noz)$ is the complex dimension of leaves of the kernel
foliation in $\cH(d'_1, \dots, d'_n)$ and in $\cH_g$\index{0H10@$\cH_g$ -- moduli space of holomorphic 1-forms}\index{Moduli space!of holomorphic 1-forms}
correspondingly.
\end{FProblem}

\paragraph{\textbf{Particular cases of classification problem}}

\begin{FProblem}[Exceptional strata of quadratic differentials;
Sec.~\ref{zorich:ss:Classification:of:Connected:Components:of:the:Strata}]

Find an invariant which would be easy to evaluate and which would
distinguish  half-translation\index{Surface!half-translation}\index{Half-translation surface}
surfaces from  different connected
components\index{Moduli space!connected components of the strata}
of   the  four   exceptional   strata\index{0M10@$\cQ$ -- moduli space of quadratic differentials}\index{Moduli space!of quadratic differentials}
$\cQ(-1,9)$,
$\cQ(-1,3,6)$, $\cQ(-1,3,3,3)$ and $\cQ(12)$.
\end{FProblem}


\begin{FProblem}[Veech surfaces; Sec.~\ref{zorich:ss:Veech:Surfaces}]

Classify primitive Veech surfaces\index{Veech!surface}\index{Surface!Veech}.
\end{FProblem}


\begin{FProblem}[Orbits of square-tiled surfaces;
Sec~\ref{zorich:ss:Teichmuller:Discs}]

Classify orbits of square-tiled surfaces\index{Surface!square-tiled}\index{Square-tiled surface}
in any stratum. Describe
their Teichm\"uller discs.

Same problem for the particular case $\cH(1,1)$.
\end{FProblem}

\index{Action on the moduli space!ofSL@of $SL(2,\R{})$|)}
\index{0SL@$SL(2,\R{})$-action on the moduli space|)}

\paragraph{\textbf{Geometry of individual flat surfaces}}

\begin{FProblem}[Quadratic asymptotics for any surface;
Sec.~\ref{zorich:ss:Counting:Closed:Geodesics:and:Saddle:Connections}]

Is it true that \emph{any} very flat surface  has exact quadratic
asymptotics  for  the number of saddle connections  and  for  the
number of regular closed geodesics?
\end{FProblem}

\begin{FProblem}[Error term for counting functions;
Sec.~\ref{zorich:ss:Counting:Closed:Geodesics:and:Saddle:Connections}]

What can  be said about the  error term in  quadratic asymptotics
for counting functions  $N(S,L)\sim c\cdot L^2$ on a generic flat
surface $S$? In particular, is it true that the limit
$$
\limsup_{L\to\infty}
\cfrac{\log |N(S,L)-c\cdot L^2|}{\log L}\overset{?}{<}
2
$$
is strictly  less than two? Is it  the same  for almost all  flat
surfaces in a given connected component of a stratum?
\end{FProblem}

One of the key properties used by C.~McMullen for the
classification of the closures of orbits of $GL(2,\R{})$ in
$\cH(1,1)$ was the knowledge that on \emph{any} flat surface in
this stratum one can find a pair of homologous saddle
connections. Cutting the surface along these saddle connections
one decomposes the surface into two tori and applies machinery
of Ratner theorem.

\begin{FProblem}[\textbf{A.~Eskin; C.~McMullen}: Decomposition of surfaces;
Sec.~\ref{zorich:ss:Multiple:Isometric:Geodesics:and:Principal:Boundary}]

Given a connected component of the stratum $\cH(d_1, \dots,
d_\noz)$ of Abelian differentials (or of quadratic differentials
$\cQ(d_1, \dots, d_\noz)$\index{0M10@$\cQ$ -- moduli space of quadratic differentials}\index{Moduli space!of quadratic differentials})
find those configurations of
homologous saddle connections (or homologous closed geodesics),
which are presented at \emph{any} very flat surface $S$ in the
stratum.
\end{FProblem}

For quadratic differentials the notion of ``homologous'' saddle
connections (closed geodesics) should be understood in terms of
homology with local coefficients,
see~\cite{zorich:Masur:Zorich}.

\paragraph{\textbf{Topological, geometric, and dynamical
properties of the strata}}

\begin{FProblem}[\textbf{M.~Kontsevich}: Topology of strata;
Sec.~\ref{zorich:s:Families:Of:Flat:Surfaces:and:Moduli:Spaces:of:Abelian:Differentials}]

Is   it    true    that    strata   $\cH(d_1,\dots,d_\noz)$\index{0H20@$\cH(d_1,\dots,d_\noz)$ -- stratum in the moduli space}\index{Stratum!in the moduli space}   and
$Q(q_1,\dots,q_n)$  are  $K(\pi,1)$-spaces (i.e.  their universal
covers are contractible)?
\end{FProblem}


\begin{FProblem}[Compactification of moduli spaces;
Sec.~\ref{zorich:ss:Simplest:Cusps:of:the:Moduli:Space}
and~\ref{zorich:ss:Multiple:Isometric:Geodesics:and:Principal:Boundary}]

Describe  natural  compactifications  of  the  moduli  spaces  of
Abelian  differentials  $\cH_g$\index{0H10@$\cH_g$ -- moduli space of holomorphic 1-forms}\index{Moduli space!of holomorphic 1-forms}  and  of  the  moduli  spaces  of
meromorphic quadratic  differentials  with  at  most simple poles
$\cQ_g$.  Describe  natural  compactifications  of  corresponding
strata    $\cH(d_1,\dots,d_\noz)$\index{0H20@$\cH(d_1,\dots,d_\noz)$ -- stratum in the moduli space}\index{Stratum!in the moduli space}    and
$\cQ(q_1,\dots,q_n)$\index{0M10@$\cQ$ -- moduli space of quadratic differentials}\index{Moduli space!of quadratic differentials}.
\end{FProblem}


\begin{FProblem}[Dynamical Hodge decomposition;
Sec.~\ref{zorich:ss:Asymptotic:Flag:and:Dynamical:Hodge:Decomposition},
\ref{zorich:ss:Teichmuller:Geodesic:Flow} and
Appendix~\ref{zorich:s:Multiplicative:ergodic:theorem}]

Study properties  of  distributions  of  Lagrangian  subspaces in
$H^1(S;\R{})$ defined  by  the
Teichm\"uller  geodesic  flow\index{Teichm\"uller!geodesic flow}, in
particular,  their   continuity.  Is  there  any  topological  or
geometric way to define them?
\end{FProblem}


\begin{FProblem}[Lyapunov exponents;
Sec.~\ref{zorich:ss:Teichmuller:Geodesic:Flow} and
Appendix~\ref{zorich:s:Multiplicative:ergodic:theorem}]

Study \emph{individual} Lyapunov exponents of the Teichm\"uller geodesic flow\index{Teichm\"uller!geodesic flow}
\newline\noindent
-- for all known $SL(2;\R{})$-invariant subvarieties;
\newline\noindent
-- for strata in large genera.

Are  they  related to  characteristic  numbers  of  some  natural
bundles  over  appropriate  compactifications   of   the  strata?
\end{FProblem}

Some other open problems can be found
in~\cite{zorich:Hubert:Masur:Schmidt:Zorich:open:problems}.


\appendix

\section{Ergodic Theorem}
\label{zorich:s:Ergodic:Theorem}

We   closely    follow    the    presentation    in   Chapter   1
of~\cite{zorich:Cornfeld:Fomin:Sinai}. However, for the sake of  brevity
we do not consider flows; for the flows the theory  is absolutely
parallel.

\paragraph{Ergodic Theorem}

Consider a  manifold $M^n$ with  a measure $\mu$. We shall assume
that the measure comes from a  volume form on $M^n$, and that the
total  volume  (total  measure)  of  $M^n$  is  finite.  We shall
consider only measurable subsets of $M^n$.

Let $T:M^n\to M^n$ be a smooth  map. We do not assume that $T$ is
a bijection  unless it is explicitly  specified. We say  that $T$
\emph{preserves measure} $\mu$  is  for any subset $U\subset M^n$
measure $\mu(T^{-1}U)$  of  the  preimage  coincides with measure
$\mu(U)$ of the set. For example the double  cover $T:S^1\to S^1$
of the circle $S^1=\R{}/\Z{}$ over itself  defined as $T:x\mapsto
2x  \mod 1$  preserves  the Lebesgue measure  on  $S^1$. In  this
section we consider only measure preserving maps.

We say  that some property  is valid \emph{for almost all} points
of $M^n$ if it  is valid for a subset $U\subset M^n$  of complete
measure $\mu(U)=\mu(M^n)$.

A subset $U\subset M^n$ is \emph{invariant} under the  map $T$ if
the preimage  $T^{-1}U$ coincides with  $U$. Thus, a notion of an
invariant  subset  is   well-defined  even  when  $T$  is  not  a
one-to-one map.  The  measure-preserving  map  $T:M^n\to  M^n$ is
\emph{ergodic}\index{Ergodic} with respect  to $\mu$ if any invariant subset has
measure $0$ or $1$. The measure-preserving  map  $T:M^n\to  M^n$ is
\emph{uniquely ergodic}\index{Ergodic!uniquely ergodic} with respect  to $\mu$ if there is no other
invariant probability measure.

Note that if $T$ has a fixed point or, more generally,
a periodic orbit (that is $T^k(x_0)=x_0$  for some $x_0\in M^n$
and some $k\in\N$),
one can consider an invariant probability measure concentrated
at the points of the orbit. Thus, such map $T$ cannot be uniquely
ergodic with respect to any Lebesgue equivalent measure.

Now we  are ready to  formulate the keystone theorem. Consider an
integrable function $f$ on $M^n$ and some point  $x$. Consider an
orbit of $T$ of  length $n$ starting at $x$. Let us  evaluate the
values of  $f$ at  the points of the orbit,  and let us calculate
the                        ``mean                         value''
$\cfrac{1}{n}\Big(f(x)+f(Tx)+\dots+f(T^{n-1}x)\Big)$ with respect
to  the  discrete  time  $k$  of  our  dynamical  system  $\dots,
T^{k-1}x, T^kx, T^{k+1}x, \dots$.

\begin{ErgodicTheorem}
Let $T:M^n\to M^n$ preserve a finite measure $\mu$ on $M^n$. Then
for any integrable function $f$ on $M^n$ and for almost all point
$x\in M^n$ there exists the time mean: there  exist the following
limit:
$$
\lim_{n\to\infty} \cfrac{1}{n} \sum_{k=0}^{n-1} f(T^kx) =
\bar{f}(x).
$$
The function $\bar f(x)$ is  integrable  and  invariant under the
map $\bar  f(Tx) = \bar f(x)$. In particular,  if $T$ is ergodic,
$\bar f$ is constant almost everywhere. Moreover,
$$
\int_{M^n} \bar f\,d\mu = \int_{M^n} f\,d\mu
$$
\end{ErgodicTheorem}

\paragraph{First Return Map}

The following theorem allows to construct numerous \emph{induced}
dynamical systems which are closely related to the initial one.

\begin{NNTheorem}[Poincar\'e Recurrence Theorem]
For any subset $U\subset M^n$ of positive measure  and for almost
any  starting  point  $x\in  U$  the  trajectory  $x,  Tx, \dots$
eventually returns to $U$, i.e. there is some $n\ge 1$  such that
$T^n x \in U$.
\end{NNTheorem}

The minimal  $n=n(x)\in\N$  as  above  is  called the \emph{first
return  time}.   According   to   Poincar\'e  Recurrence  Theorem
integer-valued function  $n(x)$  is  defined almost everywhere in
$U$. Consider the  \emph{first  return map} $T|_U:U\to U$ defined
as $T|_U: x\mapsto T^{n(x)} x$, where $x\in U$. In other words, the
map $T|_U$ maps a point  $x\in  U$ to the point where  trajectory
$Tx, T^2 x, \dots$ first meets $U$.

\begin{NNLemma}
For  any  subset  $U\subset  M^n$ of positive measure  the  first
return  map  $T|_U:U\to U$ preserves measure $\mu$ restricted  to
$U$. If $T:  M^n\to  M^n$ is ergodic than  $T|_U:U\to  U$ is also
ergodic.
\end{NNLemma}

The first return  time  induced  by an ergodic map  $T$  has  the
following geometric property.

\begin{KacLemma}
For an ergodic diffeomorphism  $T:M^n\to  M^n$ and for any subset
$U\in M^n$ of positive measure the mean value of the first return
time equals to the volume of entire space:
$$
\int_U n(x) d\mu = \mu(M^n)
$$
\end{KacLemma}

\toread{Ergodic Theory} There are numerous nice  books on ergodic
theory. I can  recommend  a classical textbook of I.~P.~Cornfeld,
S.~V.~Fomin  and  Ya.~G.~Sinai~\cite{zorich:Cornfeld:Fomin:Sinai} and  a
recent        survey        of         B.~Hasselblatt         and
A.~Katok~\cite{zorich:Hasselblatt:Katok:1A}.

\section{Multiplicative Ergodic Theorem}
\label{zorich:s:Multiplicative:ergodic:theorem}

In this section we discuss multiplicative ergodic theorem and
the notion of
\emph{Lyapunov exponents}\index{Lyapunov exponent}\index{Exponent!Lyapunov exponent}
and  then  we  present  some  basic   facts  concerning  Lyapunov
exponents. As  an  alternative  elementary  introduction  to this
subject    we    can    recommend    beautiful    lectures     of
D.~Ruelle~\cite{zorich:Ruelle}. A comprehensive information representing
the state-of-the-art in  this subject can  be found in  the  very
recent survey~\cite{zorich:Barreira:Pesin:Handbook}.

\subsection{A Crash Course of Linear Algebra}
\label{zorich:ss:A:Crash:Course:of:Linear:Algebra}

Consider a linear transformation $A:\R{n}\to\R{n}$ represented by
a  matrix  $A\in SL(n,\R{})$. Assume that $A$  has  $n$  distinct
eigenvalues   $e^{\lambda_1},   \dots   ,   e^{\lambda_n}$;   let
$\vec{v}_1, \dots, \vec{v}_n$ be the corresponding  eigenvectors.
Note that  since $\det A=1$ we get $\lambda_1+\dots+\lambda_n=0$;
in particular, $\lambda_1>0$ and $\lambda_n<0$.

Consider now a linear  transformation  $A^N:\R{n}\to\R{n}$, where
$N$  is  a very  big  positive integer.  For  almost all  vectors
$\vec{v}\in\R{n}$ the linear transformation $A^N$ acts roughly as
follows: it takes the projection $\vec{v}_{proj}$ of $\vec{v}$ to
the line  $\cV_1=Vec(\vec{v}_1)$  spanned  by the top eigenvector
$\vec{v}_1$ and then  expands  it with a factor $e^{N\lambda_1}$.
So, roughly,  $A^N$ smashes the  whole space to the straight line
$\cV_1$ and then stretches this  straight  line  with an enormous
coefficient  of   expansion  $e^{N\lambda_1}$.  (Speaking   about
projection to  $\cV_1=Vec(\vec{v}_1)$  we mean a projection along
the   hyperplane   spanned   by   the   remaining    eigenvectors
$\vec{v}_2,\dots,\vec{v}_n$.)

To be more precise,  we have to note that the image  of $\vec{v}$
would not have exactly the direction  of  $\vec{v}_1$.  A  better
approximation  of  $A^N(\vec{v})$  would  give us a vector  in  a
two-dimensional subspace $\cV_2=Vec(\vec{v}_1,\vec{v}_2)$ spanned
by  the   two   top  eigenvectors  $\vec{v}_1,  \vec{v}_2$.  When
$\lambda_2>0$ the direction of $A^N(\vec{v})$ would be very close
to the direction of $\cV_1$ though the endpoint of $A^N(\vec{v})$
might be at a distance  of  order  $e^{N\lambda_2}$ from $\cV_1$,
which  is  very   large.  However,  this  distance  is  small  in
comparison with the length of  $A^N(\vec{v})$  which  is of order
$e^{N\lambda_1}\gg e^{N\lambda_2}$.

Further terms of approximation give us  subspaces $\cV_j$ spanned
by  the  top  $j$   eigenvectors,   $\cV_j=Vec(\vec{v}_1,  \dots,
\vec{v}_j)$.   Note   that  starting  with  some  $k\le  n$   the
eigenvalues  $e^{\lambda_k}$  become  strictly smaller than  one.
This  means  that  the  image $A^N(\vec{v})$ of {\it  any}  fixed
vector  $\vec{v}$   gets  exponentially  close  to  the  subspace
$\cV_{k-1}$.

Going into details we have to admit that vectors $\vec{v}$ from a
subspace  of  measure  zero  in  the  set  of  directions  expose
different behavior. Namely, vectors $\vec{v}$ from the hyperplane
$\cL_2$ spanned  by  $\vec{v}_2,  \dots,  \vec{v}_n$ have smaller
coefficient of distortion  than the generic ones. From this point
of  view  vectors  from  linear  subspaces  $\cL_j=Vec(\vec{v}_j,
\vec{v}_{j+1},\dots,\vec{n})$   expose   more  and   more  exotic
behavior;  in   particular,   all   vectors   from  the  subspace
$\cL_k=Vec(\vec{v}_k,     \vec{v}_{j+1},\dots,\vec{v}_n)$     get
exponentially contracted.

\subsection{Multiplicative Ergodic Theorem
for a Linear Map on the Torus}\index{Multiplicative ergodic theorem}\index{Ergodic!multiplicative ergodic theorem}\index{Theorem!multiplicative ergodic}
\label{zorich:ss:Multiplicative:ergodic:theorem:torus}

Consider  a  linear transformation  $A:\R{n}\to\R{n}$  this  time
represented  by  an  {\it  integer} matrix $A\in  SL(n,\Z{})$  as
above.  Consider  the induced map $F:\T{n}\to\T{n}$ of the  torus
$\T{n}=\R{n}/\Z{n}$.  This  map  preserves  the  natural   linear
measure on  the torus. Consider the  $N$-th iterate $F^N$  of the
map $F$, where $N$ is a {\it very large} number.

Note, that the differential of $F$ in the  natural coordinates on
the torus is represented by the constant matrix $D_{x_0} F=A$ for
any $x_0\in\T{n}$. Note  also that the differential of the $N$-th
iterate of  $F$ is represented  as a product of $N$ differentials
of  $F$   along  the  trajectory  $x_0,  F(x_0),  F(F(x_0))\dots,
F^{N-1}(x_0)$ of $x_0$:
\begin{equation}
\label{zorich:eq:cocycle}
D_x(F^N)=D_{F^{N-1}(x_0)}F\circ\cdots\circ D_{F(x_0)}F\circ
D_{x_0} F
\end{equation}
Hence, in ``linear'' coordinates we get $D_x(F^N)=A^N$. Thus, the
results of the  previous section are literarily applicable to the
local  analysis  of  the  map $F^N$, where now  vector  $\vec{v}$
should be  interpreted as a  tangent vector to the torus $\T{n}$.
In particular,  these  results have the following interpretation.
If we consider the  trajectory  $x_0, F(x_0), \dots, F^N(x_0)$ of
the initial point and  the  trajectory $x, F(x), \dots$, $F^N(x)$
of a point $x$ obtained from $x_0$ by a very small deformation in
direction $\vec{v}$, then for {\it most} of the vectors $\vec{v}$
trajectories would deviate exponentially fast one from the other;
while for {\it some special} vectors  they  would  approach  each
other exponentially fast.

Namely, we get a distribution of linear subspaces  $\cL_k$ in the
tangent space to the torus such that deforming the starting point
of  trajectory  in  any  direction  in   $\cL_k$   we   get   two
exponentially    converging    trajectories.     The     subspace
$\cL_k=Vec(\vec{v}_k, \vec{v}_{k+1},\dots,\vec{v}_n)$ is  spanned
by the eigenvectors of  the  matrix $A$ having eigenvalues, which
are smaller than one. This distribution is integrable; it defines
a so-called {\it stable} foliation.

There  is   also   a   complementary   {\it  unstable}  foliation
corresponding     to    the     distribution     of     subspaces
$\cV_{k-1}=Vec(\vec{v}_1,  \dots,\vec{v}_{k-1})$  spanned by  the
eigenvectors of  the  matrix  $A$  having  eigenvalues, which are
greater than one. Passing from  the  map $F$ to the map  $F^{-1}$
the  stable  and  unstable  foliations change the  roles:  stable
foliation of $F$ becomes unstable for $F^{-1}$ and vice versa.

When matrix $A$  has an eigenvalue (or several eigenvalues) equal
to $\pm  1$, we get  also a {\it neutral} foliation corresponding
to the distribution spanned by the corresponding eigenvectors.

\begin{Exercise}
Evaluate the limit
\begin{equation}
\label{zorich:eq:limit:in:Oseledets:theorem}
\lim_{N\to\infty} \cfrac{\log\|D_x(F^N)(\vec{v}\|}{N}
\end{equation}
for a tangent  vector $\vec{v}\in T_{x_0} \T{n}$ having a generic
direction. What  are the possible values  of this limit  for {\it
any} tangent vector $\vec{v}\in T_{x_0} \T{n}$? Show that vectors
leading to  different values of this  limit are organized  into a
flag of subspaces\index{Flag!asymptotic}\index{Asymptotic!flag}
$\cL_1\supset\cL_2\supset \dots \supset \cL_n$,
where  we  assume that  all  eigenvalues of  the  matrix $A$  are
positive  and  distinct.  How  would  this  flag  change  if some
eigenvalues would have multiplicities? Would  we  have  a flag of
subspaces defined by different values of the limit  above for the
most  general  matrix $A\in  SL(n,\Z)$  (which  may  have  Jordan
blocks, complex eigenvalues, multiplicities, ...)?
\end{Exercise}

\subsection{Multiplicative Ergodic Theorem}
\label{zorich:ss:Multiplicative:ergodic:theorem}

Consider now a  smooth measure-preserving map $F:M^n\to M^n$ on a
manifold $M^n$. We consider the case when that  the total measure
of $M^n$  is finite,  and when the map $F$  is {\it ergodic} with
respect to this measure.

Consider some generic point $x_0$. Let us study,  whether we have
convergence  of  the  limit~\eqref{zorich:eq:limit:in:Oseledets:theorem}
for tangent vectors $\vec{v}\in T_{x_0} M^n$ in this more general
situation. We can always  trivialize  the tangent bundle to $M^n$
on an  open subset of  full measure. Using this trivialization we
can   reduce   our   problem   to   the   study  of  product   of
matrices~\eqref{zorich:eq:cocycle}.  This   study   is   now  much  more
difficult   than   in  the  previous  case  since  the   matrices
$D_{F^k(x_0)}F$ are  not  constant  anymore.  The  following {\it
multiplicative ergodic theorem}\index{Multiplicative ergodic theorem}\index{Ergodic!multiplicative ergodic theorem}\index{Theorem!multiplicative ergodic}
formulated  for  general mappings
of general manifolds  mimics the simplest situation with a linear
map on the torus.

\begin{NNTheorem}[Oseledets]
Let a  smooth map $F: M^n\to M^n$  be ergodic  with respect to  a
finite measure. Then, there exists a collection of numbers
$$
\lambda_1 > \lambda_2 > \dots > \lambda_k,
$$
such that  for almost any point $x\in M$  there is an equivariant
filtration\index{Flag!asymptotic}\index{Asymptotic!flag}
$$
\R{n}\simeq T_x M^n = \cL_1 \supset \cL_2 \supset \dots \supset
\cL_k \supset \cL_{k+1}=\{0\}
$$
in the  fiber $T_x M^n$  of the  tangent bundle at  $x$ with  the
following property. For every $\vec v  \in  \cL_j  -  \cL_{j+1}$,
$j=1, \dots, k$, one has
$$
\lim_{N\to+\infty} \frac{1}{N} \log\|(DF^N)_x (\vec v)\| =
\lambda_j
$$
\end{NNTheorem}

The    multiplicative    ergodic   theorem    was    proved    by
V.~I.~Oseledets~\cite{zorich:Oseledets};   a   similar   statement   for
products   of    random    matrices   was   proved   earlier   by
H.~Furstenber~\cite{zorich:Furstenberg}.

Multiplicative    ergodic    theorem    has    several    natural
generalizations. The theorem essentially  describes  the behavior
of products~\eqref{zorich:eq:cocycle} of matrices along trajectories  of
the map $F$. Actually, matrices $D_x F$ are  not distinguished by
any special property. One can consider any matrix-valued function
$A:M^n\to GL(m,\R{})$ and study the products of matrices
$$
A(F^{N-1}(x))\cdot\dots\cdot A(F(x))\cdot A(x)
$$
along  trajectories  $x,  F(x),  \dots,  F^{N-1}(x)$  of  $F$.  A
statement completely analogous to  the  above Theorem is valid in
this more general case provided the matrix-valued function $A(x)$
satisfy some very moderate requirements. Namely, we do not assume
that $A(x)$ is continuous  or  even bounded. The only requirement
is that
$$
\int_{M^n} \log_+ \|A(x)\| \mu < +\infty,
$$
where   $\log_+(y)=\max(\log(y),0)$.   When  this   condition  is
satisfied  one  says  that  $A(x)$  defines  an  {\it  integrable
cocycle}. The numbers $\lambda_1>\dots>\lambda_k$ are called  the
{\it Lyapunov exponents} of the corresponding cocycle.

\begin{Exercise}
Formulate a ``continuous-time'' version of multiplicative ergodic
theorem\index{Multiplicative ergodic theorem}\index{Ergodic!multiplicative ergodic theorem}\index{Theorem!multiplicative ergodic}
when instead of a map $F:M^n\to M^n$ we have a flow $F_t$
which is ergodic with respect to a finite measure on  $M^n$. Show
that under  the  natural normalization the corresponding Lyapunov
exponents coincide with the  Lyapunov  exponents of the map $F_1$
obtained as an action of the flow at the time $t=1$.

Consider  a  vector bundle over $M^n$; suppose  that  the  vector
bundle is endowed with a flat connection. Formulate  a version of
multiplicative ergodic Theorem for the natural action of the flow
on such vector bundle.

Note  that  in   the  latter  case  the  Lyapunov  exponents  are
responsible  for  the ``mean holonomy'' of the  fiber  along  the
flow. Namely, we take a fiber of the vector bundle  and transport
it along  a very  long piece of trajectory of  the flow. When the
trajectory  comes  close to the starting point  we  identify  the
fibers  using  the  flat  connection and we study  the  resulting
linear transformation of the fiber.
\end{Exercise}

Note that  the choice of a norm in the fibers $V_x$ is in a sense
irrelevant. Consider two norms $\|\ \|$ and $\|\ \|'$ and let
$$
c(x) = \min_{\|\vec v\|=1} \|\vec v\|' \qquad C(x) =
\max_{\|\vec v\|=1} \|\vec v\|'.
$$
If
$$
\int_M \max\left(|\log(c(x)|,|\log(C(x)|\right) \mu < +\infty,
$$
then neither the filtration $\cL_k(x)$ nor the Lyapunov exponents
$\lambda_k$ do not change when  we  replace the norm $\|\ \|$  by
the norm $\|\ \|'$. In particular, when $M$ is a compact manifold
all nonsingular norms are equivalent.

In general, even  for smooth maps  $F:M\to M$ (flows  $F_t$)  the
subspaces  defined  by  the  terms $\cL_k(x)\subset V_x$  of  the
filtration  do  not   change   continuously  with  respect  to  a
deformation  of  the  base  point $x$. However,  these  subspaces
behave  nicely  for  maps  (flows) which have  strong  hyperbolic
behavior (see~\cite{zorich:Pollicott} for a short introduction; a recent
quite   accessible   textbook~\cite{zorich:Barreira:Pesin:AMS}   and   a
survey~\cite{zorich:Barreira:Pesin:Handbook} describing the contemporary
status of
{\it Pesin theory}\index{Pesin theory}\index{Theory!Pesin theory}).

Currently there are no general methods of computation of Lyapunov
exponents other  than  numerically.  There  are  some  particular
situations, say, when the vector  bundle  has  a  one-dimensional
equivariant subspace,  or when $F_t$ is  a homogeneous flow  on a
homogeneous  space;   in   these   rather   special   cases   the
corresponding  Lyapunov  exponents can  be  computed  explicitly.
However,  in  general  it  is extremely difficult to  obtain  any
nontrivial information (positivity, simplicity of spectrum) about
Lyapunov exponents.


\index{Abelian differential|see{Holomorphic 1-form}}
\index{Cycle!asymptotic|see{Asymptotic cycle}}
\index{Conical!point|see{Conical singularity}}
\index{Connected component of a stratum|see{Moduli space: connected components of the strata}}
\index{Connection!saddle connection|see{Saddle connection}}
\index{Differential!Abelian|see{Holomorphic 1-form}}
\index{Differential!quadratic|see{Holomorphic quadratic differential}}
\index{Integer point of the moduli space|see{Lattice in the moduli space}}
\index{Interval exchange map|see{Interval exchange transformation}}
\index{Flag!of subspaces|see{Flag: asymptotic}}
\index{Flow!Teichm\"uller geodesic flow|see{Teichm\"uller geodesic flow}}
\index{Foliation!defined by a closed 1-form on a surface|see{Foliation: measured foliation}}
\index{Form!holomorphic 1-form|see{Holomorphic 1-form}}
\index{Holomorphic!differential|see{Holomorphic 1-form}}
\index{Metric!flat|see{Surface: flat; very flat}}
\index{Moduli space!integer point of|see{Lattice in the moduli space}}
\index{Moduli space!of Abelian differentials|see{Moduli space of holomorphic 1-forms}}
\index{Multiplicative cocycle|see{Cocycle: multiplicative}}
\index{Rational polygon|see{Polygon: rational}}
\index{Polygonal billiard|see{Billiard in polygon}}
\index{Principal boundary of the moduli space|see{Moduli space: principal boundary}}
\index{Quadratic differential|see{Holomorphic quadratic differential}}
\index{Saddle!point|see{Conical singularity}}
\index{Singularity!conical|see{Conical singularity}}
\index{Space of interval exchange transformations|see{Interval exchange transformation: space of}}
\index{Stratum!connected component|see{Moduli space: connected components of the strata}}
\index{Surface!translation|see{Very flat}}
\index{Translation surface|see{Surface: very flat}}
\index{Vertical direction|see{Direction: vertical}}
\index{Very flat surface|see{Surface: very flat}}
\index{Volume of a stratum|see{Moduli space: volume of the moduli space}}

\printindex
\end{document}